\documentclass[11pt]{amsart}

\usepackage{amsfonts} 
\usepackage{amssymb}
\usepackage{amsmath} 
\usepackage{amsthm} 
\usepackage{stmaryrd}
\usepackage{mathtools}
\usepackage{color}
\usepackage{xypic}
\usepackage{tikz-cd}

\usepackage{enumitem}

\usepackage{fancyhdr}
\newcommand{\lb}{\llbracket}
\newcommand{\rb}{\rrbracket}
\newcommand{\ovl}[1]{\mkern 1.5mu\overline{\mkern-1.5mu#1\mkern-1mu}\mkern 1mu}
\newcommand{\into}{\hookrightarrow}

\newcommand{\cA}{\mathcal{A}}
\newcommand{\frA}{\mathfrak{A}}
\newcommand{\Aa}{\mathfrak{a}}
\newcommand{\A}{\mathbb{A}}
\newcommand{\Bb}{\mathfrak{b}}

\newcommand{\B}{\mathcal{B}}

\newcommand{\cC}{\mathcal{C}}
\newcommand{\C}{\mathbb{C}}
\newcommand{\E}{\mathbb{E}}

\newcommand{\F}{\mathbb{F}}

\newcommand{\bL}{\mathbb{L}}

\newcommand{\m}{\mathfrak{m}}

\newcommand{\norm}{\mathrm{norm}}
\newcommand{\OK}{\mathcal{O}}

\newcommand{\proj}{\mathbb{P}}
\newcommand{\p}{\mathfrak{p}}
\newcommand{\Q}{\mathbb{Q}}

\newcommand{\Z}{\mathbb{Z}}

\DeclareMathOperator{\ad}{ad}
\newcommand{\alg}{{\rm alg}}
\DeclareMathOperator{\Aut}{Aut}

\newcommand{\di}{{\rm di}}
\DeclareMathOperator{\End}{End}

\DeclareMathOperator{\Gal}{Gal}

\DeclareMathOperator{\GL}{GL}
\DeclareMathOperator{\Hom}{Hom}

\DeclareMathOperator{\im}{Im}
\DeclareMathOperator{\Ind}{Ind}

\DeclareMathOperator{\PGL}{PGL}
\DeclareMathOperator{\PSL}{PSL}
\DeclareMathOperator{\rad}{rad}

\newcommand{\sep}{{\rm sep}}
\DeclareMathOperator{\Sl}{\mathfrak{sl}}
\DeclareMathOperator{\SL}{{SL}}

\DeclareMathOperator{\Sym}{Sym}
\newcommand{\tame}{{\rm tame}}
\DeclareMathOperator{\tr}{tr}
\DeclareMathOperator{\univ}{univ}
\newcommand{\unr}{{\rm unr}}

\newcommand{\til}{\widetilde}
\newcommand{\BB}{{\mathbb B}}

\newcommand{\rhobar}{{\ovl{\rho}}}
\newcommand{\Qbar}{{\ovl{\Q}}}
\newcommand{\cmod}[1]{\ \mathrm{mod}\ #1}
\newcommand{\bbb}{{\mathfrak b}}
\newcommand{\act}[2]{{}^{{#1}} \! {#2}} 
\newcommand{\llaur}{{(\!(}}
\newcommand{\rlaur}{{)\!)}}

\newcommand{\onto}{\twoheadrightarrow}

\mathchardef\myhyphen="2D
\newcommand{\nekovar}{{N}ekov\'a\v r}
\newcommand{\bellaiche}{{{B}ella\"iche}}
\newcommand{\teichmuller}{{T}eichm\"uller}
\newcommand{\teichmueller}{{T}eichm\"uller}

\usepackage{hyperref}
\hypersetup{
    colorlinks=true,
    linkcolor=blue,
    urlcolor=blue, 
    citecolor= blue}
\usepackage[noabbrev, capitalize, nameinlink]{cleveref}

\newtheorem{theoremA}{Theorem}

\newtheorem{theoremAprime}{Theorem}

\newtheorem{theorem}{Theorem}[section]
\newtheorem{corollary}[theorem]{Corollary}
\newtheorem{proposition}[theorem]{Proposition}
\newtheorem{lemma}[theorem]{Lemma}

\theoremstyle{definition}
\newtheorem{definition}[theorem]{Definition}

\newtheorem{questionx}[theorem]{Question}
\newenvironment{question}
  {\pushQED{\qed}\questionx}
  {\popQED\endquestionx}

\newtheorem{examplex}[theorem]{Example}
\newenvironment{example}
  {\pushQED{\qed}\examplex}
  {\popQED\endexamplex}
\Crefname{examplex}{Example}{Examples}

\theoremstyle{remark}
\newtheorem{remarkx}[theorem]{Remark}
\newenvironment{remark}
  {\pushQED{\qed}\remarkx}
  {\popQED\endremarkx}

\newtheoremstyle{dotless}{}{}{}{}{\itshape}{}{ }{}
\theoremstyle{dotless}

\usepackage[normalem]{ulem}

\begin{document}

\setlist[itemize]{labelindent=\parindent/2, itemsep = 2pt, leftmargin=*}
\setlist[enumerate]{labelindent=\parindent/2, itemsep = 2pt, leftmargin=*}

\title{Big images of two-dimensional pseudorepresentations}
\author{Andrea Conti}
\address{Universit{\'e} du Luxembourg}
\email{andrea.conti@uni.lu}

\author{Jaclyn Lang}
\address{Temple University}
\email{jaclyn.lang@temple.edu}

\author{Anna Medvedovsky}
\address{Boston University} 
\email{medved@bu.edu}


\begin{abstract}
\bellaiche\ has recently applied Pink-Lie theory to prove that, under mild conditions, the image of a continuous 2-dimensional pseudorepresentation $\rho$ of a profinite group on a local pro-$p$ domain $A$ contains a nontrivial congruence subgroup of $\SL_2(B)$ for a certain subring $B$ of $A$. We enlarge \bellaiche's ring and give this new $B$ a conceptual interpretation both in terms of {\it conjugate self-twists} of $\rho$, symmetries that constrain its image, and in terms of the \textit{adjoint trace ring} of $\rho$, which we show is both more natural and the optimal ring for these questions in general. Finally, we use our purely algebraic result to recover and extend a variety of arithmetic big-image results for $\GL_2$ Galois representations arising from elliptic, Hilbert, and Bianchi modular forms and $p$-adic Hida or Coleman families of elliptic and Hilbert modular forms.
\end{abstract}

\maketitle


Let $\rho: \Pi \to \GL_2(A)$ be a continuous representation of a profinite group $\Pi$ on a local pro-$p$ domain~$A$. How big or small can the image of $\rho$ be? Arithmetic versions of this question, with~$\rho$ the $p$-adic Galois representation attached to a cuspidal modular eigenform --- or more often a compatible family of such $\rho$ considered adelically --- have been considered since the 1960s, when Serre proved an adelic open-image result for the Galois representation on the Tate module of a non-CM 
elliptic curve \cite{Serre68}. Serre's result was 
adapted
by Ribet \cite{Ribet85}\footnote{Ribet's work appeared in a series of papers starting in 1975; \cite{Ribet85} is a convenient reference.} and then by Momose~\cite{Momose81}
to the more 
delicate setting of modular forms, where certain symmetries naturally bound the size of the image. Recently \nekovar\ \cite{Nekovar2012} generalized their work to Hilbert modular forms.  

In the 1990s Pink began a purely algebraic study of this kind of question by characterizing closed pro-$p$ subgroups of $\SL_2(A)$ in terms of associated ``Pink-Lie" algebras. Pink's investigations fueled further exploration of arithmetic big-image questions, this time for Galois representations attached to $p$-adic Hida families of modular forms, 
by Hida \cite{Hida15} and then, accounting for Ribet--Momose-type symmetries, by Lang (second author here) \cite{Lang2016}. 
Simultaneously, \bellaiche\ began adapting Pink-Lie theory to the (pseudo)representation setting, obtaining abstract big-image
results, but in a form that was difficult to compare to the symmetry formulations of Ribet, Momose, and Lang. 

In the present work we finally unite the two approaches, refining \bellaiche-Pink-Lie theory to relate $p$-adic big-image results to natural symmetry bounds in an abstract algebraic setting. We thereby recover the $p$-adic big-image results of Ribet, Momose, \nekovar, and Lang, improving the latter, under mild conditions on
$\rhobar$. 
We also obtain the first big-image results for Galois representations attached to $p$-adic Coleman families of modular forms (rather than for associated rigid-analytic Lie algebras, as in Conti(first author here)--Iovita--Tilouine \cite{CIT2016}), $p$-adic families of Hilbert modular forms, and Bianchi modular forms. Along the way we propose shifting our 
perspective 
towards 
formulating big-image results in terms of 
rings of definition of 
adjoint representations rather than in terms of rings fixed by symmetries. 
We emphasize that our results require \underline{absolutely no arithmetic input} and are \underline{provably optimal}.  
We explain how our results apply to a wide variety of modular Galois representations, and we anticipate that this framework can yield 
even more arithmetic fruit, from understanding \emph{even} Galois representations to relating reducibility/dihedrality ideals and automorphic congruence modules. Finally we hope that the algebraic nature of our results might portend similar phenomena in higher dimension.


\subsection*{Acknowledgements} 
The authors would like to thank Jo\"el Bella\"iche for helpful discussions about his work. We also thank Gebhard B\"ockle, Keith Conrad, Shaunak Deo, Mladen Dimitrov, Haruzo Hida, David Loeffler, Robert Pollack, and Jacques Tilouine for useful comments. Thanks to David Rohrlich for help with Clifford theory. 
Finally, we thank the anonymous referee, whose detailed comments improved both the content and the exposition of this article, not once but twice. 
All three authors express gratitude to the Max Planck Institute for Mathematics for creating space for us to get together; without the hospitality of the MPIM, this collaboration would not have begun.

The first author was supported by the DFG Forschergruppe 1920, the second and third authors by the National Science Foundation through awards DMS-1604148 (second author) and DMS-1703834 (third author). This article underwent extensive revisions under soft covid lockdown conditions. The third author thanks her husband and especially her nanny for months of childcare support.

\tableofcontents
\section{Introduction}\label{intro}


\subsection{The question}\label{quest}
Let $p$ be an odd prime, $A$ a local pro-$p$ domain with maximal ideal $\m$ and (finite) residue field $\F \coloneqq A/\m$, and $\Pi$ a profinite group.  Let $\rho \colon \Pi \to \GL_2(A)$ be a continuous representation\footnote{In fact we consider $2$-dimensional pseudorepresentations, but we stick with representations for the introduction.}
with the property that the residual representation $\rhobar \coloneqq \rho \cmod{\m}$ is semisimple and \emph{multiplicity free}: either absolutely irreducible, or a sum of two distinct characters to $\F^\times$.  
Roughly, the objective is to show that the image of $\rho$ is as big as possible.

Note that if $\rho$, or its restriction to an index-$2$ subgroup of $\Pi$, is reducible, then the image of $\rho$ is both well understood and not big.   Similarly, one cannot expect a big-image result when the image of $\rho$ is up to twist isomorphic to that of $\rhobar$, as happens when $\rho$ arises from a modular form of weight one. 
Let us call these three kinds of representations \emph{a priori small}. The a priori small representations are  exactly those that are not strongly absolutely irreducible.\label{saiintro}

Suppose now that $\rho$ is not a priori small. 
We cannot expect $\rho$ to be surjective: even its determinant need not be surjective. Nor can we expect the image of $\rho$ to contain all of $\SL_2(A)$, unless the image of $\rhobar$ contains all of $\SL_2(\F)$. Following ideas of Hida, we settle on the notion of \emph{fullness}. If $B$ is any ring and $\bbb \subseteq B$ is any nonzero ideal, the subgroup of $\SL_2(B)$ given by the kernel of reduction modulo $\bbb$ is a \emph{congruence subgroup} of $\SL_2(B)$ (of \emph{level} $\bbb$):  
$$\Gamma_B(\bbb) \coloneqq \ker \big(\SL_2(B) \to \SL_2(B/\bbb) \big).$$
If the image of $\rho$, up to conjugation, contains such a congruence subgroup, we say that $\rho$ is \emph{$B$-full}. A key part of the big-image game is the search for an optimal fullness ring --- or rather for an optimal equivalence class of \emph{fullness peers}, rings that each contain an ideal of the other.

Historically, the constraints on the fullness rings of a representation $\rho$ have been described in terms of certain symmetries of $\rho$. If $\sigma$ is an automorphism of $A$ and $\eta$ is a character of $\Pi$, the pair $(\sigma, \eta)$ 
is a \emph{conjugate self-twist} of $\rho$ if applying the automorphism gives the same representation as twisting by the character: $\act{\sigma}{\rho} \cong \eta \otimes \rho$.
If $\rho$ has a 
nontrivial conjugate self-twist $(\sigma, \eta)$,
then
$\rho$ cannot be $A$-full: indeed, the equation
\begin{equation}\label{cstforce}
\sigma\big(\tr \rho(g)\big) =  \tr \act{\sigma}{\rho}(g) =  \eta(g)\tr \rho(g)
\end{equation} 
means that
the trace of $\rho(g)$ is an eigenvector for $\sigma$ viewed as a linear map over the $\sigma$-invariant scalars.
But the trace of a congruence subgroup of $A$ is not so constrained. Accordingly, the known arithmetic big-image results --- of Ribet and Momose, of \nekovar, of Lang, described in \cref{hist} below --- have all established fullness with respect to $A^{\Sigma_\rho}$, the subring of $A$ fixed by the conjugate self-twists of~$\rho$.

Stepping outside the constraints of the arithmetic setting, however, reveals shortcomings of the $A^{\Sigma_\rho}$ perspective: $A^{\Sigma_\rho}$ may simply be too big for $\rho$ to be $A^{\Sigma_\rho}$-full. For one thing, $A$ itself may not see all the conjugate self-twists of~$\rho$, as in \cref{gen cst example}, a failure of normality. Enlarging $A$ may still not suffice if $A$ has inseparable elements, as in \cref{insepex}, or worse yet, transcendental ones. The limitations are exactly those of 
Galois theory: there may not be  enough automorphisms to carve down deep enough to a fullness ring with conjugate self-twists alone. 

Rather than carving down from above, we propose building a fullness ring from below. Let $A_0$ be the \emph{adjoint trace ring} of $\rho$, 
the closed subring of $A$ topologically generated 
by the elements $(\tr \rho(g))^2/\det \rho(g)$ for $g \in \Pi$,
the traces of the adjoint representation.

It is easy to see that this generating set is both twist-invariant and fixed by all conjugate self-twists. Thus $A_0$ acts as a base ring for the conjugate self-twist automorphisms, and the question of whether $A$ or its extensions have enough conjugate self-twists turns into the usual one of Galois theory: are there enough automorphisms to isolate the base?  On the other hand, $A_0$ is a potential fullness ring free from the limitations of $A^{\Sigma_\rho}$ outlined above.  Moreover, in all arithmetic settings where we recover existing fullness results, we show that $A_0$ and $A^{\Sigma_\rho}$ are always fullness peers, so that $A_0$- and $A^{\Sigma_\rho}$-fullness are equivalent. Our first theorem shows that $A_0$ is the optimal fullness ring in all cases. 
\begin{theoremA}[{Optimality theorem. See \cref{A0optimalitythm}}] \label{opthmintro}\ \\
If $\rho$ is $B$-full for some ring $B$, then a fullness peer of $B$ is contained in $A_0$. 
\end{theoremA} 
\noindent 

Accompanying the $A_0$-optimality theorem, we present the main result of this paper: $A_0$-fullness. 
We assume that the pro-$p$ part of $\Pi$ is topologically finitely generated: for more on this \emph{$p$-finiteness} condition of Mazur, always satisfied by 
characteristic-zero 
local Galois groups and global Galois groups with ramification restricted to a finite set of places, 
see \cref{mazur}.

\begin{theoremAprime}[{{Main fullness theorem}, preliminary version. See \cref{main thm intro} and \cref{nonconst det main thm}}]\label{mainthmvague}
Suppose that $\rhobar$ satisfies a mild condition. If $\rho$ is not a priori small, then $\rho$ is $A_0$-full.  
\end{theoremAprime}

\noindent For a more precise formulation and a detailed discussion of the mild condition, see \cref{main thm intro} in \cref{mainthmsec} below.  For now we merely note that the most limiting condition for applications is a \emph{regularity} assumption: we require that the image of $\rhobar$ contain a matrix whose eigenvalue ratio differs from $\pm 1$ but is contained in residue field of $A_0$. For the most general version of our theorem, which is formulated for images of \emph{pseudorepresentations} in the sense of Chenevier \cite{ChenevierDet} as studied by \bellaiche\ \cite{Bellaiche18}, see \cref{nonconst det main thm}. 

The last stand-alone result of this paper is a refinement of \cref{mainthmvague} in the case where the residual image of $\rho$ is large. Let $\E$ be the residue field of $A_0$; here we assume that $\#\E \geq 7$. 
 \setcounter{theoremA}{2}
\begin{theoremA}[{See \cref{SL2(A0)}}]\label{residuallyfullintro}
If 
$\im \rhobar \supseteq \SL_2(\E)$ then $\im \rho$ contains $\SL_2(A_0)$ up to conjugation.  
\end{theoremA}
 \setcounter{theoremA}{1}
 
\subsection{History}\label{hist}  We now survey the history of big-image results, 
both arithmetic and algebraic, using the terminology introduced above,
to situate \cref{mainthmvague} in context.  In all of the theorems stated in \cref{hist}, $A_0$-fullness is equivalent to $A^{\Sigma_\rho}$-fullness.

\subsubsection{Classical modular forms} The big-image line of inquiry began in the late 1960s, when Serre showed that if $\rho$ comes from the $p$-adic Tate module (including for $p = 2$) of a non-CM elliptic curve over a number field $F$, so that $\Pi = \Gal(\ovl F/F)$ and $A = A_0 = \Z_p$, then $\rho$ is $\Z_p$-full \cite[Theorem IV.2.2]{Serre68}.\footnote{Serre's result is better known as an open-image theorem; and in fact he shows much more: the image of all the $p$-adic Tate modules for all $p$ at once is open adelically.}

In the 80s, Ribet and Momose generalized Serre's theorem to elliptic modular forms. Let $f$ be a cuspidal non-CM eigenform of weight at least 2. Given a prime $p$ and an embedding $\iota_p \colon \ovl{\Q} \hookrightarrow \ovl{\Q}_p$, one can associate to $f$ a 2-dimensional Galois representation $\rho = \rho_{\iota_p}$ of $\Pi = \Gal(\ovl{\Q}/\Q)$ over the ring of integers $A$ of a finite extension $\Q_p$.

\begin{theorem}[{\cite[Theorem 3.1]{Ribet85}}, {\cite[Theorem 4.1]{Momose81}}; see also \cref{ribetmomose2}] 
\label{classical modular forms}\leavevmode\\
For all but finitely many primes $p$, the representation $\rho$ is $A^{\Sigma_{\rho}}$-full.\footnote{Like Serre, Ribet and Momose prove stronger adelic big-image results.
See \cref{classical applications1} for more details.}
\end{theorem}

More recently, \nekovar\ generalized \cref{classical modular forms} to representations coming from Hilbert modular forms, in which case $\Pi$ is the absolute Galois group of a totally real number field and $A$ is still a finite extension of $\Z_p$ \cite[Appendices B.3--B.6]{Nekovar2012}.

Our main theorem (\cref{mainthmvague}) recovers the at-$p$ statements of both the Ribet--Momose and the \nekovar\ results, under the assumption that residual representation satisfies our regularity condition. See \cref{classical applications1} and \cref{classical applications2} for details. 

\subsubsection{Families of $p$-adic modular forms} \label{padicfam} Although we have stated the work of Serre, Ribet, Momose, and \nekovar\ for a fixed prime $p$ to better fit our $p$-adic framework, all of these theorems are actually adelic open-image results proved using geometric methods. Much work has been done to generalize such theorems to groups other than $\GL_2$, but that is not the direction that interests us.  Rather, we are interested in fixing $p$ and deforming representations $p$-adically, which necessitates a completely different approach. There has been some progress in this direction in special cases. Recall that we are assuming throughout that $p \neq 2$.

First we suppose that $\rho$ arises from a non-CM cuspidal Hida family. In this case $\Pi = \Gal(\ovl{\Q}/\Q)$ and $A$ is a finite domain over $\Lambda \coloneqq \Z_p \lb X \rb$.  
When $A$ is a constant extension of $\Lambda$ and the image of~$\ovl{\rho}$ contains $\SL_2(\F_p)$, Boston \cite[Proposition 3]{MazurWiles86} and Fischman \cite[Theorem 4.8]{Fischman02} show that the image of $\rho$ contains $\SL_2(A^{\Sigma_\rho})$, hence $\rho$ is $A^{\Sigma_\rho}$-full.\footnote{All the works in \cref{padicfam} consider only conjugate self-twists that fix $\Lambda$; see \cref{Hida families} for details.}

More recently, Hida proved that if $\rhobar$ is locally-at-$p$ multiplicity free then ${\rho}$ is $\Lambda$-full \cite[Theorem I]{Hida15}, but his work did not relate $\Lambda$ to $A_0$ or conjugate self-twists of $\rho$. Lang 
then improved Hida's result from $\Lambda$-fullness to $A^{\Sigma_\rho}$-fullness under the assumption that $\ovl{\rho}$ is absolutely irreducible, proving the following result.

\begin{theorem}[{\cite[Theorem 2.4]{Lang2016}}; see also {\cref{jackiethesis}}]\label{Hida families intro}
If $\rho$ arises from a non-CM cuspidal Hida family, and $\ovl{\rho}$ is absolutely irreducible and satisfies additional multiplicity-freeness conditions locally at $p$ then $\rho$ is $A^{\Sigma_\rho}$-full.  
\end{theorem}

The case when $\rho$ arises from a Coleman family was studied by Conti--Iovita--Tilouine \cite{CIT2016}.  In this case we again have $\Pi = \Gal(\ovl{\Q}/\Q)$ and $A$ is a domain finite over $\Lambda$.  In \cite[Theorem 6.2]{CIT2016} it is proved that, under hypotheses similar to those in \cref{Hida families intro}, a certain rigid analytic Lie algebra attached to $\im \rho$ contains that 
of a congruence subgroup of $A^{\Sigma_\rho}$. 
This strongly suggests that $\rho$ should be $A^{\Sigma_\rho}$-full, though this statement does not follow from \cite{CIT2016}.

\subsubsection{Abstract $p$-adic representations}
Both Hida \cite{Hida15} and Lang \cite{Lang2016} rely in a key way on results of Pink \cite{Pink93} classifying, for odd $p$, pro-$p$ subgroups of $\SL_2(A)$ in terms of a correspondence with purely algebraically defined ``Pink-Lie algebras". The analogous role in \cite{CIT2016} is played by  rigid-analytic Lie theory, whence  
the different form of the conclusion in that case. Although the conclusions of the big-image theorems in all of \cite{Hida15,Lang2016,CIT2016} are stated in terms of pure algebra 
--- a feature that is most clear in the fullness results of \cite{Hida15} and \cite{Lang2016} --- nonetheless all of their proofs are  
arithmetic in nature: they rely on special information about the restriction of $\rho$ to the local Galois group at $p$, and they 
use the 
results of Ribet and Momose as input. 

In contrast, Bella\"iche in \cite{Bellaiche18} studies the image of $\rho \colon \Pi \to \GL_2(A)$ in a purely algebraic way.  More precisely, he applies Pink's theory from \cite{Pink93} to images of 2-dimensional (pseudo)representations with \emph{constant determinant}, that is, with $\det \rho$ equal to the \teichmuller\ lift of $\det \rhobar$. \bellaiche's main application is to density results for mod-$p$ modular forms, but along the way he also proves the following theorem, under the same $p$-finiteness assumption on $\Pi$ (\cref{mazur}). 

\begin{theorem}[{\cite[Theorem 7.2.3]{Bellaiche18}}; see also {\cref{Bellaiche main theorem}}]\label{Bellaiche thm intro}
Suppose that the image of $\rhobar$ contains an element with eigenvalues in $\F_p^\times$ whose ratio is not $\pm 1$.  If $\rho$ has constant determinant and is not a priori small, then there is a subring $\B_\rho(\F_p)$ of $A$ such that $\rho$ is $\B_\rho(\F_p)$-full.
\end{theorem}

Bella\"iche's ring $\B_\rho(\F_p)$ is defined as the subring of $A$ topologically generated by a $\Z_p$-module $I_1(\rho)$ created out of the Pink-Lie algebra of $\im \rho$.  See \cref{Bellaiche main theorem} and the discussion following it for the definition of the generating set $I_1(\rho)$ and the ring $\B_\rho(\F_p)$. 

Unfortunately, as Bella\"iche himself notes, it is not straightforward to relate the ring $\B_\rho(\F_p)$ to the rings $A_0$ or $A^{\Sigma_\rho}$ from previous results.  Indeed, he gives no conceptual interpretation of $\B_\rho(\F_p)$ at all. 
The goal of the present work is to refine the definition of $\B_\rho(\F_p)$ by enlarging scalars and then give it a conceptual interpretation.  Under mild assumptions we thus recover, and in the case of \mbox{$p$-adic} families improve, the results mentioned above in a uniform and purely algebraic way. We point out that prior to Bella\"iche's work, Hida's work was the only fullness result when $\ovl{\rho}$ is reducible and $\rho$ comes from a $p$-adic family of modular forms.  In the case of Coleman families, a true fullness result was not previously known. Additionally, we obtain first results in other $\GL_2$-contexts, including Galois representations attached to Bianchi modular modular forms and to $p$-adic families of Hilbert modular forms. See \cref{applications section} for all the details.

\subsection{Main theorem} \label{mainthmsec} We now state our main fullness theorem in more detail. Recall that $A$ is a local pro-$p$ domain with maximal ideal $\m$ and residue field $\F$ and $\Pi$ is a $p$-finite profinite group. Let $\rho: \Pi \to \GL_2(A)$ be a representation with mod-$\m$ reduction $\rhobar$. Let $A_0$ be the adjoint trace ring of $\rho$ and $\E$ its residue field, so that $\E \subseteq \F \subseteq \ovl\F_p$.  We say that $\ovl{\rho}$ is \textit{regular} if there is some $g_0 \in \Pi$ such that $\ovl{\rho}(g_0)$ has eigenvalues $\lambda_0, \mu_0 \in \ovl{\F}^\times$ with $\lambda_0\mu_0^{-1} \in \E^\times \setminus \{\pm 1\}$.  See  \cref{regularity remark background} for an analysis of this condition. 
For the notion of \emph{goodness}, see \cref{good}.

\begin{theoremA}[Main fullness theorem; see also {\cref{nonconst det main thm}}]\label{main thm intro}\ \\
Assume that $\ovl{\rho}$ is regular. If the projective image of $\ovl{\rho}$ is isomorphic to $S_4$, assume that $\ovl{\rho}$ is good.  
If $\rho$ is not a priori small, then $\rho$ is $A_0$-full.  
\end{theoremA}
In fact we prove something slightly more general in that we can replace $\rho$ by a pseudodeformation $(t,d) \colon \Pi \to A$ of $\ovl{\rho}$: see \cref{nonconst det main thm} for a precise statement. 
Recall that \cref{main thm intro} is provably optimal, in the sense that, if $\rho$ is $B$-full for some subring $B$ of $A$, then a fullness peer of $B$ is contained in $A_0$ (see \cref{opthmintro} or \cref{A0optimalitythm}). 

Let us point out some features of the statement of \cref{main thm intro}.  First, the group $\Pi$ can be quite general.  For that reason, representations coming from Hilbert modular forms and their $p$-adic families are no more difficult than representations coming from elliptic modular forms or their $p$-adic families.  Similarly, since $\Pi$ need not be the absolute Galois group of a number field, the notion of oddness does not play a role in the paper.  In particular, \cref{main thm intro} applies to deformations of~$\ovl{\rho}$ when $\rhobar$ is an even Galois representation. 

The proof of \cref{main thm intro} proceeds via \bellaiche-Pink-Lie theory with various refinements and improvements. As in \bellaiche\ \cite{Bellaiche18}, we linearize the pro-$p$ normal core of the image of a constant-determinant $\rho$ by considering its Pink-Lie algebra $L$. But rather than building up an $A_0$-structure on $L$ directly and interpreting this for $\im \rho$, we proceed by showing fullness in turn for a sequence of fullness peer rings, first $\B_\rho(\E)$, and then $A^{\Sigma_\rho}$, and finally $A_0$, which then proves $A_0$-fullness for all the twists of $\rho$. More precisely, the argument proceeds in steps as follows.
\begin{enumerate}[itemsep = 3pt, parsep = 3pt, topsep = 5pt]
\item \label{1twist} We show that fullness rings and adjoint trace rings are unchanged under twisting by a character (\cref{twisting does not affect fullness,a0fixedbytwist}), so that we may assume that $\rho$ has constant determinant. 
\item \label{2fp} For constant-determinant $\rho$, we show that $A_0$ and $A^{\Sigma_\rho}$ are fullness peers (\cref{A0 Afixedbytwists fullness peers}).
Therefore it suffices to prove $A^{\Sigma_\rho}$-fullness in the constant-determinant case.

\item \label{3bell} We refine \cref{Bellaiche thm intro} to show that a constant-determinant $\rho$ is $\B_\rho(\E)$-full (\cref{formal fullness corollary}).

\noindent Here $\B_\rho(\E)$ is the $W(\E)$-algebra generated by the same $\Z_p$-module $I_1(\rho)$ as \bellaiche's $\B_\rho(\F_p)$ --- see discussion following \cref{Bellaiche thm intro} above. Although a small improvement, this is crucial for the next
step, and it allows our regularity hypothesis to be weaker than that of Bella\"iche.

\item \label{4small} We show that for $\rho$ with constant determinant, $\B_\rho(\E)$ and $A^{\Sigma_\rho}$ are fullness peers \mbox{(\cref{final conclusion}),} so that $\rho$ is $A^{\Sigma_\rho}$-full.

\noindent More precisely, we show that $\B_\rho(\E)$ is always contained in $A^{\Sigma_\rho}$, and the discrepancy between them is explained on the residue field  by a failure of conjugate self-twists of $\rhobar$ to lift to $\rho$. The residue field of $\B_\rho(\E)$ is $\E = \F^{\Sigma_{\rhobar}}$, the fixed field of residual twists (\cref{residual csts cut out E}), whereas the residue field of $A^{\Sigma_\rho}$ is $\F^{\Sigma_\rho}$, the field fixed by those residual twists that lift to twists of $\rho$.  We show both that $A^{\Sigma_\rho} = \B_\rho(\F^{\Sigma_\rho})$ and that the elements of $\F^{\Sigma_\rho}$ ``missing" in $\B_\rho(\E)$ show up in its fraction field. This is the most delicate step of the argument. 
\end{enumerate}

\noindent We do not show that non-constant-determinant $\rho$ are $A^{\Sigma_\rho}$-full: \cref{gen cst example,insepex} show that such a result is impossible in general, even under favorable regularity conditions. The key point is that $A^{\Sigma_\rho}$ depends on $\rho$ itself and may chance after twisting, whereas $A_0$ depends only on the projective representation of $\rho$. 
Nor do we show that $A_0 = A^{\Sigma_\rho}$ for constant-determinant $\rho$ --- they are merely fullness peers: see \cref{a0small}. 

The circuitous nature of our argument may well be merely a matter of our historical bias: the known arithmetic big-image results are formulated in terms of conjugate self-twists, so our original motivation was to relate \bellaiche's $\B_\rho(\F_p)$ to $A^{\Sigma_\rho}$. The advantages of $A_0$ revealed themselves only later. A more direct argument for proving \cref{main thm intro} is the subject of our current investigation.

\subsection{Structure of the paper} \label{leitfaden} 
This article is informally organized into parts as follows.

\smallskip

\noindent {\bf Background:}  \cref{background section}. We review known material: pseudorepresentations, generalized matrix algebras, Pink-Lie theory, and Bella\"iche's recent 
results from \cite{Bellaiche18}.  We also introduce our notion of regularity.

\smallskip

\noindent {\bf Our philosophy:} In \crefrange{section: fullness}{optimalitysec} we present and justify our approach: our goal is to show that $\rho$ is full with respect to the adjoint trace ring $A_0$, which is both the optimal fullness result and equivalent to the historically familiar $A^{\Sigma_\rho}$-fullness of known applications. 

\begin{itemize}[itemsep = 2pt, topsep = 2pt, parsep = 2pt] 
\renewcommand\labelitemi{--}
\item \cref{section: fullness} discusses the basic properties of fullness and fullness peer rings. 
We prove that fullness is twisting-invariant (\cref{twisting does not affect fullness}).
\item \cref{cst section} studies the relationship between $A_0$ and conjugate self-twists, paying particular attention to constant-determinant, and nearly so, settings where conjugate self-twists carve out the fraction field of $A_0$ (\cref{AcstA0}), so that $A_0$ and $A^{\Sigma_\rho}$ are fullness peers. 
\item \cref{optimalitysec} proves two related optimality results: every fullness ring has a fullness peer contained in $A_0$ (\cref{A0 optimality}) and is also fixed by all conjugate self-twists (\cref{fullfixed}).
\end{itemize}

\smallskip

\noindent {\bf Technical results:} \crefrange{rings acting}{conceptual interpretation of Bellaiche ring} are the technical heart of the paper. We prove that constant-determinant $\rho$ are $\B_\rho(\E)$-full, which implies $A^{\Sigma_\rho}$-fullness, which in turn guarantees $A_0$-fullness.
\begin{itemize}[itemsep = 2pt, topsep = 2pt, parsep = 2pt] 
\renewcommand\labelitemi{--}
\item In \cref{rings acting} we show that the \bellaiche-Pink-Lie algebra $L$ attached to $\rho$, a priori only a $\Z_p$-module, is in fact a $W(\E)$-module. Under certain regularity assumptions, $\rho$ is therefore $\B_\rho(\E)$-full. This completes Step \eqref{3bell} of the proof. 
\item In \cref{section: lifting} we study residual conjugate self-twists and their lifting properties via the universal deformation ring --- how does $\Sigma_\rhobar$ compare to $\Sigma_\rho$? \item \cref{section: technical residual} explores how {the regularity of} $\ovl \rho$ {imposes structure on residual conjugate self-twists} 
in preparation for \cref{conceptual interpretation of Bellaiche ring}. We also introduce the \textit{goodness} constraint. 
\item \cref{conceptual interpretation of Bellaiche ring} is the technical heart of the technical heart of the paper. Its  main goal is to prove that $\B_\rho(\E)$ has the same field of fractions as $A^{\Sigma_\rho}$ to show that they are fullness peers (\cref{final conclusion}), completing Step \eqref{4small} of the proof. This turns out to be intimately related to the lifting properties of conjugate self-twists of $\ovl{\rho}$ to $\rho$ explored in \cref{section: lifting}.  
\end{itemize} 

\smallskip

\noindent {\bf Interpretation and applications:} \crefrange{beyond admissible}{applications section} interpret our results for general $\rho$ and explain how to apply them to various modular form contexts in detail.
\begin{itemize}[itemsep = 2pt, topsep = 2pt, parsep = 2pt]
\renewcommand\labelitemi{--}
\item In \cref{nonadmissible main thm} we derive our main fullness results: \cref{nonconst det main thm} or \cref{repok}. 
\item \cref{residual full case} is independent of the main thrust of the paper; it gives an improvement on previous very-big-image results, showing that the image of $\rho$ contains $\SL_2(A_0)$ if it does so residually.
\item \cref{applications section} explains in detail how and to what extent \cref{main thm intro} recovers and refines known big-image results about Galois representations arising from modular forms and their $p$-adic families. 
We also apply \cref{main thm intro} to obtain new results for Galois representations attached to Bianchi modular forms (\cref{Bianchi applications}), Coleman families of elliptic modular forms (\cref{Coleman families}), and $p$-adic families of Hilbert modular forms (\cref{Hilbert families}).
\end{itemize} 

\smallskip

\noindent {\bf Appendix:} \cref{appx} houses a variety of lemmas on representation theory and commutative algebra for which we failed to find convenient references in the literature, in particular the statement that semisimple representations with isomorphic adjoints are isomorphic up to twist (\cref{isomorphic ad implies differ by twist}). No claims to originality here. 

\subsection{Leitfaden suggestions}\label{therealleitfaden}

We propose different levels of  interaction with this paper for different readers. Those who are merely interested in results and applications should read the present \cref{intro} and skip to \cref{applications section}. Those also curious about our methods should additionally
skim \cref{background section}, and then read \cref{section: fullness,cst section,optimalitysec} and \cref{nonadmissible main thm}. 
Readers who have the stomach for the technical weeds should begin the same way, but then also brave  \crefrange{rings acting}{conceptual interpretation of Bellaiche ring}. Finally, any of the previous types of readers who are interested in very-big-image results ($\im \rho$ contains an $\SL_2$ if it does so residually) should peek at \cref{residual full case}.

\subsection{Notation}\label{notation}
We establish some notation and conventions.  All rings are unital.  Given any ring $R$ (not necessarily commutative), we will let $R^\times$ denote the multiplicative group of invertible elements in $R$.  For brevity, we call a commutative ring $B$ containing a subring $C$ an \emph{extension} of $C$, \emph{finite} if $B$ is module-finite as a $C$-algebra.\label{finiteextension} 
If $B$ is a domain, then $Q(B)$  denotes its field of fractions.  If $F$ is a field, then $\ovl F$ (respectively, $F^\sep$) denotes a fixed algebraic (respectively, separable) closure.

If $C \subseteq B$ is an extension of rings, write $\Aut(B)$ for the group of ring automorphisms of $B$ and $\Aut(B/C)$ for the subgroup fixing $C$ pointwise. For $\Sigma \subseteq \Aut(B)$, write $B^\Sigma$ for the subring pointwise fixed by every $\sigma \in \Sigma$. If $\sigma \in \Aut(B)$ and $f: X \to B$ is any set map, write $\act\sigma f: X \to B$ for the map $x \mapsto \sigma(f(x)).$ 

For any integer $n$, let $\zeta_n$ denote a primitive $n^{\rm th}$ root of unity.  Given a finite field $\F$ of size $q$ a power of a prime $p$, its ring of Witt vectors, isomorphic to $\Z_p[\zeta_{q-1}]$, will be denoted by $W(\F)$.  Let $s \colon \F^\times \to W(\F)^\times$ be the Teichm\"uller lift.  We extend it to $s \colon \F \to W(\F)$ by defining $s(0) \coloneqq 0$.  If $A$ is a $W(\F)$-algebra, we use $W(\F)$ to denote the image of $W(\F)$ in $A$ under the structure map.  In particular, $s$ may be viewed as being $A$-valued by composing with the structure map.  

Throughout the paper we fix a prime $p \neq 2$.  The ring $A$ will always denote a local pro-$p$ commutative ring with maximal ideal $\m$ and residue field $\F$, which is then automatically a finite extension of $\F_p$. A closed subring of such an $A$ is also automatically local and pro-$p$.\footnote{\label{closedsublocal}{Let $B$ be closed subring of $A$. Its ideal $\m_B := B \cap \m$ is maximal since $\F_B := B/\m_B$ is a subring of $\F$. Given any $\alpha \in A$ one can check that the sequence $\{\alpha^{p^n}\}_n$ converges and that its limit is $s(\ovl \alpha)$, the \teichmuller\ lift of the image of $\alpha$ in $\F$. Since $B$ is closed, it thus contains both $s(\F_B^\times)$ and inverses of elements of $1 + \m_B$. Therefore every element of $B - \m_B$ is invertible, and $B$ is local. And as a closed subgroup of a pro-$p$ group, $B$ is automatically pro-$p$.}} Since $2$ is invertible in $A$, we can always take square roots of elements $x \in 1 + \m$ via the formula
\[
\sqrt{x} \coloneqq \sum_{n = 0}^\infty {{1/2}\choose{n}}(x - 1)^n.
\]
In particular, when we write $\sqrt{x}$, we always choose the root congruent to $1$ modulo $\m$. In general, the profinite topology on $A$ is coarser than the $\m$-adic topology, but if such an $A$ is noetherian, then the profinite topology coincides with the $\m$-adic topology, so that $A$ is a complete local noetherian ring \cite[Proposition 2.4]{deSmitLenstra}.  In this case, every finite $A$-module is equipped with its natural $\m$-adic $A$-module topology, which is compatible on submodules by the Artin-Rees lemma --- see \cite[Theorem 10.11]{AtiyahMacdonald}.
\label{uniquetop} In particular, every ideal in $A$ is closed. By the Cohen structure theorem \cite[Theorem 7.7]{EisenbudCommAlg}, if $A$ is noetherian then it is a quotient of $W(\F)\lb x_1, \ldots, x_n \rb$ for some $n$.
If $A$ is additionally a domain then $A$ enjoys the so-called $N2$ (or sometimes ``japanese") property: the integral closure of $A$ in a finite extension of its field of fractions is finite over $A$ \cite[Chapter~12, proof of Corollary~2]{matsumura}, and hence a pro-$p$ local noetherian  domain in its own right.\label{n2} Note that in this case there need not be a topology on $Q(A)$ under which {$Q(A)$} is a topological field containing $A$ as a closed subring: since $A$ is compact, the existence of any such topology would mean that $Q(A)$ is a locally compact field, so that $Q(A)$ is a finite extension of $\F_p$, $\Q_p$, or $\F_p\llaur X \rlaur$ \cite[Theorem~4.12]{FANT}. But our $A$ are more general.\label{notop}

If $M$ is a subset of a $W(\F)$-module $N$, then we will write $W(\F)M$ for the $W(\F)$-linear span of $M$ in $N$.  When $R$ is a topological ring and $S$ a subring of $R$, we say that $S$ is \textit{topologically generated} by a set $X$ if $S$ is the smallest closed subring of $R$ containing $X$.  Similarly, we can talk about an additive subgroup or a $W(\F)$-algebra topologically generated by a set.   

Finally, $\Pi$ always denotes a $p$-finite profinite group (\cref{p finite condition}).  If $\ovl{\rho} \colon \Pi \to \GL_2(\F)$ is a representation over a finite field $\F$, we can compose $\ovl{\rho}$ with the natural projection \mbox{$\proj \colon \GL_2(\F) \to \PGL_2(\F)$}.  We shall refer to the image of $\Pi$ under the composition $\proj \circ \ovl{\rho}$ as \textit{the projective image of $\ovl{\rho}$}.  It is well known that the projective image of $\ovl{\rho}$ is cyclic, dihedral, or isomorphic to $A_4, S_4, A_5$ or one of $\PSL_2(\F')$ or $\PGL_2(\F')$ for some subfield $\F'$ of $\F$ (\cite[Chapter XII]{Dickson58} or see \cite[Section~3.1]{Bellaiche18})\label{FWIW}.  If $\proj \ovl{\rho}(\Pi) \cong A_4$ (respectively, $S_4, A_5$), we say that $\ovl{\rho}$ is \textit{tetrahedral} (respectively, \textit{octahedral}, \textit{icosahedral}).  If $\ovl{\rho}$ is tetrahedral, octahedral, or icosahedral, then we say that $\ovl{\rho}$ is \textit{exceptional}.  If $\proj\ovl{\rho}(\Pi)$ contains $\PSL_2(\F_p)$ and $\ovl{\rho}$ is not exceptional, then we say that $\ovl{\rho}$ is \textit{large}.  Be warned that there are exceptional isomorphisms $\PSL_2(\F_3) \cong A_4, \PGL_2(\F_3) \cong S_4, \PSL_2(\F_5) \cong A_5$.

\vspace{0.5cm}

\medskip

\section{Bella\"iche-Pink-Lie theory}\label{background section}

In this section we introduce the basic objects of study --- 2-dimensional pseudorepresentations and their associated realizations over generalized matrix algebras --- along with the primary tools we use to study them: Pink-Lie algebras and Bella\"iche's structure theorem (\cref{Bellaiche main theorem}).  The main reference for everything in this section is Bella\"iche's paper \cite{Bellaiche18}, which we refer to for most proofs.  The only exception is that our definition of regularity (\cref{regular}) is weaker than that of Bella\"iche.

\subsection{Pseudorepresentations}\label{pseudo-reps background}
In this section we summarize the definitions and notation related to two-dimensional pseudorepresentations, algebraic gadgets introduced by Chenevier in \cite{ChenevierDet} (where they are called ``determinants") to mimic the behavior of trace and determinant functions of true $2$-dimensional representation of groups. We follow  \cite[Example 1.8]{ChenevierDet} and \cite[2.1.1]{Bellaiche18} for our definitions. 

\subsubsection{Abstract pseudorepresentations}
\begin{definition}\label{pseudo-rep definition}
A ($2$-dimensional) \textit{pseudorepresentation} of a group $G$ over a commutative ring $B$ is a pair of functions $t \colon G \to B$ and $d \colon \Pi \to B^\times$ such that
\begin{enumerate}
\item $t(1) = 2$;
\item $d(gh) = d(g)d(h)$ for all $g, h \in G$;
\item \label{pseudocharid} $t(gh) = t(hg)$ for all $g, h \in G$;
\item $t(gh) + d(h)t(gh^{-1}) = t(g)t(h)$ for all $g, h \in G$.
\end{enumerate}
If $G$ is a topological group and $B$ is a topological ring, we say that a pseudorepresentation \mbox{$(t,d) \colon G \to B$} is \textit{continuous} if $t$ and $d$ are continuous maps. 
\end{definition} One can verify that if $\rho \colon G \to \GL_2(B)$ is a (continuous) representation, then $(\tr \rho, \det \rho)$ is a (continuous) $2$-dimensional pseudorepresentation. Conversely, if $B$ is an algebraically closed field, then every pseudorepresentation $(t, d): G \to B$ is carried by a unique semisimple representation $G \to \GL_2(B)$. If $2$ is invertible in $B$, then a $2$-dimensional pseudorepresentation $(t, d)$ is determined by $t$ alone: setting $h = g$ in \eqref{pseudocharid} above yields \mbox{$d(g) = \frac{t(g)^2 - t(g^2)}{2}$}. In this way, pseudorepresentations generalize earlier work of Wiles \cite{Wiles90}, Taylor \cite{Taylor91}, and Rouquier \cite{Rouquier96} on \emph{pseudocharacters}, which mimic representations by keeping track only of a trace function. 

If $(t, d) \colon \Pi \to B$ is a (continuous) pseudorepresentation and $\chi \colon \Pi \to B^\times$ is a (continuous) character, then $(\chi t, \chi^2d)$ is also a (continuous) pseudorepresentation, called the \textit{twist of $(t,d)$ by $\chi$}. 
We say $(t,d)$ is \textit{reducible} if $t = \chi_1 + \chi_2$ with $\chi_i \colon \Pi \to  B^\times$ characters.  Otherwise $(t,d)$ is \textit{irreducible}. We say $(t,d)$ is \textit{dihedral} if it is irreducible and there is a nontrivial character $\eta \colon \Pi \to  B^\times$ such that $(\eta t, \eta^2d) = (t, d)$. 

The \emph{kernel} of a pseudorepresentation $(t, d): G \to B$ is 
$$\ker (t, d):= \{g \in G: d(g) = 1 \mbox{ and } t(gx) = t(x) \mbox{ for all } x \in G\} \subseteq G.$$ This is a normal subgroup of $G$, closed if $(t,d)$ is continuous. Moreover, $(t, d)$ factors through the quotient $\Pi/\ker(t,d)$.

\subsubsection{Pseudorepresentations of profinite groups over pro-$p$ rings}\label{propprep} 

We henceforth assume that all pseudorepresentations are continuous if both the group and the ring have topologies. Recall that $A$ is a local pro-$p$ commutative ring with maximal ideal $\m$ and residue field $\F$ and $\Pi$ is a profinite group. 

If $(t,d): \Pi \to A$ is a pseudorepresentation, its \emph{residual} pseudorepresentation $(\bar t, \bar d): \Pi \to \F$ is obtained by composing $(t, d)$ with reduction modulo $\m$. Because the Brauer group of $\F$ is trivial, the semisimple representation $\rhobar$ carrying $(\bar t, \bar d)$ is always realizable over $\F$ (see also \cref{atmostquad}).

\begin{definition} \label{constantdet}
A pseudorepresentation $(t, d): \Pi \to A$  has \emph{constant determinant} if $d$ is the \teichmueller\ lift of its reduction modulo $\m$: that is, $d = s(\bar d)$. Since $A^\times \cong s(\F^\times) \times 1 + \m$, we see that $d$ is always the product of $s(\bar d)$ and a pro-$p$ character $d_1: \Pi \to 1 + \m$. The twist of $(t,d)$ by $d_1^{-1/2}$ is the unique constant-determinant pseudorepresentation with the same $\rhobar$: this is the \emph{constant determinant twist} of $(t, d)$. 
\end{definition}

If $(t,d) \colon \Pi \to A$ is a pseudorepresentation, we call the subring of $A$ topologically generated by $t(\Pi)$ the \textit{trace ring} of $(t,d)$.  Note that the residue field of the trace ring is the trace ring of $(\bar t, \bar d)$ or $\rhobar$.  Sometimes we need the residue field $\F$ of $A$ to be a quadratic extension of the trace ring of $(\bar t, \bar d)$.  Thus if $(t, d) \colon \Pi \to A$ is a pseudorepresentation and $\F$ is the residue field of $A$, we define the \textit{trace algebra} of $(t,d)$ to be the $W(\F)$-subalgebra of $A$ topologically generated by $t(\Pi)$.  Thus the residue field does not change when restricting the codomain of a pseudorepresentation to its trace algebra. Both the trace ring and the trace algebra are pro-$p$ local rings (\cref{closedsublocal}).

\begin{definition} \label{irreducible dihedral} 
We call a pseudorepresentation $(t,d): \Pi \to A$ is \emph{a priori small} if it is reducible or dihedral, or if the kernel of its constant-determinant twist is equal to the kernel of $\rhobar$.
\end{definition} 

If $A$ is a domain it turns out that the a priori small notion coincides with a certain weak kind of reducibility. Recall that a representation $\rho: G \to \GL_2(F)$ over a field $F$ is \emph{strongly (absolutely) irreducible} if $\left.\rho\right|_H$ is (absolutely) irreducible for any finite-index subgroup $H$ of $G$.

\begin{proposition} \label{apssai}
Let $(t, d): \Pi \to A$ be a pseudorepresentation to a local pro-$p$ domain $A$ with field of fractions $K$, and let $\rho: \Pi \to \GL_2(\ovl K)$ be the semisimple representation carrying $(t, d)$.
The following are equivalent: 
\begin{enumerate}
\item\label{aps} $(t, d)$ is a priori small;
\item\label{finim} $\rho$ is reducible, dihedral, or its constant-determinant twist has finite image;
\item\label{sai} $\rho$ is not strongly irreducible. 
\end{enumerate}
\end{proposition}

\noindent Note that if any twist of $\rho: \Pi \to \GL_2(\ovl K)$ has finite image, then the constant-determinant twist does; equivalently, the image of the projective representation $\proj \rho: \Pi \to \PGL_2(\ovl K)$ is finite.

\begin{proof} Since all the notions in question are twist-invariant, we may replace $(t, d)$ and $\rho$ by constant-determinant twists. Clearly \eqref{aps} $\implies$ \eqref{finim} $\implies$  \eqref{sai}. 
To see that  \eqref{sai} $\implies$ \eqref{finim}, we follow \cite[Theorem~2.3]{RibetEllAdic}. Suppose $\left.\rho\right|_H$ reducible for some finite-index subgroup $H$ of $\Pi$. Up to replacing $H$ with its normal core ({i.e., the intersection of all conjugates}, which is still of finite index in $\Pi$), we may assume that $H$ is normal in $\Pi$, so that by Clifford's theorem \cite[Theorem 7.1.1]{Craven}, $\left.\rho\right|_H$ is semisimple, and hence has abelian image. If $\rho(H)$ is not contained in the center of $\GL_2(\ovl K)$ then $H$ contains a semisimple element $h$ with distinct $\rho$-eigenvalues, and $\rho(H)$ is contained in the maximal torus centralizing $\rho(h)$. Moreover, since $H$ is normal in $\Pi$, all of $\rho(\Pi)$ is contained in the normalizer of $\rho(h)$, so that $\rho(H)$ has index $1$ or $2$ in $\rho(\Pi)$, and $\rho$ is either reducible or dihedral. On the other hand, if $\left.\rho\right|_H$ is scalar, then since its trace is $A$-valued, $\left.\rho\right|_H= \alpha \oplus \alpha$ for some character $\alpha: H \to A^\times$; since $\rho$ has constant determinant, $\alpha^2 = s(\ovl\alpha^2)$, so that $\left.\rho\right|_H$ takes values in the finite set of prime-to-$p$ roots of unity in $A$, whence $\rho$ has finite image.

For \eqref{finim} $\implies$ \eqref{aps} it suffices to consider $(t, d)$ residually exceptional or large, in which case by Rouquier-Nyssen $\rho$ is the base change of a representation $\rho_A: \Pi \to \GL_2(A)$ with $\ker \rho_A = \ker (t,d)$ (see also \cref{Bellaiche 2.4.2} and \cref{B and C in A}). We show that if $\rho$ has finite image, then reduction modulo the pro-$p$ subgroup $1 + M_2(\m)$ induces an isomorphism \mbox{$\rho(\Pi) \cong \rhobar(\Pi)$}. If $A = \F$ there is nothing to show, so we can assume that $A$ is infinite. (Indeed, if $A$ has characteristic zero, then $A$ contains $\Z_p$; otherwise, $A$ is a local $\F$-algebra with residue field $\F$, and any such ring that is also finite over $\F$ is equal to $\F$.) We claim that the projective image of $\rho$ is isomorphic to that of $\rhobar$: If $A$ has characteristic zero, then the finite subgroups of $\GL_2(K)$ are the same as those of $\GL_2(\C)$ \cite[Proposition 16]{SerrePropGalois}, hence the projective image of $\rho$ is isomorphic to $A_4, S_4$ or $A_5$. And if $A$ has characteristic $p$, then the finite subgroups of $\GL_2(A)$ are all defined over $\F$, because the eigenvalues of any finite-order element are roots of unity.  In any case, the kernel of the reduction map $\rho_A(\Pi) \to \rhobar(\Pi)$ is pro-$p$, so that it can be seen on the map $\proj \rho_A(\Pi) \to \proj \rhobar (\Pi)$. But none of $A_4$, $S_4$, $A_5$, $\PSL_2(\F')$, or $\PGL_2(\F')$ have normal subgroups of $p$-power order for $p > 2$. 
\end{proof}

\subsubsection{Pseudodeformations}

Fix a continuous semisimple representation $\ovl{\rho} \colon \Pi \to \GL_2(\F)$.

\begin{definition}\label{pseudo-def definition}
We say that a pseudorepresentation $(t, d) \colon \Pi \to A$ is a \textit{pseudodeformation} of $\ovl{\rho}$ if \mbox{$(t, d) \equiv (\tr \ovl{\rho}, \det \ovl{\rho}) \bmod \m$}. 
\end{definition}

Let $\cC$ be the category of local pro-$p$ commutative rings with residue field $\F$, which have a natural $W(\F)$-algebra structure, and with morphisms being local continuous $W(\F)$-algebra homomorphisms.  We are interested in the deformation functors
\begin{align*}
F &\colon \cC \to \underline{\mbox{SET}}\\
A &\mapsto \{(t, d) \colon \Pi \to A \text{ pseudodeformation of } \ovl{\rho}\}.
\intertext{and}  
G &\colon \cC \to \underline{\mbox{SET}}\\
A &\mapsto \{(t, d) \in F(A) \colon d = s(\det \ovl{\rho})\}.
\end{align*}

These functors are always representable.  In order for the representing ring to be noetherian, we need to impose a finiteness condition on $\Pi$ due to Mazur, which we now recall.

\begin{definition}\cite[\S 1.1]{Mazur1989}\label{p finite condition}\label{mazur}
A profinite group $\Pi$ \textit{satisfies the $p$-finiteness condition} or is \textit{$p$-finite} if, for every open subgroup $\Pi_0$ of $\Pi$, the set $\Hom(\Pi_0, \F_p)$ is finite.
\end{definition}

It is well known that $F$ is represented by a pro-$p$ local \textit{noetherian} $W(\F)$-algebra $\tilde{\cA}$ whenever $\Pi$ is a $p$-finite profinite group.  See, for example, \cite[Proposition 3.3]{ChenevierDet} or \cite[Proposition 2.3.1]{BoeckleDef}.  In particular, the trace algebra of any pseudorepresentation of a $p$-finite profinite group on a local pro-$p$ ring is a quotient of $\tilde{\cA}$ and hence noetherian.  Let $(t^{\univ}, d^{\univ}) \colon \Pi \to \tilde{\cA}$ be the universal pseudodeformation of $\ovl{\rho}$.  It is easy to see that the constant-determinant condition is a deformation condition.  Indeed, let $\Aa$ be the ideal of $\tilde{\cA}$ topologically generated by $\{d^{\univ}(g) - s(\det \ovl{\rho}(g)) \colon g \in \Pi\}$.  Then $\cA \coloneqq \tilde{\cA}/\Aa$ represents $G$.  In particular, $\cA$ is also a pro-$p$ local noetherian $W(\F)$-algebra with residue field $\F$.  We use $(T, d) \colon \Pi \to \cA$ to denote the universal constant-determinant pseudodeformation.

\begin{definition}\label{res mult free}
If $\F'$ is a subfield of $\F$, then we say that a 2-dimensional semisimple representation $\ovl{\rho}$ is \textit{multiplicity free over $\F'$} if either $\ovl{\rho}$ is absolutely irreducible or $\ovl{\rho} \cong \chi_1 \oplus \chi_2$ such that \mbox{$\chi_1, \chi_2 \colon \Pi \to \F'^\times$} are distinct characters. 
\end{definition}

The following notion of admissibility, introduced by \bellaiche, plays a central role in \cite{Bellaiche18}.

\begin{definition}\label{admissible}\cite[Section 5.2]{Bellaiche18}
A tuple $(\Pi, \ovl{\rho}, t, d)$ is an \textit{admissible pseudodeformation over $A$} if the following conditions are satisfied:
\begin{enumerate}
\item $\Pi$ is a $p$-finite profinite group;
\item $\ovl{\rho} \colon \Pi \to \GL_2(\F)$ is a continuous representation that is multiplicity free over $\F$;
\item $(t,d) \colon \Pi \to A$ is a continuous pseudodeformation of $\ovl{\rho}$;
\item $d(g) \in s(\F^\times)$ for all $g \in \Pi$, that is, $(t,d)$ has constant determinant;
\item $A$ is the trace algebra of $(t,d)$. 
\end{enumerate}
\end{definition}
\noindent A local pro-$p$ $A$ accepting an admissible pseudodeformation is a complete noetherian local ring. 

\subsection{GMAs and $(t,d)$-representations}\label{gma background}
It is natural to ask when a given pseudodeformation $(t, d) \colon \Pi \to A$ arises as the trace and determinant of an actual representation $\rho \colon \Pi \to \GL_2(A)$.  This has been studied in great generality; see the introduction of Chenevier's paper \cite{ChenevierDet} for a thorough history.  Bella\"iche and Chenevier \cite[Section 1.4]{Bellaiche-ChenevierBook} have shown that, under the residual multiplicity-free assumption, $(t,d)$ always comes from a representation if one allows something more general than matrix algebras for the target.  In \cref{gma background} we summarize Bella\"iche's \cite[Section~2]{Bellaiche18}, where he specializes his work with Chenevier to the 2-dimensional setting. All proofs that can be found in Bella\"iche's work are omitted.

\begin{definition}\label{GMA definition}
A \emph{generalized matrix algebra (GMA)} over a commutative ring $A$ is given by a tuple of data $(A, B, C, m)$, where $B$ and $C$ are $A$-modules, $m \colon B \otimes_A C \to A$ is a morphism of $A$-modules satisfying
\[
m(b\otimes c)b' = m(b'\otimes c)b \text{ and } m(b\otimes c')c = m(b\otimes c)c' \text{ for all } b,b' \in B, c, c' \in C.
\]
Given such data, define the $A$-module $R \coloneqq A \oplus B \oplus C \oplus A = \left(\begin{smallmatrix} 
A & B\\
C & A
\end{smallmatrix}\right)$ and give $R$ a multiplicative structure via
\[
\left(\begin{smallmatrix} 
a & b\\
c & d
\end{smallmatrix}\right)\left(\begin{smallmatrix} 
a' & b'\\
c' & d'
\end{smallmatrix}\right) \coloneqq \left(\begin{smallmatrix} 
aa' + m(b\otimes c') & ab' + bd'\\
a'c + dc' & dd' + m(b'\otimes c) 
\end{smallmatrix}\right) \qquad \mbox{for $a, a', d, d' \in A, b, b' \in B, c, c' \in C$},
\] 
so that $R$ has the structure of an $A$-algebra via the ring homomorphism $a \mapsto \left(\begin{smallmatrix}a &  0 \\ 0 &a \end{smallmatrix}\right)  \in R$.  We refer to the GMA given by $(A, B, C, m)$ simply by $R$.  

A \emph{morphism of GMAs} $(A, B, C, m)\to (A', B', C', m')$ (with associated $A$-algebras $R$ and $R'$) is a triple $(f_A,f_B,f_C)$ consisting of a ring morphism $f_A\colon A\to A'$ and two $A'$-module morphisms $f_B\colon B\otimes_{A,f_A}A'\to B'$, $f_C\colon C\otimes_{A,f_A}A'\to C'$ such that $f_A\circ m=m'\circ(f_B \otimes f_C)$. The data $(f_A,f_B,f_C)$ defines in a natural way an $A$-algebra morphism $\psi\colon R\to R'$; we say that $\psi$ is associated with $(f_A,f_B,f_C)$.

If $A$ is a topological ring and $B, C$ are topological $A$-modules, then $R$ inherits a natural topology, and we call $R$ a \textit{topological GMA} if $m$ is continuous.  We say that $R$ is \textit{faithful} if $m$ is nondegenerate as a pairing of $A$-modules.  As with matrix algebras, we have the notion of a trace and determinant on a GMA $R$ given by $\tr \bigl(\begin{smallmatrix} a & b\\
c & d
\end{smallmatrix}\bigr) = a + d$ and $\det \bigl(\begin{smallmatrix} a & b\\
c & d
\end{smallmatrix}\bigr) = ad - m(b \otimes c)$.
\end{definition}

The following lemma shows that when $A$ is a domain, faithful GMAs can be embedded into a matrix algebra over the field of fractions of $A$.

\begin{lemma}\cite[Lemmas 2.2.2, 2.2.3]{Bellaiche18}\label{GMAs over domains}
Assume that $A$ is a domain with field of fractions $K$ and that $R = \left(\begin{smallmatrix} 
A & B\\
C & A
\end{smallmatrix}\right)$ is a faithful GMA over $A$.  Then there exist embeddings of $A$-modules $B, C \hookrightarrow K$ such that (identifying $B, C$ with their images in $K$), $m \colon B \otimes_A C \to A$ is induced by multiplication in $K$.  In particular, if $BC \neq 0$, then $R \otimes_A K$ is isomorphic over $K$ as a GMA to $M_2(K)$.
\end{lemma}

We recall the following result of Bella\"iche, which explains that any residually multiplicity-free pseudorepresentation can be realized as the trace of a GMA-valued representation.

\begin{proposition}\label{Bellaiche 2.4.2}\cite[Proposition 2.4.2]{Bellaiche18}
Let $\ovl{\rho} \colon \Pi \to \GL_2(\F)$ be multiplicity free over $\F$.  Let $(t,d) \colon \Pi \to A$ be a pseudodeformation of $\ovl{\rho}$.  
\begin{enumerate}
\item There exists a faithful GMA $R$ over $A$ and a morphism of groups $\rho \colon \Pi \to R^\times$ such that $\tr \rho = t, \det \rho = d$, and $A\rho(\Pi) = R$.
\item If $(\rho, R)$ and $(\rho', R')$ are as in (1), then there is a unique isomorphism of $A$-algebras \mbox{$\Psi \colon R \to R'$} such that $\Psi \circ \rho = \rho'$.
\item\label{pinning} If $g_0 \in \Pi$ such that $\ovl{\rho}(g_0)$ has distinct eigenvalues $\lambda_0, \mu_0 \in \F^\times$, then there exists $(\rho, R)$ as in (1) such that $\rho(g_0)$ is diagonal and $\rho(g_0) \equiv \left(\begin{smallmatrix} 
\lambda_0 & 0\\
0 & \mu_0
\end{smallmatrix}\right) \bmod \m$. 
\item If $g_0 \in \Pi$ and $(\rho, R), (\rho', R')$ are as in (3), then the unique isomorphism of $A$-algebras \mbox{$\Psi \colon R \to R'$} such that $\Psi \circ \rho = \rho'$ is associated with an isomorphism of GMAs.
\item If $\ovl{\rho}$ is irreducible and $(\rho,R)$ is as in (1), then $R = (A, B, C, m, R)$ is isomorphic to $M_2(A)$ as a GMA over $A$.  If $\ovl{\rho}$ is reducible, then $BC \subseteq \m$.
\item If $(\rho, R)$ are as in (1), then $\ker \rho = \ker (t, d)$.
\item Assume that $A$ is noetherian and $\Pi$ is $p$-finite.  If $(t,d)$ is continuous, then for $(\rho, R)$ as in (1), $R$ is of finite type as an $A$-module.  If $R$ is given its unique topology as an $A$-algebra, then $\rho$ is continuous.
\end{enumerate}
\end{proposition}

\begin{remark}\label{B and C in A}
When $\overline{\rho}$ is absolutely irreducible, \cref{Bellaiche 2.4.2} allows us to identify the GMA $R$ with the matrix algebra $M_2(A)$.  We follow Bella\"iche in always implicitly making such an identification.  In particular, in the dihedral case, elements of $B$ and $C$ are viewed as elements of~$A$.
\end{remark}

Following Bella\"iche, we make the following definitions.

\begin{definition}\cite[Definition 2.4.3]{Bellaiche18}\label{(t,d)-rep}
A representation $\rho \colon \Pi \to R^\times$ satisfying condition (1) in \cref{Bellaiche 2.4.2} is called a \textit{$(t,d)$-representation}.  If in addition $\rho$ satisfies condition (3), then we say that $\rho$ is \textit{adapted} to $(g_0, \lambda_0, \mu_0)$.
\end{definition}

In fact, it is often useful to have the following strengthening of \cref{Bellaiche 2.4.2}(3).

\begin{lemma}\label{better adapted td-rep}
Let $\ovl{\rho} \colon \Pi \to \GL_2(\F)$ be multiplicity free over $\F$ and $\lambda_0 \neq \mu_0 \in \F^\times$ be the eigenvalues of an element in $\im \ovl{\rho}$.  Let $(t,d) \colon \Pi \to A$ be a pseudodeformation of $\ovl{\rho}$.  Then there exists $g_0 \in \Pi$ and a $(t,d)$-representation $\rho$ adapted to $(g_0, \lambda_0, \mu_0)$ such that $\rho(g_0) = \left(\begin{smallmatrix} 
s(\lambda_0) & 0\\
0 & s(\mu_0)
\end{smallmatrix}\right)$.
\end{lemma}

\begin{proof}
Let $g_0' \in \Pi$ be any element such that $\ovl{\rho}(g_0')$ has eigenvalues $\lambda_0, \mu_0$.  Then \cref{Bellaiche 2.4.2}(3) guarantees the existence of a $(t,d)$-representation $\rho \colon \Pi \to R^\times$ adapted to $(g_0', \lambda_0, \mu_0)$.  By \cite[Theorem 6.2.1]{Bellaiche18}, it follows that $\left(\begin{smallmatrix} 
s(\lambda_0) & 0\\
0 & s(\mu_0)
\end{smallmatrix}\right) \in \im \rho$. Let $g_0$ be any element in $\rho^{-1}\left(\begin{smallmatrix} 
s(\lambda_0) & 0\\
0 & s(\mu_0)
\end{smallmatrix}\right)$.  Then $\rho$ is a $(t,d)$-representation adapted to $(g_0, \lambda_0, \mu_0)$ and $\rho(g_0) = \left(\begin{smallmatrix} 
s(\lambda_0) & 0\\
0 & s(\mu_0)
\end{smallmatrix}\right)$.
\end{proof}

\subsection{Pink-Lie algebras}\label{Pink thy}
In \cref{Pink thy} we recall Pink's theory relating pro-$p$ subgroups of $\SL_2(A)$ to closed Lie subalgebras of $\Sl_2(A)$ \cite{Pink93}.  In fact, we use Bella\"iche's generalization to GMAs \cite[Section 4]{Bellaiche18}.

Recall that $A$ is a local pro-$p$ ring with $p \neq 2$.  The assumption that $p \neq 2$ is critical for Pink's theory.  We denote by $\m$ the maximal ideal of $A$.  Fix a compact topological GMA $R = \left(\begin{smallmatrix} 
A & B\\
C & A
\end{smallmatrix}\right)$ over~$A$.  (The compactness condition is satisfied, for instance, when $R$ is finite as an $A$-module.)  Write 
\[
SR^\times \coloneqq \{r \in R^\times \colon \det r = 1\}.
\]
Let $\rad R$ be the Jacobson radical of $R$, and $R^1 \coloneqq 1 + \rad R$.  We let $SR^1 \coloneqq SR^\times \cap R^1$, which is a closed normal pro-$p$ subgroup of $R^\times$.  See \cite[Remark 4.2.1]{Bellaiche18} for an explicit description of these objects. We mention here that in the case when $BC=A$ there is by \cite[Lemma 2.2.1]{Bellaiche18} an isomorphism of GMAs $R\cong M_2(A)$ that we can use to identify $\rad R$ with $\m M_2(A)$ and $R/\rad R$ with $M_2(\F)$, while if $BC\subset\m$ then $\rad R=\left(\begin{smallmatrix} 
\m & B\\
C & \m
\end{smallmatrix}\right)$
and $R/\rad R=\left(\begin{smallmatrix} 
\F & 0\\
0 & \F
\end{smallmatrix}\right)$.

Given any subset $S$ of $R$, we write
\[
S^0 \coloneqq \{s\in S\colon\tr s=0\}.
\]
Then $(\rad R)^0$ has a Lie algebra structure with bracket given by $[r_1, r_2] \coloneqq r_1r_2 - r_2r_1$.

For any topological group $G$ and closed subgroup $H$ of $G$, write $(G, H)$ for the smallest closed subgroup of $G$ containing $\{g^{-1}h^{-1}gh \colon g \in G, h \in H\}$.  Fix a closed subgroup $\Gamma \subseteq SR^1$.  Recall that the lower central series of $\Gamma$ is defined by $\Gamma_1 \coloneqq \Gamma$ and define $\Gamma_{n+1} \coloneqq (\Gamma, \Gamma_n)$.  We describe how Pink associates a filtration of Lie algebras to $\Gamma$ \cite[Section 2]{Pink93}.

Define a function
\begin{align*}
\Theta \colon R^\times &\to R^0\\
r &\mapsto r - \frac{\tr r}{2},
\end{align*}
where $(\tr r)/2$ is regarded as a scalar via the structure morphism $A \to R$.  Let $L(\Gamma) = L_1(\Gamma)$ be the (additive) subgroup of $(\rad R)^0$ topologically generated by $\Theta(\Gamma)$.  For $n \geq 2$, define $L_n(\Gamma)$ recursively as the subgroup of $(\rad R)^0$ topologically generated by the set 
\[
\{xy - yx \colon x \in L_1(\Gamma), y \in L_{n-1}(\Gamma)\}.
\]
Although the $L_n(\Gamma)$ are a priori only subgroups of $(\rad R)^0$, Pink shows that they are closed under Lie brackets and form a descending filtration, as summarized in the following proposition, which is due to Pink when $R = M_2(A)$ \cite[Proposition 3.1, Proposition 2.3]{Pink93} and to Bella\"iche in the GMA case \cite[Proposition 4.7.1]{Bellaiche18}.

\begin{proposition}
For all $n \geq 1$, we have $L_{n+1}(\Gamma) \subseteq L_n(\Gamma)$.  In particular, each $L_n(\Gamma)$ is a Lie subalgebra of $(\rad R)^0$.
\end{proposition}

We emphasize that, a priori, each $L_n(\Gamma)$ is just a $\Z_p$-module, even if the ring $A$ is very large.  The point of \cref{rings acting} is to prove that, under mild conditions, $L_n(\Gamma)$ is in fact an algebra over an (in general) much larger ring.
 
Conversely, given a closed Lie subalgebra $L$ of $(\rad R)^0$, define $H(L) \coloneqq \Theta^{-1}(L) \cap SR^1$.  Let $H_n \coloneqq H(L_n(\Gamma))$.  A priori, $H(L)$ is only a subset of $SR^1$.  However, we have the following theorem of Pink \cite[Proposition 2.4, Theorem 2.7]{Pink93}, which was generalized to GMAs by Bella\"iche \cite[Theorem 4.7.3]{Bellaiche18}.  

\begin{theorem}\label{Pink main thm}
We have that $H_n$ is a pro-$p$ subgroup of $SR^1$.  Furthermore, $\Gamma$ is a normal subgroup of $H_1$, and $H_1/\Gamma$ is abelian.  For $n \geq 2$, we have $H_n = \Gamma_n$.
\end{theorem}

\begin{remark} Pink's construction satisfies the following two important properties.
\begin{enumerate}
\item It is functorial with respect to surjective ring homomorphisms.  Namely, let $\Aa$ be a closed ideal of $A$ and $\varphi \colon R \to R/\Aa R$ the natural projection.  Then for all $n \geq 1$ we have
\[
\varphi(L_n(\Gamma)) = L_n(\varphi(\Gamma)).
\]
\item Pink's Lie algebra $L_n(\Gamma)$ is closed under conjugation by the normalizer of $\Gamma$ in $R^\times$.  This follows easily from the definitions since $\Theta$ is invariant under conjugation.\qedhere
\end{enumerate}
\end{remark} 

See \cref{congrL1} for an example calculating $L_n(\Gamma)$ when $\Gamma$ is a congruence subgroup.

\subsection{Decomposability and regularity}\label{regularity background}
In order to prove fullness theorems, it is useful to be able to decompose Pink's Lie algebra according to its entries.  In \cref{regularity background} we define this precisely and then define regularity, which will turn out to ensure that the Lie algebras of the representations we work with are decomposable.

\begin{definition}\cite[Section 4.9]{Bellaiche18}\label{decomposable}
Let $R$ be a GMA over $A$ and $L$ a closed subspace of $(\rad R)^0$.  We say that $L$ is \textit{decomposable} if 
\[
\left(\begin{smallmatrix} 
a & b\\
c & -a
\end{smallmatrix}\right) \in L\textrm{ implies that }\left(\begin{smallmatrix} 
a & 0\\
0 & -a
\end{smallmatrix}\right) \in L\textrm{ and }\left(\begin{smallmatrix} 
0 & b\\
c & 0
\end{smallmatrix}\right) \in L.
\] 
We say that $L$ is \textit{strongly decomposable} if $L$ is decomposable and 
\[
\left(\begin{smallmatrix} 
a & b\\
c & -a
\end{smallmatrix}\right) \in L\textrm{ implies that }\left(\begin{smallmatrix} 
0 & b\\
0 & 0
\end{smallmatrix}\right) \in L\textrm{ and }\left(\begin{smallmatrix} 
0 & 0\\
c & 0
\end{smallmatrix}\right) \in L.
\]
\end{definition}

If $L_n(\Gamma) \subseteq R = \left(\begin{smallmatrix} 
A & B\\
C & A
\end{smallmatrix}\right)$ is decomposable, we write \label{B1C1}
\begin{align*}
I_n(\Gamma) &\coloneqq \{a \in A \colon \left(\begin{smallmatrix} 
a & 0\\
0 & -a
\end{smallmatrix}\right) \in L_n(\Gamma)\},\\
\nabla_n(\Gamma) &\coloneqq \{\left(\begin{smallmatrix} 
0 & b\\
c & 0
\end{smallmatrix}\right) \in L_n(\Gamma)\},\\ 
B_n(\Gamma) &\coloneqq \{b \in B \colon \exists c \in C \text{ such that } \left(\begin{smallmatrix} 
0 & b\\
c & 0
\end{smallmatrix}\right) \in L_n(\Gamma)\},\\
C_n(\Gamma) &\coloneqq \{c \in C \colon \exists b \in B \text{ such that } \left(\begin{smallmatrix} 
0 & b\\
c & 0
\end{smallmatrix}\right) \in L_n(\Gamma)\}.
\end{align*}

Eventually, $L$ will be a Pink-Lie algebra associated to some admissible pseudodeformation of $\ovl{\rho} \colon \Pi \to \GL_2(\F)$.  Regularity is a condition on $\ovl{\rho}$ that will allow us to decompose $L$, as we will see in \cref{rings acting}.  

Let $\E$ be the subfield of $\F$ generated by $\{(\tr \ovl{\rho}(g))^2/\det \ovl{\rho}(g) \colon g \in \Pi\}$; equivalently, $\E$ is generated by the traces of $\ad \ovl{\rho}$.  
If $\lambda_g, \mu_g$ are the eigenvalues of $\ovl{\rho}(g)$, then we see that\\ \mbox{$(\tr \ovl{\rho}(g))^2/\det \ovl{\rho}(g) = \lambda_g\mu_g^{-1} + \lambda_g^{-1}\mu_g + 2$}.  Hence $\E$ is generated over $\F_p$ by the set
\begin{equation}\label{generating set for E}
\{\lambda\mu^{-1} + \lambda^{-1}\mu \colon \lambda, \mu \text{ are the eigenvalues of }\ovl{\rho}(g)  
\text{ for some } g \in \Pi\}.
\end{equation}
In particular, $g$ will not contribute to $\E$ if the multiplicative order of $\lambda_g\mu_g^{-1}$ is strictly less than~5.  Using this reasoning, it is straightforward to calculate $\E$ when $\ovl{\rho}$ exceptional.  Namely, if $\ovl{\rho}$ is tetrahedral or octahedral, then $\E = \F_p$.  If $\ovl{\rho}$ is icosahedral, then $\E = \F_p$ if $p = 5$ and $\E = \F_p(\zeta_5 + \zeta_5^{-1}) = \F_p(\sqrt{5})$ otherwise. 

\begin{definition}\label{regular}
Let $\ovl{\rho} \colon \Pi \to \GL_2(\F)$ be a semisimple representation.  We say that $\ovl{\rho}$ is \textit{regular} if there exists $g_0 \in \Pi$ such that $\ovl{\rho}(g_0)$ has eigenvalues $\lambda_0, \mu_0 \in \ovl{\F}_p^\times$ satisfying $\lambda_0\mu_0^{-1} \in \E^\times \setminus \{\pm 1\}$.  We call $g_0$ a \textit{regular element} for $\ovl{\rho}$.  If in addition $\lambda_0, \mu_0 \in \E^\times$, then we say that $\ovl{\rho}$ is \textit{strongly regular}.
\end{definition}

\cref{regular} is weaker than Bella\"iche's definition of regularity \cite[Definition 7.2.1]{Bellaiche18}, where the eigenvalues $\lambda_0$ and $\mu_0$ are required in addition to belong to $\F_p$.  Examining \eqref{generating set for E}, we see that the only way $\ovl{\rho}$ can fail to be regular is if, for every matrix in $\im \ovl{\rho}$ with eigenvalues $\lambda, \mu$, either $\lambda\mu^{-1} = \pm 1$ or the unique quadratic extension of $\E$ is $\E(\lambda\mu^{-1})$.

\begin{remark}\label{regularity remark background} 
Let us analyze regularity depending on the projective image of $\ovl{\rho}$.  With notation as in \cref{regular}, note that the order of $\lambda_0\mu_0^{-1}$ in $\E^\times$ corresponds to the order of $\ovl{\rho}(g_0)$ in the projective image of $\ovl{\rho}$.
\begin{enumerate}
\item If $\ovl{\rho}$ is large, then $\ovl{\rho}$ is regular.  Indeed, $\proj \ovl{\rho}(\Pi)$ contains $\PSL_2(\E)$ up to conjugation.  Since $\ovl{\rho}$ is not exceptional, $\E^\times$ contains an element $x$ such that $x^2 \neq \pm 1$.  Then the image of $\ovl{\rho}$ contains, up to conjugation, a scalar multiple of $\bigl(\begin{smallmatrix} 
x & 0\\
0 & x^{-1}
\end{smallmatrix}\bigr)$, which satisfies the regularity property.
\item If $\ovl{\rho}$ is tetrahedral and $p > 3$, then a regular element must map to a $3$-cycle in the projective image of $\ovl{\rho}$, since the other elements of $A_4$ have order at most 2.  Thus in this case regularity is equivalent to $\zeta_3 \in \E = \F_p$, which is equivalent to $p \not\equiv 2 \bmod 3$.  By a similar argument we see that if $\ovl{\rho}$ is octahedral and $p > 3$, then regularity is equivalent to one of $\zeta_3$ or $\zeta_4$ being in $\E = \F_p$, which is equivalent to $p \not \equiv 11 \bmod 12$.  If $\ovl{\rho}$ is icosahedral and $p \neq 5$, then regularity is equivalent to one of $\zeta_3$ or $\zeta_5$ being in $\E = \F_p(\sqrt{5})$, which is equivalent to $p \not\equiv 14 \bmod 15$.  
\item If $\ovl{\rho} \cong \varepsilon \oplus \delta$, then $\ovl{\rho}$ is regular if and only if $\varepsilon\delta^{-1}$ takes values in~$\E^\times$ (\cref{regular reducible reps}).
\item If $\ovl{\rho} = \Ind_{\Pi_0}^\Pi \chi$ is dihedral, then elements in $\Pi \setminus \Pi_0$ have projective order 2, and so any regular element must lie in $\Pi_0$.  Furthermore, elements in $\Pi \setminus \Pi_0$ have trace 0, and so the field $\E$ associated to $\ovl{\rho}$ is the same as the field $\E$ associated to $\ovl{\rho}|_{\Pi_0}$.  Hence we are reduced to the previous case when $\ovl{\rho}$ is reducible.
\item If the projective image of $\ovl{\rho}$ is isomorphic to $\Z/2\Z$ or $(\Z/2\Z)^2$, or if $\E = \F_3$, then $\ovl{\rho}$ is never regular.  In particular, if $p = 3$ and $\ovl{\rho}$ is tetrahedral or octahedral, then $\ovl{\rho}$ is not regular.  If $p = 5$ and $\ovl{\rho}$ is icosahedral, then $\proj \ovl{\rho}(\Pi)$ is conjugate to $\PSL_2(\F_5)$.  Thus $\E = \F_5$ and any potential regular element has eigenvalues in $(\F_5^\times)^2 = \{\pm1\}$, so $\ovl{\rho}$ is not regular in this case.\qedhere
\end{enumerate}
\end{remark}

\subsection{Bella\"iche's results}\label{bellaiche results background}
The purpose of this section is to state Bella\"iche's main results that form the basis for our work in this paper. 
We state them in slightly less generality than \cite[Section 6]{Bellaiche18}.  As before, $A$ denotes a local pro-$p$ ring with maximal ideal $\m$ and residue field $\F$.  In particular, $A$ is naturally a topological $W(\F)$-algebra.   

Let $R$ be a faithful GMA over $A$. Recall the description of $R/\rad R$ that we gave in the beginning of \cref{Pink thy}. We define $s \colon R/\rad R \to R$ by 
\[
s\begin{pmatrix} 
a & b\\
c & d
\end{pmatrix} \coloneqq \begin{cases} \left(\begin{smallmatrix} 
s(a) & s(b)\\
s(c) & s(d)
\end{smallmatrix}\right) & \textrm{if } R = M_2(A)\vspace{0.1cm}\\
\vspace{0.1cm}\left(\begin{smallmatrix} 
s(a) & 0\\
0 & s(d)
\end{smallmatrix}\right) & \text{ else.}
\end{cases}
\]
Note that in the latter case, we have a priori that $b = c = 0$.

Let us fix an admissible pseudodeformation $(\Pi, \ovl{\rho}, t, d)$ over $A$.  If $p = 3$, let us assume that $\ovl{\rho}$ is not tetrahedral.  By \cref{Bellaiche 2.4.2}, there exists a $(t,d)$-representation $\rho \colon \Pi \to R^\times$.  Given such a $(t,d)$-representation, write $G = G_\rho \coloneqq \rho(\Pi)$ and $\Gamma = \Gamma_\rho \coloneqq G \cap SR^1$.  Furthermore, let $\overline{G}$ denote the image of $G$ modulo $\rad R$.  (Note that the image of $\overline{G}$ under an embedding $R/\rad R\to\GL_2(\F)$ is a conjugate of $\ovl{\rho}(\Pi)$.)  We will write $L_n(\rho) \coloneqq L_n(\Gamma_\rho)$ and analogously for $I_n(\rho), \nabla_n(\rho), B_n(\rho), C_n(\rho)$.

Bella\"iche chooses his $(t,d)$-representations very carefully in order to give a nice description of their Pink-Lie algebras.  How this is done depends upon the projective image of $\ovl{\rho}$.  Since $\ovl{\rho}$ is multiplicity free over $\F$, we can let $\lambda_0 \neq \mu_0 \in \ovl{\F}_p^\times$ be the eigenvalues of a matrix $x_0 \in \im \ovl{\rho}$ chosen such that the following conditions are satisfied:
\begin{itemize}
\item if $\ovl{\rho}$ is large, then $(\lambda_0\mu_0^{-1})^2 \neq 1$ and $\lambda_0, \mu_0 \in \F_p^\times$;
\item if $p = 3$ and $\ovl{\rho}$ is octahedral, then $\lambda_0\mu_0^{-1}$ is a primitive fourth root of unity;
\item if $p = 5$ and $\ovl{\rho}$ is icosahedral, then $\lambda_0\mu_0^{-1}$ is a primitive third root of unity;
\item if $\ovl{\rho}$ is exceptional and does not belong to one of the previous to scenarios, then $\lambda_0\mu_0^{-1}$ is a primitive third, fourth, or fifth root of unity;
\item otherwise, the multiplicative order of $\lambda_0\mu_0^{-1}$ is equal to the maximal order of an element in the projective image of $\ovl{\rho}$.
\end{itemize}

\begin{definition}\label{well adapted}
Suppose $(\Pi, \ovl{\rho}, t, d)$ is an admissible pseudodeformation.  We say that a $(t,d)$-representation $\rho$ is \textit{well adapted} if 
\begin{enumerate}
\item $\rho$ is adapted to an element $g_0$ such that $\rho(g_0) = \left(\begin{smallmatrix} 
s(\lambda_0) & 0\\
0 & s(\mu_0)
\end{smallmatrix}\right)$, where $\lambda_0, \mu_0$ satisfy the relevant property listed above;
\item if the projective image of $\ovl{\rho}$ is dihedral and nonabelian, then $\ovl{G}$ contains a matrix of the form $\left(\begin{smallmatrix} 
0 & b\\
c & 0
\end{smallmatrix}\right)$ with $bc^{-1} \in \F_p^\times$ and $s(\ovl{G}) \subseteq \im \rho$.
\end{enumerate}
\end{definition}

Bella\"iche shows that well-adapted $(t,d)$-representations always exist, provided that one is willing to replace $\F$ by a quadratic extension in the dihedral case \cite[Proposition 6.3.2, Lemma 6.8.2]{Bellaiche18}.  

Define $\F_q$ as in the table below.  We will see in \cref{regular reducible reps} that if $\ovl{\rho}$ is regular and reducible or dihedral, then $\F_q$ can be taken to be $\E$.  If $\ovl{\rho}$ is not projectively cyclic or dihedral, then \mbox{$\F_q \subseteq \E$ by definition}.  (In the $A_5$ case, this follows from the calculation that $\E = \F_p(\sqrt{5})$ prior to \cref{regular}.)  

\begin{remark}
Our definition of $\F_q$ differs from that of Bella\"iche when $\ovl{\rho}$ is exceptional.  If $\ovl{\rho}$ is tetrahedral, then Bella\"iche defines $\F_q = \F_p(\zeta_3)$.  If $\ovl{\rho}$ is octahedral, he defines $\F_q$ to be $\F_p(\zeta_3)$ if the ratio $\lambda_0\mu_0^{-1}$ chosen prior to \cref{well adapted} has order 3 and $\F_p(\zeta_4)$ if that ratio has order 4.  If $\ovl{\rho}$ is icosahedral, then he defines $\F_q = \F_p(\zeta_5)$.  The key property that Bella\"iche needs is that $\ovl{\rho}$ can be conjugated so that its image lies in $Z\GL_2(\F_q)$ and $\bigl(\begin{smallmatrix} 
\lambda_0 & 0\\ 0 & \mu_0
\end{smallmatrix}\bigr) \in \im \ovl{\rho}$, where $Z$ is the group of scalar matrices in $\F$ (cf. \cite[Lemma 6.8.5]{Bellaiche18}).  This change of definition will be justified in \cref{Bellaiche Fq vs E}.
\end{remark}
\begin{table}[h!]
   \begin{center}
      \label{definition of Fq}
      \begin{tabular}{| l | l |}
         \hline
         the projective image of $\ovl{\rho}$ is & $\F_q$\\
     	 \hline
         cyclic of order $m$ or dihedral of order $2m$ & any subfield of $\F$ such that $(m, q - 1) > 2$\\
         exceptional & $\E(\lambda_0\mu_0^{-1})$\\
         otherwise & $\F_p$\\
				 \hline
      \end{tabular}
   \end{center}   
\end{table}

The following theorem summarizes Bella\"iche's results describing the structure of $W(\F_q)L_1(\rho)$ from \cite[Section 6]{Bellaiche18}.  We recall from \cref{B and C in A} that in the dihedral case, elements in $B$ and~$C$ can be viewed as elements of $A$.

\begin{theorem}[Bella\"iche]\label{Bellaiche main theorem}
Let $(\Pi, \ovl{\rho}, t, d)$ be an admissible pseudodeformation such that the projective image of $\ovl{\rho}$ is not isomorphic to $\Z/2\Z$ nor $(\Z/2\Z)^2$.  Then every well-adapted $(t,d)$-representation $\rho \colon \Pi \to R^\times$ with $R = \left(\begin{smallmatrix}
A & B\\
C & A
\end{smallmatrix} \right)$ has the following properties:
   \begin{enumerate} 
      \item $L_1(\rho)$ is decomposable;
      \item the ring $A$ is equal to
			\[
			\begin{cases}
      W(\F) + W(\F)I_1(\rho) + W(\F)I_1(\rho)^2 + W(\F)B_1(\rho) & \text{if } \ovl{\rho} \text{ is projectively dihedral}\\
      W(\F) + W(\F)I_1(\rho) + W(\F)I_1(\rho)^2 & \text{otherwise};
      \end{cases}
			\]
      \item $W(\F)C_1(\rho) = C$ and $W(\F)B_1(\rho) = B$;
      \item up to possibly replacing $\rho$ with its conjugate by a certain matrix $\left(\begin{smallmatrix} 
      1 & 0\\
      0 & a
      \end{smallmatrix}\right)$ with $a \in A^\times$ when $\ovl{\rho}$ is exceptional or large, $W(\F_q)L_1(\rho)$ is equal to  
         \[
         \begin{pmatrix} 
      W(\F_q)I_1(\rho) & W(\F_q)B_1(\rho)\\
      W(\F_q)C_1(\rho) & W(\F_q)I_1(\rho)
      \end{pmatrix}^0.
         \]
   \end{enumerate}
Furthermore
\begin{enumerate}[label=(\roman*)]
\item $(W(\F_q)I_1(\rho))^3 \subseteq W(\F_q)I_1(\rho)$;
\item if $\ovl{\rho}$ is not reducible, then $W(\F_q)C_1(\rho) = W(\F_q)B_1(\rho)$;
\item if $\ovl{\rho}$ is exceptional or large, then $W(\F_q)B_1(\rho) = W(\F_q)I_1(\rho)$ and $(W(\F_q)I_1(\rho))^2 \subset W(\F_q)I_1(\rho)$.
\end{enumerate}
\end{theorem}

For a subfield $\F'$ of $\F$, we shall often refer to the $W(\F')$-subalgebra of $A$ generated by $I_1(\rho)$.  We denote it by $\B_\rho(\F')$, which is simply equal to $W(\F') + W(\F')I_1(\rho) + W(\F')I_1(\rho)^2$ whenever $\F_q \subseteq \F'$.  When $\ovl{\rho}$ is not reducible or dihedral, we have $\B_\rho(\F') = W(\F') + W(\F')I_1(\rho)$ by \cref{Bellaiche main theorem}(ii, iii).  In particular, \cref{Bellaiche main theorem} says that $A = \B_\rho(\F)$ unless $\ovl \rho$ is projectively dihedral, in which case $A = \B_\rho(\F) + W(\F)B_1(\rho)$.

Bella\"iche uses \cref{Bellaiche main theorem} to deduce that, under certain hypotheses, the representation $\rho$ is $\B_\rho(\F_p)$-full.  See \cref{Bellaiche thm intro} or \cite[Theorem 7.2.3]{Bellaiche18} for a precise statement of his result.  The first step in our main theorem is to improve this to $\B_\rho(\E)$-fullness in \cref{rings acting}.   
But first we discuss fullness, conjugate self-twists, and the connections between them in detail.

\medskip


\section{Fullness}\label{section: fullness}

\label{fullness background} 

In this section we explore the notion of \emph{fullness}, our measure of the size of the image of a continuous (pseudo)representation on a noetherian local pro-$p$ ring $A$. Here we additionally assume that $A$ is a domain, with field of fractions $K$. 
\subsection{Fullness for (pseudo)representations}
Let $B$ be {any ring}.
For any nonzero $B$-ideal $\Bb$, let
\[
\Gamma_{B}(\Bb) \coloneqq \ker(\SL_2(B) \to \SL_2(B/\Bb)) = \left\{\left(\begin{smallmatrix}
1 + a & b\\
c & 1 + d
\end{smallmatrix}\right) \in \SL_2(B) \colon a, b, c, d \in \Bb\right\}
\]
be the \emph{congruence subgroup} of $\SL_2(B)$ of \emph{level} $\Bb$.

\begin{definition}\label{fullness}
Let $G$ be a subgroup of $\GL_2(K)$.  {For a subring $B$ of $K$} we say that $G$ is $B$-\textit{full} if there exists a nonzero $B$-ideal $\Bb$ and $x \in \GL_2(K)$ such that 
\[
x^{-1}Gx \supseteq \Gamma_{B}(\Bb). 
\]
A $\GL_2(K)$-valued representation is \textit{$B$-full} if its image is $B$-full.  If $(t,d) \colon \Pi \to A$ is a pseudorepresentation, we say $(t,d)$ is \textit{$B$-full} if there exists a $(t,d)$-representation $\rho \colon \Pi \to R^\times$ such that $\iota \circ \rho$ is $B$-full, where $\iota$ is an embedding of $R$ into $M_2(K)$. Such an $\iota$
exists by \cref{GMAs over domains}; {note that by replacing $\iota$ by a conjugate embedding we may insist that $\Gamma_B(\Bb) \subset \GL_2(K)$ is contained in $\iota(R^\times)$ on the nose.}
We will say $B$ is a $(t,d)$-\textit{fullness ring} if $(t,d)$ is $B$-full.
\end{definition}

{The notion of fullness, which has appeared in earlier incarnations in \cite[last introductory paragraph]{Hida15} and \cite[Definition 2.2]{Lang2016}, is analogous to Bella\"iche's notion of  ``congruence large-image" \cite[Definition 7.2.1]{Bellaiche18}.}
We now show that 
fullness is well defined for pseudorepresentations and gives compatible notions for representations and pseudorepresentations.

\begin{lemma}\label{fullness well defined}
Let $(t,d) \colon \Pi \to A$ be a pseudorepresentation and $\rho \colon \Pi \to \GL_2(K)$ a representation whose trace takes values in $A$. {Let $B$ be any subring of $K$.}
\begin{enumerate}
\item \label{fullness for pseudoreps well defined} If there exists a $(t,d)$-representation that is $B$-full, then every $(t,d)$-representation is $B$-full.  
\item \label{fullness for reps well defined} The representation $\rho$ is $B$-full if and only if its pseudorepresentation $(\tr \rho, \det \rho)$ is $B$-full.
\end{enumerate}
\end{lemma}

\begin{proof}
To prove \eqref{fullness for pseudoreps well defined}, let $\rho \colon \Pi \to R^\times$ and $\rho' \colon \Pi \to R'^\times$ be two $(t,d)$-representations.  We just have to verify that the $A$-algebra isomorphism $\Psi \colon R \to R'$ such that $\rho' = \Psi \circ \rho$ from \cref{Bellaiche 2.4.2} is given by conjugation by an element of $\GL_2(K)$.  Consider $\Psi \otimes 1 \colon R \otimes_A K \to R' \otimes_A K$ and recall that $R \otimes_A K \cong M_2(K) \cong R' \otimes_A K$ by  \cref{GMAs over domains}.  By the Skolem-Noether theorem, it follows that $\Psi \otimes 1$ (and hence $\Psi$) is conjugation by an element of $\GL_2(K)$.

For \eqref{fullness for reps well defined}, let $(t,d) = (\tr \rho, \det \rho)$, which is a pseudorepresentation over $A$ by assumption.  Let $r \colon \Pi \to R^\times$ be a $(t,d)$-representation, and embed $R$ into $\GL_2(K)$ by \cref{GMAs over domains}, thus viewing~$r$ as valued in $\GL_2(K)$.  Note that fullness of $\rho$ (respectively, $r$) implies that $\rho$ (respectively, $r$) is absolutely irreducible since this is true for the inclusion representation of a congruence subgroup in $\GL_2(K)$.  As $(\tr \rho, \det \rho) = (t,d) = (\tr r, \det r)$ and either $\rho$ or $r$ is absolutely irreducible, by the Brauer-Nesbitt Theorem it follows that $\rho$ and $r$ are conjugate by a matrix in $\GL_2(K)$.  Since fullness is defined up to conjugation in $\GL_2(K)$, it follows that $\rho$ is $B$-full if and only if $r$, and hence $(t,d)$, is $B$-full.
\end{proof}

The next two propositions suggest that we can restrict our attention to fullness rings that are closed subrings of the trace algebra. These ideas are made precise \cref{fullclosedprop} and \cref{closedinAt} below. Fix a continuous pseudorepresentation $(t, d): \Pi \to A$.

\begin{proposition}\label{closednormal2} If $(t, d)$ is $B$-full for some subring $B$ of K, then the $\Z$-linear span of $t(\Pi)$ contains a nonzero $B$-ideal. In particular, the trace algebra $A_t$ of $t$ contains a nonzero $B$-ideal.
\end{proposition}

\begin{proof} 
Let $0 \neq \Bb$ be an $B$-ideal such that $\Gamma_{B}(\Bb)$ is contained in the image of some $(t, d)$-representation. We claim that pairwise products of elements of $\Bb$ are all in the set $\{t(g) - 2: g \in \Pi\}$, so that $\Bb^2$ is contained in the trace algebra $A_t$. Indeed, an element of $\Gamma_{B}(\Bb)$ is of the form $\bigl(\begin{smallmatrix} 
1+a & b\\
c & 1 + d
\end{smallmatrix}\bigr)$ with $a, b, c, d \in \Bb$ such that $a + d + ad - bc = 0$.  In particular, for any $b, c \in \Bb$, taking $d = bc$ and $a = 0$ shows that 
\[
\tr\bigl(\begin{smallmatrix} 
1 & b\\
c & 1 + bc
\end{smallmatrix}\bigr) = 2 + bc \in t(\Pi).
\]  
Since $2 = t(1) \in t(\Pi)$ and elements of the form $bc$ generate $\Bb^2$, we see that the $\Z$-span of $t(\Pi)$ contains $\Bb^2$, which is nonzero since it is the square of a nonzero ideal of a domain.   
\end{proof}

\begin{proposition}\label{fullERprop}\label{fullclosedprop}
Let $(t, d): \Pi \to A$ be a continuous pseudorepresentation. 
If $(t, d)$ is $B$-full for a subring $B$ of $A$, then $(t, d)$ is also full for the closure $\ovl B$ of $B$ in $A$. Conversely, suppose that $(t, d)$ is $\ovl B$-full. 
If the image of a $(t, d)$-representation contains a congruence subgroup of $\SL_2(\ovl B)$ whose level is the closure $\ovl \Bb$ of an ideal $\Bb$ of $B$, then $(t, d)$ is also $B$-full. 

\end{proposition}

\begin{proof} 
We show that for a closed subgroup $G$ of the unit group $R^\times$ of a faithful 
GMA $R$ over $A$ equipped with an embedding $\iota: R \into M_2(K)$, if $\iota(G)$ contains $\Gamma_B(\Bb)$ for some nonzero ideal $\Bb$ of $B$, then $\iota(G)$ also contains $\Gamma_{\ovl B}(\ovl \Bb)$ for the closure $\ovl \Bb$ of $\Bb$. Since $A$ is noetherian, and both $R$ and $M_2(A)$ are finite $A$-algebras, the preimage $S \coloneqq\iota^{-1}\big(M_2(A)\big)$ is a finite, hence closed, 
$A$-subalgebra of $R$. Moreover, the induced map $\left.\iota\right|_S: S \to M_2(A)$ is a homeomorphism onto its image by the compatibility of topologies on finite $A$-modules (see p.~\pageref{uniquetop}).
In particular, any closed subset of $R$ containing $\iota^{-1}\big(\Gamma_B(\Bb)\big)$ will also contain its closure $\iota^{-1}\big(\Gamma_{\ovl B}(\ovl \Bb)\big).$ The converse claim is clear since $\ovl \Bb$ is nonzero only if $\Bb$ is. 
\end{proof}

\subsection{Fullness peers}\label{fpsec}
A pseudorepresentation {may be} $B$-full for more than one choice of ring $B$, even if $B$ is a closed subring of $A$. For example, any pseudorepresentation that happens to be full for $\Z_{p^2} = W(\F_{p^2})$ is also full for the order $\Z_p + p \Z_{p^2} \subset \Z_{p^2}$. 
We say two subrings $B_1, B_2$ of $K$ are \textit{fullness peers} if every nonzero ideal of $B_1$ contains a nonzero ideal of $B_2$ and vice versa, in which case $B_1$-fullness is equivalent to $B_2$-fullness. One easily checks that fullness peerage is an equivalence relation on subrings of $K$. The next lemma gives a criterion for establishing when nested domains are fullness peers.

\begin{lemma}\label{fg plus same quotient field implies ideal containment}\label{fgqf}
Let $B_1 \subseteq B_2$ be domains.  The following conditions are equivalent:
\begin{enumerate}
\item \label{f1} $B_1$ contains a nonzero ideal of $B_2$;
\item \label{f2} there exists $y \in B_1 \setminus \{0\}$ such that $yB_2 \subseteq B_1$;
\item \label{f3} $B_1$ and $B_2$ are fullness peers.
\end{enumerate}
These equivalent conditions imply that $B_2$ and $B_1$ have the same field of fractions.  
If moreover $B_1$ is noetherian, then conditions (\ref{f1},\ref{f2},\ref{f3}) are equivalent to:
\begin{enumerate}[resume]
\item \label{f4} $Q(B_2) = Q(B_1)$ and $B_2$ is a finite $B_1$-algebra.
\end{enumerate}
\end{lemma}

\begin{proof}
For \eqref{f1} implies \eqref{f2}, take $y$ to be any nonzero element of the nonzero ideal of $B_2$ contained in $A_1$. If \eqref{f2} holds, then an arbitrary nonzero ideal $\Bb$ of $A_1$ contains $(yB_2)\Bb$, which is a nonzero ideal of $B_2$, implying  \eqref{f3}.  Clearly \eqref{f3} implies \eqref{f1}.

To see that $Q(B_2) = Q(B_1)$ under any of (\ref{f1},\ref{f2},\ref{f3}), note that any $x \in B_2$ can be written as $(yx)/y \in Q(B_1)$ with $y$ as in \eqref{f2}.

For the rest of the proof, assume that $B_1$ is noetherian. Suppose first that any of  (\ref{f1},\ref{f2},\ref{f3}) holds. If $J$ is a non-zero ideal of $B_2$ contained in $B_1$, then $J$ is a finitely generated $B_1$-module. By replacing $J$ with a smaller $B_2$-ideal, we can assume that $J$ is principal in $B_2$, that is, $J = bB_2$ for some $b\in B_1$. Now choose a finite set of generators $\{bx_1, \ldots, bx_n\}$ of $bB_2$ as an $B_1$-module, with $x_1, \ldots, x_n$ in $B_2$. Then, for every $y$ in $B_1$, $by$ is a linear combination $\sum_i a_ibx_i$ for some $a_i\in B_1$, which means that $y = \sum_i a_i x_i$, so the set $\{x_1, \ldots, x_n\}$ generates $B_2$ as an $B_1$-module.

Conversely, suppose that \eqref{f4} is satisfied.  Let $x_1, \ldots, x_n$ be generators for $B_2$ as an $B_1$-module.  Write $x_i = b_{i1}/b_{i2}$ with $b_{ij} \in B_1 \setminus \{0\}$.  Set $b = \prod_{i = 1}^n b_{i2} \in B_1 \setminus \{0\}$.  Then $bx_i \in B_1$ for all $i$, and it follows that $bB_2 \subseteq A_1$, proving \eqref{f2}.
\end{proof}

\begin{question}\label{transporter} 
Note that $\Z_p \lb X \rb$  and its 
non-noetherian subring $\Z_p + p \Z_p \lb X \rb$ and are fullness peers even though the extension is not finite. Could profiniteness 
substitute for finiteness in \eqref{f4} above? That is, if $B_2$ is a local noetherian pro-$p$ domain, and $B_1 \subset B_2$ is a closed subring with the same field of fractions, are $B_1$ and $B_2$ necessarily fullness peers? 
\end{question}

\begin{corollary}\label{closedinAt}
Let $(t, d): \Pi \to A$ be a pseudorepresentation. If $(t, d)$ is $B$-full for some subring $B$ of $K$, then $B \cap A$ is a fullness peer of $B$. 
\end{corollary}

\begin{proof}
 By \cref{closednormal2}, the trace algebra of $(t, d)$, and hence $A$, contains a nonzero ideal of $B$. Therefore so does $B \cap A$, and by \cref{fg plus same quotient field implies ideal containment}, $B \cap A \subseteq B$ is an extension of fullness peers. 
\end{proof} 

\cref{closedinAt} together with \cref{fullERprop} 
 allows us to restrict our attention to fullness for closed subrings of $A$, though we continue to point out features of the general case for completeness. 

\begin{corollary}\label{normalpuppies}
Let $B$ be a complete 
{local noetherian} domain. Then all the extensions of $B$ contained in the normalization $B^{\norm}$ are fullness peers of $B$.
\end{corollary}

\begin{proof}
Since $B^{\norm}$ is the integral closure of $B$ in its field of fractions, $B^{\norm}$ is finite over $B$ by the $N2$ property~(see p.~\pageref{n2}), hence a noetherian $B$-module. Therefore any ring $C$ with {$B \subseteq C \subseteq B^\norm$} 
is finite over $B$. Fullness peerage then follows from \cref{fg plus same quotient field implies ideal containment}. 
\end{proof}

\subsection{Fullness and twisting}

A key property of fullness, shown in \cref{twisting does not affect fullness}, is that it does not change when we twist a pseudorepresentation by a character.  This will allow us to do all of our technical work in the setting of constant determinant pseudorepresentations where we have Bella\"iche's \cref{Bellaiche main theorem} available. The proof of the twist invariance of fullness relies on a calculation of the Pink-Lie algebras of a congruence subgroup.  
Let $B$ be a local pro-$p$ domain. For a {closed} nonzero $B$-ideal $\Bb$, define
\[
\Sl_2(\Bb) \coloneqq \left\{\left(\begin{smallmatrix} 
a & b\\
c & -a
\end{smallmatrix}\right) \colon a, b, c \in \Bb\right\} {\subset M_2(\Bb)}.
\]

\begin{lemma}[cf. {\cite[Proposition 4.8.2]{Bellaiche18}}]\label{congrL1}
Let $\Bb$ be a closed ideal of $B$.  Then $\Gamma_B(\Bb)$ is a closed pro-$p$ subgroup of $\GL_2(B)$ and $L_n(\Gamma_B(\Bb)) = \Sl_2(\Bb^n)$.
\end{lemma}
\begin{proof}
For $x = \left(\begin{smallmatrix} 1+a & b \\ c & 1+d \end{smallmatrix}\right) \in \Gamma_B(\Bb)$ one has $\Theta(x)=\left(\begin{smallmatrix}\frac{a-d}{2} & b \\ c & \frac{d-a}{2} \end{smallmatrix}\right) \in \Sl_2(\Bb)$, so $L_1(\Gamma_B(\Bb)) \subseteq \Sl_2(\Bb)$.  In particular, for any $b, c \in \Bb$, we have $\Theta\left(\begin{smallmatrix} 
1 & b\\
0 & 1
\end{smallmatrix}\right) = \left(\begin{smallmatrix} 
0 & b\\
0 & 0
\end{smallmatrix}\right)$ and $\Theta\left(\begin{smallmatrix} 
1 & 0\\
c & 1
\end{smallmatrix}\right) = \left(\begin{smallmatrix} 
0 & 0\\
c & 0
\end{smallmatrix}\right)$.  For $a \in \Bb$ we have $\left(\begin{smallmatrix} 
1 + 2a & -2a\\
2a & 1 - 2a
\end{smallmatrix}\right) \in \Gamma_B(\Bb)$, and so
\[
\Theta\left(\begin{smallmatrix} 
1 + 2a & -2a\\
2a & 1 - 2a
\end{smallmatrix}\right) = \left(\begin{smallmatrix} 
a & -2a\\
2a & -a
\end{smallmatrix}\right) = \left(\begin{smallmatrix} 
a & 0\\
0 & -a
\end{smallmatrix}\right) + \Theta\left(\begin{smallmatrix} 
1 & -2a\\
0 & 1
\end{smallmatrix}\right) + \Theta\left(\begin{smallmatrix} 
1 & 0\\
2a & 1
\end{smallmatrix}\right).
\]
It follows that $\Sl_2(\Bb)$ is contained in the additive subgroup generated by $\Theta(\Gamma_B(\Bb))$.  Since $\Sl_2(\Bb)$ is closed in $\Sl_2(B)$, it follows that $\Sl_2(\Bb) = L_1(\Gamma_B(\Bb))$. 

It is straightforward to calculate by induction on $n$ that the subgroup topologically generated by 
\[
\{xy - yx \colon x \in \Sl_2(\Bb), y \in \Sl_2(\Bb^n)\}
\]
is $\Sl_2(\Bb^{n+1})$.  That is, $L_n(\Gamma_B(\Bb)) = \Sl_2(\Bb^n)$ for all $n \geq 1$.
\end{proof}

\begin{corollary}\label{congrcomm}
Let $\Bb$ be a closed $B$-ideal different from $B$.  Then $\big(\Gamma_B(\Bb), \Gamma_B(\Bb)\big) = \Gamma_B(\Bb^2)$. If $\#\BB>3$, this also holds for $\Bb = B$. 

\end{corollary}

\begin{proof}
First assume that $\Bb \neq B$.  By \cref{Pink main thm},
\[
\big(\Gamma_B(\Bb), \Gamma_B(\Bb)\big) = \Theta^{-1}(L_2(\Gamma_B(\Bb))) \cap \Gamma_B(\m).
\]
By \cref{congrL1}, $L_2(\Gamma_B(\Bb)) = \Sl_2(\Bb^2)$.  

Clearly $\Gamma_B(\Bb^2) \subseteq \Theta^{-1}(\Sl_2(\Bb^2)) \cap \Gamma_BA(\m)$.  We compute $\Theta^{-1}\left(\begin{smallmatrix} 
a & b\\
c & -a
\end{smallmatrix}\right) \cap \Gamma_B(\m)$ for $\left(\begin{smallmatrix} 
a & b\\
c & -a
\end{smallmatrix}\right) \in \Sl_2(\Bb^2)$.  If $\left(\begin{smallmatrix} 
\alpha & \beta\\
\gamma & \delta
\end{smallmatrix}\right) \in \Theta^{-1}\left(\begin{smallmatrix} 
a & b\\
c & -a
\end{smallmatrix}\right) \cap \Gamma_B(\m)$ then we must have $\beta = b, \gamma = c, \alpha - \delta = 2a$, and $1 = \alpha\delta - \beta\gamma$.  From this one calculates that $\alpha = a \pm \sqrt{1 + a^2 + bc}$ and $\delta = -a \pm \sqrt{1 + a^2 + bc}$.  But only one of these possibilities has $\alpha \equiv 1 \equiv \delta \bmod \m$ and thus is in $\Gamma_B(\m)$.  That is, there is a unique element in $\Theta^{-1}\left(\begin{smallmatrix} 
a & b\\
c & -a
\end{smallmatrix}\right) \cap \Gamma_B(\m)$.  It follows that $\Theta^{-1}(\Sl_2(\Bb^2)) \cap \Gamma_B(\m) = \Gamma_B(\Bb^2)$, as desired. 

We now prove that $\big(\SL_2(B), \SL_2(A)\big) = \SL_2(A)$ when $\#\BB>3$.  By the first statement in the corollary, we know that $\Gamma_B(\m^2) \subseteq \big(\SL_2(B), \SL_2(B)\big)$, so we may assume that $\m^2 = 0$.  Furthermore, the residual image of $\big(\SL_2(B), \SL_2(B)\big)$ is $\big(\SL_2(\BB), \SL_2(\BB)\big)$, which is equal to $\SL_2(\BB)$.  Therefore, it suffices to show that $\left(\begin{smallmatrix} 
1 + a & b\\
c & 1-a
\end{smallmatrix}\right) \in  \big(\SL_2(B), \SL_2(B)\big)$ for any with $a, b, c \in \m$.  Since $\m^2 = 0$, we can decompose
\[
\left(\begin{smallmatrix} 
1 + a & b\\
c & 1-a
\end{smallmatrix}\right) = \left(\begin{smallmatrix} 
1 + a & 0\\
0 & 1-a
\end{smallmatrix}\right)\left(\begin{smallmatrix} 
1 & b\\
c & 1
\end{smallmatrix}\right).
\]  

Let $x \in B^\times$ such that $x^2 \not\equiv 1 \bmod \m$, which exists since $\#\BB > 3$.  Note that for any $\beta, \gamma \in \m$ we have
\begin{align*}
\left(\begin{smallmatrix} 
1 & b\\
c & 1
\end{smallmatrix} \right) 
&= \left(\begin{smallmatrix} 
1 & b(1 - x^2)^{-1}\\
c(1 - x^{-2})^{-1} & 1
\end{smallmatrix}\right) \left(\begin{smallmatrix} 
x & 0\\
0 & x^{-1}
\end{smallmatrix}\right)  \left(\begin{smallmatrix} 
1 & -b(1 - x^2)^{-1}\\
-c(1 - x^{-2})^{-1} & 1
\end{smallmatrix}\right) \left(\begin{smallmatrix} 
x^{-1} & 0\\
0 & x
\end{smallmatrix} \right) \\
&\in \big(\SL_2(B), \SL_2(B)\big)
\end{align*}
and
\[
\left(\begin{smallmatrix} 
1 + a & 0\\
0 & 1 - a
\end{smallmatrix} \right) = \left(\begin{smallmatrix} 
1 + \frac{a}{2} & 0\\
0 & 1 - \frac{a}{2}
\end{smallmatrix} \right)\left(\begin{smallmatrix} 
0 & 1\\
-1 & 0
\end{smallmatrix} \right)
\left(\begin{smallmatrix} 
1 - \frac{a}{2} & 0\\
0 & 1 + \frac{a}{2}
\end{smallmatrix} \right)\left(\begin{smallmatrix} 
0 & -1\\
1 & 0
\end{smallmatrix} \right) \in \big(\SL_2(B), \SL_2(B)\big).
\]
It follows that $\Gamma_B(\m) \subseteq \big(\SL_2(B), \SL_2(B)\big)$ and hence that $\SL_2(B)$ is its own topological derived subgroup.
\end{proof}

Having calculated the derived subgroup of a congruence subgroup, we can now prove that fullness for closed subrings of $A$ is inherited by restrictions to finite-index and coabelian subgroups.

\begin{proposition}\label{closednormal1}\label{finiteindex}
Suppose that the pseudorepresentation $(t, d): \Pi \to A$ is $B$-full for a subring $B$ of $K$. Let $\Pi_0$ be a closed normal subgroup of $\Pi$ so that $\Pi/\Pi_0$ is abelian. Then $\left(\left. t \right|_{\Pi_0}, \left. d \right|_{\Pi_0} \right)$ is also $B$-full. The same is true if $\Pi_0$ is a closed finite-index subgroup of $\Pi$ so long as $B$ is not finite. 
\end{proposition}

\begin{proof}
First, assume that $B$ is a closed subring of $A$. 

Let $\rho \colon \Pi \to R^\times$ be a $(t,d)$-representation such that $\rho(\Pi)$ contains $\Gamma_{B}(\Bb)$ for some nonzero $B$-ideal $\Bb$. 
Write \mbox{$G \coloneqq \rho(\Pi)$} and let \mbox{$G_0:=\rho(\Pi_0)$}. If $\Pi_0$ is coabelian, then $G/G_0$ is abelian, so that $G_0$ contains the derived subgroup $[G, G]$. In particular, $G_0$ contains $\left[\Gamma_{B}(\Bb), \Gamma_{B}(\Bb)\right]$, which is $\Gamma_{B}(\Bb^2)$ by \cref{congrcomm}. Since $A$ is a domain, $\Bb^2$ is nonzero if $\Bb$ is.

Suppose alternatively that $\Pi_0$ is finite index in $\Pi$ and $B$ is not finite. Replacing $\Bb$ by $\m_A \cap \Bb$ (that is, 
if
$\Bb = B$, we replace $\Bb$ by the  maximal ideal $\m_A \cap B$ of $B$, nonzero by the assumption on $B$), we note that $\Gamma_B(\Bb)$ is contained in $\rho(\Gamma)$ for $\Gamma := \ker \rhobar \subseteq \Pi$. Let $\Gamma_0$ be the normal core 
of $\Gamma \cap \Pi_0$ inside $\Gamma$, so that $\Gamma_0$ is a finite-index normal subgroup of $\Gamma$ contained in~$\Pi_0$. Since $\Gamma$ is pro-$p$, and hence pro-solvable, and $\Gamma/\Gamma_0$ is finite, $\Gamma_0$ must contain the $n^{\rm th}$ derived subgroup of $\Gamma$ for some $n \geq 1$. Therefore $\rho(\Gamma_0)$, and hence $\rho(\Pi_0)$, contains the $n^{\rm th}$ closed derived subgroup of $\Gamma_B(\Bb)$, namely $\Gamma_B(\Bb^{2^n})$ (\cref{congrcomm}). This last is again a nontrivial congruence subgroup of~$\SL_2(B)$. 

For arbitrary $B$, \cref{closedinAt} allows us to assume that $B$ is a subring of $A$. If $(t, d)$ is $B$-full of level $\Bb$ for some ideal $\Bb$ of $B$, then the argument above and the first part of \cref{fullERprop} tell us that $\left(\left. t \right|_{\Pi_0}, \left. d \right|_{\Pi_0} \right)$ is $\ovl B$-full of level $\ovl \Bb^{2^n}$ for some $n \geq 1$. This ideal is the closure of the ideal $\Bb^{2^n}$ of~$B$, so we are done by the second part of \cref{fullclosedprop}. 
\end{proof}

We conclude that fullness is unchanged under twisting. 
\begin{corollary}\label{twisting does not affect fullness}
Let $(t, d): \Pi \to A$ be a pseudorepresentation and $\chi: \Pi \to A^\times$ a continuous character. If $(t, d)$ is $B$-full for some subring $B$ of $K$, then $(\chi t, \chi^2 d)$ is also $B$-full. 
\end{corollary}

\begin{proof}
Follows from \cref{closednormal1} by setting $\Pi_0 := \ker \chi$. 
\end{proof}

\medskip
\section{Adjoint trace rings and conjugate self-twists}\label{cst section}


Having established the terms of the investigation --- finding congruence subgroups for \emph{fullness rings} (\cref{fullness}) contained in images of GMA-valued representations (\cref{Bellaiche 2.4.2}) --- in this section we search for optimal (fullness peerage equivalence classes of) fullness rings. Our first stop is the ring fixed by the conjugate self-twists of $(t, d)$ (\cref{cstadjcst} below), symmetries that naturally limit its image.  Although this fixed-by-twist-automorphisms subring is a fullness ring for the historical big-image results that serve as our inspiration, one cannot expect fullness with respect to the ring fixed by conjugate self-twists in the general setting.   
Indeed, there may not be enough automorphisms to carve down to a fullness ring  as illustrated in \cref{gen cst example,insepex}, reflecting the limits of Galois theory.

Instead of trying to cut out a fullness ring from above, we build one from below by considering the trace ring $A_0$ of the adjoint pseudorepresentation.  This \emph{adjoint trace ring} (\cref{A0 defn} below) does not change when twisting $(t, d)$ by a character, and moreover is morally expected to be pointwise fixed by all conjugate self-twists. Because of topological considerations, we do not actually show that $A_0$ is fixed by all conjugate self-twists (\cref{A0fullsimpleCST}) until after we prove our main fullness result, so this idea is merely a guiding principle --- except for $(t, d)$ whose determinant is \emph{$A_0$-constant} (\cref{a0constantdetdef}), a condition expected to be satisfied by all intended applications.

The main result of this section is \cref{AcstA0}: 
if $(t, d)$ has $A_0$-constant--determinant, then $A_0$ and the ring fixed by conjugate self-twists are fullness peers.  We crucially use this fullness peerage result when deriving our $A_0$-fullness results for certain constant-determinant pseudorepresentations satisfying our mild conditions (\cref{final goal}), which we then propagate to all such pseudorepresentations using the twist-invariance of $A_0$ (\cref{nonconst det main thm}). 

\subsection{Conjugate self-twists}\label{csts intro} 

Recall that $A$ is a local pro-$p$ noetherian ring and $\Pi$ is a $p$-finite profinite group. Fix a continuous pseudorepresentation \mbox{$(t, d): \Pi \to A$} with trace algebra~$A_t$. 

\begin{definition}\label{cstadjcst}\label{a0def}\
If $A$ is a domain, let $B$ be a domain extending it; otherwise let $B = A$.  
A \emph{($B$-valued) conjugate self-twist} of $(t, d)$ is a pair $(\sigma, \eta)$, where $\sigma$ is an automorphism of $B$ as a $\Z_p$-algebra and $\eta: \Pi \to B^\times$ is a character.  We also consider $\til\Sigma_t(B/C)$, the conjugate self-twists whose automorphisms $\Sigma_t(B/C)$ fix a subring $C$ of $B$. 
\end{definition}

The set of all $B$-valued conjugate self-twists of a pseudorepresentation forms a group $\til \Sigma_t(B)$, with composition law $(\sigma_1, \eta_1) \circ (\sigma_2, \eta_2) = (\sigma_1 \sigma_2,\eta_1 \act{\sigma_1}{\eta_2})$
and inverse $(\sigma, \eta)^{-1} = (\sigma^{-1}, \act{\sigma^{-1}}{\eta^{-1}}).$ 
Forgetting the character is therefore is a group homomorphism \mbox{$\til \Sigma_t(B) \to \Aut_{\Z_p}(B)$} whose image we denote $\Sigma_t(B)$. If $A$ is a domain and \mbox{$\rho: \Pi \to \GL_2(A)$} is a semisimple representation, we use the notation $\Sigma_\rho(B) := \Sigma_{\tr \rho}(B)$ and $\til\Sigma_\rho(B)\coloneqq \til\Sigma_{\tr \rho}(B)$.

The kernel of the forget-the-character map are the \emph{dihedral} conjugate self-twists $\Sigma_t^{\di}$: 
$$1 \to \Sigma_t^{\di} \to \til \Sigma_t(B) \to \Sigma_t(B) \to 1.$$
If $(1, \eta)$ is a nontrivial dihedral conjugate self-twist, then one can check that $H := \ker \eta$ is an index-$2$ subgroup of $G$, that $(\left.t\right|_H, \left.d\right|_H)$ is a sum of two characters, and that $(t, d)$ is carried by the induction of either of them. 
In particular $\eta$ is quadratic so that $\Sigma_t^{\di}$ does not depend on $B$.  Moreover, when $B$ is a field, $H$, and hence $\eta$, is uniquely defined by $(t, d)$ unless the projective image of the semisimple $\rho$ carrying $(t, d)$ is the Klein-$4$ group, in which case there are three possibilities for $H$. In other words, when $B$ is a field, $\Sigma_t^{\di}$ is abelian dihedral if $\rho$ is projectively dihedral, cyclic of order $2$ if $\rho$ is reducible of projective order $2$, and trivial in all other cases.  See \cref{cm cst implies eta quadratic,resolve definitions of dihedral}, and the proof of \cref{precise structure of Sigma} for details.

\begin{proposition}\label{cstaftertwist} If $\chi: \Pi \to A^\times$ is a continuous character and $B$ is an extension of $A$, then $$\til \Sigma_t(B) \cong \til \Sigma_{\chi t}(B).$$
\end{proposition}

\begin{proof}
The map $(\sigma, \eta) \mapsto (\sigma, \act{\sigma}{\chi} \chi^{-1} \eta)$ realizes the isomorphism.
\end{proof}

If a conjugate self-twist $(\sigma, \eta)$ happens to be $A$-valued, then the automorphism $\sigma$ is automatically continuous: since $A$ is local, algebraic automorphisms automatically send the maximal ideal to itself, and since $A$ is noetherian, the maximal ideal defines the topology.\label{localautcont}
It turns out that an $A$-valued conjugate self-twist character $\eta$ must also be continuous --- see \cref{a0all}\eqref{etafinitecont} below. But in general we cannot expect that all conjugate self-twists can be restricted to ones defined over $A$, as seen in \cref{gen cst example} below illustrating a failure of normality in the sense of Galois theory. 

\begin{example}\label{gen cst example}
For any odd prime $p$, let $\Pi$ be the subgroup of $\GL_2(\Z_p[\sqrt[p]{p}])$ generated by $\GL_2(\Z_p)$ and the scalar matrix $1 + \sqrt[p]{p}$.
Let $\rho \colon \Pi \to \GL_2(\Z_p[\sqrt[p]{p}])$ be the inclusion, a representation with trace algebra $A = \Z_p[\sqrt[p]{p}]$. Note that $\Sigma_\rho(A)$ is trivial, since $A$ has no automorphisms --- but that's not the whole story. Let $B = \Z_p[\sqrt[p]{p}, \zeta_p]$ and consider the automorphism $\sigma$ in $\Aut(B)$ sending $\sqrt[p]{p}$ to $\zeta_p\sqrt[p]{p}$, together with the character $\eta: \Pi \to B^\times$ with kernel $\GL_2(\Z_p)$ sending $1 + \sqrt[p]{p}$ to $(1 + \zeta_p\sqrt[p]{p})(1 + \sqrt[p]{p})^{-1}$. Then $(\sigma, \eta)$ is a $B$-valued conjugate self-twist of $\rho$. 
One can show that $\rho$ is not $A$-full ($\tr \rho$ does not span an ideal of $A$; see \cref{closednormal2}). On the other hand, $\rho$ is visibly $\Z_p$-full, corresponding to the fact that $\Z_p$ is the ring fixed by~$\Sigma_\rho(B)$. 
\end{example}

Because it is not immediately clear that we can demand that relevant extensions of $A$ be endowed with a sensible topology (if $A$ is a domain we already cannot expect a topology on $Q(A)$: see p.\,\pageref{notop}), there is no way to require conjugate self-twists to be continuous.  On the other hand, since trace algebras are \textit{topologically} generated rings, they only behave well under conjugate self-twists satisfying continuity conditions. As a result, much of this section is devoted to finding settings where conjugate self-twists are $A$-valued and hence continuous.  This will eventually allow us to show that all conjugate self-twists of constant-determinant pseudorepresentations satisfying the mild conditions of \cref{main thm intro} are continuous: see 
\cref{A0fullsimpleCST}.
Write $A_t^\alg$ for the $W(\F)$-subalgebra of $A$ algebraically generated by $t(\Pi)$, so that the trace algebra $A_t$ is the closure of  $A_t^\alg$ in $A$.

\begin{proposition}\label{a0all} Let $(\sigma, \eta)$ be a $B$-valued conjugate self-twist.
\begin{enumerate}
\item \label{a0etafinite} If $(t,d)$ has constant determinant, then $\eta$ is $W(\F)$-valued and finite order. 
\item \label{a0sigmalikecont} If $\eta$ is ${A_t^\alg}$-valued (for example, if $\eta$ is $W(\F)$-valued), then $\sigma$ restricts to an automorphism of $A_t^\alg$, and there is an $A_t$-valued conjugate self-twist $(\sigma', \eta)$ satisfying $\left.\sigma'\right|_{A_t^\alg} =  \left.\sigma\right|_{A_t^\alg}$. 
\item \label{etafinitecont} If $\eta$ is $A$-valued, then $\eta$ is continuous.

\end{enumerate}
\end{proposition}

\begin{proof}
\begin{enumerate} 
\item We follow \cite[1.5]{Momose81}. Recall that $d = s(\ovl{d})$ has finite order by the assumption. 
Since
$\act{\sigma}{d} = \eta^2 d$ 
we have 
$\eta^2 = \act{\sigma}{d}d^{-1}$.  Since $d$ is finite order, $\act{\sigma}{d}$ must be a power of $d$. We now claim that $\act{\sigma}{d}d^{-1}$ has a power of $d$ as a square root, so that $\eta$ is differs from a power of $d$ by at most a quadratic character and hence takes values in $W(\F)$. 
Indeed, if $d$ has odd order, then $d$ itself has a power-of-$d$ square root, so that any power of $d$ has the same property. And if $d$ has even order, then since $\sigma$ preserves orders, $\act\sigma d$ must be an odd power of $d$; which means that $\act{\sigma}{d}d^{-1}$ is an even power of $d$ and hence has a power-of-$d$ square root.

\item Given that $\eta$ and $t$ are both $A_t^\alg$-valued, 
$\act\sigma t = \eta t$ is $A_t^\alg$-valued as well. Since $B$ is either a domain or a local ring and $\sigma$ is a $\Z_p$-algebra homomorphism, it follows that $\sigma$ permutes $W(\F)$, and hence it permutes $A_t^\alg$ as well. Moreover, the action of $\sigma$ is continuous on $A_t^\alg$ in the topology from $A_t$: with $\m^\alg: = \m \cap A_t^\alg$, we have $\sigma(\m^\alg) \subseteq \eta(\Pi) \m^\alg \subseteq \m^\alg$.  
Therefore $\left. \sigma\right|_{A_t^\alg}$ extends uniquely to a continuous  automorphism $\sigma'$ of the closure $A_t$ of $A_t^\alg$. Finally, since $\sigma$ and $\sigma'$ agree on $t(\Pi)$, the pair $(\sigma', \eta)$ is still a conjugate self-twist, this time $A_t$-valued.

\item First, we claim that $\ker (t, d) \subseteq \ker \eta$, so that $\eta$ factors through $G:=\Pi/\ker (t, d)$. Indeed, for $g \in \ker (t, d)$ we have in particular $t(g) = 2$  so that \mbox{$\eta(g) = \act\sigma t(g) / t(g)  = 1$}.
Since $\ker (t, d)$ is closed, $G = \Pi/\ker (t, d)$ is still $p$-finite profinite, and we check continuity of $\eta$ as a character on $G$. Moreover $\Gamma := \ker \bar \rho/\ker(t, d) \subseteq G$ is a finite-index pro-$p$ subgroup of $G$ \cite[Lemma~3.8]{ChenevierDet}. 
By $p$-finiteness, $\Gamma$ is topologically finitely generated, so that any finite-index subgroup of $\Gamma$ is open in $\Gamma$ \cite[{$\S$4.2~exercise~6(d)}]{SerreGalCoh}, and hence in $G$. Now use the fact that $A^\times \cong \F^\times \times (1 + \m)$ to write $\eta = \eta^{(p)} \eta_p$, where $\eta_p: G \to 1 + \m$ is a pro-$p$ character and $\eta^{(p)}:G\to \F^\times$ has prime-to-$p$ order. 
Then $\eta^{(p)}$ is continuous because its kernel contains the open pro-$p$ subgroup $\Gamma$.
And $\eta_p$ is continuous because the preimage $U \subseteq \Gamma$ of any open subgroup of $1 + \m$ along $\left.\eta_p\right|_\Gamma$ is finite index in $\Gamma$, hence open in $G$. Therefore $\eta$ is continuous.\qedhere
\end{enumerate}
\end{proof}

\begin{corollary}\label{atisbigenough}
If $(t, d)$ has constant determinant, then any conjugate self-twist $(\sigma, \eta)$ {of $(t,d)$} 
has~$\eta$ continuous, finite-order, and $W(\F)$-valued. Moreover, 
there exists an automorphism $\sigma'$ of $A_t$ agreeing with $\sigma$ on $A_t^\alg$ so that $(\sigma', \eta)$ is an $A_t$-valued conjugate self-twist of {$(t,d)$}. 
\end{corollary}

\noindent {A posteriori after our main fullness results, we will deduce that, at least for most constant-determinant $(t, d)$ {that are not a priori small}, $\sigma$ restricts to an automorphism of $A_t$, necessarily continuous, and consequently $\sigma' = \left.\sigma\right|_{A_t}$: see \cref{nonconst det main thm,fullfixed}.}

\subsection{The adjoint trace ring}\label{toptroub}
Quite generally, if $V$ is a $2$-dimensional vector space over a field~$F$, then the adjoint action of matrices $m \in \GL(V)$ on $\End_F(V)^0$ --- that is, the conjugation action on the vector-space of trace-zero endomorphisms --- has trace $(\tr m)^2/\det m$ and factors through $\PGL(V)$ since scalars act trivially. 
By analogy, we define the adjoint-trace elements and the adjoint-trace ring attached to $(t, d)$. 

\begin{definition}\label{A0 defn}
The \emph{adjoint trace ring} of $(t, d)$, denoted $A_{0, t}$ or simply $A_0$, is the closed subring of $A$ topologically generated by the \emph{adjoint-trace elements} $\tr \ad t(g):=t(g)^2/d(g)$ for $g \in \Pi$.
\end{definition}

\begin{proposition}\label{a0n2}
The adjoint trace ring $A_0$ of any pseudorepresentation of a $p$-finite profinite group is a local noetherian pro-$p$ ring. In particular, it is N2 (see p.\:\pageref{n2}). 
\end{proposition}

\begin{proof}
By definition, $A_0$ is a closed subring of the local pro-$p$ ring $A$, so that $A_0$ is also local and pro-$p$. Moreover, by construction $A_0$ is the trace ring of a pseudodeformation of $\ad^0 \rhobar$, where $\rhobar$ is the semisimple residual representation carrying $(t, d)$. Since $\Pi$ is $p$-finite, the universal deformation ring of the pseudorepresentation associated to $\ad^0 \rhobar$ is noetherian, hence so is its quotient $A_0$. 
\end{proof}

Adjoint-trace elements don't change when $\big(t(g), d(g)\big)$ is replaced by $\big(\alpha t(g), \alpha^2 d(g)\big)$ for nonzero scalars $\alpha$. This has two consequences. First, the adjoint-trace ring is unchanged under twisting.

\begin{proposition} \label{a0fixedbytwist} \label{a0indepoftwist}
If $\chi: \Pi \to A$ is a continuous character, then $A_{0, t} = A_{0, \chi t}$.
\end{proposition}

\begin{proof}
The adjoint-trace elements of $(t,d)$ and of $(\chi t, \chi^2 t)$ are the same: $\frac{t(g)^2}{d(g)} = \frac{\left( \chi(g) t(g)\right)^2}{ \chi^2(g) d(g)}$.
\end{proof} 

\noindent To state the second, we define $A_{0, t}^\alg \subseteq A_{0, t}$ as the subring generated algebraically rather than topologically by the values of $t^2/d$, so that by definition $A_{0, t}$ is the closure of $A_{0, t}^\alg$ in $A$.

\begin{proposition} \label{a0algfix} If $(\sigma, \eta)$ is an arbitrary conjugate self-twist of $(t,d)$, then $\sigma$ fixes $A_{0, t}^\alg$ pointwise. 
\end{proposition}

\begin{proof}
The adjoint-trace elements are fixed by any $\sigma$. Indeed, for $g \in \Pi$, 
\[
\act{\sigma}{\left(t(g)^2/d(g)\right)} = \act \sigma t (g)^2/\act\sigma d(g) =  \big( \eta(g) t(g)\big)^2/\eta(g)^2 d(g) = t(g)^2/d(g).
\qedhere
\]
\end{proof}

In spite of \cref{a0algfix}, there may not be enough automorphisms to cut out $A_0$ or $A_0^\alg$, as the following example illustrates.

\begin{example}\label{insepex} Let $G = \GL_2(\F_p\lb x \rb)$ and let $\rho:G \to \GL_2\big(\F_p \lb x \rb\big)$ be the inclusion representation. Define $\chi: G \to \F_p \lb x^{1/p}\rb^\times$ as follows: for $g \in G$ let $d_g := \det g \in \F_p \lb x \rb$, and set $\chi(g) \coloneqq d_g(x^{1/p})$. Then $\chi$ is a character to $A \coloneqq \F_p \lb x^{1/p} \rb$, and we consider the representation $\rho \otimes \chi: G \to \GL_2\big(A)$. Since $\rho$ has no nontrivial conjugate self-twists --- or, said more precisely, for any extension $B$ of $A_0 := \F_p \lb x \rb$, we have $\til\Sigma_\rho(B/A_0) \cong \Aut(B/A_0)$, each appearing with trivial character --- the conjugate self-twists of $\rho \otimes \chi$ defined over an extension $B$ of $A$ are all of the form $(\sigma, \act{\sigma}\chi/\chi)$ for $\sigma \in \Aut(B/A_0)$ (\cref{cstaftertwist}). Since any automorphism that fixes $A_0$ also fixes $A$, the ring fixed by all conjugate self-twists of $\rho \otimes \chi$ is $A$ itself. But in this case it's clear by inspection that $A$ is not a fullness ring (see, for example \cref{closednormal2}). On the other hand $A_0$ is a fullness ring --- the image of $\rho \otimes \chi$ contains $\SL_2(A_0)$. 
\end{example}

In the next section we show that, working over fields rather than rings, one can always cut out the analogue of $A_0$ with conjugate self-twists provided the usual conditions from Galois theory hold.

\subsection{Interlude: the theory over abstract fields}
In this subsection we switch gears to a topology-free setting, where we show that under Galois-theoretically favorable conditions, the \emph{adjoint trace field} (defined below) is exactly the fixed field of all conjugate self-twists. This result allows us to give an interpretation of the residue field of $A_0$ as the fixed field of residual conjugate self-twists, and is both inspiration and key in deducing an analogous result for (a generalization of) constant-determinant pseudorepresentations in our pro-$p$ setting in the next subsection --- a lodestar and a tool.

Let $F$ be an arbitrary field whose characteristic is not $2$, $G$ an abstract group, and $(t, d): G \to F$ a pseudorepresentation.  For an extension $L$ of $F$, define an $L$-valued conjugate self-twist $(\sigma, \eta)$ of $(t, d)$ as in \cref{cstadjcst}, but where $\sigma$ is simply a field automorphism rather than a morphism of $\Z_p$-algebras.  Write $\til \Sigma_{t}(L)$ for the set of $L$-valued conjugate self-twists of $(t,d)$, and $\Sigma_t(L)$ for the group of automorphisms appearing in $\til \Sigma_t(L)$. Define the \emph{adjoint-trace field} of $(t, d)$ to be the subfield $F_0 := F_{0, t}$ generated by the adjoint-trace elements $t(g)^2/d(g)$ over all $g \in G$. Then as in \cref{a0algfix}, the adjoint trace field $F_{0, t}$ is fixed by all conjugate self-twist automorphisms, so that $F_{0, t} \subseteq L^{\Sigma_t(L)}$. The following theorem gives conditions that guarantee reverse containment.

\begin{theorem}\label{fautscsts} 
Let $L$ be a separably closed extension of $F$, and $E \subseteq L$ an extension of $F_{0, t}$. Then $$\Sigma_t(L/E) = \Aut(L/E).$$
In particular, $\Sigma_t(F^\sep) = \Aut(F^\sep/F_{0, t}).$ 
\end{theorem}

\begin{proof}
To show that $\Aut(L/E) \subseteq \Sigma_t(L/E)$, we start with $\sigma \in \Aut(L/E)$ and produce a twist character. Let $\rho$ be a semisimple representation over a finite extension $F'$ of $F$ carrying $(t, d)$. Since $\sigma$ fixes the adjoint trace algebra $F_0 \subseteq E$, we have $\tr \ad \act{\sigma}{\rho} = \tr \ad \rho$. Moreover, $\ad \rho$ is semisimple if $\rho$ is. By Brauer-Nesbitt, therefore, the multiplicities of irreducible representations of $G$ inside $\ad \rho$ and inside  $\ad (\act{\sigma}{\rho})$ are equal in (the prime field of) $F$.  If $F$ has characteristic $0$ or characteristic $p \geq 5$, then we can already conclude that $\ad \rho \cong \ad (\act\sigma\rho)$. If $F$ has characteristic~$3$, then we split off a trivial character acting on the center using $\ad \rho = 1 \oplus \ad^0 \rho$ and need only eliminate the case that $\ad^0 \rho \cong \phi^{\oplus 3}$ and $\ad^0 \act\sigma\rho \cong \act{\sigma}{\phi}^{\oplus 3}$ for some character $\phi$ with $\phi \neq \act\sigma\phi$. But since at least one of the eigenvalues of $\ad^0 \rho(g)$ is $1$ for any $g \in G$ (see  \cref{trace of adjoint}), this is impossible. Finally,  \cref{isomorphic ad implies differ by twist} says that, since $\ad \rho \cong \ad \act\sigma\rho$ as representations over $L$, there is a character $\eta \colon G \to L^\times$ such that $\act{\sigma}{\rho} \cong \eta \otimes \rho$. Therefore $(\sigma, \eta)$ is a conjugate self-twist! 

The second statement is a special case of the first since every conjugate self-twist fixes $F_{0, t}$.
\end{proof}

\begin{corollary} \label{f0ifsep}
\begin{enumerate}
\item \label{basicff0} If $F/F_{0,t}$ is a separable extension, then $(F^\sep)^{\Sigma_t(F^\sep)} = F_{0, t}$. 
\item \label{sick}  
Suppose that $L$ is a separably closed extension of $F$, and $E \subseteq L$ is an extension of $F_{0,t}$.\\ If $L / E$ is separable, then $L^{\Sigma_t(L/E)} = E.$
\end{enumerate}
\end{corollary}

\noindent 
{Both} statements follow from \cref{fautscsts}. The conditions of \cref{f0ifsep}\eqref{sick} are nontrivially satisfied, for example, if $F_{0,t} = \Q$, $F = \Q(\sqrt[3]{2})$, $E = \Q_p$ and $L = \ovl\Q_p$. Or see \cref{repok}\eqref{cstcarvek0end}.

In particular \cref{f0ifsep}\eqref{basicff0}  for $F = \F$ finite applies in our pro-$p$ setting. In this case, every pseudorepresentation $(t,d): \Pi \to F$ is carried by a unique semisimple $\rhobar: \Pi \to \GL_2(\F)$. Write $\E$ for the adjoint trace field $\F_{0, \rhobar}$.

\begin{corollary}\label{residual csts cut out E}
$\E =  \ovl \F^{\Sigma_{\ovl \rho}(\ovl \F)} = \F^{\Sigma_{\ovl \rho}(\F)}$ 
\end{corollary}

\begin{proof}
The first equality is \cref{f0ifsep}\eqref{basicff0}. The second follows from \cref{atisbigenough}, which implies that every conjugate self-twist of $\rhobar$ restricts to one defined over $\F$. 
\end{proof}

Our next goal is to limn a setting where the pro-$p$ analogue of the extension $F/F_{0,t}$ is Galois, so that we can obtain analogues of \cref{fautscsts} and \cref{f0ifsep}\eqref{basicff0}.

\subsection{$A_0$-constant determinant and simple conjugate self-twists} We now return to our \mbox{pro-$p$} setting and consider a restriction on the determinant of $(t, d)$ that shares many properties with the constant-determinant case, but is mild enough to be expected to be satisfied by all our arithmetic applications. This constraint allows us to restrict our attention to conjugate self-twists valued in the trace algebra.

\begin{definition}\label{a0constantdetdef}\
We say that $(t,d)$ has \emph{$A_0$-constant determinant} if $d$ is the product of an $A_0$-valued character and a character of finite prime-to-$p$ order. Said another way, $d$ is $A_0$-constant if its pro-$p$ part $d_1$ takes values in $A_0$. If $(t, d)$ has $A_0$-constant determinant,  we call the conjugate self-twists in $\til \Sigma_t(A_t)$ \emph{simple}. For brevity we write $\til \Sigma_t$ and $\Sigma_t$ for these simple conjugate self-twists in Sections \ref{section: lifting} through \ref{conceptual interpretation of Bellaiche ring}.
\end{definition}

\noindent By automatic continuity, all automorphisms in $\Sigma_t(A_t)$ fix $A_0$, so that $\Sigma_t(A_t) = \Sigma_t(A_t/A_0)$.   The following lemma shows that twist characters are also  continuous. 

\begin{lemma}\label{nearly const det}
Let $(t, d): \Pi \to A$ be an $A_0$-constant-determinant pseudorepresentation, and let $(t',d')$ be its constant-determinant twist obtained by twisting off the pro-$p$ part $d_1$ of $d$. Then: \begin{enumerate}
\item \label{sameas} $A_t = A_{t'}$
\item \label{samesigs} $\til \Sigma_t(A_t) = \til \Sigma_{t'}(A_{t'})$
\item \label{etagood} If $(\sigma, \eta)$ in $\til \Sigma_t(A_t)$ is a simple conjugate self-twist, then $\eta$ is $W(\F)^\times$-valued and continuous. 
\end{enumerate}
\end{lemma}

\begin{proof}
First note that $(t, d)$ and $(t', d')$ have the same adjoint trace algebra $A_0$, since the latter is twist invariant (\cref{a0indepoftwist}).  Write $\chi$ for $d_1^{-1/2}$, a continuous $A_0$-valued character. Note that $A_0 \subset A_{t} \cap A_{t'}$. Since $\chi$ and $\chi^{-1}$ are both valued in $A_0 \subset A_{t} \cap A_{t'}$, we can move back and forth using $t' = \chi t$ and $t = \chi^{-1} t$. That is $A_{t'} = A_{\chi t} \subseteq A_t$ and $A_t = A_{\chi^{-1} t'} \subseteq A_{t'}$. Part~\eqref{sameas} follows. For~\eqref{samesigs}, because $\chi$ is fixed by $\sigma$, the map $\til \Sigma_t(A_t) \to \til \Sigma_{\chi t'}(A_t)$ from \cref{cstaftertwist} sending $(\sigma, \eta)$ to $(\sigma, \act \sigma \chi \chi^{-1} \eta)$ is the identity. For~\eqref{etagood}, use \eqref{samesigs} and \cref{atisbigenough}, or redo \cref{a0all}\eqref{a0etafinite} directly and use \cref{a0all}\eqref{etafinitecont}. 
\end{proof}

\subsection{$K/K_0$ as a Galois extension}

Finally we assume that $A$ is a (local pro-$p$) \underline{domain} and fix an \underline{$A_0$-constant-determinant} pseudorepresentation $(t, d): \Pi \to A$. By replacing $A$ by $A_t$, we may assume that $A$ is the trace algebra of $(t, d)$; let $K$ be the fraction field of $A$. Recall that $\rhobar$ is the semisimple representation $\Pi \to \GL_2(\F)$ carrying $(\bar t, \ovl d)$. 
Let $A_0$ be the adjoint trace ring of $(t, d)$ with fraction field $K_0$; write $\E$ for the adjoint trace ring of $\rhobar$. Here $\F$ is the residue field of $A$ and $\E$ is the residue field of $A_0$. 
	Let $\Gamma := \ker \rhobar \subset \Pi$. This is a finite-index normal subgroup of $\Pi$, and we can restrict $(t, d)$ to $\Gamma$ to obtain $(\left.t \right|_\Gamma, d_1)$.

\begin{lemma}\label{tgammainA0}
We have $t(\Gamma) \subset A_0$.
\end{lemma}

\begin{proof}
For $\gamma \in \Gamma$, we have $t(\gamma) = 2 + m$ for some $m \in \m$. The corresponding adjoint-trace element is $t(\gamma)^2/ d_1(\gamma) = d_1(\gamma)^{-1} (4 + 4m + m^2) \in A_0$. Since $d_1(\gamma)^{-1} \in A_0$ by assumption, so is $4 + 4m + m^2$. But any closed subring containing $4 + 4m + m^2$ contains $2 + m = t(\gamma)$ as well since $2$ is invertible. 
\end{proof}

\begin{proposition} \label{A0finiteA} $A$ is a finite $A_0$-algebra. 
\end{proposition} 

\begin{proof}
The subgroup $\Gamma := \ker \bar\rho$ is a finite-index normal subgroup of $\Pi$. Moreover, $A_{\left. t\right|_{\Gamma}} \subseteq A_0$ (\cref{tgammainA0}). Thus $K_0$ contains the trace field of $\left. t\right|_\Gamma$, so that by \cref{atmostquad} there is an at-most quadratic extension $L_0$ over $K_0$ over which we can define a representation carrying $(\left. t\right|_\Gamma, d_1)$. By the proof of \cref{tracefieldextension}, there is a finite extension $L$ over $L_0$ containing $t(\Pi)$. In fact, letting $A^\alg$ be the $W(\F) A_0$-submodule of $A$ (algebraically, not topologically) generated by $t(\Pi)$, we have $A^\alg \subseteq L$. Note that $A$ is the closure of $A^\alg$.

We claim that $A^\alg$ is integral over $A_0$. For every $g \in \Pi$, the element $t(g)$ is a square root of~$t(g)^2$, which is in $A_0[d] \subseteq W(\F) A_0$. Since $A^\alg$ is generated by the $t(g)$ over $W(\F)A_0$, the former is integral over the latter. The claim follows. 

Finally, let $B$ be the integral closure of $A_0$ in $L$. Since $A_0$ is noetherian and $N2$ (\cref{a0n2}), $B$ is a finite, hence noetherian, $A_0$-algebra. Because $B$ contains $A^\alg$ by integrality, the latter is also finite over $A_0$. But this means that $A^\alg$ is a sum of finitely many compact submodules of $A$, so compact itself, hence closed. In other words, $A^\alg = A$. Therefore $A$ is finite over $A_0$. 
\end{proof} 

\begin{corollary}\label{AoverA0explicit}
$A$ is a multiquadratic extension of $W(\F)A_0$. 
More precisely, there is an integer $n$ prime to $p$, elements $a_1, \ldots, a_r$ of $A_0$, and integers $k_1, \ldots, k_r$ modulo $n$ so that for $\zeta = \zeta_n$ a generator of~$s(\F^\times)$ we have 
\begin{equation} \label{explicitA}
A = A_0[\zeta, \sqrt{\zeta^{k_1} a_1}, \ldots, \sqrt{\zeta^{k_r} a_r}].
\end{equation}
\end{corollary} 

\begin{proof}
Let $n$ be the order of $\F^\times$, so that $\zeta$ generates $W(\F)$ and $A_0[\zeta] = W(\F)A_0$. 
Since $A_0$ contains $$a(g) := d_1(g) t^2(g)/d(g) = t^2(g)/s(\ovl d(g)) = t^2(g) / \zeta^{k(g)}$$ for every $g \in \Pi$ and $k(g)$ depending on $g$, we see that  $A$ is topologically generated over $W(\F)A_0$ by $t(g) =  \sqrt{\zeta^{k(g)} a(g)}$. By \cref{A0finiteA} and compactness only finitely many of these are needed. 
\end{proof}

\begin{theorem}\label{KgaloisK0}
The extension $K$ over $K_0$ is finite abelian, with $$\Gal(K/K_0) = \Aut(A/A_0) = \Sigma_t(A).$$ 
\end{theorem} 

We give two different arguments both fundamentally rooted in \cref{fautscsts}.

\begin{proof}[First proof of \cref{KgaloisK0}]
Since $K = K_0(\zeta, \sqrt{\zeta^{k_1} a_1}, \ldots, \sqrt{\zeta^{k_r} a_r})$ by \cref{AoverA0explicit} (and using the same notation), it is clear that $\Aut(K/K_0) = \Aut(A/A_0)$. To see that $K/K_0$ is Galois, note that it is a subextension of $K_0(\zeta_{2n}, \sqrt{a_1}, \ldots, \sqrt{a_r})/K_0$, which is a compositum of abelian extensions and hence abelian as well. Here we use the fact that the characteristic of (the prime field of) $K$ is either zero or an odd prime $p$, so that  $K_0(\zeta_{2n})/K_0$ for $n$ prime to $p$ and the $K_0(\sqrt{a_i})/K_0$ are all separable and hence abelian. Therefore $K/K_0$ is abelian with $\Gal(K/K_0)$ a subgroup of $(\Z/2\Z)^r \times (\Z/n\Z)^\times$. 

To see that $\Aut(A/A_0) = \Sigma_t(A)$, start with $\sigma \in \Aut(A/A_0)$ and proceed as in the proof of \cref{fautscsts} to produce a character $\eta: \Pi \to K^\times$ with $(\sigma, \eta)$ a conjugate self-twist of $(t, d)$. By \cref{nearly const det}, $\eta$ is $A$-valued, so that $(\sigma, \eta)$ is in $\Sigma_t(A)$. 
\end{proof}

\begin{proof}[Second proof of \cref{KgaloisK0}]
As in the first proof, it is clear that $K$ is separable over $K_0$. First suppose that $(t,d)$ has constant determinant.  To show that $K/K_0$ is normal, take an embedding $\sigma$ of $K$ into $\ovl K$ over $K_0$. Proceed as in the proof of \cref{fautscsts} to get a character $\eta: \Pi \to \ovl K^\times$ so that $(\sigma, \eta)$ is a $\ovl K$-valued conjugate self-twist. Use \cref{atisbigenough} 
to find $\sigma'$ so that $(\sigma', \eta)$ is an $A$-valued conjugate self-twist. Note that $\left.\sigma'\right|_{W(\F) A_0} = \left.\sigma\right|_{W(\F)A_0}$ since both $\sigma$ and $\sigma'$ fix $A_0$ and act the same way on $W(\F)$ by construction.  Therefore $\sigma'\sigma^{-1}$ is an embedding of $K$ that fixes $Q(W(\F)A_0) = K_0(\zeta)$.  Since $K/K_0(\zeta)$ is a multiquadratic extension, it is Galois, so $\sigma'\sigma^{-1}$ sends $K$ to itself.  But $\sigma'$ also sends $K$ to itself by construction, so $\sigma$ must be an automorphism of $K$. In fact, using that $\sigma'\sigma^{-1}$ fixes $K_0$ and $t(\Pi)$, we see that $\sigma|{_A} = \sigma'$ and thus $\Gal(K/K_0) = \Sigma_t(A)$.  The fact that $\Sigma_t(A)$ is abelian in the constant-determinant case follows from \cref{tildeSigma rho-bar finite} and the diagram following \cref{lifting arbitrary csts}. 
In general when $(t,d)$ has $A_0$-constant determinant, the result follows from the constant determinant-case and \cref{nearly const det}.  
\end{proof}

\begin{corollary}\label{AcstA0}
We have $K_0 = K^{\Sigma_t(A)}$. Moreover, $A^{\Sigma_t(A)}$ is a finite extension of $A_0$ with the same field of fractions and the same normalization. 
\end{corollary} 

\begin{proof}
The first statement follows from \cref{KgaloisK0}. Since $A_0 \subseteq A^{\Sigma_t(A)} \subseteq K^{\Sigma_t(A)} = K_0$, the field of fractions of $A^{\Sigma_t(A)}$ is $K_0$. Finally, since $A$ is integral over $A_0$, so is $A^{\Sigma_t(A)}$.  
\end{proof}

\cref{AcstA0,normalpuppies} imply the following fullness comparison result for two key subrings of $A$.

\begin{corollary} \label{A0 Afixedbytwists fullness peers} 
The rings $A_{0}$ and $A^{\Sigma_t(A)}$, as well as their normalizations, are all fullness peers. 
\end{corollary}

The example below shows that proper containment $A_0 \subsetneq A^{\Sigma_t(A)}$ is possible: $A_0$ and $A^{\Sigma_t(A)}$ need not have the same residue field. 

\begin{example}\label{a0small}
Define $$\ovl G := \left\{\begin{pmatrix} a & b \\ b^p &  a^p \end{pmatrix} \right\} \cup \left\{\begin{pmatrix} a & b \\ - b^p & -a^p \end{pmatrix} \right\} \subset \GL_2(\F_{p^2}),$$
so that in each set above $a$ and $b$ are in $\F_{p^2}$ satisfying $N(a) \neq N(b)$. Let $\widetilde G \subset\GL_2(\F_{p^2} \lb X \rb)$ be the set of matrices whose determinants are in $\F_{p^2}$, and let $G \subset \widetilde G$ be the set of invertible matrices that are residually in $\ovl G$. Finally, let $\Pi := G$ and set \mbox{$\rho: \Pi \to \GL_2(\F_{p^2} \lb X \rb)$} to be the inclusion map, a constant-determinant representation. Then the trace algebra $A$ of $\rho$ is \mbox{$A = \F_{p^2} \lb X \rb$} --- indeed, fix a generator $\beta$ of $\F_{p^2}^\times$, so that $N(\beta) = \beta^{p+1} \neq 1$. Then for every $a \in \F_{p^2}^\times$, $G$ contains both 
$$g_a := \begin{pmatrix} 1 & \beta \\ \beta^p  + \beta^{-1} a X & 1 + aX \end{pmatrix} \qquad \mbox{and}\qquad h_a := \begin{pmatrix} a &0 \\ 0  & -a^p \end{pmatrix},$$ 
so that $A$ contains both $\tr g_a = 2 + aX$ and $\tr h_a = a - a^p$. A similar computation shows that the adjoint trace ring is $A_0 = \F_p + X \F_{p^2}\lb X \rb$: on one hand, $A_0$ contains 
$$(\tr g_a)^2/\det g_a = (4 + 4aX + a^2 X^2)(1 - N(\beta))^{-1}$$
for every $a \in \F_{p^2}^\times,$ and on the other hand every element of $A_0$ is {residually} in $\F_p$, reflecting the fact that $\ovl\rho$  admits a conjugate self-twist and illustrating \cref{residual csts cut out E}. Since $\rho$ itself has no conjugate self-twist --- any automorphism of $A$ {appearing in a conjugate self-twist} must fix $A_0$ and hence extends to the identity on $K_0$, which contains $A$ --- we have $A^{\Sigma_\rho(A)} = A$.  
\end{example}

We close with a technical observation necessary for our a posteriori justification for focusing on $A$-valued conjugate self-twists for $A_0$-constant--determinant pseudorepresentations.

\begin{proposition}\label{a0detstablefixed}
Let $(\sigma, \eta)$ be any conjugate self-twist of an $A_0$-constant--determinant pseudorepresenation $(t, d)$. Its trace algebra $A$ is $\sigma$-stable if and only if $\sigma$ fixes $A_0$ pointwise. Both of these conditions are satisified if $(t, d)$ is $A_0$-full.
\end{proposition}

\begin{proof}
If $\sigma$ restricts to an automorphism of $A$, then $\left.\sigma\right|_{A}$ is automatically continuous, so that $\sigma$ fixing $A_0^\alg$ pointwise per \cref{a0algfix} is enough to conclude that $\sigma$ fixes all of $A_0$. Conversely, if $A_0$ is fixed by $\sigma$, then since $A$ is finite over $A_0$ in this setting (\cref{A0finiteA}), $A^\alg$ is dense in $A$, and $A$ is noetherian, we can express $A = a_1 A_0 + \cdots + a_n A_0$ for $a_i \in A^\alg$, so that $\sigma(A) \subseteq A$. 

If $(t, d)$ is $A_0$-full, then $A^\alg$ contains some nonzero $A_0$-ideal $\Bb_0$ (\cref{closednormal2}). The expression $A = a_1 A_0 + \cdots + a_n A_0$ for $a_i \in A^\alg$ implies that $\Bb_0 A \subseteq A^\alg$, whence $A^\alg$ and $A$ are fullness peers, and in particular have the same field of fractions (\cref{fgqf}). Since $A \subseteq Q(A^\alg)$, the values of $\sigma$ on $A$ are entirely determined by the restriction of $\sigma$ to $A^\alg$. Therefore the action of $\sigma$ on $A$ coincides with the action on $A$ by the continuous extension of $\left.\sigma\right|_{A^\alg}$ guaranteed by \cref{atisbigenough}. Thus $\sigma$ restricts to an automorphism of $A$, as claimed. 
\end{proof} 

\begin{corollary}\label{A0fullsimpleCST} Let $(t, d): \Pi \to A$ be any pseudorepresentation. If $(t, d)$ is full for its adjoint trace ring $A_0$, then every arbitrary conjugate self-twist of $(t, d)$ fixes $A_0$ pointwise.
\end{corollary}
\begin{proof} Let $(\sigma, \eta)$ be a conjugate self-twist of $(t, d)$, and let $(t', d')$
 be its constant-determinant twist. Then $(\sigma, \eta')$ is a conjugate self-twist of $(t', d')$ for some character $\eta'$ (\cref{cstaftertwist}). Since neither fullness nor $A_0$ is changed under twisting (\cref{twisting does not affect fullness,a0fixedbytwist}), the constant-determinant twist $(t', d')$ is still $A_0$-full, and hence by \cref{a0detstablefixed} $\sigma$ fixes $A_0$ pointwise.
\end{proof}

\noindent Combined with our main $A_0$-fullness result (\cref{nonconst det main thm}), \cref{A0fullsimpleCST} allows us to conclude a posteriori that under mild conditions on {non--a-priori-small} $(t, d)$, all conjugate self-twists of $(t, d)$ fix all of $A_0$, resolving the topological concerns we have danced around in \cref{csts intro,toptroub}. Of course one continues to hope for a less circuitous argument that works for all $(t, d)$.

\section{Optimality}\label{optimalitysec}

Armed with the notions of conjugate self-twists and the adjoint trace ring from \cref{cst section}, we prove two related optimality results. We show that the adjoint trace ring of a pseudorepresentation $(t,d)$ contains a fullness peer of any fullness ring for $(t,d)$ (\cref{A0optimalitythm}). In this way we establish the adjoint trace ring as the optimal fullness ring up to fullness peerage. We also show that every fullness ring is fixed by all conjugate self-twists (\cref{fullfixed}), generalizing \cref{A0fullsimpleCST}. 

As usual, $A$ is a local noetherian pro-$p$ domain carrying a continuous pseudorepresentation \mbox{$(t, d): \Pi \to A$} of a $p$-finite profinite group $\Pi$. We assume that $A$ is a domain with field of fractions $K$. Also let $A_0$ and $K_0$ be the adjoint trace ring of $(t, d)$ (\cref{a0def}) and its fraction field, respectively. Recall that $(t, d)$ is $B$-full for any subring $B$ of $K$ if there exists a $(t,d)$-representation $\rho:\Pi \to R^\times$ and an embedding $R \into M_2(K)$ such that via this embedding $\im \rho \supseteq \Gamma_B(\Bb)$ for some nonzero ideal $\Bb$ of $B$ (\cref{fullness}).

We first show that fullness rings are fixed by conjugate self-twist automorphisms whose characters are continuous, or at least nearly so.
\begin{proposition}\label{conditional fixed by cst}\label{closednormal3}
Suppose $(t, d)$ is $B$-full for some subring $B$ of $A$.  Let~$(\sigma, \eta)$ be a conjugate self-twist of $(t, d)$ valued in any domain extending $A$. If $\ker \eta$ contains a subgroup $H$, closed and normal in $\Pi$, such that $\Pi/H$ is abelian, then $\sigma$ fixes $B$ pointwise.

In particular, if $(\sigma, \eta)$ is $A$-valued, then $\sigma$ always fixes $B$ pointwise.
\end{proposition}

\begin{proof}
By \cref{fullclosedprop} we may replace $B$ by its closure in $A$. On one hand, since $H \subseteq \ker \eta$, we have that $t(H)$ is fixed by $\sigma$. On the other hand, by \cref{closednormal1}, $(\left.t\right|_H, \left.d\right|_H)$ is still $B$-full, so that by \cref{closednormal2} applied to $(\left.t\right|_H, \left.d\right|_H)\colon H \to A$, the $\Z$-span of $t(H)$ contains some nonzero $B$-ideal $\Bb$. Therefore every element of $\Bb$ is fixed by $\sigma$. Since any element of $B$ is a ratio of elements of $\Bb$ (if $b \neq 0$ is in $\Bb$, then  $x = \frac{xb}{b}$ for any $x \in B$), every element of $B$ is fixed by $\sigma$. 

If $(\sigma, \eta)$ is $A$-valued, then $\eta$ is continuous (\cref{a0all}\eqref{etafinitecont}) so that we can take $H = \ker \eta$.
\end{proof}

\begin{corollary}\label{A0 optimality}
If $(t, d)$ is $B$-full for some subring $B$ of $K$, then $B \subseteq K_0$.
\end{corollary} 
\begin{proof}
Replace $(t, d)$ with its constant-determinant twist (\cref{twisting does not affect fullness,a0fixedbytwist}) and let $A_t$ be its trace algebra. Replace $B$ by its fullness peer $B \cap A_t$ (\cref{closedinAt,fgqf}). The result now follows from \cref{closednormal3} using \cref{KgaloisK0}. \end{proof}

Our main optimality theorem is an improvement on \cref{A0 optimality}: any fullness ring is not only contained in $K_0$, it in fact has a fullness peer subring contained in $A_0$.

\begin{theorem}\label{A0optimalitythm}
Suppose $(t, d)$ is $B$-full for some subring $B$ of $K$. Then $B \cap A_0$ is a fullness peer of $B$ contained in $A_0$. 
\end{theorem}

\begin{proof}
First suppose that $B$ is finite. Then $B$ is a finite field; its only nonzero ideal is itself, and the only congruence subgroup of $\SL_2(B)$ is $\SL_2(B)$ itself. In this case, $B$ has no proper subring fullness peers and we show directly that $B$ is contained in $A_0$. Indeed, $K$, and hence $A$, is now an $\F_p$-algebra, and it follows that $\rhobar$ is also $B$-full. Now \cref{closednormal3} applied to $\rhobar$ tells us that $B$ is fixed by all conjugate self-twist automorphisms of $\rhobar$: in other words, $B \subseteq \E$ (\cref{residual csts cut out E}), which is the residue field of, and hence here contained in, $A_0$.

Now assume that $B$ is not finite. Use \cref{twisting does not affect fullness} to replace $(t,d)$ by its constant-determinant twist and let $A_t \subseteq A$ be its trace algebra. Replace $B$ by its fullness peer $B \cap A_t$ (\cref{twisting does not affect fullness,a0fixedbytwist}). Let $\Gamma = \ker \rhobar$, a finite-index closed subgroup of $\Pi$. By \cref{finiteindex}, the restriction $\left(\left.t \right|_{\Gamma}, \left.d \right|_\Gamma\right)$ is $B$-full. By \cref{closednormal2}, $A_0$, which contains the trace algebra of $\Gamma$ (\cref{tgammainA0}), contains a nonzero ideal of~$B$. Therefore so does $B \cap A_0$, making $B \cap A_0 \subseteq B$ an extension of fullness peers by \cref{fg plus same quotient field implies ideal containment}.
\end{proof}
\noindent We do not know whether we can sharpen \cref{A0optimalitythm} to conclude that any closed fullness ring $B$ contained in $A_0$ with $Q(B) = K_0$ must in fact be a fullness peer of $A_0$. This would follow from an affirmative answer to \cref{transporter}.

We close with an observation that fullness rings are always fixed by conjugate self-twist automorphisms. \cref{A0fullsimpleCST} has already established this for $A_0$; \cref{fullfixed} below is a generalization.

\begin{theorem}  \label{fullfixed}
Let $(t, d)\colon \Pi \to A$ be a pseudorepresentation. If $(t, d)$ is $B$-full for some subring $B$ of 
$K$, then $B$ is fixed by any conjugate self-twists valued in any 
domain extending $A$ containing~$B$. 
\end{theorem} 

\noindent Note that any automorphism valued in a domain $E$ containing $A$ extends uniquely to an automorphism of $Q(E)$, which contains $K$ and hence $B$.

\begin{proof} 
Let $(t', d')$ be the constant determinant twist associated to $(t,d)$, obtained by twisting off the pro-$p$ part $d_1$ of $d$. As in \cref{cstaftertwist}, if $(\sigma, \eta)$ is an $E$-valued conjugate self-twists $(t, d)$ for some domain $E$ extending $A$, 
then $\big(\sigma,\ \sigma(d_1^{-1/2})d_1^{1/2}\eta\big)$ is a $E$-valued conjugate self-twist of $(t',d')$.  
By \cref{atisbigenough}, $\chi = \sigma(d_1^{-1/2})d_1^{1/2}\eta$ is continuous and $W(\F)$-valued. 
Note that \mbox{$\ker d_1 = \ker d_1^{1/2} = \ker \sigma(d_1^{-1/2})$} is closed since $d_1$ is also continuous.  Thus $\ker d_1 \cap \ker \chi$ is a closed subgroup of $\ker \eta$.  Moreover, since the order of $\chi$ is prime to $p$ and $d_1$ is pro-$p$ we get that
\[
\Pi/(\ker d_1 \cap \ker \chi) \cong \Pi/\ker d_1 \times \Pi/\ker \chi
\]
is abelian.  Thus we can take $H = \ker d_1 \cap \ker \chi$ in \cref{conditional fixed by cst} to deduce that $B \cap A$ 
is fixed by~$\sigma$. As $B \cap A$ and $B$ are fullness peers, $\sigma$ thus fixes an ideal of $B$, hence all of $B$. 
\end{proof}

\medskip


\section{Fullness for $\B_\rho(\E)$}
\label{rings acting}

Throughout this section we fix a local pro-$p$ ring $A$, not necessarily a domain, with residue field~$\F$.  Fix an admissible pseudodeformation $(\Pi, \ovl{\rho}, t, d)$ over $A$.  Let $A_0$ be the adjoint trace ring of $(t, d)$  and $\E$ the residue field of $A_0$. 

\subsection{$L_2(\rho)$ is a $W(\E)$-module}\label{section action of E}
Recall that, a priori, Pink's construction only gives Lie algebras that are $\Z_p$-modules.  The goal of \cref{section action of E} is to show that if $\rho$ is a $(t,d)$-representation, then in fact its associated Lie algebras are modules over $W(\E)$ (\cref{extending constant scalars}).  Although this is a minor improvement on $\Z_p$ (indeed, it is no improvement at all if $W(\E) = \Z_p$), it is an essential input for proving the results of \cref{conceptual interpretation of Bellaiche ring}.

We assume throughout \cref{section action of E} that the eigenvalues of $\ovl{\rho}(g)$ are in $\F^\times$ for all $g \in \Pi$.  This requires at most replacing $\F$ by its unique quadratic extension.

Let $\lambda \neq \mu \in \F^\times$ be the eigenvalues of a matrix in $\im \ovl{\rho}$.  By \cref{better adapted td-rep}, there is a $(t,d)$-representation $\rho_{\lambda, \mu} \colon \Pi \to R_{\lambda, \mu}^\times$ and $g_{\lambda, \mu} \in \Pi$ such that 
\[
\rho_{\lambda, \mu}(g_{\lambda, \mu}) = \begin{pmatrix} 
s(\lambda) & 0\\
0 & s(\mu)
\end{pmatrix}.
\]  
Recall that $G_{\rho_{\lambda, \mu}} \coloneqq \im \rho_{\lambda, \mu}, \Gamma_{\rho_{\lambda, \mu}} \coloneqq G_{\rho_{\lambda, \mu}} \cap SR_{\lambda, \mu}^1$, and $L_n(\rho_{\lambda, \mu}) \coloneqq L_n(\Gamma_{\rho_{\lambda, \mu}})$.  Since $\lambda \neq \mu$, the Lie algebra $L_1(\rho_{\lambda, \mu})$ is decomposable \cite[Corollary 6.2.2]{Bellaiche18}.  Note that although the Teichm\"uller map $s$ is not additive, it is easy to check that $W(\F_p(\lambda\mu^{-1} + \lambda^{-1}\mu)) = \Z_p[s(\lambda)s(\mu)^{-1} + s(\lambda)^{-1}s(\mu)]$, a fact that we make frequent use of in \cref{nabla is a module} below.

\begin{lemma}\label{nabla is a module} 
With the notation introduced at the beginning of \cref{bellaiche results background}, we have
\begin{enumerate}
\item $\nabla_1(\rho_{\lambda, \mu}), B_1(\rho_{\lambda, \mu}), C_1(\rho_{\lambda, \mu}),$ and $L_2(\rho_{\lambda, \mu})$ are $W(\F_p(\lambda\mu^{-1} + \lambda^{-1}\mu))$-modules;
\item if the projective image of $\ovl{\rho}$ contains $\PSL_2(\F_p)$ and $p \geq 7$, then $I_1(\rho_{\lambda, \mu})$ is a \mbox{$W(\F_p(\lambda\mu^{-1} + \lambda^{-1}\mu))$}-module; after possibly replacing $\rho_{\lambda, \mu}$ with its conjugate by a certain $\left(\begin{matrix} 
1 & 0\\
0 & a
\end{matrix}\right)$ with $a \in A^\times$, one has that $L_1(\rho_{\lambda, \mu})$ is a $W(\F_p(\lambda\mu^{-1} + \lambda^{-1}\mu))$-module.
\end{enumerate}
\end{lemma}

\begin{proof}
Note that
\[
L_2(\rho_{\lambda, \mu}) = \left[I_1(\rho_{\lambda, \mu})\left(\begin{smallmatrix} 
1 & 0\\
0 & -1
\end{smallmatrix}\right), \nabla_1(\rho_{\lambda, \mu})\right] + [\nabla_1(\rho_{\lambda, \mu}), \nabla_1(\rho_{\lambda, \mu})].
\]
Furthermore, from their definitions on page~\pageref{B1C1}, it's clear that $B_1(\rho_{\lambda, \mu})$ and $C_1(\rho_{\lambda, \mu})$ inherit whatever structure $\nabla_1(\rho_{\lambda, \mu})$ has.
  Therefore it suffices to show that $\nabla_1(\rho_{\lambda, \mu})$ is a $W(\F_p(\lambda\mu^{-1} + \lambda^{-1}\mu))$-module.

To prove that
\[
(s(\lambda)s(\mu)^{-1} + s(\lambda)^{-1}s(\mu))\nabla_1(\rho_{\lambda, \mu}) \subseteq \nabla_1(\rho_{\lambda, \mu}),
\]
recall that $L_1(\rho_{\lambda, \mu})$ is closed under conjugation by $G_{\rho_{\lambda, \mu}}$ (in fact, by any element in the normalizer of $\Gamma_{\rho_{\lambda, \mu}}$).  In particular, it is closed under conjugation by $\left(\begin{smallmatrix} 
s(\lambda) & 0\\
0 & s(\mu)
\end{smallmatrix}\right)$ and $\left(\begin{smallmatrix} 
s(\lambda)^{-1} & 0\\
0 & s(\mu)^{-1}
\end{smallmatrix}\right)$.  Using this, a short matrix calculation shows that if $\bigl(\begin{smallmatrix} 
0 & b\\
c & 0
\end{smallmatrix}\bigr) \in \nabla_1(\rho_{\lambda, \mu})$ then 
\[
\left(\begin{smallmatrix} 
0 & s(\lambda)s(\mu)^{-1}b\\
s(\lambda)^{-1}s(\mu)c & 0
\end{smallmatrix}\right), \left(\begin{smallmatrix} 
0 & s(\lambda)^{-1}s(\mu)b\\
s(\lambda)s(\mu)^{-1}c & 0
\end{smallmatrix}\right) \in \nabla_1(\rho_{\lambda, \mu}).
\]
Therefore $(s(\lambda)s(\mu)^{-1} + s(\lambda)^{-1}s(\mu))\nabla_1(\rho_{\lambda, \mu}) \subseteq \nabla_1(\rho_{\lambda, \mu})$.

Finally, if the projective image of $\ovl{\rho}$ contains $\PSL_2(\F_p)$ and $p \geq 7$, then by \cref{Bellaiche main theorem}, up to replacing $\rho_{\lambda, \mu}$ with its conjugate by a certain $\left(\begin{smallmatrix} 
1 & 0\\
0 & a
\end{smallmatrix}\right)$ with $a \in A^\times$, we have
\[
L_1(\rho_{\lambda, \mu}) = \left(\begin{smallmatrix} 
I_1(\rho_{\lambda, \mu}) & I_1(\rho_{\lambda, \mu})\\
I_1(\rho_{\lambda, \mu}) & I_1(\rho_{\lambda, \mu})
\end{smallmatrix}\right)^0,
\]
and thus $B_1(\rho_{\lambda, \mu}) = I_1(\rho_{\lambda, \mu}) = C_1(\rho_{\lambda, \mu})$.  (Note that conjugation by $\left(\begin{smallmatrix} 
1 & 0\\
0 & a
\end{smallmatrix}\right)$ with $a \in A^\times$ does not change $I_1(\rho_{\lambda, \mu})$.)  In particular, $L_1(\rho_{\lambda, \mu})$ is strongly decomposable.  By the second statement of the lemma, we see that $I_1(\rho_{\lambda, \mu})$ is a $W(\F_p(\lambda\mu^{-1} + \lambda^{-1}\mu))$-module.  The above description of $L_1(\rho_{\lambda, \mu})$ shows that it is also a $W(\F_p(\lambda\mu^{-1} + \lambda^{-1}\mu))$-module.
\end{proof}

\begin{proposition}\label{extending constant scalars}
Let $\rho \colon \Pi \to R^\times$ be a $(t, d)$-representation.  Then $L_n(\rho)$ is a $W(\E)$-module for all $n \geq 2$.  If the projective image of $\ovl{\rho}$ contains $\PSL_2(\F_p)$ and $p \geq 7$, then $L_1(\rho)$ is a $W(\E)$-module.
\end{proposition}

\begin{proof}
By the generating set for $\E$ given in \eqref{generating set for E} in \cref{regularity background}, it suffices to show that $L_n(\rho)$ is closed under multiplication by $s(\lambda)s(\mu)^{-1} + s(\lambda)^{-1}s(\mu)$ for all $\lambda, \mu \in \ovl{\F}_p^\times$ that are distinct eigenvalues of an element in $\im \ovl{\rho}$.  Fix such $\lambda, \mu$.  Let $\rho_{\lambda, \mu} \colon \Pi \to R_{\lambda, \mu}^\times$ be the $(t,d)$-representation over $A$ described prior to \cref{nabla is a module}.  Let us assume furthermore that, in the case when $\ovl{\rho}$ is not projectively cyclic or dihedral, that we have already replaced $\rho_{\lambda, \mu}$ by its relevant diagonal conjugate so that the description of $W(\F_q)L_1(\rho_{\lambda, \mu})$ from \cref{Bellaiche main theorem} applies to $\rho_{\lambda, \mu}$.

Since $\rho \colon \Pi \to R^\times$ and $\rho_{\lambda, \mu} \colon \Pi \to R_{\lambda, \mu}^\times$ are both $(t,d)$-representations over $A$, it follows from \cref{Bellaiche 2.4.2} that there is a unique $A$-algebra isomorphism $\Psi \colon R_{\lambda, \mu} \to R$ such that $\rho = \Psi \circ \rho_{\lambda, \mu}$.  We claim that this implies that $L_n(\rho) = \Psi(L_n(\lambda, \mu))$ for all $n \geq 1$.  If this is true, then $L_2(\rho)$ is closed under multiplication by $s(\lambda)s(\mu)^{-1} + s(\lambda)^{-1}s(\mu) \in A$ since $L_2(\rho_{\lambda, \mu})$ is by \cref{nabla is a module} and $\Psi$ is an $A$-algebra homomorphism.  Since $L_2(\rho)$ is a $W(\E)$-module, it follows immediately from the definition that $L_n(\rho)$ is a $W(\E)$-module for all $n \geq 2$.  Furthermore, the argument in this paragraph applies to $L_1(\rho)$ under the assumption that the projective image of $\ovl{\rho}$ contains $\PSL_2(\F_p)$ for $p \geq 7$.

To see that $L_n(\rho) = \Psi(L_n(\rho_{\lambda, \mu}))$, note that $G_\rho = G_{\Psi\circ \rho_{\lambda, \mu}} = \Psi(G_{\rho_{\lambda, \mu}})$.  Since $\Psi$ is an algebra morphism, it follows that $\Psi(\rad R_{\lambda, \mu}) = \rad R$.  Furthermore, since $\rho$ and $\rho_{\lambda, \mu}$ are both $(t, d)$-representations, it follows that $\Psi$ preserves determinants.  Therefore $\Psi(SR_{\lambda, \mu}^1) \supset SR^1$.  Since $\Psi$ is a continuous algebra homomorphism, it follows directly from the definition of $\Theta$ that $\Psi(L_1(\rho_{\lambda, \mu})) = L_1(\rho)$ and hence $\Psi(L_n(\rho_{\lambda, \mu})) = L_n(\rho)$ for all $n \geq 1$.
\end{proof}

\subsection{$L_2(\rho)$ is a $\B_\rho(\E)$-module}\label{L1 ABmodule}
In \cref{L1 ABmodule} we use Bella\"iche's work to show that, for any well-adapted $(t,d)$-representation $\rho$, $L_n(\rho)$ is a module over a ring comparable to $A$.  This is the key input into \cref{formal fullness corollary}, which is our improvement on Bella\"iche's fullness theorem.

\begin{proposition}\label{the biggest ring for L2}
Let $\rho$ be a $(t,d)$-representation adapted to $(g_0, \lambda_0, \mu_0)$.  Then
\begin{enumerate}
\item $L_2(\rho)$ is a module over $W(\E)[I_1(\rho)^2] \coloneqq W(\E) + W(\E)I_1(\rho)^2$;  
\item if $n \geq 1$ and $L_n(\rho)$ is strongly decomposable, then $L_{n+1}(\rho)$ is a module over $\B_\rho(\E)$;
\item if the projective image of $\ovl{\rho}$ contains $\PSL_2(\F_p)$ for $p \geq 7$, then up to replacing $\rho$ with its conjugate by some $\left(\begin{smallmatrix} 
1 & 0\\
0 & a
\end{smallmatrix}\right)$ with $a \in A^\times$, $L_1(\rho)$ is a module over $\B_\rho(\E)$.
\end{enumerate}
\end{proposition}

\begin{proof}
Since $\rho$ is adapted to $(g_0, \lambda_0, \mu_0)$, it follows that $L_1(\rho)$ is decomposable \cite[Corollary 6.2.2]{Bellaiche18}.  Note that $[I_1(\rho)\bigl(\begin{smallmatrix} 1 & 0\\ 0 & -1
\end{smallmatrix}\bigr), \nabla_1(\rho)] \subset \nabla_1(\rho)$ since $L_1(\rho)$ is a Lie algebra.  That is, for all $a \in I_1(\rho)$ and $\bigl(\begin{smallmatrix} 
0 & b\\
c & 0
\end{smallmatrix}\bigr) \in \nabla_1(\rho)$, we have $2a\bigl(\begin{smallmatrix} 
0 & b\\
-c & 0
\end{smallmatrix}\bigr) \in \nabla_1(\rho)$.  To prove the first statement, we can apply this fact a second time to $\alpha \in I_1(\rho)$ and $2a\bigl(\begin{smallmatrix}
0 & b\\
-c & 0
\end{smallmatrix} \bigr)$ to see that $4a\alpha\bigl(\begin{smallmatrix} 
0 & b\\
c & 0
\end{smallmatrix}\bigr) \in \nabla_1(\rho)$.  Therefore $\nabla_1(\rho)$ is closed under multiplication by $I_1(\rho)^2$.  Since
\[
L_2(\rho) = \left[I_1(\rho)\left(\begin{smallmatrix} 
1 & 0\\
0 & -1
\end{smallmatrix}\right), \nabla_1(\rho)\right] + [\nabla_1(\rho), \nabla_1(\rho)],
\]
we see that $L_2(\rho)$ is closed under multiplication by $I_1(\rho)^2$.

For the second statement, if $L_n(\rho)$ is strongly decomposable, then we can write
\[
L_n(\rho) = I_n(\rho)\bigl(\begin{smallmatrix} 
1 & 0\\
0 & -1
\end{smallmatrix}\bigr) \oplus B_n(\rho)\bigl(\begin{smallmatrix} 
0 & 1\\
0 & 0
\end{smallmatrix}\bigr) \oplus C_n(\rho)\bigl(\begin{smallmatrix} 
0 & 0\\
1 & 0
\end{smallmatrix}\bigr).
\]
By calculating $\left[\left(\begin{smallmatrix} 
1 & 0\\
0 & -1
\end{smallmatrix}\right), \left(\begin{smallmatrix} 
0 & 1\\
0 & 0
\end{smallmatrix}\right)\right]$ and $\left[\left(\begin{smallmatrix} 
1 & 0\\
0 & -1
\end{smallmatrix}\right), \left(\begin{smallmatrix} 
0 & 0\\
1 & 0
\end{smallmatrix}\right) \right]$, we find that $I_1(\rho)B_n(\rho) = B_{n+1}(\rho) \subset B_n(\rho)$ and $I_1(\rho)C_n(\rho) = C_{n+1}(\rho) \subset C_n(\rho)$.  Therefore $B_n(\rho), C_n(\rho)$ are closed under multiplication by $I_1(\rho)$.  Since
\begin{gather*}
L_{n+1}(\rho) = [B_n(\rho)\bigl(\begin{smallmatrix} 
0 & 1\\
0 & 0
\end{smallmatrix}\bigr), C_n(\rho)\bigl(\begin{smallmatrix} 
0 & 0\\
1 & 0
\end{smallmatrix}\bigr)] + [I_1(\rho)\bigl(\begin{smallmatrix} 
1 & 0\\
0 & -1
\end{smallmatrix}\bigr), B_n(\rho)\bigl(\begin{smallmatrix} 
0 & 1\\
0 & 0
\end{smallmatrix}\bigr)] + \\
+ [I_1(\rho)\bigl(\begin{smallmatrix} 
1 & 0\\
0 & -1
\end{smallmatrix}\bigr), C_n(\rho)\bigl(\begin{smallmatrix} 
0 & 0\\
1 & 0
\end{smallmatrix}\bigr)],
\end{gather*}
it follows that $L_{n+1}(\rho)$ is closed under multiplication by $I_1(\rho)$.

The first two results now follow from \cref{extending constant scalars}.  The last statement follows from \cref{Bellaiche main theorem} and \cref{extending constant scalars}.
\end{proof}

\begin{remark}
It would be nice to remove the assumption that $L_1(\rho)$ is strongly decomposable and still conclude that $L_2(\rho)$ is a $\B_\rho(\E)$-module, but we do not see a way to do this. 
\end{remark}

\subsection{Regularity implies $\B_\rho(\E)$-fullness}\label{formal fullness section}
The goal of \cref{formal fullness section} is to establish a slightly stronger version of \cite[Theorem 7.2.3]{Bellaiche18}, which is Bella\"iche's \cref{Bellaiche thm intro} of the introduction.  We do so in \cref{formal fullness corollary} below.  Our result is different from that of Bella\"iche mainly in that we can weaken his definition of regularity and enlarge his ring $\B_\rho(\F_p)$ to $\B_\rho(\E)$.   

Throughout \cref{formal fullness section} the ring $A$ will be a local pro-$p$ domain with residue field $\F$ and field of fractions $K$.  We fix an admissible pseudodeformation $(\Pi, \ovl{\rho}, t, d)$ over $A$ throughout this section.  If $\rho$ is a $(t,d)$-representation that is adapted to some $(g_0, \lambda_0, \mu_0)$, then $L_1(\rho)$ is decomposable by \cite[Corollary 6.2.2]{Bellaiche18}.  Thus $I_1(\rho)$ is defined.  We write $K_1$ for the field of fractions of $\B_\rho(\E)$.

\begin{proposition}\label{formal fullness}
Assume that $\ovl{\rho}$ is regular.  Let $\rho \colon \Pi \to R^\times$ be a $(t,d)$-representation adapted to $(g_0, \lambda_0, \mu_0)$ for a regular element $g_0$ such that $\rho(g_0) = \left(\begin{smallmatrix} 
s(\lambda_0) & 0\\
0 & s(\mu_0)
\end{smallmatrix}\right)$.  If $B_1(\rho), C_1(\rho) \neq 0$, then $\rho$ is $\B_\rho(\E)$-full.
\end{proposition}

\begin{proof}
It is easy to see that the eigenvalues of $\rho(g_0) = \left(\begin{smallmatrix} 
s(\lambda_0) & 0\\
0 & s(\mu_0)
\end{smallmatrix}\right)$ acting on $L_n(\rho)$ by conjugation are $1, s(\lambda_0)s(\mu_0^{-1}), s(\lambda_0^{-1})s(\mu_0)$, which are distinct elements of $W(\E)^\times$ since $g_0$ is a regular element.  Since $L_n(\rho)$ is a $W(\E)$-module for $n \geq 2$ by \cref{extending constant scalars}, it follows that $L_n(\rho)$ is the direct sum of the eigenspaces for the conjugation action of $\rho(g_0)$.  Thus, $L_n(\rho)$ is strongly decomposable for $n \geq 2$.  By \cref{the biggest ring for L2}, it follows that $L_n(\rho)$ is an $\B_\rho(\E)$-module for $n \geq 3$.

Since $A$ is a domain, we may view $R$ inside of $M_2(K)$ by \cref{GMAs over domains}.  Note that if \mbox{$B_1(\rho), C_1(\rho) \neq 0$}, then since $I_n(\rho), B_n(\rho), C_n(\rho) \subset K$, it follows that $I_n(\rho), B_n(\rho)$, and $C_n(\rho)$ are nonzero for all $n \geq 1$.  In particular, $I_3(\rho), B_3(\rho), C_3(\rho)$ are nonzero $\B_\rho(\E)$-modules.   

Define
\[ 
R_1 \coloneqq \left(\begin{smallmatrix} 
\B_\rho(\E) & B_3(\rho)\\
C_3(\rho) & \B_\rho(\E)
\end{smallmatrix}\right).
\]
Then $R_1$ is a faithful GMA over $\B_\rho(\E)$.  By the proof of \cite[Lemma 2.2.2]{Bellaiche18}, if $0 \neq b_0 \in B_3(\rho)$ and $x = \bigl(\begin{smallmatrix} 
1 & 0\\
0 & b_0
\end{smallmatrix}\bigr)$, it follows that $xR_1x^{-1} \subseteq \GL_2(K_1)$.  Thus, by replacing $\rho$ with $x\rho x^{-1}$, which is still a $(t,d)$-representation adapted to $(g_0, \lambda_0, \mu_0)$ that sends $g_0$ to $\left(\begin{smallmatrix}
s(\lambda_0) & 0\\
0 & s(\mu_0)
\end{smallmatrix}\right)$, we may assume that $B_3(\rho), C_3(\rho) \subseteq K_1$.  (Note that $I_1(\rho) = I_1(x\rho x^{-1})$.)  

Note that any nonzero $\B_\rho(\E)$-submodule of $K_1$ contains a nonzero element of $\B_\rho(\E)$ and thus contains a non-zero $\B_\rho(\E)$-ideal.  Therefore there exists a nonzero $\B_\rho(\E)$-ideal $\Bb$ contained in $I_3(\rho) \cap B_3(\rho) \cap C_3(\rho)$.  Hence $\Sl_2(\Bb) \subseteq L_3(\rho)$.  Using \cref{Pink main thm} we deduce that $\Gamma_{\B_\rho(\E)}(\Bb) \subset \im \rho$ and $\rho$ is $\B_\rho(\E)$-full.
\end{proof}

\begin{corollary}\label{formal fullness corollary}
Assume that $\ovl{\rho}$ is regular.  Let $\rho$ be a well-adapted $(t,d)$-representation adapted to $(g_0, \lambda_0, \mu_0)$ for a regular element $g_0$ such that $\rho(g_0) = \bigl(\begin{smallmatrix}s(\lambda_0) & 0\\
0 & s(\mu_0) \end{smallmatrix} \bigr)$.  If $(t, d)$ is not a priori small, then $(t, d)$ is $\B_\rho(\E)$-full.
\end{corollary}

\begin{proof}
By \cref{formal fullness} it suffices to show that $B_1(\rho), C_1(\rho) \neq 0$.  We do this by analyzing the different possibilities for $\ovl{\rho}$.  
By \cref{dihedral or ad irreducible}, we see that either $\ovl{\rho}$ is reducible, dihedral, or $\ad^0 \ovl{\rho}$ is irreducible. Assume first that we are in the last case. The group $\Gamma$ is equipped with a decreasing normal filtration 
\[
\Gamma_n \coloneqq \Gamma \cap \Gamma_A(\m^n)
\]
whose quotients $\Gamma_n/\Gamma_{n+1}$ have the structure of $\F_p$-vector spaces with an action of $\ovl G$ by conjugation.  Moreover, the map $\gamma \mapsto \gamma - 1$ gives a $\ovl G$-equivariant embedding $\Gamma_n/\Gamma_{n+1} \hookrightarrow \Sl_2(\m^n/\m^{n+1})$; write $V_n$ for its image. The assumption that $(t, d)$ is not a priori small implies that $\Gamma$ is nontrivial, so that $\Gamma_n/\Gamma_{n +1}$, and hence $V_n$, is nontrivial for some $n \geq 1$.  To see that $B_1(\rho)$ (respectively, $C_1(\rho)$) is nonzero, it suffices to show that this $V_n$ contains an element whose upper right (respectively, lower left) entry is nonzero. This can be checked on the $\F$-span of $V_n$. Choosing an $\F$-basis $x_1, \ldots, x_d$ of $\m^n/\m^{n+1}$ gives a $\ovl{G}$-equivariant splitting $\Sl_2(\m^n/\m^{n+1}) = \oplus_{i = 1}^d \Sl_2(\F)x_i$.  Since $V_n$ is nontrivial, there is some $i$ such that the projection $W$ of $\F V_n$ to $\Sl_2(\F)x_i$ is nonzero.  Then $W$ is a stable subspace of $\Sl_2(\F)x_i$, which is simple since $\ad^0 \ovl{\rho}$ is irreducible.  Thus $W = \Sl_2(\F)x_i$ and $W$, and hence $\F V_n$, contains an element whose upper right (respectively, lower left) entry is nonzero.\footnote{Together with the observation that $B_n(\rho), C_n(\rho) \neq 0$ for all $n$ if it is true for $n = 1$ from \cref{formal fullness}, the argument here gives an independent proof of the \eqref{finim} $\implies$ \eqref{aps} implication in \cref{apssai}.}

Now suppose that $\ovl{\rho}$ is reducible.  Since $\rho$ is well adapted by assumption, it follows that $\rho$ is adapted to $(g_0, \lambda_0, \mu_0)$, where $\ovl{\rho}(g_0)$ generates the projective image of $\ovl{\rho}$.  In particular, $\rho$ is automatically adapted to a regular element.  Suppose for contradiction that $C_1(\rho) = 0$ (respectively, $B_1(\rho) = 0$).  Then $\Gamma_\rho$ is contained in the upper (respectively, lower) triangular matrices.  By \cite[Theorem 6.2.1]{Bellaiche18}, we know that $s(\ovl{G}) \subset G_\rho$ since $\rho$ is well adapted.  Thus $G_\rho = s(\ovl{G})\Gamma_\rho$. But then $G_\rho$ is contained in the upper (respectively, lower) triangular matrices, and hence $\rho$ is reducible.  Therefore $t$ is the sum of two continuous characters $\Pi \to A^\times$, which contradicts our assumption that it is not a priori small.  Thus $B_1(\rho), C_1(\rho) \neq 0$ if $\ovl{\rho}$ is reducible.

Finally suppose that $\ovl{\rho}$ is dihedral.  By \cref{resolve definitions of dihedral} there is a unique subgroup $\Pi_0$ of index 2 in $\Pi$ such that $\ovl{\rho} \cong \Ind_{\Pi_0}^\Pi \chi$ for some character $\chi \colon \Pi_0 \to \F^\times$.  Applying the reducible case to $\ovl{\rho}|_{\Pi_0}$, we see that either $\rho|_{\Pi_0}$ is reducible or $B_1(\rho\vert_{\Pi_0}), C_1(\rho\vert_{\Pi_0}) \neq 0$.  The first possibility is not allowed by hypothesis, so we must have $B_1(\rho|_{\Pi_0}), C_1(\rho|_{\Pi_0}) \neq 0$.  But $B_1(\rho|_{\Pi_0}) \subseteq B_1(\rho)$ and $C_1(\rho|_{\Pi_0}) \subseteq C_1(\rho)$, which proves the desired result when $\ovl{\rho}$ is dihedral.
\end{proof}

\medskip


\section{Lifting residual conjugate self-twists}\label{section: lifting}

In \cref{section: lifting} we study the ($A_t$-valued) conjugate self-twists of constant-determinant pseudorepresentations.  In particular, we show in \cref{section lifting csts} that they are all controlled by those of $\ovl \rho$, and all residual conjugate self-twists lift to the universal constant-determinant pseudodeformation ring.  Having shown in \cref{section: res csts} that the group of residual conjugate self-twists is finite and abelian, we deduce the same for any constant determinant pseudorepresentation.  Finally, in \cref{residual trivial csts} we study special phenomena that arise when $\ovl \rho$ is dihedral and hence there may be conjugate self-twists that are residually dihedral.

Throughout \cref{section: lifting} we only consider simple conjugate self-twists.  Thus we write $\til \Sigma_t = \til \Sigma_t(A_t)$ and $\til \Sigma_{\ovl \rho}$ for $\til \Sigma_{\ovl \rho}(\F)$; similarly for $\Sigma_t$ and $\Sigma_{\ovl \rho}$. 

\subsection{Residual conjugate self-twists}\label{section: res csts}
Fix a semisimple representation $\ovl{\rho} \colon \Pi \to \GL_2(\F)$, and recall that $\proj \colon \GL_2(\F) \to \PGL_2(\F)$ denotes the natural projection.  For \cref{section: res csts} only, we do not require $\ovl \rho$ to be residually multiplicity-free.  
We begin by studying $\Sigma_{\ovl{\rho}}^{\di}$ and use that to show that $\til \Sigma_{\ovl{\rho}}$ is finite and abelian. 

\begin{lemma}\label{precise structure of Sigma}\mbox{}
\begin{enumerate}
\item If $\im \proj \ovl{\rho}$ is not dihedral or cyclic of order 2, then $\Sigma_{\ovl{\rho}}^{\di}$ is trivial.  
\item If $\im \proj \ovl{\rho}$ is either a nonabelian dihedral group or has order 2, then $\Sigma_{\ovl{\rho}}^{\di}$ has order 2. 
\item If $\im \proj \ovl{\rho}$ is isomorphic to $(\Z/2\Z)^2$, then $\Sigma_{\ovl{\rho}}^{\di}$ is isomorphic to $(\Z/2\Z)^2$.  
\end{enumerate}
\end{lemma}

\begin{proof}
We claim that if $\ovl{\rho}$ is irreducible, then the following sets are in bijection:
\begin{enumerate}[label=(\alph{enumi})]
\item\label{Sigma0} $\Sigma_{\ovl{\rho}}^{\di} \setminus \{(1,1)\}$;
\item\label{subs of Pi} subgroups $\Pi_0$ of $\Pi$ such that $[\Pi \colon \Pi_0] = 2$ and $\ovl{\rho}(\Pi_0)$ is abelian;
\item\label{subs of Prhobar} subgroups $H$ of $\proj\ovl{\rho}(\Pi)$ such that $[\proj\ovl{\rho}(\Pi) \colon H] = 2$ and $H$ is abelian.
\end{enumerate}
Indeed, the maps between them can be described as follows.  Given $(1, \eta) \in \Sigma_{\ovl{\rho}}^{\di} \setminus \{(1,1)\}$, let $\Pi_0 \coloneqq \ker \eta$.  The fact that $[\Pi \colon \Pi_0] = 2$ follows from \cref{resolve definitions of dihedral}, and \cref{cm cst implies eta quadratic} shows that $\ovl{\rho}(\Pi_0)$ is abelian.  Conversely, given $\Pi_0$ as in \ref{subs of Pi}, let $\eta_{\Pi_0} \colon \Pi \to \Pi/\Pi_0 \cong \{\pm 1\}$ be the natural projection.  Note that $\ovl{\rho}|_{\Pi_0}$ is reducible since $\ovl{\rho}(\Pi_0)$ is abelian.  Let $\chi \colon \Pi_0 \to \F^\times$ be a constituent of $\ovl{\rho}|_{\Pi_0}$.  Then $\ovl{\rho} \cong \Ind_{\Pi_0}^\Pi \chi$ by Frobenius reciprocity since $\ovl{\rho}$ is irreducible.  Thus $(1, \eta_{\Pi_0}) \in \Sigma_{\ovl{\rho}}^{\di} \setminus \{(1,1)\}$ by \cref{resolve definitions of dihedral}.

Given $\Pi_0$ as in \ref{subs of Pi}, let $H \coloneqq \proj\ovl{\rho}(\Pi_0)$.  Given $H$ as in \ref{subs of Prhobar}, let $\Pi_0 \coloneqq \proj\ovl{\rho}^{-1}(H)$.  It is clear that $[\Pi \colon \Pi_0] = 2$.  That $\ovl{\rho}(\Pi_0)$ is abelian follows from the fact that $H$ is abelian and scalar matrices commute with everything.

When $\ovl{\rho}$ is irreducible, the lemma now follows from counting subgroups as in \ref{subs of Prhobar} in each of the possible projective images of $\ovl{\rho}$.  (The fact that elements in $\Sigma_{\ovl{\rho}}^{\di}$ have order at most 2 follows from the fact that $\det \ovl{\rho} = \eta^2 \det \ovl{\rho}$ and so $\eta^2 = 1$.)

Finally, suppose that $\ovl{\rho} = \varepsilon \oplus \delta$.  If $(1, \eta) \in \Sigma_{\ovl{\rho}}^{\di}$ and $\eta$ is nontrivial, then we must have $\eta\varepsilon = \delta$ and $\eta\delta = \varepsilon$.  Thus
\[
\varepsilon\delta^{-1} = \eta = \delta\varepsilon^{-1},
\]
which implies that $\varepsilon\delta^{-1}$ has order 2.  But the projective image of $\ovl{\rho}$ is isomorphic to the image of $\varepsilon\delta^{-1}$.  Thus $\Sigma_{\ovl{\rho}}^{\di}$ is trivial unless the projective image of $\ovl{\rho}$ has order 2, in which case there is one nontrivial element.
\end{proof}

\begin{corollary}\label{tildeSigma rho-bar finite}
The group $\til \Sigma_{\ovl{\rho}}$ is finite and abelian. 
\end{corollary}

\begin{proof}
Since $\F$ is a finite field, there are only finitely many automorphisms of $\F$.  Any two elements in $\til \Sigma_{\ovl \rho}$ with the same automorphism differ by an element of $\Sigma_{\ovl \rho}^{\di}$, which is finite by \cref{precise structure of Sigma}.  

To see that $\til \Sigma_{\ovl \rho}$ is abelian, fix a generator $\sigma$ of the cyclic group $\Sigma_{\ovl \rho} = \Gal(\F/\E)$.  Let $\eta$ be a character such that $(\sigma, \eta) \in \til \Sigma_{\ovl \rho}$.  Then $\til \Sigma_{\ovl \rho}$ is generated by $\{(\sigma, \eta\eta') \colon (1, \eta') \in \Sigma_{\ovl \rho}^{\di}\}$.  Since $\eta'$ is at most quadratic by \cref{precise structure of Sigma}, the action of $\sigma$ on $\eta'$ is trivial and hence one easily checks that any two of the generators commute.
\end{proof}

\subsection{Lifting conjugate self-twists}\label{section lifting csts}
Let $\Pi$ be a profinite group satisfying the $p$-finiteness condition.  Fix a multiplicity-free representation $\ovl{\rho} \colon \Pi \to \GL_2(\F)$.  Recall from \cref{pseudo-reps background} that there is a local pro-$p$ noetherian $W(\F)$-algebra $\cA$ with maximal ideal $\m_{\cA}$ and residue field $\F$ and a pseudodeformation $(T, d) \colon \Pi \to \cA$ that is universal among all constant-determinant pseudodeformations of $\ovl{\rho}$.  The purpose of \cref{section lifting csts} is to show that every conjugate self-twist of $\ovl\rho$, and in fact of every constant-determinant pseudodeformation of $\ovl{\rho}$, can be lifted to a conjugate self-twist of $(T,d)$ (see \cref{lifting residual csts} and \cref{lifting arbitrary csts} below).  

Since we are working only with constant-determinant pseudodeformations, we shall identify any $\F$-valued character $\eta$ with the $W(\F)$-valued character $s(\eta)$.  Furthermore, we will consider $\eta$ as being valued in any $W(\F)$-algebra via the structure map.  If $\sigma$ is an automorphism of $\F$, we write $W(\sigma)$ the automorphism of $W(\F)$ induced by $\sigma$. 

We introduce some notation that will be used in the proof of \cref{lifting residual csts}.  For any $W(\F)$-algebra $A$, let $A^{\sigma} \coloneqq A \otimes_{W(\F),W(\sigma)} W(\F)$, where $W(\F)$ is considered as a $W(\F)$-algebra via $W(\sigma)$.  We can equip $A^{\sigma}$ with two different $W(\F)$-algebra structures by letting $W(\F)$ act either on the first or second factor of the tensor product.  In what follows, we refer to these actions respectively as the \textit{first} or \textit{second $W(\F)$-algebra structure} on $A^{\sigma}$.  Let $\iota(\sigma, A) \colon A \to A^{\sigma}$ be the natural map given by $\iota(\sigma, A)(a) = a \otimes 1$.  It is an isomorphism of rings with inverse given by $\iota(\sigma^{-1}, A^{\sigma})$ since $(A^{\sigma})^{\sigma^{-1}}$ can be naturally identified with $A$ as a $W(\F)$-algebra.  Furthermore, $\iota(\sigma, A)$ is a morphism of $W(\F)$-algebras with respect to the first structure on $A^{\sigma}$.  Note that if we view $\cA^{\sigma}$ with respect to its second $W(\F)$-algebra structure, its residue field is $\F \otimes_{\F, \sigma} \F$, which is identified with $\F$ via $x \otimes y \mapsto \sigma(x)y$.  The proof of the following proposition is a more streamlined treatment of the arguments in \cite[Section 2]{Lang2016}.  

\begin{proposition}\label{lifting residual csts}
Let $(\sigma, \eta) \in \til \Sigma_{\ovl{\rho}}$.  Then there is an automorphism $\tilde{\sigma}$ of $\cA$ such that \mbox{$(\tilde{\sigma}, \eta) \in \til \Sigma_T$} and $\tilde{\sigma}$ induces $\sigma$ modulo $\m_{\cA}$.  Furthermore, for any $w$ in the image of $W(\F)$ in $\cA$, we have \mbox{$\tilde{\sigma}(w) = W(\sigma)(w)$}.
\end{proposition}

Note that any such $\tilde{\sigma}$ is necessarily unique, because it is determined by the character $\eta$.

\begin{proof}
Note that $\eta T \colon \Pi \to \cA$ is the universal constant-determinant pseudodeformation of \mbox{$\eta \otimes \ovl{\rho} \cong \act{\sigma}{\ovl{\rho}}$}.  We claim that, considering $\cA^{\sigma^{-1}}$ with its second $W(\F)$-algebra structure, $\iota(\sigma^{-1}, \cA) \circ (\eta T)$ is a constant-determinant pseudodeformation of $\ovl{\rho}$.  Indeed, reducing $\iota(\sigma^{-1}, \cA) \circ (\eta T)$ modulo the maximal of ideal of $\cA^{\sigma^{-1}}$ gives
\[
\iota(\sigma^{-1}, \F) \circ (\eta \tr \ovl{\rho}) = \act{\sigma}{\tr \ovl{\rho}} \otimes_{\F, \sigma^{-1}} 1 = 1 \otimes \tr \ovl{\rho},
\]
which is identified with $\tr \ovl{\rho}$ under the identification of $\F \otimes_{\F, \sigma^{-1}} \F$ with $\F$ discussed prior to the proposition.  

By universality, there is a unique $W(\F)$-algebra homomorphism $\alpha \colon \cA \to \cA^{\sigma^{-1}}$, where $\cA^{\sigma^{-1}}$ is given its second $W(\F)$-algebra structure, such that 
\[
\alpha \circ T = \iota(\sigma^{-1}, \cA) \circ (\eta T).
\]
Since $\iota(\sigma, \cA^{\sigma^{-1}})$ is the inverse of $\iota(\sigma^{-1}, \cA)$, we have that
\begin{equation}\label{cst relation}
\iota(\sigma, \cA^{\sigma^{-1}}) \circ \alpha \circ T = \eta T.
\end{equation}
Define $\tilde{\sigma} \coloneqq \iota(\sigma, \cA^{\sigma^{-1}}) \circ \alpha$, which is a ring endomorphism of $\cA$.  The relation \eqref{cst relation} implies that $\tilde{\sigma}$ is an automorphism of $\cA$ since the image of $T$ topologically generates $\cA$ as a $W(\F)$-module \cite[Proposition 5.3.3]{Bellaiche18} and $\eta$ takes values in $W(\F)$.  The relation \eqref{cst relation} also shows that $(\tilde{\sigma}, \eta) \in \til \Sigma_T$.

Finally, let $w \in W(\F)$.  Since $\alpha$ is a $W(\F)$-algebra homomorphism with respect to the second $W(\F)$-algebra structure on $\cA^{\sigma^{-1}}$, we have that
\begin{align*}
\tilde{\sigma}(w) &= \iota(\sigma, \cA^{\sigma^{-1}})\circ\alpha(w) = \iota(\sigma, \cA^{\sigma^{-1}})(1 \otimes w)\\
&= \iota(\sigma^{-1}, \cA)^{-1}(W(\sigma)(w) \otimes 1) = W(\sigma)(w).\qedhere\end{align*}\end{proof}

For the rest of \cref{section lifting csts}, fix a local pro-$p$ $W(\F)$-algebra $A$ with residue field $\F$ and a constant-determinant pseudodeformation $(t, d) \colon \Pi \to A$ of $\ovl{\rho}$.  Assume that $A$ is the $W(\F)$-algebra generated by $t(\Pi)$.  Let $\alpha_t \colon \cA \to A$ be the unique $W(\F)$-algebra homomorphism such that $\alpha \circ T = t$ given by universality.  The following corollary shows that conjugate self-twists of $(t,d)$ also lift to conjugate self twists of $(T,d)$.

\begin{corollary}\label{lifting arbitrary csts}
Given $(\sigma, \eta) \in \til \Sigma_t$, there is a unique $(\tilde{\sigma}, \eta) \in \til \Sigma_T$ such that $\alpha_t \circ \tilde{\sigma} = \sigma \circ \alpha_t$.
\end{corollary}

\begin{proof}
Let $\ovl{\sigma}$ denote the automorphism of $\F$ induced by $\sigma$.  Let $\tilde{\sigma}$ be the automorphism of $\cA$ given by \cref{lifting residual csts} lifting $\ovl{\sigma}$ to $\cA$.  Then we just have to show that $\alpha_t \circ \tilde{\sigma} = \sigma \circ \alpha_t$.  Note that $\sigma$ acts by $W(\ovl{\sigma})$ on the image of $W(\F)$ in $A$, so $\sigma^{-1} \circ \alpha \circ \tilde{\sigma}$ is a $W(\F)$-algebra homomorphism.  Thus by universality, it suffices to show that $t = \sigma^{-1} \circ \alpha_t \circ \tilde{\sigma} \circ T$.  Since $\eta$ takes values in $W(\F)$ and $\alpha_t$ is a $W(\F)$-algebra homomorphism, we have that
\[
\sigma^{-1} \circ \alpha_t \circ \tilde{\sigma} \circ T = \sigma^{-1} \circ \alpha_t(\eta T) = \sigma^{-1}(\eta t) = \sigma^{-1}(\sigma(t)) = t.\qedhere
\]
\end{proof}

We end \cref{section lifting csts} with some observations about the consequences of \cref{lifting residual csts} and \cref{lifting arbitrary csts}.  They give the following commutative diagram with exact rows.
\[
\xymatrix@=20pt@R=20pt{
	1\ar@{->}[r] & \Sigma_t^{\di}\ar@{->}[r]\ar@{^{(}->}[d] & \til \Sigma_t \ar@{->}[r]\ar@{^{(}->}[d] & \Sigma_t\ar@{->}[r]\ar@{^{(}->}[d] & 1\\
	1\ar@{->}[r] & \Sigma_T^{\di}\ar@{->}[r]\ar@{^{(}->}[d] & \til \Sigma_T \ar@{->}[r]\ar@{->}[d]^*[@]{\cong} & \Sigma_T\ar@{->}[r]\ar@{->>}[d] & 1\\
	1\ar@{->}[r] & \Sigma_{\ovl{\rho}}^{\di}\ar@{->}[r] & \til \Sigma_{\ovl{\rho}}\ar@{->}[r] & \Sigma_{\ovl{\rho}}\ar@{->}[r] & 1
	}
\]
We write 
\[
\beta_t \colon \Sigma_t \to \Sigma_{\ovl{\rho}}
\] 
for the composition of the vertical maps on the right in the above diagram.  It is induced by the composition $\tilde{\beta}_t \colon \til \Sigma_t \to \til \Sigma_{\ovl{\rho}}$ of the middle maps, which reflects the fact that every conjugate self-twist of $(t,d)$ induces a conjugate self-twist of $\ovl{\rho}$.  Combining \cref{lifting arbitrary csts} with \cref{tildeSigma rho-bar finite}, we see that $\til\Sigma_t$ is a finite abelian group for any constant-determinant pseudodeformation $(t,d)$ of $\ovl{\rho}$.

In this paper, we will only be concerned with pseudodeformations $(t, d)$ of $\ovl{\rho}$ that are not a priori small.  Under this assumption, if $t \neq \tr \ovl{\rho}$ then $\Sigma_t^{\di} = 1$ and $\Sigma_T^{\di} = 1$ by \cref{precise structure of Sigma}(1).  In particular, $\til \Sigma_t = \Sigma_t$ and $\til \Sigma_T = \Sigma_T$, so (except for $\ovl{\rho}$) a conjugate self-twist $(\sigma, \eta)$ is determined uniquely by the automorphism $\sigma$.  

\subsection{The dihedral case}\label{residual trivial csts}
As usual, let $\Pi$ be a $p$-finite profinite group, $A$ a local pro-$p$ ring and $(t,d) \colon \Pi \to A$ a constant determinant pseudorepresentation with trace algebra $A$.  Throughout this section we assume that its associated residual representation $\overline{\rho}$ is dihedral with nonabelian projective image.  We also fix any well-adapted $(t,d)$-representation $\rho$, which can be taken to be valued in $\GL_2(A)$ since $\overline{\rho}$ is absolutely irreducible \cref{Bellaiche 2.4.2}(5).  This case requires special care for two related reasons.  First, it is the only case when $A$ is not generated simply by $I_1(\rho)$ as a $W(\F)$-algebra; one also needs to include $B_1(\rho)$ in the generating set by \cref{Bellaiche main theorem}. (As explained in \cref{B and C in A}, it makes sense to view $B_1(\rho)$ as a subset of $A$ in this case.)  Second, this is the only case when $\Sigma_{\overline{\rho}}^{\di}$ is nontrivial and hence $\ker \beta_t$ can be nontrivial.  As we will see, a nontrivial element in $\ker \beta_t$ necessarily behaves quite differently from all other conjugate self-twists, because its action cannot be seen on the residue field.  In \cref{plus/minus parts} we explain how a nontrivial element in $\ker \beta_t$ interacts with $I_1(\rho)$ and $B_1(\rho)$.  In contrast, when $\ker \beta_t=1$ we show that $\B_\rho(\F)$ has the same fraction field as $A = \B_\rho(\F)+ W(\F)B_1(\rho)$.  

We use $\tau$ to refer to the nontrivial element of $\ker \beta_t$, should it exist.  For $\varepsilon \in \{+, -\}$, let
\[
A^\varepsilon \coloneqq \{a \in A \colon \act{\tau}{a} = \varepsilon a\}.
\]

\begin{proposition}\label{plus/minus parts}
Assume that $\ovl{\rho}$ is projectively dihedral and nonabelian.  Suppose there exists $1 \neq \tau \in \ker \beta_t$.  If $\rho \colon \Pi \to \GL_2(A)$ is a well-adapted $(t,d)$-representation, then $A^+ = \B_\rho(\F)$ and $A^- = W(\F)B_1(\rho)$.
\end{proposition}

\begin{proof}
Since $A = A^+ \oplus A^-$ and $A = \B_\rho(\F) + W(\F)B_1(\rho)$ by Bella\"iche's \cref{Bellaiche main theorem}, it suffices to show that $\tau$ acts trivially on $I_1(\rho)$ and by $-1$ on $B_1(\rho)$.  Let $\eta \colon \Pi \to \{\pm 1\}$ be the unique quadratic character such that $\ovl{\rho} \cong \ovl{\rho} \otimes \eta$.  (It is unique by \cref{precise structure of Sigma} since the projective image of $\ovl{\rho}$ is not isomorphic to $(\Z/2\Z)^2$).  Since $\tau \in \ker \beta_t$, it follows that $\eta$ must be the character such that $(\tau, \eta) \in \til \Sigma_t$.

We first prove that $I_1(\rho)$ is fixed by $\tau$.  As usual, let $\Gamma \coloneqq \im \rho \cap \Gamma_A(\m)$.  Recall that by definition $I_1(\rho)$ is the $\Z_p$-module topologically generated by $\{\alpha - \delta \colon \left(\begin{smallmatrix} 
1 + \alpha & b\\
c & 1 + \delta
\end{smallmatrix}\right) \in \Gamma\}$.  Let $g \in \ker \eta$.  Since $\rho$ is well adapted, we can write
\[
\rho(g) = \gamma \cdot \begin{pmatrix} 
s(\lambda) & 0\\
0 & s(\mu)
\end{pmatrix}
\]
for some $\gamma \in \Gamma$ and $\lambda, \mu \in \F^\times$.  Write $\gamma = \left(\begin{smallmatrix} 
1 + \alpha & b\\
c & 1 + \delta
\end{smallmatrix}\right)$ with $\alpha - \delta = 2a$, $\left(\begin{smallmatrix} 
a & b\\
c & -a
\end{smallmatrix}\right) \in L_1(\rho)$ and $0 = \alpha + \delta + \alpha\delta - bc$.  Then we have 
\[
\rho(g) = \begin{pmatrix} 
s(\lambda)(1 + \alpha) & s(\mu)b\\
s(\lambda)c & s(\mu)(1 + \delta)
\end{pmatrix}.
\]
Since $g \in \ker \eta$, it follows that
\[
s(\lambda)(1 + \alpha) + s(\mu)(1 + \delta) = \tr \rho(g) = \act{\tau}{(\tr \rho(g))} = s(\lambda)(1 + \act{\tau}{\alpha}) + s(\mu)(1 + \act{\tau}{\delta}),
\]
since $\tau$ acts trivially on $W(\F)$.  Thus, we obtain
\[
s(\lambda_0\mu_0^{-1})(\alpha - \act{\tau}{\alpha}) = \act{\tau}{\delta} - \delta
\]
for all $\lambda,\mu$ such that $\left(\begin{smallmatrix} 
\lambda & 0\\
0 & \mu
\end{smallmatrix}\right) \in \im \ovl{\rho}$.  As the projective image of $\ovl{\rho}$ is not isomorphic to $\Z/2\Z$ or $(\Z/2\Z)^2$, it follows that when $g$ varies $\lambda\mu^{-1}$ takes at least two distinct values in $\F^\times$.  Thus it follows that $\alpha - \act{\tau}{\alpha} = 0 = \act{\tau}{\delta} - \delta$.  Since $I_1(\rho)$ is generated by $\alpha - \delta$ with $\alpha,\delta$ as above, it follows that $I_1(\rho)$ is fixed by $\tau$.

The proof that $\tau$ acts by $-1$ on $B_1(\rho)$ is similar.  Namely, recall that $B_1(\rho)$ is topologically generated by $\{b, c \in A \colon \left(\begin{smallmatrix} 
1 + \alpha & b\\
c & 1 + \delta
\end{smallmatrix}\right) \in \Gamma\}$.  Let $g \in \Pi \setminus \ker \eta$.  Again since $\rho$ is well adapted, we can write 
\[
\rho(g) = \gamma \cdot \begin{pmatrix} 
0 & s(\lambda)\\
s(\mu) & 0
\end{pmatrix}
\]
for some $\gamma \in \Gamma, \lambda, \mu \in \F^\times$.  As above, write $\gamma = \left(\begin{smallmatrix} 
1 + \alpha & b\\
c & 1 + \delta
\end{smallmatrix}\right)$.  Then we have
\[
\rho(g) = \begin{pmatrix} 
s(\mu)b & s(\lambda)(1 + \alpha)\\
s(\mu)(1 + \delta) & s(\lambda)c
\end{pmatrix}.
\]
Since $g \not\in \ker \eta$, it follows that
\[
s(\mu)b + s(\lambda)c = \tr \rho(g) = -\act{\tau}{(\tr \rho(g))} = -s(\mu)\act{\tau}{b} - s(\lambda)\act{\tau}{c}.
\]
Thus 
\[
s(\mu\lambda^{-1})(b + \act{\tau}{b}) = -(c + \act{\tau}{c})
\]
for all $\left(\begin{smallmatrix} 
0 & \lambda\\
\mu & 0
\end{smallmatrix}\right) \in \im \ovl{\rho}$.    Once again, since the projective image of $\ovl{\rho}$ is not isomorphic to $\Z/2\Z$ or $(\Z/2\Z)^2$, it follows that $\lambda\mu^{-1}$ takes at least two distinct values in $\F^\times$.  Therefore $b + \act{\tau}{b} = 0 = c + \act{\tau}{c}$.  Since $B_1(\rho)$ is generated by such $b$ and $c$, it follows that $\tau$ acts on $B_1(\rho)$ by $-1$.
\end{proof}

In particular, we can always apply \cref{plus/minus parts} to the universal pseudorepresentation \mbox{$(T, d) \colon \Pi \to \cA$} of $\overline{\rho}$ since $\Sigma_T \cong \til \Sigma_{\overline{\rho}}$ and thus $\ker \beta_T$ is always nontrivial in the dihedral case whenever a nondihedral deformation exists.  Fix a well-adapted $(T,d)$-representation $\rho^{\univ} \colon \Pi \to \GL_2(\cA)$.  By universality, we have a $W(\F)$-algebra homomorphism $\alpha_t \colon \cA \to A$ such that $t = \alpha_t \circ T$.  Let $\rho = \rho_t$ be the well-adapted $(t,d)$-representation obtained by composing $\rho^{\univ}$ with the map $\GL_2(\cA) \to \GL_2(A)$ induced by $\alpha_t$.  Since the Pink-Lie algebra is functorial with respect to surjective ring homomorphisms, we see that $I_1(\rho^{\univ})$ (respectively, $B_1(\rho^{\univ})$) surjects onto $I_1(\rho)$  (respectively, $B_1(\rho)$).

For any subfield $\F' \subseteq \F$, let $\cA' = \B_{\rho^{\univ}}(\F') + W(\F')B_1(\rho^{\univ})$. We have $\cA' = (\cA')^+ \oplus (\cA')^-$ with $(\cA')^+ = \B_{\rho^{\univ}}(\F')$ and $(\cA')^- = W(\F')B_1(\rho^{\univ})$ by \cref{plus/minus parts}.

\begin{proposition}\label{field of fractions}
Suppose that $A$ is a local pro-$p$ domain and $(t,d) \colon \Pi \to A$ is a constant-determinant pseudodeformation of a dihedral $\overline{\rho}$ such that $\ker \beta_t = 1$.  Let $\rho$ be a well-adapted $(t,d)$-representation obtained from a universal one as described above.  Then for any subfield $\F' \subseteq \F$:
\begin{enumerate} 
\item\label{noetherian rings} $\cA'$ and $\B_\rho(\F')$ are noetherian rings;
\item\label{noetherian modules} $W(\F')B_1(\rho)$ and hence $\B_\rho(\F') + W(\F')B_1(\rho)$ are noetherian $\B_\rho(\F')$-modules;
\item\label{same fraction field} $\B_\rho(\F')+W(\F')B_1(\rho)$ and $\B_\rho(\F')$ have the same field of fractions.
\end{enumerate}
\end{proposition}

\begin{proof}
For \eqref{noetherian rings} note that $\B_\rho(\F')$ is the image of $(\cA')^+$ under $\alpha_t$, which is noetherian if $\cA'$ is by \cref{R+ noetherian}.  To see that $\cA'$ is noetherian when $\F' = \F$, we have \mbox{$\cA' = \cA = \B_{\rho^{\univ}}(\F) + W(\F)B_1(\rho^{\univ})$} by \cref{Bellaiche main theorem}, and hence $\cA$ is noetherian by the $p$-finiteness of $\Pi$.  Note that $\cA$ is finite and integral over $\B_{\rho^{\univ}}(\F') + W(\F')B_1(\rho^{\univ})$ since this is true of $W(\F)$ over $W(\F')$.  Then\\ \mbox{$\cA' = \B_{\rho^{\univ}}(\F') + W(\F')B_1(\rho^{\univ})$} is noetherian by \cite[Theorem 2]{Eakin}.

The statements in \eqref{noetherian modules} follow from the corresponding statements for $\rho^{\univ}$, which in turn follow from \cref{R fg as R+ module} since $\cA'$ is noetherian.

As $\B_\rho(\F')$ is the image of $(\cA')^+$ under $\alpha_t$ while $W(\F')B_1(\rho)$ is he image of $(\cA')^-$ and $\cA'$ is noetherian, \eqref{same fraction field} follows from \cref{nongraded quotients}.  
\end{proof}

\medskip


\section{Regularity and residual conjugate self-twists}
\label{section: technical residual}

In this section we prepare the groundwork for \cref{conceptual interpretation of Bellaiche ring} by studying conjugate self-twists of $\ovl \rho$, particularly how they interact with regularity (\cref{regular}).  In particular, $\ovl \rho \colon \Pi \to \GL_2(\F)$ is a semisimple regular representation throughout this section, and after \cref{red/reg} it is always absolutely irreducible.  We only consider simple conjugate self-twists in this section and thus write $\til \Sigma_{\ovl \rho}$ for $\til \Sigma_{\ovl \rho}(\F)$ and similarly for $\Sigma_{\ovl \rho}$.  In \cref{red/reg} we see that when $\ovl \rho$ is regular we may assume that it has no conjugate self-twists.  We consider the restriction of $\ovl \rho$ to the kernels of twist characters in \cref{res to Pi0}, followed by some technical lemmas about the extension $\F/\E$ in \cref{F over E}.  In \cref{choose basis for rhobar section} we define the condition of \textit{goodness} when $\ovl \rho$ is octahedral, which weighs on our main theorem.  \cref{choose basis for rhobar section} culminates in \cref{choice of basis} where we choose the basis of $\ovl \rho$ that will be used throughout \cref{conceptual interpretation of Bellaiche ring}.

Recall that $\ovl{\rho}$ is \textit{regular} if $\im \ovl{\rho}$ contains an element with eigenvalues $\lambda_0, \mu_0 \in \ovl{\F}^\times$ such that $\lambda_0\mu_0^{-1} \in \E^\times \setminus \{\pm1\}$.  

\subsection{Reducible regular representations}\label{red/reg} 
We show that if $\ovl{\rho}$ is reducible and regular, then one can eliminate the conjugate self-twists of $\ovl{\rho}$ by twisting $\ovl{\rho}$ by a character, making the proof of our main theorem especially easy in that case (cf. \cref{the reducible case}).

\begin{lemma}\label{regular reducible reps}
Suppose that $\ovl{\rho} = \varepsilon \oplus \delta$ and $\ovl{\rho}$ is regular.  If $(\sigma, \eta) \in \til \Sigma_{\ovl{\rho}}$, then $\act{\sigma}{\varepsilon} = \eta\varepsilon$ and $\act{\sigma}{\delta} = \eta\delta$.  In particular, $\varepsilon\delta^{-1}$ takes values in $\E$.
\end{lemma}

\begin{proof}
It suffices to show that if $\ovl{\rho}$ is regular then we cannot have $\act{\sigma}{\varepsilon} = \eta\delta$ and $\act{\sigma}{\delta} = \eta\varepsilon$.  If this were true, then we would have $\act{\sigma}{\varepsilon}\delta^{-1} = \eta = \act{\sigma}{\delta}\varepsilon^{-1}$, which implies that
\begin{equation}\label{sigma inverts}
\act{\sigma}{(\varepsilon\delta^{-1})} = \delta\varepsilon^{-1}.
\end{equation}

Since $\ovl{\rho}$ is regular, there is some $g \in \Pi$ such that $\varepsilon(g)\delta(g)^{-1} \in \E \setminus \{\pm 1\}$.  As $\E$ is fixed by $\sigma$ by \cref{a0algfix}, it follows from \eqref{sigma inverts} that $\varepsilon(g)\delta(g)^{-1} =  \pm 1$, a contradiction.

The last sentence in the statement of the lemma follows from the fact that, for any \mbox{$\sigma \in \Sigma_{\ovl{\rho}} = \Gal(\F/\E)$}, we have $\act{\sigma}{\varepsilon}\varepsilon^{-1} = \eta = \act{\sigma}{\delta}\delta^{-1}$ and hence $\varepsilon\delta^{-1}$ is fixed by $\Gal(\F/\E)$.
\end{proof}

\begin{corollary}\label{no csts reducible case}
Suppose $\ovl{\rho} = \varepsilon \oplus \delta$ and $\ovl{\rho}$ is regular.  Then $\ovl{\rho}' \coloneqq \ovl{\rho} \otimes \delta^{-1}$ has no conjugate self-twists.  
\end{corollary}

\begin{proof}
Since $\ovl{\rho}$ is regular, its projective image cannot have order 2.  Therefore it suffices to show that $\Sigma_{\ovl{\rho}'}$ is trivial by \cref{precise structure of Sigma}.  Let $\F'$ be the extension of $\F_p$ generated by the trace of $\ovl{\rho}'$.  Then $\Sigma_{\ovl{\rho}'} = \Gal(\F'/\E)$, so it suffices to show that $\F' \subseteq \E$.  But $\F'$ is generated by the values of $\varepsilon\delta^{-1}$, which takes values in $\E$ by \cref{regular reducible reps}.
\end{proof}

\subsection{Kernels of twist characters and regularity}\label{res to Pi0}
In this section we introduce the subgroup of $\Pi$ given by intersecting the kernels of all twist characters.  It is often useful to restrict to this subgroup because doing so kills the conjugate self-twists and but retains fullness.  In this section we study how this restriction affect the residual representation and regularity, first in the exceptional/large image case, then when $\ovl \rho$ is dihedral. Define 
\begin{equation}\label{Pi0 defn}
\Pi_0(\ovl{\rho}) \coloneqq \bigcap_{(\sigma, \eta) \in \til \Sigma_{\ovl{\rho}}} \ker \eta.
\end{equation}

\begin{remark}\label{Pi over Pi0 abelian}
The quotient $\Pi/\Pi_0(\ovl \rho)$ is abelian. Indeed, $\Pi_0(\ovl \rho)$ is the kernel of the natural map diagonal map of $\Pi$ to $\prod_{(\sigma, \eta) \in \til\Sigma_{\ovl{\rho}}}\Pi/\ker\eta$, 
so $\Pi/\Pi_0(\ovl \rho)$ can be embedded into an abelian group.
\end{remark}

As we will now see, it is easiest to control $\ovl \rho|_{\Pi_0(\ovl \rho)}$ when the order of $\det \ovl \rho$ is a power of $2$, which is an assumption we are forced to make in \cref{conceptual interpretation of Bellaiche ring}.  We may always twist $\ovl{\rho}$ by a character to assume that the order of $\det \ovl{\rho}$ is a power of 2:

\begin{lemma}\label{detorder2}
Let $d \colon \Pi \to \F^\times$ be a character.  Then there is a character $\chi \colon \Pi \to \F^\times$ such that the order of $d\chi^2$ is a power of 2.
\end{lemma}

\begin{proof}
The odd-order part of $d$ has a square root $\psi$.  Take $\chi = \psi^{-1}$.
\end{proof}

\begin{lemma}\label{big image mult free}
Assume that $\ovl{\rho}$ is exceptional or large.  If the order of $\det \ovl{\rho}$ is a power of 2, then $\ovl{\rho}|_{\Pi_0(\ovl \rho)}$ is absolutely irreducible.
\end{lemma}

\begin{proof}
If $(\sigma, \eta) \in \til \Sigma_{\ovl{\rho}}$, then $\eta^2$ is equal to a power of $\det \ovl{\rho}$ and hence the order of $\eta$ is a power of 2.  Thus $[\Pi \colon \Pi_0(\ovl \rho)]$ is a power of 2, and $\Pi/\Pi_0(\ovl \rho)$ is abelian by \cref{Pi over Pi0 abelian}.  

By hypothesis, the projective image of $\ovl{\rho}$ is isomorphic to one of $A_4, S_4, A_5$, $\PSL_2(\E), \PGL_2(\E)$.  None of $A_4, A_5, \PSL_2(\E)$ has a subgroup of 2-power index with abelian quotient.  Both $S_4$ and $\PGL_2(\E)$ have a unique proper 2-power index subgroup with abelian quotient, namely $A_4$ and $\PSL_2(\E)$, respectively.  Therefore the possible projective images of $\ovl{\rho}|_{\Pi_0(\ovl \rho)}$ are the same as for $\ovl \rho$, so $\ovl \rho|_{\Pi_0(\ovl \rho)}$ is absolutely irreducible.  
\end{proof}

\begin{proposition}\label{dihedral regular implies mult free}
Assume that $\ovl{\rho}$ is regular dihedral, say $\ovl{\rho} = \Ind_{\Pi_0}^\Pi \chi$.  Then $\ovl{\rho}|_{\Pi_0(\ovl{\rho})}$ is multiplicity free over $\E$.  Furthermore, given $g \in \Pi_0$, we have $g \in \Pi_0(\ovl{\rho})$ if and only if $\chi(g) \in \E^\times$.
\end{proposition}

\begin{proof}
Since $\ovl{\rho}$ is regular, it follows from \cref{resolve definitions of dihedral} that there is a unique subgroup $\Pi_0$ of $\Pi$ of index 2 such that $\ovl{\rho} \cong \Ind_{\Pi_0}^\Pi \chi$ for some character $\chi \colon \Pi_0 \to \F^\times$.  For any $h \in \Pi$, define $\chi^h \colon \Pi_0 \to \F^\times$ by $\chi^h(g) \coloneqq \chi(h^{-1}gh)$.  The character $\chi^h$ only depends on the class of $h$ in $\Pi/\Pi_0$.  Fix an element $c \in \Pi \setminus \Pi_0$.  Fix a generator $\sigma \in \Sigma_{\ovl{\rho}} = \Gal(\F/\E)$, and choose $\eta$ such that $(\sigma, \eta) \in \til \Sigma_{\ovl{\rho}}$.  (Note that there are two choices for $\eta$, and they differ by the character $\eta_0 \colon \Pi \twoheadrightarrow \Pi/\Pi_0 \cong \{\pm 1\}$.)  Then $\Pi_0(\ovl{\rho}) = \ker \eta_0 \cap \ker \eta$ since $\sigma$ generates $\Sigma_{\ovl{\rho}}$.  Therefore $\Pi_0(\ovl{\rho}) = \ker \eta|_{\Pi_0}$.

Note that any regular element for $\ovl{\rho}$ must be in $\Pi_0$ since elements in $\Pi \setminus \Pi_0$ have projective order~2.  By applying \cref{regular reducible reps} to $\ovl{\rho}|_{\Pi_0}$, we find that $\act{\sigma}{\chi} = \eta\chi$ and $\act{\sigma}{\chi^c} = \eta\chi^c$.  Hence $\eta|_{\Pi_0} = \act{\sigma}{\chi}\chi^{-1}$, so $g \in \Pi_0(\ovl{\rho})$ if and only if $\chi(g) \in \E^\times$.  In particular, $\ker \chi \subseteq \Pi_0(\ovl{\rho})$.  Furthermore, using the fact that $\act{\sigma}{\chi}\chi^{-1} = \eta|_{\Pi_0} = \act{\sigma}{\chi^c}(\chi^c)^{-1}$, we find that the character $\chi/\chi^c$ takes values in $\E^\times$.

We know that $\ovl{\rho}|_{\Pi_{0}(\ovl{\rho})}$ is multiplicity free over $\E$ if and only if there is some $g \in \Pi_0(\ovl{\rho})$ such that $\chi(g) \neq \chi^c(g)$.  If $\ker \chi \neq \ker \chi^c$, then we can choose $g \in \ker \chi \setminus \ker \chi^c$.  Then $\chi(g) = 1 \neq \chi^c(g)$ and $g \in \Pi_0(\ovl{\rho})$ by the previous paragraph.  Therefore we may assume that $\ker \chi = \ker \chi^c$.

Let $n$ denote the order of $\chi$.  Since $\ker \chi = \ker \chi^c$, we have that $\chi^c = \chi^a$ for some $a \in (\Z/n\Z)^\times$.  Note that $\chi^{c^2} = \chi$ since $c^2 \in \Pi_0$.  Therefore
\[
\chi = \chi^{c^2} = (\chi^c)^c = (\chi^a)^c = (\chi^c)^a = (\chi^a)^a = \chi^{a^2}.
\]
Fix $g_0 \in \Pi_0$ such that $\ovl{\rho}(g_0)$ generates the projective image of $\ovl{\rho}(\Pi_0)$.  We will show that $h \coloneqq g_0^{a-1}$ is in $\Pi_0(\ovl{\rho})$ and $\chi(h) \in \E^\times$ with $\chi(h) \neq \chi^c(h)$.  First we calculate, using the fact that $\chi^{a^2} = \chi$,
\[
\chi^c(h) = \chi^a(g_0^{a-1}) = \chi^{a^2}(g_0)\chi^{-1}(g_0) = 1.
\]
Hence $h \in \ker \chi^c = \ker \chi \subseteq \Pi_0(\ovl{\rho})$.  On the other hand, 
\[
\chi(h) = \chi(g_0^{a-1}) = \chi^c(g_0)\chi^{-1}(g_0).
\]
We saw in the second paragraph that $\chi/\chi^c$ is an $\E$-valued character.  Furthermore, $\chi^c(g_0)/\chi(g_0) \neq 1$ since $g_0$ was chosen has a generator of the projective image of $\ovl{\rho}(\Pi_0)$, which is isomorphic to the image of $\chi/\chi^c$.
\end{proof}

\begin{corollary}\label{rhobar irred over complement to kernel betat}
Assume that $\ovl{\rho}$ is regular and dihedral and that the order of $\det \ovl{\rho}$ is a power of 2.  Let $\sigma$ be a generator of $\Sigma_{\ovl{\rho}}$ and $\eta \colon \Pi \to \F^\times$ a character such that $(\sigma, \eta) \in \til \Sigma_{\ovl{\rho}}$.  Then either $\Sigma_{\ovl{\rho}}$ is trivial or $\ovl{\rho}|_{\ker \eta}$ is absolutely irreducible.
\end{corollary}

\begin{proof}
Write $\ovl{\rho} = \Ind_{\Pi_0}^\Pi \chi$ and fix $c \in \Pi \setminus \Pi_0$.  We shall make frequent use of \cref{resolve definitions of dihedral} in this proof without referencing it every time.  We saw in the proof of \cref{dihedral regular implies mult free} that $\Pi_0(\ovl{\rho}) = \Pi_0 \cap \ker \eta$.  Thus \cref{dihedral regular implies mult free} implies that $\chi|_{\Pi_0 \cap \ker \eta} \neq \chi^c|_{\Pi_0 \cap \ker \eta}$.  

If $\ker \eta \neq \Pi_0 \cap \ker \eta$, then $[\ker \eta \colon \Pi_0 \cap \ker \eta] = 2$ since $[\Pi \colon \Pi_0] = 2$.  Thus $\ovl{\rho}|_{\ker \eta} \cong \Ind_{\Pi_0 \cap \ker \eta}^{\ker \eta} \chi|_{\Pi_0 \cap \ker \eta}$.  Since $\chi|_{\Pi_0 \cap \ker \eta} \neq \chi^c|_{\Pi_0 \cap \ker \eta}$ it follows that $\ovl{\rho}|_{\ker \eta}$ is irreducible.

If $\ker \eta = \Pi_0 \cap \ker \eta$ then $\Pi_0 \supseteq \ker \eta$ and $\Pi/\ker \eta$ is a cyclic group whose order is a power of~2 since $\eta^2$ is a power of $\det \ovl{\rho}$.  If $\Pi_0 \neq \ker \eta$, then there is a subgroup $\ker \eta \subseteq \Pi' \subset \Pi_0$ such that $[\Pi_0 \colon \Pi'] = 2$.  Note that $\chi|_{\Pi'} \neq \chi^c|_{\Pi'}$ since $\chi|_{\ker \eta} \neq \chi^c|_{\ker \eta}$.  Then $\ovl{\rho}|_{\Pi_0} \cong \Ind_{\Pi'}^{\Pi_0} \chi|_{\Pi'}$ is irreducible, a contradiction since $\ovl{\rho} \cong \Ind_{\Pi_0}^\Pi \chi$.  Thus we must have $\Pi_0 = \ker \eta$.  Therefore $\ovl{\rho} \cong \ovl{\rho} \otimes \eta$ and so $\sigma$, and hence $\Sigma_{\ovl{\rho}}$, is trivial.  
\end{proof}

\subsection{$\F/\E$ when $\det \ovl \rho$ is a power of $2$}\label{F over E}

In \cref{conceptual interpretation of Bellaiche ring} we will assume that the order of $\det \ovl{\rho}$ is a power of 2.  A large part of the reason for that assumption is that it guarantees that $[\F \colon \E]$ can be taken to be a power of 2 as well, as the next lemma shows.  We need this in an induction argument in \cref{conceptual interpretation of Bellaiche ring}.  Given any $\F$-valued function $f$ and any subfield $\F'$ of $\F$, let us write $\F'(f)$ for the subfield of $\F$ generated over $\F'$ by the values of $f$.  

\begin{lemma}\label{structure of Sigma power of 2}
Assume that the order of $\det \ovl{\rho}$ is a power of 2.  Then the degree of $\F_p(\tr \ovl{\rho})$ over $\E$ is a power of~2.
\end{lemma}

\begin{proof}
Let $d \coloneqq \det \ovl{\rho}$.  Since the order of $d$ is a power of 2, the degree of $\E(d)$ over $\E$ is a power of~2.  But, for an arbitrary $g'\in\Pi$, the extension $\E(\tr \ovl{\rho}(g'))$ is at most quadratic over $\E(d)$ because $\tr \ovl{\rho}(g')$ satisfies
\[
d(g')x^2 - (\tr \ovl{\rho}(g'))^2/d(g') \in \E(d)[x].
\]
The field $\F_p(\tr \ovl{\rho})$ is obtained from $\E(d)$ by adding finitely many values of $\tr \ovl{\rho}$.
\end{proof}

In \cref{conceptual interpretation of Bellaiche ring} we will be interested in gradings coming from conjugate self-twists.  To be able to apply \cref{automorphisms to gradings} in those situations, we now verify one of the hypotheses.

\begin{lemma}\label{condition * satisfied}
Assume that both the order of $\det \ovl{\rho}$ and $[\F \colon \E]$ are powers of 2.  If $n = \#\Sigma_{\ovl{\rho}}$, then $\F$ contains a primitive $n^{\rm th}$ root of unity.  In particular, condition $(\ast)$ from \cref{autom grad} is satisfied.
\end{lemma}

\begin{proof}
Let $d \coloneqq \det \ovl{\rho}$, and write $2^s$ for the order of $d$.  We have that $\E(d)$ contains a primitive $(2^s)^{\rm th}$ root of unity.  If $[\F \colon \E(d)] = 2^r$, then $\F$ contains a primitive $(2^{r+s})^{\rm th}$ root of unity.  On the other hand, 
\[
n = \#\Sigma_{\ovl{\rho}} = [\F \colon \E] = 2^r[\E(d) \colon \E].
\]
Since $d$ has order $2^s$, it follows that $[\E(d) \colon \E]$ divides $2^{s-1}$.  Thus $n$ divides $2^{r+s-1}$, and so $\F$ contains a primitive $n^{\rm th}$ root of unity.
\end{proof}

\subsection{A good basis for $\ovl{\rho}$}\label{choose basis for rhobar section}
We need to carefully choose a basis for $\ovl{\rho}$ that has many good properties and will allow us to choose a good $(t,d)$-representation in \cref{choosing a basis}.  In this section we explain how to find this basis when $\ovl{\rho}$ is exceptional or large.  Let us first define an extra condition on octahedral representations.

\begin{definition}\label{good}
We say a regular octahedral representation $\ovl{\rho}$ is \textit{good} if at least one of the following properties is satisfied:
\begin{enumerate}
\item $p \equiv 1 \bmod 3$;
\item $\ovl{\rho}$ is strongly regular;
\item there is a regular element $g_0 \in \Pi$ such that $g_0^2 \in \Pi_0(\ovl{\rho})$.
\end{enumerate}
\end{definition}

We shall need to know that if $\ovl{\rho}$ is good, then twisting away the odd part of the determinant of $\ovl{\rho}$ gives a representation that is also good.

\begin{lemma}\label{good under twisting}
Let $\ovl{\rho} \colon \Pi \to \GL_2(\F)$ be a good representation.  Let $\chi \colon \Pi \to \F^\times$ be the unique odd-order character such that the order of $\chi^2\det \ovl{\rho}$ is a power of 2.  Then $\ovl{\rho} \otimes \chi$ is good.
\end{lemma}

\begin{proof}
First note that twisting by any character does not change the projective image, so $\ovl{\rho} \otimes \chi$ is octahedral.  Regularity is also invariant under twisting.  The claim is clear if $p \equiv 1 \bmod 3$, so we assume that $p \equiv 2 \bmod 3$.  The regularity assumption then implies that $\zeta_4 \in \F_p$ by \cref{regularity remark background}.  As in the proof of \cref{detorder2}, decompose $\det \ovl{\rho} = d_1d_2$, where $d_i \colon \Pi \to \F^\times$ are characters such that the order of $d_1$ is odd and the order of $d_2$ is a power of 2.

First suppose that $\ovl{\rho}$ is strongly regular.  Then there is a matrix $g_0 \in \Pi$ such that $\ovl{\rho}(g_0)$ has eigenvalues $\lambda_0, \mu_0 \in \E^\times$ such that $\lambda_0\mu_0^{-1} = \zeta_4$.  We have $\lambda_0\mu_0 = \det \ovl{\rho}(g_0) = d_1(g_0)d_2(g_0)$.  Note that any $\sigma \in \Gal(\F/\E)$ fixes $\lambda_0\mu_0$ since $\lambda_0,\mu_0 \in \E$.  Therefore $\sigma(d_1(g_0)d_2(g_0)) = d_1(g_0)d_2(g_0)$.  But since $d_1(g_0)$ is an odd order root of unity and $d_2(g_0)$ is a 2-power order root of unity, it follows that $\sigma$ must fix both $d_1(g_0)$ and $d_2(g_0)$.  Write $a$ for the order of $d_1$.  Then $\chi = d_1^{-(a+1)/2}$ by the proof of \cref{detorder2}.  In particular, $\chi(g_0) \in \E^\times$.  Thus the eigenvalues $\chi(g_0)\lambda_0$ and $\chi(g_0)\mu_0$ of $(\ovl{\rho} \otimes \chi)(g_0)$ are in $\E$.  Thus $g_0$ is a strongly regular element for $\ovl{\rho} \otimes \chi$, as desired.

Finally, suppose that there is a regular element $g_0 \in \Pi$ such that $g_0^2 \in \Pi_0(\ovl{\rho})$.  Let $\sigma$ be a generator for $\Gal(\F/\E)$ and let $\eta \colon \Pi \to \F^\times$ such that $(\sigma, \eta) \in \til \Sigma_{\ovl{\rho}}$.  Then $\Pi_0(\ovl{\rho}) = \ker \eta$ and $\Pi_0(\ovl{\rho} \otimes \chi) = \ker \act{\sigma}{\chi}\chi^{-1}\eta$.  Since $g_0^2 \in \Pi_0(\ovl{\rho})$ and $\act{\sigma}{\det \ovl{\rho}} = \eta^2 \det \ovl{\rho}$, it follows that $\det \ovl{\rho}(g_0) \in \E$.  But $\det \ovl{\rho}(g_0) = d_1(g_0)d_2(g_0)$, and since $d_1(g_0)$ is an odd order root of unity and $d_2(g_0)$ has 2-power order, it follows that both $d_1(g_0)$ and $d_2(g_0)$ are in $\E$.  Therefore $\chi(g_0) = d_1^{-(a+1)/2}(g_0) \in \E$.  Thus $g_0^2 \in \ker \act{\sigma}{\chi}\chi^{-1}\eta = \Pi_0(\ovl{\rho} \otimes \chi)$.
\end{proof}

Finally we describe the basis of $\ovl{\rho}$ that we shall work with in \cref{conceptual interpretation of Bellaiche ring}.  Let $Z$ denote the group of scalar matrices in $\GL_2(\F)$.  The following lemma justifies our definition of $\F_q$ in \cref{bellaiche results background} for exceptional representations.

\begin{lemma}\label{Bellaiche Fq vs E}
Up to conjugation, the image of $\ovl{\rho}$ is contained in $Z\GL_2(\E)$.  If $\F_q$ is an extension of $\E$ and  $\lambda_0, \mu_0 \in \ovl{\F}_p^\times$ are eigenvalues of a matrix in the image of $\ovl{\rho}$ such that $\lambda_0\mu_0^{-1} \in \F_q$, then we may further conjugate $\ovl{\rho}$ to assume that $\bigl(\begin{smallmatrix} 
\lambda_0 & 0\\
0 & \mu_0
\end{smallmatrix}\bigr) \in \im \ovl{\rho}$ and $\im \ovl{\rho} \subseteq Z\GL_2(\F_q)$.
\end{lemma}

\begin{proof}
By \cref{residual csts cut out E}, $\E = \F^{\Sigma_{\ovl{\rho}}}$.  First we show that $\ovl{\rho}$ can be conjugated to land in $Z\GL_2(\E)$.  Let $\sigma \in \Gal(\F/\E)$ be a generator and $\eta$ a character such that $(\sigma, \eta) \in \Sigma_{\ovl{\rho}}$.  Then there is some $x \in \GL_2(\F)$ such that for all $g \in \Pi$, we have $\act{\sigma}{\ovl{\rho}}(g) = x^{-1}\eta(g)\ovl{\rho}(g)x$.  By a theorem of Serge Lang \cite[Corollary to Theorem 1]{Lang56}, it follows that there is some $y \in \GL_2(\F)$ such that $x = \act{\sigma}{y}y^{-1}$.  Thus $\act{\sigma}{(y^{-1}\ovl{\rho}(g)y)} = \eta(g)(y^{-1}\ovl{\rho}(g)y)$.  Replacing $\ovl{\rho}$ by its conjugate by $y$, we have that the projective image of $\ovl{\rho}$ is fixed by $\Gal(\F/\E)$, and hence the image of $\ovl{\rho}$ lands in $Z\GL_2(\E)$, as desired.  

If $\F_q, \lambda_0, \mu_0$ are as in the statement of the lemma, then $\ovl{\rho}$ can be further conjugated such that $\bigl(\begin{smallmatrix} 
\lambda_0 & 0\\
0 & \mu_0
\end{smallmatrix}\bigr) \in \im \ovl{\rho}$ while preserving the property that the image of $\ovl{\rho}$ is in $Z\GL_2(\F_q)$.
\end{proof}

\begin{proposition}\label{choice of basis}
Let $\ovl{\rho} \colon \Pi \to \GL_2(\F)$ be regular and either exceptional or large.  If $\ovl{\rho}$ is octahedral, assume further that $\ovl{\rho}$ is good.  Assume that the order of $\det \ovl{\rho}$ is a power of $2$.  Then there is a regular element $g_0 \in \Pi$ and a basis for $\ovl{\rho}$ such that the following are simultaneously true:
\begin{enumerate}
\item $\im \ovl{\rho} \subseteq Z\GL_2(\E)$;
\item $\ovl{\rho}(g_0) = \bigl(\begin{smallmatrix} 
\lambda_0 & 0\\
0 & \mu_0
\end{smallmatrix}\bigr)$ for some $\lambda_0,\mu_0\in\F$;
\item if $p \geq 7$ and $\ovl{\rho}$ is large, then $\lambda_0, \mu_0 \in \F_p^\times$;
\item there is a positive integer $n$ such that $g_0^n \in \Pi_0(\ovl{\rho})$ and $\ovl{\rho}(g_0^n)$ is not scalar.
\end{enumerate}
\end{proposition}

\begin{proof}
By \cref{Bellaiche Fq vs E} we can always conjugate $\ovl{\rho}$ so that $\im \ovl{\rho} \subseteq Z\GL_2(\E)$.  If $g_0 \in \Pi$ is a regular element and $\lambda_0$ and $\mu_0$ are the eigenvalues of $\ovl{\rho}(g_0)$, then we may assume further that $\ovl{\rho}(g_0) = \bigl(\begin{smallmatrix} 
\lambda_0 & 0\\
0 & \mu_0
\end{smallmatrix}\bigr)$.

If $\ovl{\rho}$ is large, then up to conjugation, $\im \ovl{\rho} \supseteq \SL_2(\F_p)$.  Indeed, up to conjugation we may assume that the projective image of $\ovl{\rho}$ contains $\PSL_2(\E)$.  Therefore there is some $\lambda \in \F^\times$ such that $\lambda \bigl(\begin{smallmatrix} 
1 & 1\\
0 & 1
\end{smallmatrix}\bigr) \in \im \ovl{\rho}$.  Note that the $n^{\rm th}$ power of this matrix is $\lambda^n\bigl(\begin{smallmatrix} 
1 & n\\
0 & 1
\end{smallmatrix}\bigr)$.  Since $\lambda \in \F^\times$, its order $m$ is prime to $p$.  Therefore we can write $1 = am + bp \equiv am \bmod p$ for some $a, b \in \Z$.  Thus 
\[
\bigl(\begin{smallmatrix} 
1 & 1\\
0 & 1
\end{smallmatrix}\bigr) =  \lambda^{am}\bigl(\begin{smallmatrix} 
1 & am\\
0 & 1
\end{smallmatrix}\bigr) = (\lambda\bigl(\begin{smallmatrix} 
1 & 1\\
0 & 1
\end{smallmatrix}\bigr))^{am} \in \im \ovl{\rho}.
\]
Similarly, $\bigl(\begin{smallmatrix} 
1 & 0\\
1 & 1
\end{smallmatrix}\bigr) \in \im \ovl{\rho}$.  Since $\bigl(\begin{smallmatrix} 
1 & 1\\
0 & 1
\end{smallmatrix}\bigr)$ and $\bigl(\begin{smallmatrix} 
1 & 0\\
1 & 1
\end{smallmatrix}\bigr)$ generate $\SL_2(\F_p)$, it follows that $\SL_2(\F_p) \subseteq \im \ovl{\rho}$.

If $p \geq 7$, then we can choose $\alpha \in \F_p^\times$ such that $\alpha^2 \neq \pm 1$.  Then any $g_0 \in \Pi$ such that $\ovl{\rho}(g_0)$ has eigenvalues $\alpha, \alpha^{-1}$ satisfies the first three conditions.  Note that $\proj\ovl{\rho}(g_0) \in \PSL_2(\E)$.  Recall that $\Pi_0(\ovl{\rho})$ is a normal subgroup of 2-power index in $\Pi$ since the order of $\det \ovl{\rho}$ is a power of 2.  Furthermore, $\Pi/\Pi_0(\ovl{\rho})$ is abelian.  Therefore $\proj \ovl{\rho}(\Pi_0(\ovl{\rho}))$ is either $\PGL_2(\E)$ or $\PSL_2(\E)$.  In either case, we can find $g_0 \in \Pi_0(\ovl{\rho})$ such that $\ovl{\rho}(g_0)$ has eigenvalues $\alpha, \alpha^{-1}$.  Thus all of the properties of the proposition are satisfied for this choice of $g_0$.

Next suppose that $\ovl{\rho}$ is either tetrahedral or icosahedral.  Once again, $\proj \ovl{\rho}(\Pi_0(\ovl{\rho}))$ is a normal subgroup of $\proj \ovl{\rho}(\Pi)$ with $2$-power index and abelian quotient.  Since $\proj \ovl{\rho}(\Pi)$ is isomorphic to one of $A_4$ or $A_5$, it follows that $\proj \ovl{\rho}(\Pi_0(\ovl{\rho})) = \proj \ovl{\rho}(\Pi)$.  In particular, one can choose the regular element $g_0$ to be in $\Pi_0(\ovl{\rho})$, and the resulting representation satisfies all of the desired conditions.

Finally, suppose that $\ovl{\rho}$ is octahedral and good.  If $p \equiv 1 \bmod 3$ then any $g_0 \in \Pi$ such that $\proj\ovl{\rho}(g_0)$ has order 3 is a regular element.  Since $\proj \ovl{\rho}(\Pi_0(\ovl{\rho}))$ is a normal subgroup of $\proj \ovl{\rho}(\Pi)$ with $2$-power index and abelian quotient, it follows that $\Pi_0(\ovl{\rho})$ contains an element $g_0$ such that $\proj\ovl{\rho}(g_0)$ has order~3.  Such a $g_0$ satisfies all of the necessary conditions.

Next suppose that $p \equiv 2 \bmod 3$ and that $\ovl{\rho}$ is strongly regular.  Let $g_0 \in \Pi$ be a strongly regular element.  Then $\ovl{\rho}(g_0) = \bigl(\begin{smallmatrix} 
\lambda\zeta_4 & 0\\
0 & \lambda
\end{smallmatrix}\bigr)$ for some $\lambda \in \E^\times$.  (Note that $\zeta_4 \in \F_p$ since $\ovl{\rho}$ is regular and $p \equiv 2 \bmod 3$.)  We claim that $g_0 \in \Pi_0(\ovl{\rho})$.  Indeed, let $\sigma$ be a generator of $\Gal(\F/\E)$ and $\eta$ a character such that $(\sigma, \eta) \in \til \Sigma_{\ovl{\rho}}$.  Then $\Pi_0(\ovl{\rho}) = \ker \eta$.  Since $\lambda, \zeta_4 \in \E^\times$ we have
\[
\lambda(\zeta_4 + 1) = \act{\sigma}{(\lambda(\zeta_4 + 1))} = \act{\sigma}{\tr \ovl{\rho}(g_0)} = \eta(g_0)\tr\ovl{\rho}(g_0) = \eta(g_0)\lambda(\zeta_4 + 1).
\] 
As $\zeta_4 + 1 \neq 0$ it follows that $\eta(g_0) = 1$, and so $g_0 \in \Pi_0(\ovl{\rho})$, as claimed.  Therefore $g_0$ satisfies all of the necessary conditions.

Finally suppose that $p \equiv 2 \bmod 3$ and there is a regular element $g_0 \in \Pi$ such that $g_0^2 \in \Pi_0(\ovl{\rho})$.  Note that $\proj\ovl{\rho}(g_0)$ has order 4 since $p \not\equiv 1 \bmod 3$.  Therefore $\proj\ovl{\rho}(g_0^2)$ is nontrivial, so $\ovl{\rho}(g_0^2)$ is not scalar.  Therefore $g_0$ satisfies all of the conditions of the proposition.
\end{proof}

Note that the $g_0$ chosen in \cref{choice of basis} satisfies all of the conditions prior to \cref{well adapted}.  In particular, if $(t, d) \colon \Pi \to A$ is any admissible pseudodeformation of $\ovl{\rho}$, then any $(t,d)$-representation that is adapted to the element $g_0$ from \cref{choice of basis} is well adapted.

\medskip


\section{Fullness peers: $\B_\rho(\E)$ and $A^{\Sigma_\rho}$}
\label{conceptual interpretation of Bellaiche ring}

Throughout \cref{conceptual interpretation of Bellaiche ring} we fix a local pro-$p$ domain $A$ and an admissible pseudodeformation $(\Pi, \ovl{\rho}, t, d)$ over $A$.  We only consider $A$-valued conjugate self-twists throughout this section and thus write $\Sigma_t$ for $\Sigma_t(A)$ and similarly for $\til \Sigma_t$.  The goal of \cref{conceptual interpretation of Bellaiche ring} is to prove that $(t,d)$ is $A_0$-full whenever $(t,d)$ is not a priori small and regularity is satisfied.  In view of \cref{A0 Afixedbytwists fullness peers} and \cref{formal fullness corollary}, it suffices to show that $A^{\Sigma_t}$ and $\B_\rho(\E)$ are fullness peers for some well chosen $(t,d)$-representation $\rho$.  Let us point out an easy case when this is possible.  If $\ovl{\rho}$ has no conjugate self-twists, then $\E = \F$ by \cref{residual csts cut out E} and $\Sigma_t = 1$ by the diagram following \cref{lifting arbitrary csts}.  Furthermore, the assumption that $\til \Sigma_{\ovl{\rho}} = 1$ implies that $\ovl{\rho}$ is not dihedral and so $A = \B_\rho(\F)$ by \cref{Bellaiche main theorem}.  Therefore we have
\[
\B_\rho(\E) = \B_\rho(\F) = A = A^{\Sigma_t}.
\]  

In general, the proof that $(t,d)$ is $A^{\Sigma_t}$-full, hence $A_0$-full, is structured as follows.  The case when $\ovl{\rho}$ is reducible is easily done in \cref{reducible ASigmat fullness}, so from \cref{choosing a basis} onwards we always assume that $\ovl{\rho}$ is irreducible.  In light of \cref{formal fullness corollary}, the strategy is to prove that, under certain conditions on $\ovl{\rho}$ and a good choice of a $(t,d)$-representation $\rho$, the two rings $\B_\rho(\E)$ and $A^{\Sigma_t}$ have the same fields of fractions and $A^{\Sigma_t}$ is finitely generated as a $\B_\rho(\E)$-module.  This is done in \cref{final conclusion}, although key parts of it are proved in \cref{small non domain case} and \cref{induction step}.  \cref{fg plus same quotient field implies ideal containment} then implies that $A^{\Sigma_t}$ and $\B_\rho(\E)$ are fullness peers.  In \cref{final goal}, we combine \cref{formal fullness corollary}, which established $\B_\rho(\E)$-fullness, with \cref{final conclusion} to show that $(t,d)$ is $A^{\Sigma_t}$-full under mild assumptions on $\ovl{\rho}$.  Since $A^{\Sigma_t}$ and $A_0$ are fullness peers in the constant-determinant setting, we conclude that our admissible pseudodeformation $(t,d)$ is $A_0$-full.  

Let us now establish some assumptions on our fixed residual representation $\ovl{\rho} \colon \Pi \to \GL_2(\F)$.  Assume that $\ovl{\rho}$ is regular and, after \cref{the reducible case}, absolutely irreducible.  Whenever $\ovl{\rho}$ is absolutely irreducible, assume further that $\det \ovl{\rho}$ is a power of 2, which can always be achieved by twisting $\ovl{\rho}$ by a character by \cref{detorder2}.  Furthermore, the twisting operation does not change the field $\E$ by \cref{a0fixedbytwist}. 
Assume that $[\F \colon \E]$ is a power of 2, which is possible by \cref{structure of Sigma power of 2}.  Note that we do not require $\F$ to be the trace ring of $\ovl{\rho}$ since one may need to make a quadratic extension of the trace ring in order to make representations well-adapted in the dihedral case.

\subsection{The reducible case}\label{the reducible case}
When $\ovl{\rho}$ is reducible, we can use \cref{formal fullness corollary} to show that $(t,d)$ is $A^{\Sigma_t}$- and $A_0$-full.

\begin{proposition}\label{reducible ASigmat fullness}
Suppose that $\ovl{\rho} = \varepsilon \oplus \delta$ and that $\ovl{\rho}$ is regular.  If $(t,d)$ is not a priori small then $(t,d)$ is $A^{\Sigma_t}$- and $A_0$-full.
\end{proposition}

\begin{proof}
Let $(t', d') = (s(\delta^{-1})t, s(\delta^{-1})^2d)$, which is a pseudodeformation of $\ovl r \coloneqq \ovl{\rho} \otimes \delta^{-1}$.  Let $A_{t'}$ be the trace ring of $(t', d')$; its  residue field $\E$ since $\ovl r$ has no conjugate self-twists by \cref{no csts reducible case}.  Then $(\Pi, \ovl r, t',d')$ is an admissible pseudorepresentation over $A_{t'}$.  Note that $(t', d')$ is not a priori small since $(t,d)$ is not.  By \cref{formal fullness corollary}, there is a well-adapted $(t', d')$-representation $r$ such that $(t', d')$ is $\B_r(\E)$-full.  But $A_{t'} = \B_r(\E)$ by \cref{Bellaiche main theorem}.  Since $\ovl r$ has no conjugate self-twists by \cref{no csts reducible case}, it follows that $\Sigma_{t'}$ is trivial.  Thus $A_{t'}^{\Sigma_{t'}} = A_{t'} = \B_r(\E)$.  

By \cref{twisting does not affect fullness} it follows that $(t,d)$ is $A_{t'}^{\Sigma_{t'}}$-full.  We know that $A^{\Sigma_t}$ and $A_{t'}^{\Sigma_{t'}}$ have the same fields of fractions, namely $K_0$, by \cref{AcstA0}.  Furthermore $A$ is obtained by adjoining the values of $s(\delta)$ and $W(\F)$ to $A_{t'}$.  Therefore $A$, and hence $A^{\Sigma_t}$, is finitely generated over $A'$, so $(t, d)$ is $A^{\Sigma_t}$-full by \cref{fg plus same quotient field implies ideal containment}.  Finally, $A_0$ and $A^{\Sigma_t}$ are fullness peers in this setting by \cref{A0 Afixedbytwists fullness peers}.
\end{proof}

\subsection{Choosing a good $(t,d)$-representation}\label{choosing a basis}
Throughout \cref{choosing a basis} through \cref{bigness subsection} we fix an absolutely irreducible regular representation $\ovl{\rho} \colon \Pi \to \GL_2(\F)$ such that the order of $\det \ovl{\rho}$ is a power of 2.  We assume that $[\F \colon \E]$ is a power of 2 by \cref{structure of Sigma power of 2}.  If $\ovl{\rho}$ is octahedral, we assume further that $\ovl{\rho}$ is good.  Furthermore, we fix a good basis for $\ovl{\rho}$ as follows.  If $\ovl{\rho}$ is exceptional or large, choose a basis and a regular element $g_0 \in \Pi$ such that \cref{choice of basis} holds.  If $\ovl{\rho} = \Ind_{\Pi_0}^\Pi \chi$ is dihedral, assume that $\ovl{\rho}(\Pi_0)$ is diagonal and $\im \ovl{\rho}$ contains a matrix $\bigl(\begin{smallmatrix} 
0 & b\\
c & 0
\end{smallmatrix}\bigr)$ such that $bc^{-1} \in \F_p$, which is possible by \cite[Proposition 6.3.2]{Bellaiche18}.

Recall that $(T, d) \colon \Pi \to \cA$ is the universal constant-determinant pseudorepresentation.  Part of our arguments will require appealing to a universal $(T,d)$-representation.  This requires choosing a good $(T,d)$-representation $\rho^{\univ}$ and also choosing our $(t,d)$-representation to be compatible with $\rho^{\univ}$.  In particular, we want $I_1(\rho^{\univ})$ to be fixed by all conjugate self-twists of $(T,d)$.  In \cref{choosing a basis} we make these choices and compatibilities precise.  Since we need \cref{A' sub GMA,I1 fixed by twists} for the universal ring $\cA$, in \cref{choosing a basis} we do not require $A$ to be a domain, only a local pro-$p$ ring.

Fix a generator $\sigma_1$ of $\Sigma_{\ovl{\rho}} = \Gal(\F/\E)$.  We want to choose a character $\eta_1 \colon \Pi \to \F^\times$ such that $(\sigma_1, \eta_1) \in \til \Sigma_{\ovl{\rho}}$.  There is a unique choice for $\eta_1$ when $\ovl{\rho}$ is not dihedral.  If $\ovl{\rho}$ is dihedral and $\Sigma_t = 1$, choose $\eta_1$ to be the trivial character.  Recall from the end of \cref{section lifting csts} that $\beta_t \colon \Sigma_t \to \Sigma_{\ovl{\rho}}$ is given by reducing automorphisms of $A$ modulo $\m$.  If $\ovl{\rho}$ is dihedral and $\ker \beta_t = 1$ but $\Sigma_t \neq 1$, then there is a unique complement to $\Sigma_{\ovl{\rho}}^{\di}$ in $\til \Sigma_{\ovl{\rho}}$ that contains $\tilde{\beta}_t(\til \Sigma_t)$.  Choose $\eta_1$ such that $(\sigma_1, \eta_1)$ generates that complement.  Otherwise, when $\ovl{\rho}$ is dihedral, we may take $\eta_1$ to be either of the two characters such that $(\sigma_1, \eta_1) \in \til \Sigma_{\ovl{\rho}}$.  Recall from \eqref{Pi0 defn} in \cref{section: res csts} that $\Pi_0(\ovl{\rho})$ is the intersection of the kernels of all twist characters of $\ovl{\rho}$.  Define
\[
\Pi_1 \coloneqq \begin{cases}
\ker \eta_1 & \text{if } \ovl{\rho} \text{ is dihedral}\\
\Pi_0(\ovl{\rho}) & \text{else}.
\end{cases}
\]
Let $A_1$ be the subring of $A$ topologically generated by $t(\Pi_1)$.  Note that $\ovl{\rho}|_{\Pi_1}$ is absolutely irreducible by \cref{big image mult free} and \cref{rhobar irred over complement to kernel betat}.

\begin{proposition}\label{A' sub GMA}
There exists a well-adapted $(t, d)$-representation $\rho \colon \Pi \to \GL_2(A)$ such that $\rho|_{\Pi_1}$ takes values in $\GL_2(A_1)$ and such that $\rho$ is adapted to a regular element.    
\end{proposition}

\begin{proof}
With the exception of the well-adaptedness statement, the proof is well known since $\ovl{\rho}|_{\Pi_1}$ is absolutely irreducible.  Indeed, a theorem of Rouquier \cite[Theorem 5.1]{Rouquier96} and Nyssen \cite{Nyssen96} tells us that there are representations $\rho \colon \Pi \to \GL_2(A)$ and $\rho_1 \colon \Pi_1 \to \GL_2(A_1)$ such that $\tr \rho = t$ and $\tr \rho_1 = t|_{\Pi_1}$.  By a theorem of Carayol and Serre, $\rho|_{\Pi_1}$ and $\rho_1$ are conjugate by a matrix in $\GL_2(A)$ \cite[Th\'eor\`eme 1]{Carayol94}.  

For the well-adaptedness statement, let us first assume that $\ovl{\rho}$ is not dihedral.  Choose $\rho$ adapted to $g_0$ and $\rho_1$ adapted to $g_0^n$ with $g_0$ and $n$ as in \cref{choice of basis}.  Then the matrix $M \in \GL_2(A)$ such that $M^{-1}\rho|_{\Pi_1}M = \rho_1$ commutes with $\rho(g_0^n) = \bigl(\begin{smallmatrix} s(\lambda_0^n) & 0\\
0 & s(\mu_0^n)
\end{smallmatrix}\bigr) = \rho_1(g_0^n)$.  Since $\lambda_0^n \neq \mu_0^n$ by \cref{choice of basis}, it follows that $M$ must be diagonal.  In particular, $M$ commutes with $\rho(g_0)$.  Hence $M^{-1}\rho M$ is still adapted to $g_0$ and satisfies the properties in the statement of the proposition.    

The idea is similar when $\ovl{\rho}$ is dihedral, except we can no longer assume that $\rho_1$ is adapted to the $g_0 \in \Pi_0$ such that $\ovl{\rho}(g_0)$ generates the unique index-2 subgroup of the projective image of $\ovl{\rho}$, because $g_0$ may not be in $\Pi_1$.  Let $\rho$ be a well-adapted $(t,d)$-representation, say adapted to $g_0$ with $\rho(g_0) = \bigl(\begin{smallmatrix} 
s(\lambda_0) & 0\\
0 & s(\mu_0)
\end{smallmatrix}\bigr)$.  

Since $\ovl{\rho}$ is regular, it follows that $\ovl{\rho}|_{\Pi_0(\ovl{\rho})}$ is multiplicity free over $\E$ by \cref{dihedral regular implies mult free}.  Therefore, since $\rho$ is well adapted, the image of $\rho$ contains a matrix of the form $\bigl(\begin{smallmatrix} 
s(\lambda) & 0\\
0 & s(\mu)
\end{smallmatrix}\bigr)$ with $\lambda \neq \mu$ and $\lambda, \mu \in \E^\times$.  Let $h \in \Pi$ such that $\rho(h) = \bigl(\begin{smallmatrix}
s(\lambda) & 0\\
0 & s(\mu)
\end{smallmatrix} \bigr)$.  

We claim that $h \in \Pi_1$.  It suffices to prove that $h \in \Pi_0(\ovl{\rho})$ since $\Pi_0(\ovl{\rho}) = \Pi_0 \cap \ker \eta_1 \subset \ker \eta_1 = \Pi_1$.  Note that $h \in \Pi_0$ as $\ovl{\rho}(h)$ is diagonal.  By \cref{dihedral regular implies mult free}, $h \in \Pi_0(\ovl{\rho})$ if and only if the eigenvalues of $\ovl{\rho}(h)$ are in $\E^\times$.  But the eigenvalues  $\lambda, \mu$ of $\ovl{\rho}(h)$  were chosen to be in $\E^\times$.  Therefore $h \in \Pi_1$.

By \cite[Proposition 2.4.2]{Bellaiche18} we may assume that $\rho_1$ in the first paragraph of this proof is adapted to $h$.  Therefore the matrix $M \in \GL_2(A)$ such that $M^{-1}\rho|_{\Pi_1}M = \rho_1$ commutes with $\rho(h) = \bigl(\begin{smallmatrix} 
s(\lambda) & 0\\
0 & s(\mu)
\end{smallmatrix}\bigr) = \rho_1(h)$.  Since $\lambda \neq \mu$, it follows that $M$ is a diagonal matrix.  Note that the second property in \cref{well adapted} is unchanged by conjugation by a diagonal matrix.  Therefore $M^{-1}\rho M$ is still well adapted and satisfies the statement of the proposition.
\end{proof}

We recall that when $\overline{\rho}$ is dihedral, we may view elements in $B$ as elements of~$A$ (\cref{B and C in A}).

\begin{corollary}\label{I1 fixed by twists}
There exists a well-adapted $(t,d)$-representation $\rho \colon \Pi \to R^\times$ such that $I_1(\rho) \subseteq A^{\Sigma_t}$.  If $\ovl{\rho}$ is dihedral and $\sigma \in \Sigma_t$ such that $\sigma$ and $\ker \beta_t$ generate $\Sigma_t$, then we may assume furthermore that $B_1(\rho)$ is pointwise fixed by $\sigma$.
\end{corollary}

\begin{proof}
Let $\rho$ be the $(t,d)$-representation from \cref{A' sub GMA}.  Since the order of $\det \ovl{\rho}$ is a power of 2, it follows that $[\Pi \colon \Pi_0(\ovl{\rho})]$ is a power of 2.  Since $\Gamma$ is pro-$p$ and $p \neq 2$, it follows that $\Gamma \subseteq \rho(\Pi_0(\ovl{\rho})) \subseteq \GL_2(A_1)$.  Therefore $L_1(\rho) \subseteq \Sl_2(A_1)$, and so $I_1(\rho), B_1(\rho) \subseteq A_1$.

Let $(\sigma, \eta) \in \til \Sigma_t$ such that $\Pi_1 \subseteq \ker \eta$.  Then for all $g \in \Pi_1$ we have
\[
\act{\sigma}{t}(g) = \eta(g)t(g) = t(g),
\]  
and thus $A_1$ is contained in the subring of $A$ fixed by $\sigma$.  

If $\ker \beta_t = 1$, then every $(\sigma, \eta) \in \til \Sigma_t$ satisfies $\Pi_1 \subseteq \ker \eta$ by definition of $\Pi_1$.  Thus if $\ker \beta_t = 1$, then $A_1 \subseteq A^{\Sigma_t}$, and hence $I_1(\rho), B_1(\rho) \subseteq A^{\Sigma_t}$.

Now suppose that $\ovl{\rho}$ is dihedral and $\ker \beta_t \neq 1$.  Then half of the elements $(\sigma, \eta) \in \til \Sigma_t$ satisfy $\ker \eta \subseteq \ker \eta_1 = \Pi_1$, namely all those in the preimage under $\beta_t$ of the subgroup generated by $(\sigma_1, \eta_1)$ in $\til \Sigma_{\ovl{\rho}}$.  This proves the statement about $B_1(\rho)$ in the dihedral case.  To see that $I_1(\rho)$ is fixed by all conjugate self-twists, it remains to show that $I_1(\rho)$ is fixed by the nontrivial element in $\ker \beta_t$.  This follows from \cref{plus/minus parts}.  
\end{proof}

In light of \cref{I1 fixed by twists}, let us fix a well-adapted $(T, d)$-representation $\rho^{\univ} \colon \Pi \to \GL_2(\cA)$ such that $I_1(\rho^{\univ}) \subseteq \cA^{\Sigma_T}$.  Assume furthermore in the case when the projective image of $\ovl{\rho}$ is not dihedral that we have conjugated $\rho^{\univ}$ by the relevant diagonal element so that \cref{Bellaiche main theorem} applies to $\rho^{\univ}$, and thus to any quotient of $\rho^{\univ}$.  Recall from the commutative diagram following \cref{lifting arbitrary csts} that we have a natural reduction map $\beta_T \colon \Sigma_T \to \Sigma_{\overline{\rho}}$ that sends an automorphism $\sigma$ of $\mathcal{A}$ to the automorphism it induces on $\F = \mathcal{A}/\m_{\mathcal{A}}$.  In the case when $\ovl{\rho}$ is dihedral, we need to choose a complement to $\ker \beta_T$ in $\Sigma_T$, whose generator we will denote by $\nu$. 
We choose $\nu$ such that $(\nu, \eta_1) \in \til \Sigma_T$, where $\eta_1$ is the character fixed prior to \cref{A' sub GMA}.  By \cref{I1 fixed by twists}, we may and do assume that $B_1(\rho^{\univ})$ is fixed by $\nu$.    

The universal property of $(\cA,(T,d))$ gives a surjective $W(\F)$-algebra homomorphism $\alpha_t \colon \cA \to A$.  Let $\rho_t \coloneqq \alpha_t \circ \rho^{\univ} \colon \Pi \to \GL_2(A)$.  It is a $(t,d)$-representation such that $I_1(\rho_t) \subseteq A^{\Sigma_t}$ by the diagram following \cref{lifting arbitrary csts}.  Furthermore, if $\ovl{\rho}$ is dihedral and $\ker \beta_t = 1$, then $B_1(\rho_t) \subseteq A^{\Sigma_t}$ as well.  By the functoriality of Pink-Lie algebras with respect to quotient maps, we have that $\alpha_t(I_1(\rho^{\univ})) = I_1(\rho_t)$ and $\alpha_t(B_1(\rho^{\univ})) = B_1(\rho_t)$.  All of our theorems below will be specifically for this well-chosen representation $\rho_t$.  To ease notation, write $\rho = \rho_t$.

Recall that by \cref{Bellaiche main theorem}
\[
A = \begin{cases} 
\B_\rho(\F) + W(\F)B_1(\rho) & \text{if } \ovl{\rho} \text{ is dihedral}\\
\B_\rho(\F) & \text{else.}
\end{cases}
\] 
By \cref{plus/minus parts} and the fact that $B_1(\rho) \subseteq A^{\Sigma_t}$ if $\ovl{\rho}$ is dihedral and $\ker \beta_t = 1$,  it follows that 
\[
A^{\Sigma_t} = \begin{cases}
\B_\rho(\F^{\beta_t(\Sigma_t)}) + W(\F^{\beta_t(\Sigma_t)})B_1(\rho) & \text{if } \ovl{\rho} \text{ dihedral and} \ker \beta_t = 1\\
\B_\rho(\F^{\beta_t(\Sigma_t)}) & \text{else.}
\end{cases}
\]  
We therefore define 
\[
J = J(\rho) \coloneqq \begin{cases}
W(\E)I_1(\rho) + W(\E)I_1(\rho)^2 + W(\E)B_1(\rho) & \text{if } \ovl{\rho} \text{ is dihedral and} \ker \beta_t = 1\\
W(\E)I_1(\rho) + W(\E)I_1(\rho)^2 & \text{else.}
\end{cases}
\]
We claim that $J \subset \m$ is a multiplicatively closed $W(\E)$-module by \cref{Bellaiche main theorem}.  The key is to note that, since $\ovl{\rho}$ is regular and $\rho$ is well adapted, it follows that Bella\"iche's field $\F_q$ from p.\,\pageref{definition of Fq} is contained in $\E$.  Therefore it follows from \cref{Bellaiche main theorem} that $(W(\E)I_1(\rho))^3 \subseteq W(\E)I_1(\rho)$ and $W(\E)I_1(\rho)B_1(\rho) \subseteq W(\E)B_1(\rho)$ and $(W(\E)B_1(\rho))^2 \subseteq W(\E)I_1(\rho)$, which proves that $J$ is multiplicatively closed.  Define
\[
\frA \coloneqq W(\F) + W(\F)J.
\]
We have $\frA = A$ unless $1 \neq \ker \beta_t$, in which case $\frA = A^+$ by \cref{plus/minus parts}. 

\begin{remark}\label{the strategy}
The rings $W(\E) + J$ and $A^{\Sigma_t}$ differ only in their constants, $W(\E)$ versus $W(\F^{\beta_t(\Sigma_t)})$.  Furthermore, $W(\E) + J$ is often equal to $\B_\rho(\E)$, and the goal of this section is to relate $\B_\rho(\E)$ with $A^{\Sigma_t}$.  Assume for a moment that $W(\E) + J = \B_\rho(\E)$.  Then the difference between $\B_\rho(\E)$ and $A^{\Sigma_t}$ is entirely governed by understanding which elements of $\Sigma_{\ovl{\rho}}$ lift to elements in $\Sigma_t$ under $\beta_t$.  In particular, when there are elements in $\Sigma_{\ovl{\rho}}$ that do not lift to $\Sigma_t$, we will be interested in writing the extra elements in $W(\F^{\beta_t(\Sigma_t)})$ as quotients of elements in $J$ to show that $Q(\B_\rho(\E)) = Q(A^{\Sigma_t})$.
\end{remark}

\subsection{Lifting conjugate self-twists to $\frA$}\label{lifting csts to frA}
In \cref{lifting csts to frA} we study a condition on $J$, called \textit{smallness} (\cref{J small}), that dictates which conjugate self-twists of $\ovl{\rho}$ lift to conjugate self-twists of $(t,d)$.  This study culminates in \cref{small implies surjective}.  As a consequence, we prove in \cref{small non domain case} that under such a smallness assumption, $A^{\Sigma_t} = \B_\rho(\E)$.  The reader is advised that, with the exception of the motivational remark following \cref{J small}, the assumption that $A$ is a domain is never used in \cref{lifting csts to frA}.

Throughout \cref{lifting csts to frA}, fix a subgroup $\Sigma \subseteq \Sigma_{\ovl{\rho}}$, and let $\F' \coloneqq \F^\Sigma$.  Write $W \coloneqq W(\F)$ and $W' \coloneqq W(\F')$.  For an arbitrary ring $\mathcal{R}$ and a finite group $X$ of ring automorphisms of $\mathcal{R}$, for any $\varphi \in \Hom(X, \mathcal{R}^\times)$, we write
\[
\mathcal{R}^\varphi \coloneqq \{s \in \mathcal{R} \colon \act{\sigma}{s} = \varphi(\sigma)s, \forall \sigma \in X\}.  
\]  
As explained at the beginning of \cref{conceptual interpretation of Bellaiche ring}, we assume that $[\F \colon \E]$ is a power of 2.  By \cref{condition * satisfied} we may apply \cref{automorphisms to gradings} to conclude that $\F = \oplus_{\varphi \in \Sigma^*} \F^\varphi$, where $\Sigma^* \coloneqq \Hom(\Sigma, \F^\times)$.  Note that since $\Sigma = \Gal(W/W')$, it follows that this decomposition lifts to $W$.  More precisely, viewing elements of $\Sigma$ as automorphisms of $W$ and elements of $\Sigma^*$ as valued in $W^\times$ by composing with the Teichm\"uller map, we can define $W^\varphi$ for each $\varphi \in \Sigma^*$.  Then \cref{automorphisms to gradings} gives $W = \oplus_{\varphi \in \Sigma^*} W^\varphi$. 

For all $\varphi \in \Sigma$, define 
\[
\frA(\varphi) \coloneqq W^\varphi + W^{\varphi}J,
\]
where $W^{\varphi}J \coloneqq \{\sum_i \alpha_i j_i | \alpha_i \in W^{\varphi}, j_i \in J\}$.  Since $\frA = W + WJ$ it follows immediately that $\frA = \sum_{\varphi \in \Sigma^*} \frA(\varphi)$.  We will be interested in understanding when this sum is direct, because in that case we will show that it is possible to find lifts of elements of $\Sigma$ in $\Sigma_t$.  If $\Aa$ is an ideal of $\frA$ and $\varphi \in \Sigma^*$, let $\Aa(\varphi) \coloneqq \frA(\varphi) \cap \Aa$ and let $(\frA/\Aa)(\varphi) \subset \frA/\Aa$ be the image of $\frA(\varphi)$ under the natural projection $\frA \to \frA/\Aa$.

\begin{lemma}\label{sub and quotient direct sums}
The following are equivalent:
\begin{enumerate}
\item $\frA = \oplus_{\varphi \in \Sigma^*} \frA(\varphi)$.
\item For every $\frA$-ideal $\Aa$ such that $\Aa = \oplus_{\varphi \in \Sigma^*} \Aa(\varphi)$, we have $\frA/\Aa = \oplus_{\varphi \in \Sigma^*} (\frA/\Aa)(\varphi)$. Furthermore, there exists at least one such ideal $\Aa$.
\item There exists an $\frA$-ideal $\Aa$ such that $\Aa = \oplus_{\varphi \in \Sigma^*} \Aa(\varphi)$ and $\frA/\Aa = \oplus_{\varphi \in \Sigma^*} (\frA/\Aa)(\varphi)$.
\end{enumerate}
\end{lemma}

\begin{proof}
First we show that (1) implies (2).  We can take $\Aa = 0$ for the existence statement in (2).  Now suppose that $\Aa$ is an $\frA$-ideal such that $\Aa = \oplus_{\varphi \in \Sigma^*} \Aa(\varphi)$.  If $\sum_{\varphi \in \Sigma^*} \ovl{a}_\varphi = 0 \in \frA/\Aa$ with each $\ovl{a}_\varphi \in (\frA/\Aa)(\varphi)$, then letting $a_\varphi \in \frA(\varphi)$ be a lift of $\ovl{a}_\varphi$, we see that $\sum_{\varphi \in \Sigma^*} a_\varphi \in \Aa = \oplus_{\varphi \in \Sigma^*} \Aa(\varphi)$.  Thus, there are $\alpha_\varphi \in \Aa(\varphi)$ such that $\sum_{\varphi \in \Sigma^*} a_\varphi = \sum_{\varphi \in \Sigma^*} \alpha_\varphi$.  Since $\frA = \oplus_{\varphi \in \Sigma^*} \frA(\varphi)$, it follows that $a_\varphi = \alpha_\varphi$ for all $\varphi \in \Sigma^*$.  Thus $\ovl{a}_\varphi = 0 \in \frA/\Aa$ for all $\varphi \in \Sigma^*$ and hence $\frA/\Aa = \oplus_{\varphi \in \Sigma^*} (\frA/\Aa)(\varphi)$.

The fact that (2) implies (3) is trivial.

To see that (3) implies (1), suppose that $\Aa$ is an $\frA$-ideal such that $\Aa = \oplus_{\varphi \in \Sigma^*} \Aa(\varphi)$ and \mbox{$\frA/\Aa = \oplus_{\varphi \in \Sigma^*} (\frA/\Aa)(\varphi)$}.  For each $\varphi \in \Sigma^*$ fix a set $S_\varphi \subset \frA(\varphi)$ of representatives of $(\frA/\Aa)(\varphi)$ such that \mbox{$0 \in S_\varphi$}.  Suppose that $\sum_{\varphi \in \Sigma^*} a_\varphi = 0$ with $a_\varphi \in \frA(\varphi)$.  Then there is a unique way \mbox{to write each $a_\varphi$ as}
\[
a_\varphi = s_\varphi + \alpha_\varphi
\]
with $s_\varphi \in S_\varphi$ and $\alpha_\varphi \in \Aa(\varphi)$.  Modulo $\Aa$, we see that
\[
\sum_{\varphi \in \Sigma^*} \ovl{s}_\varphi = 0.
\]
Since $\frA/\Aa = \oplus_{\varphi \in \Sigma^*} (\frA/\Aa)(\varphi)$, it follows that $\ovl{s}_\varphi = 0$ for all $\varphi$.  As $0 \in S_\varphi$, it follows that $s_\varphi = 0$ for all $\varphi \in \Sigma^*$.  Therefore $a_\varphi = \alpha_\varphi \in \Aa(\varphi)$.  As $\Aa = \oplus_{\varphi \in \Sigma^*} \Aa(\varphi)$, it follows that each $a_\varphi = 0$.  Thus $\frA = \oplus_{\varphi \in \Sigma^*} \frA(\varphi)$.
\end{proof}

\begin{definition}\label{J small}
Let $\bL_2\subset\bL_1$ be subfields of $\F$.  We say that $J$ is \textit{small with respect to} $\bL_1/\bL_2$ if
\[
\ker(W(\bL_1) \otimes_{W(\bL_2)} W(\bL_2)J \to W(\bL_1)J) = 0,
\]
where the map is  multiplication inside $\frA$.  Otherwise, we say that $J$ is \textit{big with respect to} $\bL_1/\bL_2$.
\end{definition}

To motivate \cref{J small}, recall from \cref{the strategy} that we need to be able to write elements of $W(\F^{\beta_t(\Sigma_t)})$ as quotients of elements in $J$ whenever $\F^{\beta_t(\Sigma_t)} \neq \E$.  Suppose that $\bL_2 = \E, \bL_1 = \F^{\beta_t(\Sigma_t)}$, and $[\bL_1 \colon \bL_2] = 2$.  Write $\bL_1 = \bL_2(\alpha)$.  Then $W(\bL_1) = W(\bL_2) \oplus s(\alpha)W(\bL_2)$ and so 
\[
W(\bL_1) \otimes_{W(\bL_2)} W(\bL_2)J = W(\bL_2)J \oplus (s(\alpha)W(\bL_2) \otimes_{W(\bL_2)} W(\bL_2)J).
\]
If $J$ is big with respect to $\bL_1/\bL_2$, then we can find $x, y \in W(\bL_2)J \setminus \{0\}$ such that $x + s(\alpha)y = 0$.  Thus $s(\alpha) = x/y$, and hence $W(\bL_1)$ is in the field of fractions of any domain containing $W(\bL_2)J$.  In contrast, the following proposition shows that when $J$ is small with respect to $\F/\F'$, elements of $\Sigma$ can be lifted to automorphisms of $\frA$.

\begin{proposition}\label{A is direct sum}
If $J$ is small with respect to $\F/\F'$, then $\frA = \oplus_{\varphi \in \Sigma^*} \frA(\varphi)$.  In this case, every $\ovl{\sigma} \in \Sigma$ can be lifted to an automorphism $\sigma$ of $\frA$ such that $\sigma$ acts trivially on $J$, and a lift with this property is unique.
\end{proposition}

\begin{proof}
Note that $\Aa \coloneqq WJ$ is an $\frA$-ideal since $\frA = W + WJ$ and $J$ is multiplicatively closed as discussed prior to \cref{the strategy}.  The assumption that $J$ is small with respect to $\F/\F'$ implies that $WJ = \oplus_{\varphi \in \Sigma^*} W^{\varphi}J$.  Indeed, 
\[
\bigoplus_{\varphi \in \Sigma^*}(W^\varphi \otimes_{W'} W'J) = \bigl(\bigoplus_{\varphi \in \Sigma^*} W^\varphi \bigr) \otimes_{W'} W'J = W \otimes_{W'} W'J \hookrightarrow WJ,
\]
and the image of $W^\varphi \otimes_{W'} W'J$ is exactly $W^\varphi J$.  Since $WJ = \sum_{\varphi \in \Sigma^*} W^\varphi J$, it follows that the multiplication map is an isomorphism and thus $WJ$ is graded by $\Sigma^*$.  Note that $\Aa(\varphi) = W^\varphi J$, so $\Aa = \oplus_{\varphi \in \Sigma^*} \Aa(\varphi)$.

By \cref{sub and quotient direct sums}, for the first statement of the proposition it suffices to show that \mbox{$\frA/WJ = \oplus_{\varphi \in \Sigma^*} (\frA/WJ)(\varphi)$}.  Note that
\[
\frA/WJ = (W + WJ)/WJ \cong W/(W \cap WJ)
\]
and $W \cap WJ$ is a closed $W$-submodule of $pW$ since $J \subseteq \m_{\frA}$.  Thus we have $W \cap WJ = p^nW$ and $\frA/WJ \cong W/p^nW$ for some $1 \leq n \leq \infty$, where $p^\infty W \coloneqq \{\infty\}$.  Since $W$ is graded by $\Sigma^*$, it follows from \cref{sub and quotient direct sums} that $W/p^nW$ is graded by $\Sigma^*$ as well.  Therefore $\frA = \oplus_{\varphi \in \Sigma^*} \frA(\varphi)$.

For the second statement, we let $\sigma$ act by $W(\ovl{\sigma})$ on $W$ and trivially on $J$.  The only question is to verify that this is well defined.  Since $\frA = \oplus_{\varphi \in \Sigma^*} \frA(\varphi)$, it suffices to show that $\sigma$ is well defined on each $\frA(\varphi)$.  That is, we must show
\[
\sum_{i = 1}^n \alpha_i j_i = 0 \implies \sum_{i = 1}^n W(\ovl{\sigma})(\alpha_i)j_i = 0,
\]
where $\alpha_i \in W^\varphi, j_i \in J$.  Since $J$ is small with respect to $\F/\F'$, $\sum_{i = 1}^n \alpha_i j_i = 0$ implies that $\sum_{i = 1}^n \alpha_i \otimes j_i = 0\in W\otimes_{W'}W'J$.  Since $\alpha_i \in W^{\varphi}$, we know that $W(\ovl{\sigma})(\alpha_i) = s(\varphi(\ovl{\sigma}))\alpha_i$ for all $i$.  Hence
\[
0 = \sum_{i = 1}^n \alpha_i \otimes j_i \implies 0 = s(\varphi(\ovl{\sigma})) \sum_{i = 1}^n \alpha_i \otimes j_i.
\]
Therefore $\sum_{i = 1}^n W(\ovl{\sigma})(\alpha_i)j_i = 0$, since it is the image of $s(\varphi(\ovl{\sigma})) \sum_{i = 1}^n \alpha_i \otimes j_i$ under \mbox{$W\otimes_{W'}W'J\to WJ$}.
\end{proof}

Now that we have lifted elements of $\Sigma$ to automorphisms of $\frA$ under the smallness assumption, we would like to verify that the lifts are conjugate self-twists of $(t,d)$ when $\frA = A$ and that they come from conjugate self-twists when $\frA = A^+$.  (Recall that $A^+$ is only defined when $\ker \beta_t$ is nontrivial; see \cref {plus/minus parts}.)

\begin{theorem}\label{small implies surjective}
If $J$ is small with respect to $\F/\F'$ then $\Sigma$ is contained in the image of $\beta_t \colon \Sigma_t \to \Sigma_{\ovl{\rho}}$. Furthermore, every lift of $\ovl{\sigma} \in \Sigma$ to $\Sigma_t$ acts trivially on $J$. 
\end{theorem}

\begin{proof}
Fix $\ovl{\sigma} \in \Sigma$.  By \cref{A is direct sum}, there is a unique $\sigma \in \Aut \frA$ that acts as $W(\ovl{\sigma})$ on $W$ and fixes $J$.  If $\ker \beta_t = 1$, then $\frA = A$.  If $\ker \beta_t \neq 1$, then $\frA = A^+$ and we need to extend $\sigma$ to $A = A^+ \oplus A^-$.  We do this by declaring that $\sigma$ fixes $A^-$; we will still denote the automorphism of $A$ by $\sigma$.  We already know that $\sigma$ acts trivially on $J$, so it is enough to prove that $\sigma\in\Sigma_t$.

Our strategy is to show that $\sigma$ comes from an appropriate element of $\Sigma_T$. More precisely, we claim that there is some $(\tilde{\sigma},\eta) \in \til \Sigma_T$ such that $\sigma \circ \alpha_t = \alpha_t \circ \tilde{\sigma}$, where $\alpha_t \colon \cA \to A$ is the $W$-algebra homomorphism given by universality. If this is true, then for all $g \in \Pi$ we have
\[
\act{\sigma}{t}(g) = \sigma \circ \alpha_t(T(g)) = \alpha_t \circ \tilde{\sigma}(T(g)) = \alpha_t(\eta(g)T(g)) = \eta(g)\alpha_t(T(g)) = \eta(g)t(g)
\]
since $\alpha_t$ is a $W$-algebra homomorphism and $\eta(g) \in W$.  Thus $\sigma \in \Sigma_t$.

First suppose that $\ovl{\rho}$ is not dihedral.  Then there is a unique lift $\tilde{\sigma}$ of $\ovl{\sigma}$ in $\Sigma_T$ by \cref{precise structure of Sigma} and \cref{lifting residual csts}.  By \cref{lifting residual csts}, we know that $\tilde{\sigma}$ acts as $W(\ovl{\sigma})$ on the image of $W$ in $\cA$.  Furthermore, $\tilde{\sigma}$ acts trivially on $I_1(\rho^{\univ})$ by our fixed choice of $\rho^{\univ}$ after \cref{I1 fixed by twists}.  Since $\ovl{\rho}$ is not dihedral, it follows that $\cA = W + WJ(\rho^{\univ})$.  By the construction of $\rho^{\univ}$ and $\rho_t$, we have $\alpha_t(I_1(\rho^{\univ})) = I_1(\rho)$ and thus $\alpha_t(J(\rho^{\univ})) = J(\rho) = J$.  Recall that $\sigma$ acts trivially on $J$.  Both $\tilde{\sigma}$ and $\sigma$ act on $W$ by $W(\ovl{\sigma})$.  Thus for any $\sum_{i = 1}^n a_ix_i$ with $a_i \in W, x_i \in J(\rho^{\univ}) \cup \{1\}$, we have 
\[
\sigma \circ \alpha_t\left(\sum_{i = 1}^n a_i x_i\right) = \sum_{i = 1}^n W(\ovl{\sigma})(a_i)\alpha_t(x_i) = \alpha_t \circ \tilde{\sigma}\left(\sum_{i=1}^n a_ix_i\right).  
\]

If $\ovl{\rho}$ is dihedral, then there are two lifts of $\ovl{\sigma}$ in $\Sigma_T$ by \cref{precise structure of Sigma} and \cref{lifting residual csts}.  One acts on $\cA^-$ by $+1$ and the other acts by $-1$ by \cref{plus/minus parts} and since we chose $\rho^{\univ}$ such that $B_1(\rho^{\univ})$ is fixed by $\nu$, which generates a complement of $\ker \beta_T$.  Let $\tilde{\sigma} \in \Sigma_T$ be the lift of $\ovl{\sigma}$ that is in $\langle \nu \rangle$.  Thus $\tilde{\sigma}$ acts trivially on $J(\rho^{\univ})$ and $B_1(\rho^{\univ})$.  Then an argument similar to that in the previous paragraph shows that $\sigma \circ \alpha_t = \alpha_t \circ \tilde{\sigma}$.
\end{proof}

\begin{corollary}\label{small non domain case}
If $J$ is small with respect to $\F/\E$ then $A^{\Sigma_t} = W(\E) + J$.  Suppose furthermore that either $\ovl{\rho}$ is not dihedral or $\ker \beta_t \neq 1$.  Then $A^{\Sigma_t} = \B_\rho(\E)$.  
\end{corollary}

\begin{proof}
By \cref{small implies surjective} applied to $\Sigma=\Sigma_t$, the map $\beta_t$ is a surjection and $\Sigma_t$ acts trivially on $J$. If $\ker \beta_t = 1$, then $A = \frA = W + WJ$, so $A^{\Sigma_t} = W(\E) + J$.

If $\ker \beta_t \neq 1$, then $A = W + WI_1(\rho) + WI_1(\rho)^2 + WB_1(\rho)$ and $J = W(\E)I_1(\rho) + W(\E)I_1(\rho)^2$.  Note that $A^{\Sigma_t} \subseteq A^+ = W + WJ$ since the nontrivial element in $\ker \beta_t$ acts by $-1$ on $B_1(\rho)$ by \cref{plus/minus parts}.  As above, we have that 
\[
A^{\Sigma_t} = (W + WJ)^{\Sigma_t} = W(\E) + J.
\]

The last sentence in the statement of the corollary follows from the definition of $J$.
\end{proof}

\begin{remark}
Note that none of the arguments in \cref{lifting csts to frA} require that $A$ is a domain.  In particular, when $J$ is small with respect to $\F/\E$ and either $\ovl{\rho}$ is not dihedral or $\ker \beta_t \neq 1$, \cref{small non domain case} gives a conceptual interpretation of the ring $\B_\rho(\E)$.
\end{remark}

\subsection{When $J$ is big with respect to $\F/\E$}\label{bigness subsection}
\cref{small non domain case} requires the assumption that $J$ is small with respect to $\F/\E$.  We do not always expect this to be true.  The purpose of \cref{bigness subsection} is to show that  $A^{\Sigma_t}$ and $\B_\rho(\E)$ have the same fraction field and $A^{\Sigma_t}$ is a finite type $\B_\rho(\E)$-module even without the smallness assumption.  This is done in \cref{final conclusion}, although the two key inputs to that theorem are \cref{no more lifting} and \cref{induction step}.  Then we can apply \cref{fg plus same quotient field implies ideal containment} and \cref{formal fullness corollary} to conclude that $\rho_t$, and thus $(t,d)$, is $A^{\Sigma_t}$-full in \cref{final goal}.

The discussion following \cref{J small} shows why one may expect to get $Q(A^{\Sigma_t}) = Q(\B_\rho(\E))$ when smallness fails and $[\F^{\beta_t(\Sigma_t)} \colon \E] = 2$.  Unfortunately, the assumption that $[\F^{\beta_t(\Sigma_t)} \colon \E] = 2$ is rather critical to that argument.  This is the primary reason we insist that $[\F \colon \E]$ be a power of 2 throughout this section.  It allows us to split up the extension $\F^{\beta_t(\Sigma_t)}/\E$ into a series of quadratic extensions, and thus we can apply the argument following \cref{J small} inductively.  This is the essential idea of the argument; we now prepare some notation to formalize it.

Write $[\F^{\beta_t(\Sigma_t)} \colon \E] = 2^n$ for some $n \geq 0$.  
For integers $0 \leq i \leq n$, let $\E_i$ be the unique extension of~$\E$ of degree $2^i$.  In particular, $\E_0 = \E$ and $\E_n = \F^{\beta_t(\Sigma_t)}$, and $[\E_i \colon \E_{i - 1}] = 2$ for all $1 \leq i \leq n$.  For $0 \leq i \leq n$, let $W_i$ denote the image of $W(\E_i)$ in $\frA$.  Define 
\[
\frA_i \coloneqq W_i + W_iJ \subseteq \frA.
\]
In particular, $\frA_0 = W(\E) + J$ and $\frA_n = W(\F^{\beta_t(\Sigma_t)}) + W(\F^{\beta_t(\Sigma_t)})J$.  Since $A$ is a domain, so are all of the $\frA_i$, and we write $Q(\frA_i)$ for the field of fractions of $\frA_i$.

In the case when $\frA = A^+$, there is a 2-to-1 group homomorphism $\Sigma_t \to \Aut A^+$ given by restricting elements of $\Sigma_t$ to $A^+$.  Let $\Sigma_t(\frA)$ denote the image of this map when $\frA = A^+$, and otherwise (that is, whenever $\ker \beta_t = 1$) let $\Sigma_t(\frA) = \Sigma_t$.  In either case we can identify $\Sigma_t(\frA)$ with a subgroup of $\Sigma_{\ovl{\rho}}$ via $\beta_t$, and we have $\E_n = \F^{\beta_t(\Sigma_t(\frA))}$.  We write $\Sigma_t(\frA)^*\coloneqq \Hom(\Sigma_t(\frA), \frA^\times)$.

We begin with two preliminary lemmas about the relationship between smallness and the $\E_i$. 

\begin{lemma}\label{descending to FSigmat}
We have that $J$ is small with respect to $\F/\E_n$; that is, $\ker(W \otimes_{W_n} W_nJ \to WJ) = 0$.
\end{lemma}

\begin{proof}
Recall that $WJ$ is an $\frA$-ideal that is stable under the action of $\Sigma_t(\frA)$ since $\Sigma_t$ fixes $J$ by the construction of $\rho_t$.  By \cref{condition * satisfied}, we can apply \cref{automorphisms to gradings} with $X = \Sigma_t(\frA)$.  Therefore
\[
WJ = \bigoplus_{\varphi \in \Sigma_t(\frA)^*} (WJ)^\varphi,
\]
where $(WJ)^\varphi \coloneqq \{x \in WJ \colon \act{\sigma}{x} = \varphi(\sigma)x, \forall \sigma \in \Sigma_t(\frA)\}$.  Recall that $W^\varphi J$ was defined prior to \cref{sub and quotient direct sums}.  We claim that
\[
(WJ)^\varphi = W^\varphi J.
\]
Clearly $(WJ)^\varphi \supseteq W^\varphi J$ since $\Sigma_t(\frA)$ acts trivially on $J$.  On the other hand,
\[
WJ = \bigoplus_{\varphi \in \Sigma_t(\frA)^*} (WJ)^\varphi = \sum_{\varphi \in \Sigma_t(\frA)^*} (WJ)^\varphi \supseteq \sum_{\varphi \in \Sigma_t(\frA)^*} W^\varphi J = WJ,
\]
so we must have equality.

For each $\varphi \in \Sigma_t(\frA)^*$, choose $x_\varphi \in \F^\varphi \setminus \{0\}$.  Then $\{s(x_\varphi) \colon \varphi \in \Sigma_t(\frA)^*\}$ is a $W_n$-basis for $W$.  Thus we have
\[
W \otimes_{W_n} W_nJ = \bigoplus_{\varphi \in \Sigma_t(\frA)^*} W_n s(x_\varphi) \otimes_{W_n} W_nJ.
\]

If $x \in \ker(W \otimes_{W_n} W_nJ \to WJ)$, then we can write
\[
x = \sum_{\varphi \in \Sigma_t(\frA)^*} s(x_\varphi) \otimes y_\varphi
\]
for some $y_\varphi \in W_nJ$.  Then we have
\[
0 = \sum_{\varphi \in \Sigma_t(\frA)^*} s(x_\varphi)y_\varphi
\]
and $s(x_\varphi)y_\varphi \in W^\varphi J$.  Since $WJ = \oplus_{\varphi \in \Sigma_t(\frA)^*} W^\varphi J$, it follows that each $s(x_\varphi)y_\varphi = 0$.  As $\frA$ is a domain and $s(x_\varphi) \neq 0$, it follows that $y_\varphi = 0$ for all $\varphi \in \Sigma_t(\frA)^*$.
\end{proof}

\begin{lemma}\label{behavior of kernels under tensor product}
We have
\[
\ker(W \otimes_{W_{n-1}} W_{n-1}J \to WJ) = W \otimes_{W_n} \ker(W_n \otimes_{W_{n-1}} W_{n-1}J \to W_nJ).
\]
\end{lemma}

\begin{proof}
Let $K \coloneqq \ker(W_n \otimes_{W_{n-1}} W_{n-1}J \to W_nJ)$.  We have an exact sequence of $W_n$-modules
\[
0 \to K \to W_n \otimes_{W_{n-1}} W_{n-1}J \to W_nJ \to 0.
\]
Since $W$ is free over $W_n$, tensoring with $W$ over $W_n$ gives an exact sequence
\[
0 \to W \otimes_{W_n} K \to W \otimes_{W_{n-1}} W_{n-1}J \to W \otimes_{W_n} W_nJ \to 0.
\]
We can identify the last nonzero term in this sequence with $WJ$ by \cref{descending to FSigmat}.\\  Thus \mbox{$W \otimes_{W_n} K = \ker(W \otimes_{W_{n-1}} W_{n-1}J \to WJ)$}.
\end{proof}

\begin{proposition}\label{no more lifting}
If $J$ is big with respect to $\F/\E$ then $Q(\frA_n) = Q(\frA_{n-1})$.
\end{proposition}

\begin{proof}
We claim that $J$ is big with respect to $\E_n/\E_{n-1}$.  Indeed, if $J$ were small with respect to $\E_n/\E_{n-1}$, then $J$ would be small with respect to $\F/\E_{n-1}$ by \cref{behavior of kernels under tensor product}.  Therefore we could apply \cref{small implies surjective} with $\Sigma = \Gal(\F/\E_{n-1})$, which implies that $\E_n \subseteq \E_{n-1}$, a contradiction.  

Let $\{1, \alpha\}$ be an $\E_{n-1}$-basis for $\E_n$.  Then $\{1, s(\alpha)\}$ is a $W_{n-1}$-basis for $W_n$ and so
\[
W_n \otimes_{W_{n-1}} W_{n-1}J = (W_{n-1} \otimes_{W_{n-1}} W_{n-1}J) \oplus (W_{n-1} s(\alpha) \otimes_{W_{n-1}} W_{n-1}J).
\]
Since $J$ is big with respect to $\E_n/\E_{n-1}$, there exist $x, y \in W_{n-1}J \setminus \{0\}$ such that
\[
x + s(\alpha)y = 0.
\]
Thus, $s(\alpha) = -x/y \in Q(\frA_{n-1})$.  It follows that $W_n \subset Q(\frA_{n-1})$ and hence $Q(\frA_n) = Q(\frA_{n-1})$.
\end{proof}

Finally, we descend from $Q(\frA_n)$ to $Q(\frA_0)$ by induction on $n$.  

\begin{proposition}\label{induction step}
For all $2 \leq k \leq n$, if $Q(\frA_k) = Q(\frA_{k-1})$ then $Q(\frA_{k-1}) = Q(\frA_{k-2})$.  In particular, if $J$ is big with respect to $\F/\E$, then $Q(A^{\Sigma_t}) = Q(W(\E) + J)$.
\end{proposition}

\begin{proof}
Note that for any $k \geq 1$ we have $Q(\frA_k) = Q(\frA_{k-1})$ if and only if $W_k \subseteq Q(\frA_{k-1})$.  Assume that $Q(\frA_k) = Q(\frA_{k-1})$ for some $k$, $2 \leq k \leq n$.  Choose $\ovl{\alpha} \in \E_{k-2}, \ovl{\beta} \in \E_{k-1}$ such that \mbox{$\E_{k-1} = \E_{k-2}(\sqrt{\ovl{\alpha}})$} and $\E_k = \E_{k-1}(\sqrt{\ovl{\beta}})$.  Define $\alpha \coloneqq s(\ovl{\alpha})$ and $\beta \coloneqq s(\ovl{\beta})$, so $W_{k-1} = W_{k-2}(\sqrt{\alpha})$ and $W_k = W_{k-1}(\sqrt{\beta})$.  It suffices to show that $\sqrt{\alpha} \in Q(\frA_{k-2})$.  

Since $Q(\frA_k) = Q(\frA_{k-1})$, we can write $\sqrt{\beta} = x/y$ with $x,y \in \frA_{k-1} \setminus \{0\}$.  By multiplying $x$ and~$y$ by any nonzero element of $W_{k-1}J$, we may assume that $x, y \in W_{k-1}J \setminus \{0\}$.  

Note that we can write $y = i_1 + \sqrt{\alpha}i_2$ with $i_1, i_2 \in W_{k-2}J$.  If $i_1 - \sqrt{\alpha}i_2 \neq 0$, then by multiplying $x$ and $y$ by $i_1 - \sqrt{\alpha}i_2$, we may assume that $y \in W_{k-2}J \setminus \{0\}$.  If $i_1 - \sqrt{\alpha}i_2 = 0$ and $y\notin W_{k-2}J$, then we must have $i_2 \neq 0$ since $y \neq 0$ and $\sqrt{\alpha} = i_1/i_2 \in Q(\frA_{k-2})$, as desired.  We assume henceforth that $y \in W_{k-2}J$.

Write $x = a + b\sqrt{\alpha}$ for some $a, b \in W_{k-2}J$.  Then we have $y\sqrt{\beta} = a + b\sqrt{\alpha}$ and thus
\begin{equation}\label{squaring beta}
y^2\beta = a^2 + \alpha b^2 + 2ab\sqrt{\alpha}.
\end{equation}
Since $\beta \in W_{k-1}$, we may write $\beta = e + f\sqrt{\alpha}$ for some $e, f \in W_{k-2}$.  Note that $f \not\equiv 0 \bmod p$ since $[\E_k \colon \E_{k-1}] = 2$ and $\E_k = \E_{k-1}(\sqrt{\ovl{\beta}})$.  Substituting this into equation \eqref{squaring beta}, we see that
\[
(y^2f - 2ab)\sqrt{\alpha} = a^2 + \alpha b^2 - y^2e \in W_{k-2}J.
\]
Note that $y^2f - 2ab \in W_{k-2}J$ since all of $y, f, a, b \in W_{k-2}J$.  If $y^2f - 2ab \neq 0$, then we can conclude that $\sqrt{\alpha} \in Q(\frA_{k-2})$ as desired.

Henceforth, assume that $y^2f = 2ab$.  Then we also have $y^2e = a^2+\alpha b^2$.  Thus $2f^{-1}eab = a^2 + \alpha b^2$.  Note that $a, b \neq 0$ since $2ab = y^2f$ and we know $y, f \neq 0$.  Then we have
\[
2f^{-1}e = \frac{a}{b} + \alpha\frac{b}{a}.
\]
Therefore $\frac{a}{b}$ is a root of $t^2 - 2f^{-1}et + \alpha \in W_{k-2}[t]$.  The discriminant of this polynomial is \mbox{$4(f^{-2}e^2 - \alpha)$}.  We claim that $W_{k-1} \supseteq W_{k-2}(\sqrt{f^{-2}e^2 - \alpha})$.  Note that since $\beta = s(\ovl{\beta})$ and $s$ is multiplicative we have that $e^2 - f^2\alpha = s(N_{\E_{k-1}/\E_{k-2}}(\ovl{\beta}))$, and therefore \mbox{$\sqrt{e^2 - f^2\alpha} = s(\sqrt{N_{\E_{k-1}/\E_{k-2}}(\ovl{\beta})}) \in s(\E_{k-1}^\times)$}.  Therefore $\frac{a}{b} \in W_{k-1}$.  Thus we can write $\sqrt{\beta} = \frac{x/b}{y/b} = \frac{\frac{a}{b} + \sqrt{\alpha}}{y/b}$, and so
\[
\frac{y}{b} = \left(\frac{a}{b} + \sqrt{\alpha}\right)\sqrt{\beta}^{-1}.
\]
Note that $\frac{a}{b} + \sqrt{\alpha} \neq 0$ since $y \neq 0$.  It follows that $(\frac{a}{b} + \sqrt{\alpha})\sqrt{\beta}^{-1}$ generates $W_k$ over $W_{k-1}$ since $\frac{a}{b} + \sqrt{\alpha} \in W_{k-1}$ and $\sqrt{\beta}$ generates $W_k$ over $W_{k-1}$.  Thus $(\frac{a}{b} + \sqrt{\alpha})\sqrt{\beta}^{-1} = \frac{y}{b} \in Q(\frA_{k-2})$ and so $W_k \subset Q(\frA_{k-2})$.  Therefore $Q(\frA_k) = Q(\frA_{k-2})$.

For the second statement of the proposition, note that by \cref{no more lifting} we have \mbox{$Q(\frA_n) = Q(\frA_0)$} for all $0 \leq k \leq n$.  We have $\frA_0 = W(\E) +  J$ by definition.  Since $\Sigma_t$ acts trivially on $J$, it follows that 
\[
A^{\Sigma_t} = \frA^{\Sigma_t(\frA)} = W_n + W_nJ = \frA_n.\qedhere
\]
\end{proof}

\begin{corollary}\label{final conclusion}
We have 
\begin{enumerate}
\item\label{containment} $A^{\Sigma_t} \supseteq \B_\rho(\E)$;
\item\label{fg} $A^{\Sigma_t}$ is a finitely generated $\B_\rho(\E)$-module;
\item\label{ff} $A^{\Sigma_t}$ has the same field of fractions as $\B_\rho(\E)$.
\end{enumerate}
In particular, $A^{\Sigma_t}$ and $\B_\rho(\E)$ are fullness peers.
\end{corollary}

\begin{proof}
The definition of $\E$ and \cref{I1 fixed by twists} imply \eqref{containment}.

For \eqref{fg}, if either $\bar{\rho}$ is not dihedral or $\ker \beta_t \neq 1$, then $A^{\Sigma_t} = \B_\rho(\F^{\beta_t(\Sigma_t)})$, which is finitely generated over $\B_\rho(\E)$ since $W(\F)^{\Sigma_t}$ is finitely generated over $W(\E)$.  In the case when $\bar{\rho}$ is dihedral and $\ker \beta_t = 1$, recall that $\B_\rho(\E)$ is a noetherian ring (\cref{field of fractions}\eqref{noetherian rings}) and $A = \B_\rho(\F)+W(\F)B_1(\rho)$ is a noetherian $\B_\rho(\E)$-module by \cref{field of fractions}\eqref{noetherian modules} and the fact that $W(\F)$ is noetherian over $W(\E)$.  Thus the $\B_\rho(\E)$-submodule $A^{\Sigma_t}$ of $A$ is necessarily noetherian and hence finitely generated.

The third point has largely been established already.  When $J$ is small with respect to $\F/\E$, it follows from \cref{small non domain case}.  When $J$ is big with respect to $\F/\E$ and either $\ovl{\rho}$ is not dihedral or $\ker \beta_t \neq 1$, this follows from \cref{induction step} since in those cases $W(\E) + J = \B_\rho(\E)$.  Finally, when $\ovl{\rho}$ is dihedral, $\ker \beta_t = 1$, and $J$ is big with respect to $\F/\E$ we have 
\[
Q(A^{\Sigma_t}) = Q(W(\E) + J) = Q(\B_\rho(\E)),
\]
where the first equality follows from \cref{induction step} and they second from \cref{field of fractions}\eqref{same fraction field}.  

The final statement now follows from \cref{fg plus same quotient field implies ideal containment}.
\end{proof}

We now have the following corollary, which summarizes the most general theorem we have for images of admissible pseudodeformations with $2$-power determinant.

\begin{corollary}\label{final goal}
Let $\ovl{\rho} \colon \Pi \to \GL_2(\F)$ be a regular representation such that the order of $\det \ovl{\rho}$ is a power of 2.  If $\ovl{\rho}$ is octahedral, assume furthermore that $\ovl{\rho}$ is good.  Let $A$ be a domain and $(t, d) \colon \Pi \to A$ an admissible pseudodeformation of $\ovl{\rho}$.  If $(t,d)$ is not a priori small, then $(t,d)$ is $A^{\Sigma_t}$-full, hence $A_0$-full.
\end{corollary}

\begin{proof}
By \cref{formal fullness corollary}, $\rho_t$ is $B_\rho(\E)$-full.   \cref{final conclusion} implies that $\rho_t$ is $A^{\Sigma_t}$-full, hence $A_0$-full by \cref{A0 Afixedbytwists fullness peers}. 
\end{proof}

\medskip


\section{Main fullness results}\label{beyond admissible}
\label{nonadmissible main thm} 
In \cref{beyond admissible} we draw conclusions from \cref{final goal} that are useful in applications and when comparing our work with previous results in the literature.  Although the constant determinant assumption is important to be able to use Bella\"iche's work in \cref{conceptual interpretation of Bellaiche ring}, in practice one rarely works in a constant-determinant setting. Here we give the most general fullness result we can prove; in particular, we remove the constant-determinant assumption present in \cref{final goal}.  We then recast our main theorem and other highlights of the theory of fullness in the language of representation theory rather than pseudorepresentations.  

To ensure that our main result can be read independent of much of the rest of the paper, we briefly recall our notation and terminology.  Let $p$ be an odd prime, $A$ a local pro-$p$ domain with maximal ideal $\m$ and residue field $\F$, and $\Pi$ a $p$-finite profinite group (\cref{p finite condition}).  We are interested in a continuous pseudodeformation $(t,d) \colon \Pi \to A$ of a semisimple representation $\ovl \rho \colon \Pi \to \GL_2(\F)$.  Unlike much of the paper, in \cref{nonadmissible main thm} we never require $A$ to be the trace algebra of $(t,d)$.  We say $(t,d)$ is not a priori small if it is not reducible, dihedral, or equal to a twist of its Teichm\"uller lift (\cref{irreducible dihedral}), and it is $A_0$-full if the image of some representation carrying $(t,d)$ contains, up to conjugation, an $A_0$-congruence subgroup (\cref{fullness}).

Recall that $A_0$ is the adjoint trace ring of $(t,d)$ (\cref{A0 defn}); its residue field is the trace field of $\rhobar$ and is denoted $\E$.  We say that $\ovl \rho$ is \textit{regular} if its image contains a matrix whose eigenvalue ratio is in $\E^\times \setminus \{\pm1\}$: see \cref{regular} and \cref{regularity remark background} immediately following.  When $\ovl \rho$ is octahedral (that is, projective image $S_4$), see \cref{good} for the definition of \textit{goodness}.  

\begin{theorem}\label{nonconst det main thm}
Let $p > 2$ be prime, $A$ a local pro-$p$ domain with residue field $\F$, and $\Pi$ be a $p$-finite profinite group.  Let $\ovl{\rho} \colon \Pi \to \GL_2(\F)$ be a regular semisimple representation that is good if $\ovl{\rho}$ is octahedral.  If $(t,d) \colon \Pi \to A$ is a pseudodeformation of $\ovl{\rho}$ that is not a priori small, then $(t,d)$ is $A_0$-full.
\end{theorem}

\begin{proof}
Let $\chi \colon \Pi \to A^\times$ be a character such that $(t', d') \coloneqq (\chi t, \chi^2d)$ is a constant-determinant pseudorepresentation, and write $\ovl{\rho}' \coloneqq \ovl{\chi} \otimes \ovl{\rho}$, where $\ovl{\chi} \colon \Pi \to \F^\times$ is the reduction of $\chi$ modulo $\m$.  Assume that $\chi$ is chosen such that $\ovl{\rho}'$ has no conjugate self-twists if $\ovl{\rho}$ is reducible and the order of $\det \ovl{\rho}'$ is a power of 2 if $\ovl{\rho}$ is absolutely irreducible.  This is possible by \cref{no csts reducible case} in the reducible case and \cref{detorder2} in the absolutely irreducible case.  Furthermore, note that if $\ovl{\rho}$ is octahedral and good, then so is $\ovl{\rho}'$ by \cref{good under twisting}.  Let $A'$ be the subring of $A$ topologically generated by $t'(\Pi)$.  We have seen in \cref{reducible ASigmat fullness} and \cref{final goal} that if $\ovl{\rho}'$ is regular (and under the further assumption that $\ovl{\rho}$ is good when $\ovl{\rho}$ is octahedral) and $(t', d')$ is not a priori small, then $(t', d')$ is $A_0$-full.  This is sufficient by \cref{twisting does not affect fullness}.  
\end{proof}

\begin{remark}\label{tentwo}
Let $(t, d)$ be as in \cref{nonconst det main thm}. 
Then its constant-determinant twist $(t', d'): \Pi \to A$ satisfies the conditions of \cref{final goal}, so that it is $A_{t'}^{\Sigma_{t'}(A_{t'})}$-full. By \cref{twisting does not affect fullness}, $(t, d)$ itself is also $A_{t'}^{\Sigma_{t'}(A_{t'})}$-full. However it does not follow that that $(t, d)$ is $A_{t}^{\Sigma_{t}(A)}$-full. 
Indeed, $(t,d)$ may be affected by the pathologies pointed out in \cref{gen cst example,insepex}, or even worse, the pro-$p$-part of $d$ may be transcendental over $A_0$. Such obstacles are not faced by well-behaved $A_0$.
\end{remark}

We end this section by recasting all our headline results, including our main result (\cref{nonconst det main thm}) in the language of representation theory, convenient for comparing applications to results in the literature in \cref{applications section}.  
 
Building on the notation recalled above, let $\rho:\Pi \to \GL_2\big(Q(A)\big)$ be a representation whose trace lands in $A$ and is continuous, and let $(t,d) = (\tr \rho, \det \rho)$.  Recall that $\rho$ is not a priori small if $\rho$ is strongly absolutely irreducible (\cref{apssai}). Write $A_\rho$ for the trace algebra of $\rho$ --- the $W(\F)$-subalgebra of $A$ topologically generated by $t(\Pi)$ --- and $K$ for its field of fractions. Let $K_0$ denote the field of fractions of $A_0$, the adjoint trace ring of $\rho$. There is a unique semisimple representation $\rhobar$ whose trace is equal to the reduction of $t$ modulo $\m$. 
Finally, recall that we say that the determinant of $\rho$ is $A_0$-constant if the pro-$p$ part of $\det \rho$ is $A_0$-valued. In this case $K/K_0$ is Galois (\cref{KgaloisK0}).

\begin{theorem}\label{repok} 
Let $\rho$ be not a priori small. Suppose $\rhobar$ is regular, and good if octahedral. 
\begin{enumerate}
\item\label{a0fullend}  {\bf $A_0$-fullness:} $\rho$ is $A_0$-full.
\item\label{optend} {\bf Optimality:} If $\rho$ is $A'$-full for a subring $A' \subseteq A$, then $A_0$ contains a fullness peer of $A'$. 
\item\label{cstfixa0end} {\bf All CSTs fix $A_0$:} For any extension $B$ of $A$ we have $\Sigma_\rho(B) = \Sigma_\rho(B/A_0)$.
\end{enumerate}
If further $K$ is a separable extension of $K_0$, then:
\begin{enumerate}[resume]
\item \label{cstcarvek0end} {\bf CSTs carve out $K_0$:} $(K^\sep)^{\Sigma_\rho(K^\sep)} = K_0.$
\end{enumerate}
If further still $\rho$ has $A_0$-constant determinant, then:
\begin{enumerate}[resume]
\item\label{cstsimpleend} {\bf All CSTs are simple:} Every $K^\sep$-valued conjugate self-twist restricts to a simple $A_\rho$-valued one: $\Sigma_\rho(K^\sep) \onto \Sigma_\rho(A_\rho)$. Moreover,
$K^{\Sigma_\rho(A_\rho)} = K_0$.
\item \label{cstfixfullend} {\bf CST-invariants fullness:} $\rho$ is $A_\rho^{\Sigma_\rho(A_\rho)}$-full.  
\end{enumerate}
\end{theorem}

\begin{proof} 
Recall that $\rho$ is $A'$-full if and only if $(\tr \rho,\det \rho)$ is $A'$-full by \cref{fullness well defined}\eqref{fullness for reps well defined}.  Thus the first statement follows from \cref{nonconst det main thm} and the second from \cref{A0optimalitythm}.  
The third statement follows from \eqref{a0fullend} by \cref{A0fullsimpleCST}. The fourth statement follows from \cref{f0ifsep}\eqref{sick} with \mbox{$L = K^\sep$}, \mbox{$E = K_0$}, \mbox{$F = K$} and $F_0$ the (nontopological) adjoint trace field  of $(t, d)$ viewed as valued in $K$. For the fifth statement, use  \cref{A0fullsimpleCST} to obtain the restriction map \mbox{$\Sigma_t(K^\sep/A_0) \to \Sigma_t(A_\rho)$}, surjective since $K^\sep$ is normal over $K$. Alternatively, $K$ is Galois over $K_0$ with $\Sigma_t(A_\rho) = \Gal(K/K_0)$ by \cref{KgaloisK0} and then \mbox{$\Sigma_t(K^\sep) = \Gal(K^\sep/K_0)$} by \cref{fautscsts}. The last statement follows from \eqref{a0fullend} and \cref{A0 Afixedbytwists fullness peers}. \qedhere
\end{proof}


\section{Residually large representations}
\label{residual full case} 

Here we show that by imposing stronger conditions on the residual image when $\ovl \rho$ is large, we obtain a more precise description of the image of $\rho$.  That is, we assume that $\rho \colon \Pi \to \GL_2(A)$ is a continuous representation such that $\im \ovl{\rho} \supseteq \SL_2(\E)$ and $\ovl\rho$ is large. 
Under this assumption we have a more precise understanding of the image of $\rho$ than simply fullness.  Unlike our main fullness result, in this section $A$ is any local pro-$p$ ring; it need not be a domain.

Historically, this is the case that has been studied the most, starting with the work of Boston in \cite[Appendix]{MazurWiles86}.  Boston shows that if $A$ is a complete local noetherian ring and $H$ is a closed subgroup of $\SL_2(A)$ that projects onto $\SL_2(A/\m^2)$, then $H = \SL_2(A)$ \cite[Appendix Proposition 2]{MazurWiles86}.  Bella\"iche has pointed out that Boston's result follows from his work \cite[Remark~6.8.4]{Bellaiche18}. As an application of the description of the image found in \cref{SL2 case}, we show in \cref{variation on boston} that one can replace the hypothesis that $H$ projects onto $\SL_2(A/\m^2)$ with the hypothesis that $H$ projects onto $\SL_2(A/\m)$ and $A$ is the trace ring of $H$ to obtain the same conclusion.  Theorems of this form have been obtained in special cases, for instance for the Galois representation attached to the mod-$p$ Hecke algebra \cite[Theorem, Introduction]{Amaros20}.  Note that assuming $H$ projects onto $\SL_2(A/\m)$, the hypothesis about $A$ being the trace algebra can always be arranged while Boston's hypothesis about projecting onto $\SL_2(A/\m^2)$ may fail.  Moreover, our description of the image in \cref{SL2 case} does not require that the residual image contain all of $\SL_2(\F)$, rather only $\SL_2(\E)$, and is thus more general that the setting of Boston's work.

Let $\rho \colon \Pi \to \GL_2(A)$ be a continuous representation such that $\im \ovl{\rho} \supseteq \SL_2(\E)$.  Note that such a $\ovl\rho$ is large if and only if $\E \neq \F_3, \F_5$.  Let $\chi \colon \Pi \to A^\times$ be the character described in the proof of \cref{nonconst det main thm}, and let $r \coloneqq \chi \otimes \rho$.  Assume that $\rho$ is conjugated in such a way that \cref{Bellaiche main theorem} applies to $r$.  Recall that $\B_r(\E) = W(\E)[I_1(r)]$.

\begin{proposition}\label{SL2 case}
Assume $\#\E \geq 7$.  Then
\begin{enumerate} 
\item $\im r \supseteq \SL_2(\B_r(\E))$ as a finite index subgroup;
\item $\im \rho \supseteq \SL_2(\B_r(\E))$;  
\item $\B_r(\E)$ is the largest subring $B$ of $A$ for which $\im \rho \supseteq \SL_2(B)$. 
\end{enumerate}
\end{proposition}

\begin{proof}
For ease of notation, write 
$\m_r$ for the maximal ideal of $\B_r(\E)$.

Since $\im \ovl r \supseteq \SL_2(\E)$, it follows that $\im r \supseteq \SL_2(W(\E))$ by \cite[Main Theorem]{Manohar15}.  In particular, $p \in I_1(r)$ and so $\m_r = I_1(r)$ by \cref{Bellaiche main theorem}.  Let $H$ be the subgroup of $G \coloneqq \im r$ generated by $\Gamma = \Gamma(r)$ and $\SL_2(W(\E))$.  Then $H$ is a finite index subgroup of $G$ since $\Gamma$ is.  

We claim that $H = \SL_2(\B_r(\E))$.  Indeed, note that $\Gamma = \Gamma_{\B_r(\E)}(\m_r)$ by \cite[Corollary 6.8.3]{Bellaiche18} and the fact that $\m_r = I_1(r)$.  In particular, this shows that $H \subseteq \SL_2(\B_r(\E))$.  In fact, $H$ is a subgroup of $\SL_2(\B_r(\E))$ such that $H/\Gamma = \SL_2(\E) = \SL_2(\B_r(\E))/\Gamma$.  Thus we must have equality.

Now suppose that $\im r \supseteq \SL_2(B)$ for some subring $B$ of $A$.  Without loss of generality, we may assume $B$ is closed, hence local.  Then $\Gamma \supseteq \Gamma_{B}(\m_B)$, which implies that $I_1(r) \supseteq \m_B$.  On the other hand, if $\im r \supseteq \SL_2(B)$ then $\im \ovl r \supseteq \SL_2(B/\m_B)$.  By definition of $\E$, we know that $\E$ is the largest subfield of $\F$ such that $\im \ovl r \supseteq \SL_2(\E)$.  Thus we must have $B/\m_B \subseteq \E$.  It follows that $B \subseteq \B_r(\E)$.

As for $\rho$, note that there is a character $\tilde{\chi} \colon \im \rho \to A^\times$ such that \mbox{$\im r = \{x\tilde{\chi}(x) \colon x \in \im \rho\}$}.  Now~$\tilde{\chi}$ must be trivial on $\SL_2(B)$ for any ring $B$ whose residue field has more than three elements by \cref{congrcomm}.  Therefore $\im \rho$ and $\im r$ contain the same copies of~$\SL_2$.
\end{proof}

\begin{corollary}\label{SL2(A0)}
If $\#\E \geq 7$ and $\im \ovl\rho \supseteq \SL_2(\E)$, then $\im \rho$ contains $\SL_2(A_0)$ up to conjugation. 
\end{corollary}

\begin{proof}
By \cref{SL2 case} it suffices to show that $\B_r(\E) \supseteq A_0$, and this only needs to be shown when $r = \rho$ is the constant-determinant universal representation.  Without loss of generality we may twist to assume that the order of $\det \rho$ is a power of $2$.  By \cref{Bellaiche main theorem} we have $\cA = \B_r(\F)$ since $\ovl\rho$ is large.  By \cref{I1 fixed by twists} we may conjugate $\rho$ to assume that $I_1(\rho)$ is fixed by all ($\cA$-valued) conjugate-self twists of $\rho$, and by \cref{lifting arbitrary csts} all conjugate-self twists of $\ovl\rho$ lift to conjugate self-twists of $\rho$ since $\rho$ is universal.  Therefore 
\[
\B_r(\E) = W(\F)^{\Sigma_\rho}[I_1(\rho)] = \cA^{\Sigma_\rho} \supseteq \cA_0.\qedhere
\]
\end{proof}

\begin{theorem}\label{variation on boston}
Let $A$ be a pro-$p$ local noetherian ring with $p \neq 2$ and residue field $\F$.  Let $H$ be a closed subgroup of $\SL_2(A)$ that projects onto $\SL_2(\F)$.  If $A$ is the trace ring of $H$, that is, $A$ is topologically generated by $\tr H$, then $H = \SL_2(A)$.
\end{theorem}

\begin{proof}
Let $\rho \colon H \to \SL_2(A)$ be the natural inclusion.  Note that $\rho$ has constant determinant.  Since $A$ is the trace ring of $\rho$, \cref{Bellaiche main theorem} implies that $A = \B_\rho(\F)$.  Thus by \cref{SL2 case}, it suffices to show that $\E = \F$.  Since $\im \overline{\rho} = \SL_2(\F)$, it follows that $\E$ is the subfield of $\F$ generated by the squares of traces of $\SL_2(\F)$.  A straightforward matrix calculation shows that $\E = \F$.  
\end{proof}

\begin{remark}
In the preprint \cite{BoeckleArias}, Aryas-de-Reina and B\"ockle prove a large image result for a residually full representation $\Pi\to G(A)$, where $G$ is an adjoint group and $A$ is the ring of definition of the representation. It seems possible to recover  \cref{variation on boston} by applying their result to the projective representation ${\mathbb P}\rho\colon\Pi\to\PGL_2(A)$ attached to $\rho$ and using the fact that the ring of definition of $\mathbb P\rho$ is the ring fixed by the conjugate self-twists of $\rho$.
\end{remark}


\section{Applications to Galois representations}\label{applications section}

In this section we specialize \cref{nonconst det main thm} to some arithmetic settings, more specifically to representations coming from elliptic, Hilbert, and Bianchi cuspidal eigenforms (\cref{classical applications1} through \cref{Bianchi applications}) and cuspidal $p$-adic families of elliptic and Hilbert eigenforms (\cref{Hida families} through \cref{Hilbert families}). We explain how to recover, and in some cases improve, the results already present in the literature. 
In particular, since our methods are entirely agnostic about the group $\Pi$, they reveal that many of the classical big-image results are fundamentally \emph{algebraic} in nature: they do not rely on the arithmetic input, such as local information at the places where a Galois representation is ramified, that went into the original proof.

\subsection{Classical modular eigenforms}\label{classical applications1}
Let $f$ be a non-CM cuspidal modular eigenform of some level and some weight $k \geq 2$  defined over a number field $K$. Fix an algebraic closure $\Qbar$ of $\Q$ and let $G_\Q := \Gal(\Qbar/\Q)$. For any prime $\p$ of $K$ lying over a rational prime $p$, let $K_\p$ be the completion of $K$ at $\p$ and $\OK_\p$ its ring of integers. A construction of Deligne attaches to this data an irreducible continuous representation $\rho_{f, \p}: G_{\Q} \to \GL_2(K_\p)$, unramified almost everywhere and hence factoring through a $p$-finite extension, whose traces of Frobenius elements at unramified primes 
correspond to Hecke eigenvalues of $f$. 
Because $G_\Q$ is compact we may view $\rho_{f, \p}$ as taking values in $\GL_2(\OK_\p)$.

The following result about the image of $\rho_{f, \p}$ was proved by Ribet and Momose in the 1980s, generalizing an earlier theorem of Serre about Tate modules of elliptic curves. Let $K_{\p, 0}$ be the subfield of $K_\p$ fixed by all the generalized conjugate self-twists of $\rho_{f, \p}$, and $\OK_{\p, 0}$ its ring of integers.

\begin{theorem}[Ribet, Momose  at $\p$: first version\ {\cite{Ribet85,Momose81}}]  \label{ribetmomose} \leavevmode\\
For all but finitely many primes $\p$ of $K$, the representation $\rho_{f, \p}$ is $\OK_{\p, 0}$-full.
\end{theorem} 

To show the extent to which our work recovers the result of this theorem at $\p$, we first make more explicit Ribet and Momose's condition on $\p$. 
Let $K_0$ be the subfield of $K$ fixed by the conjugate self-twists of $f$; the $p$-adic field $K_{\p, 0}$ defined above is the completion of $K_0$ at the prime under $\p$. Let $H \subseteq G_\Q$ be the intersection of the kernels of all the conjugate self-twist characters of~$f$:

\begin{equation}\label{H}
H \coloneqq \bigcap_{(\sigma, \eta) \in \til \Sigma_f} \ker \eta.
\end{equation}
Then $H$ is a finite-index normal subgroup of $G_\Q$. 
Because all the conjugate self-twist characters of $\rho_{f, \p}$ are trivial on $H$, the trace of $\left.\rho_{f, \p}\right|_H$ lands in $K_{\p, 0}$. 
As described just after \cref{tracefielddef}, there is therefore a $K_{\p, 0}$-quaternion algebra $D_\p$ splitting over $K_\p$ with $\rho_{f, \p}(H) \subset D_\p^\times \subset \GL_2(K_\p).$

Ribet and Momose describe a global analogue to this picture. They define a global $K_0$-quaternion algebra $D$ split over $K$ with $D(K_{\p, 0}) \cong D_\p$: that is, for each prime $\p$, the restriction to $H$ of $\rho_{f, \p}$ can be viewed as taking values in $D(K_{\p, 0})^\times$. By compactness again we may view $\rho_{f, \p}(H)$ as a subgroup of the units of a maximal order $\OK_{D, \p}$ of $D(K_{\p, 0})$. Ribet and Momose's adelic open-image theorems say that the image of $\rho_{f, \p}$ always contains an open subgroup of the norm-$1$ units of $\OK_{D, \p}^{\times}$, and for all but finitely many $\p$ it contains all of those norm-$1$ units.
In particular, if $D(K_{\p, 0})$ is split, then up to conjugation $\OK_{D, \p}^\times = \GL_2(\OK_{\p, 0})$; and we can therefore make \cref{ribetmomose} more precise.

\begin{theorem}[Ribet, Momose  at $\p$: second version\ {\cite{Ribet85,Momose81}}]  \label{ribetmomose2} \leavevmode\\
If $D(K_{\p, 0})$ is split, then the representation $\rho_{f, \p}$ is $\OK_{\p, 0}$-full.
\end{theorem} 

\begin{remark} \label{iff} In fact, the statement in \cref{ribetmomose2} is an if-and-only-if. Indeed, no element of $D(K_{\p, 0})$ can have distinct $\GL_2$-eigenvalues in $K_{\p, 0}$. A matrix $g$ with eigenvalues $\alpha, \beta$ satisfies $(g - \alpha)(g - \beta) = 0$. If $\alpha, \beta$ are in the center of a division algebra containing $g$, and at the same time are  eigenvalues of $g$ in any matrix setting, then the Cayley-Hamilton equation $(g - \alpha)(g - \beta) = 0$ means that either $g = \alpha$ or $g = \beta$. Thus no embedding of $D(K_{\p, 0})$ into $\GL_2(K_\p)$ can contain any congruence subgroup of 
$\SL_2(\OK_{\p, 0})$, or even of {$\SL_2(\Z_p)$}. The same is true for all of $\rho(\Pi)$: see \cref{KgaloisK0} and \cref{closednormal1} or \cref{divred} below. In other words, no nonsplit $\rho_{f, \p}$ can ever satisfy our present definition of fullness. As we've defined it, fullness is fundamentally a $\GL_2$ property; 
a fitting notion for more general algebraic groups generalizing Ribet and Momose's openness beyond Krull dimension $1$ is outside the scope of this investigation.
\end{remark}

We now show that our results recover Ribet and Momose's theorem at $\p$ in most cases. Let $\F$ be the residue field of $\OK_\p$, a finite extension of $\F_p$, and let $\ovl{\rho} \colon G_\Q \to \GL_2(\F)$ be the semisimplification $\rho_{f, \p}$ modulo the maximal ideal of $\OK_\p$.  
Recall that $\ovl{\rho}$ is \emph{regular} if its image contains a {matrix whose eigenvalue ratio} 
is not $\pm 1$ but is contained in the trace algebra $\E$ of $\ad \ovl \rho$ (which is a subfield of the residue field of $\OK_{\p, 0}$): see \cref{regular} and \cref{regularity remark background} immediately following. 
If $\ovl{\rho}$ is octahedral (that is, projective image $S_4$), see \cref{good} for the notion of \emph{goodness}.

\begin{theorem}[Our results recovering Ribet and Momose at $\p$]\label{classical application}Assume that $p$ is odd and that $\ovl{\rho}$ is regular;  if $\ovl{\rho}$ is octahedral, assume further that $\ovl{\rho}$ is good.  Then $\rho_{f, \p}$ is $\OK_{\p, 0}$-full.
\end{theorem}

\begin{proof}
Since $f$ is cuspidal, non-CM, and has weight $k \geq 2$, its associated representation $\rho_{f, \p}$ is strongly absolutely irreducible (\cite[Proposition~4.4]{Ribet77}) and hence not a priori small.  By \cref{repok}, $\rho_{f, \p}$ is $A_0$-full, where as usual $A_0 \subseteq \OK_{\p}$ is the adjoint trace ring.  By \cref{A0 Afixedbytwists fullness peers}, $A_0$ and $\OK_{\p, 0}$ are fullness peers, and $\rho_{f, \p}$ is $\OK_{\p, 0}$-full. 
\end{proof}

The regularity assumption in \cref{classical application} a posteriori forces $D(K_{\p, 0})$ to split (\cref{iff}). We can also see that a nonsplit $D(K_{\p, 0})$ means an irregular $\rhobar$ directly:

\begin{proposition}\label{divred}
If $D(K_{\p, 0})$ is a division algebra, then $\rhobar$ is reducible and not regular. 
\end{proposition}

\begin{proof}
We first show that $\left.\rhobar\right|_H$ is reducible and not regular. Let $L/K_{\p, 0}$ be the unique quadratic unramified extension, $\pi$ is a uniformizer of either, and $\sigma$ the nontrivial element of $\Gal(L/K_{\p, 0})$. Note that we do not assume that $L$ is a subfield of $K_\p$. Write $\ell$, $k$ for the residue fields of $L$, $K_{\p, 0}$, respectively; then $[\ell: k] = 2$ and $\E \subseteq k$. 
By compactness $\rho(H)$ can be viewed as a subgroup of 
$$\left\{ \begin{psmallmatrix} \alpha & \pi \beta \\ \sigma(\beta) & \sigma(\alpha)\end{psmallmatrix} \colon \alpha \in \OK_L^\times, \beta \in \OK_L \right\} \subset \GL_2(\OK_L),$$ 
the maximal order of $D(K_{\p, 0})$ as viewed inside $\GL_2(L)$,
which gives a representation $\rho'$ of $H$ over~$L$. By inspection, it is clear that its residual representation $\rhobar'$ is, up to semisimplification, a sum of two characters to $\ell$ conjugate over $k$. This means that the eigenvalue ratio $r$ of any element in $\rhobar'(H)$ 
is in the form $r = a^{\# k -1}$ for some $a \in \ell^\times$. Such an element is in $k$ if and only if $r^{\#k - 1} = 1$; in other words, if and only if $$1 = a^{(\#k)^2 - 2(\# k) + 1} = a^{-2(\#k) - 2} = r^{-2}.$$ But this last is only possible if $r = \pm 1$. In other words, $\rho'$ is residually neither absolutely irreducible nor regular. Since $\rho'$ is isomorphic to $\left.\rho\right|_H$ over $\ovl\Q_p$, the same is true for $\left.\rho\right|_H$ as well.  

We now follow \cite[B.4.8(1)]{Nekovar2012} to claim that the same is true for $\rho$ on all of $\Pi$. Change notation to let $L$ be any subextension of $K_\p$ that is quadratic over $K_{\p, 0}$ and hence splits $D(K_{\p, 0})$, with $\sigma$ a generator of $\Gal(L/K_{\p,0})$ and $\pi$ a uniformizer of $K_{\p, 0}$ that is not a norm from $L$. Since $H$ is normal in $\Pi$ and $\rho(H)$ spans  $D(K_{\p, 0})$ over $K_{\p, 0}$, the image of $\rho$, up to conjugation, will be contained in the normalizer of the subgroup $Q=\left\{ \begin{psmallmatrix} \alpha & \pi \beta \\ \sigma(\beta) & \sigma(\alpha)\end{psmallmatrix}  : (\alpha, \beta) \in L^2 - \{(0,0)\} \right\}$, isomorphic to $D(K_{\p, 0})^\times$, in $\GL_2(K_\p)$. One can show that this normalizer is just $K_\p^\times Q$\footnote{Certainly, $K_\p^\times Q$ normalizes. Conversely, if $\begin{psmallmatrix} a& b \\ c& d \end{psmallmatrix} \in \GL_2(K_\p)$ normalizes $\begin{psmallmatrix} 0 & \pi \\ 1 & 0\end{psmallmatrix}$, then the off-diagonal relation gives us $\sigma((\pi a^2 - \pi^2 c^2)\delta^{-1}) = (\pi d^2 -b^2)\delta^{-1}
,$
where $\delta = ad- bc$; if it normalizes $\begin{psmallmatrix} 0 & \pi t \\ -t & 0\end{psmallmatrix}$, where $t \in L$ with $\sigma(t) = -t$, then the off-diagonal gives 
$\sigma((\pi a^2 + \pi^2 c^2)\delta^{-1}) = (\pi d^2 + b^2)\delta^{-1};$
if it normalizes $\begin{psmallmatrix} t & 0 \\ & -t\end{psmallmatrix}$ then the off-diagonal is 
$\sigma(a (\pi c) \delta^{-1}) = d b \delta^{-1}$.
Combining the first two, we obtain $\sigma(a^2\delta^{-1}) = d^2 \delta^{-1}$ and $\sigma((\pi c)^2 \delta^{-1}) = b^2 \delta^{-1}$; adding the third gives us, for example if $a$ is invertible, $\sigma(\pi c/a) = b/d$. By considering the diagonal relations, we also get $\sigma(b/a) = \pi c/d$ and $\sigma(a/d) = d/a$. In other words any normalizing element with a nonzero entry in the upper left looks like $\begin{psmallmatrix} a & u a\, \sigma(\pi c/a)  \\ c & ua \end{psmallmatrix}$, with $a \in K_\p^\times$ arbitrary, $c \in aL$, and $u:=d/a \in L$ of norm $1$. From Hilbert 90, any norm-$1$ $u$ is $\sigma(\alpha)/\alpha$ for some $\alpha \in L$; letting $x := a/\alpha$ and $\beta := c \alpha/a$ puts our matrix in the desired form $x \begin{psmallmatrix} \alpha &  \pi \sigma(\beta)  \\ \beta & \sigma(\alpha) \end{psmallmatrix}\in  K_\p^\times Q$.
Or see \cite[B.1.6]{Nekovar2012} for a more conceptual argument. 
}, so that passing from $H$ to $\Pi$ does not affect the projective image of $\rhobar$. Therefore $\rhobar$ on all of $\Pi$ remains reducible and not regular.
\end{proof}

In other words, the regularity assumption in \cref{classical application} eliminates the division algebra case. 
One might hope for a converse, so that our methods could recover all of \cref{ribetmomose2}. But alas this is not so: there are certainly cases where $D(K_{\p, 0})$ is a matrix algebra but regularity is not satisfied, so that our methods do not apply.  In addition to $p = 2$, we do not conclude fullness if $\F = \F_3$, even if $f$ has no conjugate self-twists and hence $D$ is globally split, as is the case for a non-CM elliptic curve over $\Q$. If the image of $\rhobar$ is too small to accommodate the regularity assumption, our methods cannot handle it.

\subsection{Hilbert modular eigenforms}\label{classical applications2} 
Everything in \cref{classical applications1} has been generalized to Hilbert modular forms. In particular, our results recover, in much the same manner and to much the same extent, the big-image results of \nekovar\ generalizing Ribet and Momose's work over $\Q$. We summarize the situation very briefly.

Let $F$ be a totally real field and $f$ a non-CM cuspidal Hilbert modular eigenform over $F$ all of whose weights are at least $2$.  Fix an algebraic closure $\ovl{F}$ of $F$, and let $G_F \coloneqq \Gal(\overline{F}/F)$.  Fix a prime $p$ and an embedding $\iota_p \colon \overline{\Q} \hookrightarrow \overline{\Q}_p$. Let $\rho_{f, \iota_p} \colon G_F \to \GL_2(\overline{\Q}_p)$ be the Galois representation attached to $f$ in the usual way, which we may view as having coefficients in the ring of integers $\OK$ of some finite extension of $\Q_p$.  Let $\OK_0 \subseteq \OK$ be the ring of integers of the fixed field of all the conjugate self-twists of $\rho_{f, \iota_p}$. Like Ribet and Momose, \nekovar\ constructs a division algebra $D$ over the fixed field $K_0$ of $\Sigma_f$ and proves an adelic open-image result, which implies $\OK_0$-fullness when $D$ splits.

\begin{theorem}[{\cite{Nekovar2012}}]  \label{nekovar} 
For all but finitely many $\iota_p$, the representation $\rho_{f, \iota_p}$ is $\OK_{0}$-full.
\end{theorem}
\noindent Our results depend on hypotheses on the reduction $\rhobar_{f, \iota_p}$ of $\rho_{f, \iota_p}$ modulo the maximal ideal of $\OK$. 

\begin{theorem}[Our results recovering {\nekovar} at $\iota_p$] \label{ournekovar} Suppose that $p$ is odd and $\rhobar_{f, \iota_p}$ is regular; if $\rhobar_{f, \iota_p}$ is octahedral suppose further that it is good. Then $\rho_{f, \iota_p}$ is $\OK_{0}$-full.  
\end{theorem}
\noindent The proof is analogous to that of \cref{classical application}. In particular, the fact that the weight of $f$ is at least $2$ at each infinite place means that $\rho_{f, \iota_p}$ has distinct Hodge-Tate weights, which implies that no twist of it has finite image (see \cref{apssai} for context). Or apply \cite[Lemma 3.2.12]{CalegariEmertonGee}.

\subsection{Bianchi modular forms and generalizations}\label{Bianchi applications}
Unlike Hilbert modular forms, which are automorphic forms on $\GL_2$ over totally real fields, Galois representations associated to automorphic forms of $\GL_2$ over CM fields have only been constructed relatively recently, and even then only under some technical assumptions.  We briefly summarize how our results can be applied to that context.  

Let $E$ be a CM field with maximal totally real subfield $F$.  Fix an algebraic closure $\ovl E$, and let $G_E \coloneqq \Gal(\ovl E/E)$.  Let $\pi$ be a cuspidal automorphic representation of $\GL_2(\A_E)$, where $\A_E$ denotes the adeles of $E$; when $E$ is imaginary quadratic, $\pi$ is called a Bianchi modular form.  Assume that $\pi$ is of cohomological type with central character $\omega$.  Following Mok \cite{Mok14}, assume moreover that $\omega$ arises from an algebraic idele class character $\tilde{\omega}$ on $\A_F^\times$ via the norm map and that $\tilde{\omega} = \bigotimes_v \tilde{\omega}_v$ such that $\tilde{\omega}_v(-1)$ takes the same value for all archimedean places of $F$.  (When $F = \Q$, this is simply the condition that $\omega$ is invariant under complex conjugation.)  Suppose there is no nontrivial quadratic character $\delta$ of $E$ such that $\pi \cong \pi \otimes \delta$; this is analogous to the non-CM assumption present in \cref{classical applications1} and \cref{classical applications2}.

For each rational prime $p$ and fixed embedding $\iota_p \colon \ovl \Q \hookrightarrow \ovl \Q_p$, associated to $\pi$ there is a continuous irreducible representation $\rho_{\pi, \iota_p} \colon G_E \to \GL_2(\ovl \Q_p)$, which we may view as having coefficients in the ring of integers $\OK$ of some finite extension of $\Q_p$.  In this generality, the existence of $\rho_{\pi, \iota_p}$ is due to Mok \cite{Mok14}, who generalized the construction of Taylor in the imaginary quadratic case \cite{Taylor94}.  Mok also shows, building on work of Berger and Harcos in the imaginary quadratic case \cite{BergerHarcos07}, that $\rho_{\pi, \iota_p}$ is unramified outside a finite set of places and hence factors through a $p$-finite group $\Pi$.  

Let $\OK_0$ be the ring of integers of the fixed field of all the conjugate self-twists of $\rho_{\pi, \iota_p}$.

\begin{theorem}\label{bianchi theorem}
Suppose $p$ is odd and $\ovl \rho_{\pi, \iota_p}$ is regular and good if octahedral.  Then $\rho_{\pi, \iota_p}$ is $\OK_0$-full.
\end{theorem}

\noindent The proof is analogous to that of \cref{classical application} and \cref{ournekovar}. The Hodge-Tate weights of $\rho_{\pi, \iota_p}$ are distinct by \cite[Theorem 5.17]{Mok14}.

To our knowledge, \cref{bianchi theorem} is the most general fullness result in the literature in this context, though Taylor proves the weaker theorem that the image of $\rho_{\pi, \iota_p}$ is Zariski dense in the imaginary quadratic case \cite[Corollary 2]{Taylor94}. Our \cref{bianchi theorem} may be well known to experts. 

\begin{remark}\label{bianchifam}
In contrast to the case of elliptic or Hilbert modular forms that can be $p$-adically interpolated in families with dense classical points and thus have associated ``big" Galois {pseudo}representations (see \crefrange{Hida families}{Hilbert families}), a $p$-adic family of Bianchi modular forms often has only finitely many classical points \cite[Theorem 8.9]{CalegariMazur09}, \cite[Theorem 1.1]{Serban19}. Therefore no Galois pseudorepresentations have been attached to Bianchi families by conventional methods. 
\end{remark}

\subsection{Hida $p$-adic families of modular forms}\label{Hida families} 
In this section, we explore the extent to which our methods recover known big-image results for Galois representations attached to ordinary $p$-adic families of modular forms, often called Hida families. 

Fix $p > 2$, and let $A$ be the ring corresponding to a primitive non-CM irreducible component of Hida's cuspidal shallow Hecke algebra parametrizing $p$-ordinary cuspforms of some fixed tame level and weight $k$ such that $k-1 \equiv i \bmod p-1$ for some fixed $0 \leq i \leq p-1$ (for $i = 0$ this is the ring $\mathbb I'$ on \cite[p.~158]{Lang2016}).
Then $A$ is a finite extension of the Iwasawa algebra $\Lambda := \Z_p\lb 1 + p\Z_p \rb \cong \Z_p \lb T \rb$, which parametrizes the corresponding component of weight space. Let $\F$ be the residue field of~$A$ and $K$ the fraction field of $A$.

Let $(t, d) \colon G_\Q \to A$ be the pseudorepresentation obtained from gluing together those attached to the classical cuspforms in the family, all of which have the same semisimplifed residual representation $\rhobar: G_\Q \to \GL_2(\F)$.  Let $\varepsilon_p \colon G_\Q \to \Z_p^\times$ be the $p$-adic cyclotomic character and \mbox{$\langle \cdot \rangle \colon \Z_p^\times \to 1 + p\Z_p$} the projection onto the pro-$p$ part of $\Z_p^\times$.  The weight character $\kappa \colon G_\Q \to \Lambda^\times$ is given by \mbox{$\kappa(g) = (1 + T)^{\langle \varepsilon_p(g) \rangle}$}.  Let $\chi$ be the tame Dirichlet character associated to the family.  The determinant of $\rho$ is given by 
\begin{equation}\label{dkappachi} 
d = \kappa \chi s(\ovl \varepsilon_p)^i,
\end{equation}
Both $\rhobar$ and $(t, d)$ factor through $\Pi$, the Galois group of the maximal extension of $\Q$ unramified outside $p$ and the level, a $p$-finite profinite group. 
Let $\Pi_p \subset \Pi$ be a decomposition group at $p$ and $I_p \subset \Pi_p$ its inertia subgroup. The ordinary condition guarantees that there exists a $(t,d)$-representation $\rho$, which we view as $\GL_2(K)$-valued by \cref{GMAs over domains}, with $\rho|_{\Pi_p} = \begin{psmallmatrix} 
\epsilon & *\\
0 & \delta
\end{psmallmatrix}$, where $\delta$ an unramified character and $\epsilon$ coincides with $\kappa$ on wild inertia and therefore surjects onto \mbox{$(1 + T)^{1 + p\Z_p}$} \cite[Theorem 4.3.2]{HidaGMF}.

The image of $\rho$ has been studied by Boston \cite[Appendix]{MazurWiles86}, Fischman~\cite{Fischman02}, Hida~\cite{Hida15}, and Lang~\cite{Lang2016}. The latter two are the more recent and most general results, so we focus there.

\begin{theorem}[$\Lambda$-fullness for Hida families, Hida \cite{Hida15}]
If $\rhobar$ restricted to $\Pi_p$ is multiplicity free and $\rho$ is realizable by a representation over $A$ then $\rho$ is $\Lambda$-full.
\end{theorem} 

\noindent Hida's $\Lambda$-fullness strongly suggested that every conjugate self-twist of $(t,d)$ should fix $\Lambda$ (see \cref{fullfixed} here for a proof of this fact). Following Hida, Lang analyzes how the image of $\rho$ is constrained by $\Sigma_t(A/\Lambda)$, the conjugate self-twists that fix $\Lambda$ pointwise. 
To state her result, we let $H_\Lambda \subseteq \Pi$ be the intersection of all the $\ker \eta$ for $(\sigma, \eta)$ in $\Sigma_t(A/\Lambda)$.  

\begin{theorem}[{Big image for Hida families, \cite[Theorem 2.4]{Lang2016}}]\label{jackiethesis} 
Suppose that $\F \neq \F_3$, that $\rhobar$ is absolutely irreducible, and that 
there is an element in $\rhobar(H_{\Lambda} \cap \Pi_p)$ whose eigenvalue ratio is in $\E^\times \setminus\{1\}$.\footnote{\label{langerror}There is a small error in \cite{Lang2016}, which we correct here. Theorem 2.4 as stated {\it loc.~cit.} requires merely that $\rhobar$ restricted to $H_\Lambda \cap \Pi_p$ be multiplicity free, but in fact the result relies on the stronger regularity condition given here. Indeed, on \cite[p. 174]{Lang2016} the definition of $\overline{L}[\lambda]$ only makes sense if one knows $\overline{L}$ is closed under multiplication by $\lambda$, where $\lambda$ is an adjoint eigenvalue of the regular element.} 
 Then $\rho$ is $A^{\Sigma_\rho(A/\Lambda)}$-full. 
\end{theorem}

\noindent A posteriori Lang's fullness result by itself justifies considering only those conjugate self-twists that fix $\Lambda$, even without Hida's $\Lambda$-fullness: see \cref{fullfixed}. Our work both recovers virtually all of Lang's result (exception: \cite{Lang2016} is able to handle some cases where $\mathbb P \rhobar$ is the Klein-4 group) and extends it to include residually reducible $\rhobar$.
 Let $A_0$ be the adjoint trace ring of $\rho$; see \cref{A0 defn}. For the notion of $A_0$-constant determinant, see \cref{a0constantdetdef}; for good octahedral~$\rhobar$ see \cref{good}.

\begin{theorem}[{Our work recovering \cite{Lang2016}}]\label{ourhida}~ 
\begin{enumerate} 
\item  \label{hidafull} If $\rhobar$ is regular and good if octahedral, then $\rho$ is $A_0$-full. 
\item  \label{lambdaina0}If $\Pi_p$ contains a regular element for $\ovl \rho$, then $A_0$ contains $\Lambda$. 
\end{enumerate}

Consequently, if $\Pi_p$ contains a regular element for $\ovl \rho$ and $\rhobar$ is further good if octahedral, then
\begin{enumerate}[resume]
\item \label{lambdafull} $\rho$ is $\Lambda$-full;
\item \label{asigmafull} $\rho$ is $A^{\Sigma_\rho(A)}$-full; 
\item \label{lambdafixed} every conjugate self-twist of $\rho$ fixes $\Lambda$, so that $\Sigma_\rho(A) = \Sigma_\rho(A/\Lambda)$ and $\rho$ is $A^{\Sigma_\rho(A/\Lambda)}$-full. 
\end{enumerate}  
\end{theorem}

\begin{proof}
For \eqref{hidafull}, the representation $\rho$ is not reducible since the Hida family is cuspidal, and it is not dihedral since the Hida family is not CM. The fact that $\rho(\Pi) \not \cong \rhobar(\Pi)$  follows from the fact that a Hida family has classical specializations of weight at least 2. Therefore we know that $\rho$ is $A_0$-full by \cref{repok}.

For \eqref{lambdaina0}, let $d_1 \colon \Pi \to A^\times$ be the pro-$p$ part of $d = \det \rho$, and let $\rho' = d_1^{-1/2}\otimes \rho$ with \mbox{$(t', d') = (\tr \rho', \det \rho')$} the constant-determinant (pseudo)representation of $\rho$.  Note that $\rho'|_{\Pi_p}$ is still upper triangular since $\rho|_{\Pi_p}$ is.  Let $g_0 \in \Pi_p$ be a regular element with residual eigenvalues $\lambda_0, \mu_0$, and let $r$ be a $(t', d')$-representation adapted to $(g_0, \lambda_0, \mu_0)$. By the proof of \cite[Theorem 6.2.1]{Bellaiche18}, we see that, up to replacing $g_0$ with the limit of a sequence of its powers, we may assume that $r(g_0) = \begin{psmallmatrix} 
s(\lambda_0) & 0\\
0 & s(\mu_0)
\end{psmallmatrix}$. 
Viewing both $\rho'$ and $r$ as $\GL_2(K)$-valued by \cref{GMAs over domains}, we see that they are isomorphic since they have the same trace and are irreducible.  In particular, $\rho'(\Pi_p)$ contains an element with eigenvalues $s(\lambda_0), s(\mu_0)$, which (up to swapping $\lambda_0$ and $\mu_0$) is necessarily of the form  $M \coloneqq \begin{psmallmatrix} 
s(\lambda_0) & *\\
0 & s(\mu_0)
\end{psmallmatrix}$.  
On the other hand, using the description of $\epsilon$ and $\delta$ above, we see that 
$\rho'(\Pi_p)$ contains $J \coloneqq \begin{psmallmatrix} 
(1 + T)^{1/2} & *\\
0 & (1+T)^{-1/2}
\end{psmallmatrix}$. 

We compute adjoint-trace elements. Both \mbox{$a \coloneqq \frac{(\tr M)^2}{\det M} 
= \frac{(s(\lambda_0) + s(\mu_0))^2}{s(\lambda_0)s(\mu_0)} 
= 2 + \frac{s(\lambda_0)}{s(\mu_0)} + \frac{s(\mu_0)}{s(\lambda_0)}$} 
and
\[
b := \frac{(\tr MJ)^2}{\det MJ} =  
\frac{\big(s(\lambda_0)(1+T) +  s(\mu_0)\big)^2}{ s(\lambda_0) s(\mu_0)(1 + T)} = \left(a + \Big(2 + 2 \frac{s(\lambda_0)}{s(\mu_0)} \Big)T + \frac{s(\lambda_0)}{s(\mu_0)} T^2\right)(1+T)^{-1}.
\]
are in $A_0$ by construction.
The last expression shows that $b$ is in $W(\E)\lb T\rb$, since $\lambda_0\mu_0^{-1} \in \E$ by the regularity assumption. Moreover, the $T$-coefficient of $b$ is $2 + 2 \frac{s(\lambda_0)}{s(\mu_0)} - a=  \frac{s(\lambda_0)}{s(\mu_0)} - \frac{s(\mu_0)}{s(\lambda_0)},$ which is in $W(\E)^\times$ since $\lambda_0\mu_0^{-1} \neq \pm 1$. It follows that the closed $W(\E)$-algebra generated by $b$ in $A_0$ is all of $W(\E)\lb T \rb$. In other words $\Lambda \subseteq W(\E)\lb T \rb \subseteq A_0$, as claimed. 

For \eqref{lambdafull}, combine \eqref{hidafull} and \eqref{lambdaina0}. 
For \eqref{asigmafull}, from \eqref{lambdaina0} and the expression for $d$ in \eqref{dkappachi} $(t, d)$ has $A_0$-constant determinant. Now use \cref{A0 Afixedbytwists fullness peers}.
For \eqref{lambdafixed}, combine \eqref{lambdafull}, \cref{fullfixed}, and \eqref{asigmafull}.
\end{proof}

\begin{remark}
The regularity-on-$\Pi_p$ hypothesis in \cref{ourhida}\eqref{lambdaina0} can easily be check in terms of the data of the tame Nebentypus character $\chi$ and the mod-$p$ eigenvalue $a_p$ of $U_p$.  Indeed, since $\delta$ sends Frobenius to $U_p$, to verify regularity we need to check whether $\ovl \epsilon \ovl\delta^{-1}$ takes on a value in $\E^\times \setminus \{\pm 1\}$.  Writing $\chi|_{\Pi_p} = \chi^{\unr}\chi^{\tame}$ with $\chi^{\unr}$ unramified and $\chi^{\tame}$ a character on $\mu_{p-1}$ by local class field theory, we see that the tame part of $\ovl \epsilon\ovl\delta^{-1}$ is $\chi^{\tame}\ovl{\varepsilon}_p^i$ and the unramified part is $\chi^{\unr}\ovl\delta^{-2}$.  Note that the tame part necessarily takes values in $\F_p^\times$, so if $\chi^{\tame}\ovl\varepsilon_p^i$ has order greater than 2, then the regularity-on-$\Pi_p$ hypothesis is automatically satisfied.  Otherwise, one must look to the unramified part and check whether some power of $\chi^{\unr}(p)a_p^{-2}$ lies in $\E^\times \setminus \{\pm 1\}$. 
\end{remark}

\begin{remark}\label{CSTs fix wt space?} {\bf Do all conjugate self-twists of $p$-adic families fix weight space?} In an abstract algebraic setting, given a representation of a profinite group $\Pi$ over a ring $A$ that is finite over $\Lambda = \Z_p\lb T \rb$, we cannot expect to prove that $\rho$ is $\Lambda$-full  --- equivalently, that every conjugate self-twist of $\rho$ fixes~$\Lambda$ --- because it is simply not true.

On the other hand, one intuitively expects $\Lambda$, which parametrizes weight space, to be preserved by any conjugate self-twist of a $p$-adic family. Indeed, if $(\sigma, \chi)$ is a conjugate self-twist that doesn't fix $\Lambda$, then modulo any prime ideal of $\Lambda$, the character $\chi$ will relate pairs of eigenforms of different $p$-adic weights --- an implausible scenario. Therefore, although our $A_0$-fullness result in \cref{ourhida}\eqref{hidafull} is purely algebraic, in order to recover the full strength of \cite{Lang2016}, and to match our intuition of how conjugate self-twists in $p$-adic families behave, we necessarily need to use some modular-form-theoretic input. For Hida families the ordinary condition suffices: see the proof of \cref{ourhida} above, especially part \eqref{lambdaina0}.~
We have not extended this result to Coleman families, which is a drawback on our big-image result in \cref{Coleman family theorem}. Ultimately one hopes to formalize the geometric intuition alluded to above.
\end{remark}

\subsection{Coleman $p$-adic families of classical modular forms}\label{Coleman families}
We now relax the ordinary assumption present in \cref{Hida families} and derive the consequences of our main theorem in the context of pseudorepresentations arising from the Coleman-Mazur eigencurve, comparing with known results in the literature.

Let $X$ be a cuspidal irreducible component of the $p$-adic Coleman-Mazur eigencurve of some fixed tame level (\cite{Buzzard07}, \cite{ColemanMazur98}); having dealt with the ordinary case in \cref{Hida families}, we assume that $X$ is nonordinary.  Let $A$ be the ring of analytic functions on $X$ bounded by $1$.  It is a compact $\Z_p$-algebra (hence pro-$p$) since $X$ is nested \cite[Lemma 7.2.11(ii), Corollary 7.2.12]{Bellaiche-ChenevierBook}.  In fact, the map from $X$ to weight space endows $A$ with a $\Lambda$-algebra structure.  As usual, $A$ is a local domain since $X$ is irreducible.   
As in the Hida family setting, one obtains a $2$-dimensional pseudorepresentation $(t,d) \colon G_\Q \to A$ by gluing together those attached to classical cuspforms in the family, all of which have the same semisimplified residual representation $\ovl \rho \colon G_\Q \to \GL_2(\F)$.  Then $(t,d)$ is unramified outside of a finite set of primes, namely $p$ and the tame level, and thus factors through a $p$-finite quotient $\Pi$ of $G_\Q$.  Unlike the theorems presented in \cref{classical applications1} through \cref{Hida families}, there are no previous fullness results known for $(t,d)$ (though see \cref{compare to CIT} below), so we proceed directly to a statement of our result in this setting.

\begin{theorem}\label{Coleman family theorem}
If $\ovl \rho$ is regular and good if octahedral, then $(t,d)$ is $A_0$-full. 
\end{theorem}

\begin{proof}
By \cref{nonconst det main thm} it suffices to show that $(t,d)$ is not a priori small.  It is not reducible since it is cuspidal.  Any $(t,d)$-representation does not have finite image since $X$ admits classical specializations of weight at least 2.  Finally, since $X$ is nonordinary its CM points are isolated and hence $X$ necessarily admits a classical non-CM positive slope specialization: see \cite[Corollary 3.6]{CIT2016}.  Thus $(t,d)$ is not dihedral. 
\end{proof}

\begin{remark}\label{compare to CIT}
Besides \bellaiche's work \cite{Bellaiche18}, the only previous work on images of Galois representations of finite slope $p$-adic families of modular forms was done in \cite{CIT2016}.  We briefly compare \cref{Coleman family theorem} to their main result \cite[Theorem 1.3]{CIT2016}. 
\begin{itemize}
\item \textbf{Setup:} Rather than working with $A$ as above, they restrict to an irreducible component $\mathbb I^\circ$ of what they call the ``adapted slope $\leq h$ Hecke algebra" --- essentially a bounded-slope piece of $X$ as above: see \cite[\S 3.1]{CIT2016}.   
Note that one can replace $A$ above by $\mathbb I^\circ$ and retain the veracity of \cref{Coleman family theorem}.  
\item \textbf{Assumptions on $\ovl \rho$:} In \cite[Theorem 1.3]{CIT2016} the authors assume that $\ovl \rho$ is absolutely irreducible, even when restricted to the intersection of the kernels of twist characters.  Moreover, their regularity assumption is stronger than ours in that it requires the mod-$p$ eigenvalues of the regular element to be in $\F_p^\times$ rather than requiring their ratio to lie in $\E^\times$.
\item \textbf{Conclusion:} As mentioned above, \cite[Theorem 1.3]{CIT2016} is not a true fullness result.  Rather, it shows ``rigid-Lie fullness" --- a certain rigid analytic Lie algebra attached to the image of a $(t,d)$-representation contains the rigid analytic Lie algebra of a congruence subgroup.  While highly suggestive, one does not know how to recover an actual congruence subgroup in the image from this result.  Following \cite{Lang2016}, Conti and his coauthors show rigid-Lie fullness with respect to the ring fixed by $\Sigma_t(\mathbb I^\circ/\Lambda)$.\qedhere
\end{itemize}
\end{remark}

\begin{remark}
Although the determinant $d$ has a form similar to \eqref{dkappachi} --- the universal character $\kappa$ times a finite-order character --- in this setting we have not proved that $d$ is $A_0$-constant: we do not know whether the image of $\Lambda$ is contained in $A_0$.  See \cref{CSTs fix wt space?} for why one expects this to be true nonetheless.
\end{remark}

\subsection{$p$-Adic families of Hilbert modular forms}\label{Hilbert families}
Since our methods are agnostic about the group $\Pi$, one can proceed with a similar analysis in the context of $p$-adic families of Hilbert modular forms, which we briefly outline here.  We believe these are the first big image results in this context.

As in \cref{classical applications2}, fix a totally real field $F$.  Let $X$ be a cuspidal irreducible component of a $p$-adic eigenvariety interpolating classical Hilbert modular forms over $F$ of a fixed tame level; there are several possible constructions, for instance \cite{Urban} or \cite{AndreattaIovitaPilloniHilbert}.  Let $A$ be the ring of analytic functions on $X$ bounded by $1$, which is again a pro-$p$ local domain.  Gluing together the pseudorepresentations attached to the classical cusp forms parametrized by $X$, all of which have the same semisimplified residual representation $\ovl \rho \colon G_F \to \GL_2(\F)$, yields a $2$-dimensional pseudorepresentation $(t,d) \colon G_F \to A$.  It is unramified outside the tame level and $p$ and hence factors through a $p$-finite quotient $\Pi$ of $G_F$.

\begin{theorem}\label{Hilbert families theorem}
Suppose that $\ovl \rho$ is regular and good if octahedral.  If $X$ admits a non-CM classical specialization, then $(t,d)$ is $A_0$-full.
\end{theorem}

\begin{proof}
By \cref{nonconst det main thm} it suffices to check that $(t,d)$ is not a priori small, which follows from the fact that $X$ admits classical specializations that are cuspidal, not CM, and whose weights are at least $2$.
\end{proof}

\begin{remark}
As in \cref{Coleman families}, we do not know whether $(t,d)$ has $A_0$-constant determinant in this case and hence lack an $A_t^{\Sigma_t}$-fullness result.  As in \cref{CSTs fix wt space?}, we expect the image in $A$ of the ring $\Lambda_F$ of analytic functions on weight space to be contained in $A_0$.  Note that in this case, $\Lambda_F$ is a power series ring over $\Z_p$; the number of variables depends on the totally real field $F$.
\end{remark}

\medskip
\appendix
\section{Algebraic sundries}\label{appx}

\subsection{Representations with isomorphic adjoint differ by a character}\label{adjoint rep section}
Throughout \cref{adjoint rep section}, let $G$ be a group and $F$ a separably closed field of odd characteristic. All representations are assumed to be finite-dimensional.  
Let $\Sl_n(F)$ denote the $F$-vector space of $n \times n$-matrices of trace 0 and $\ad^0 \colon \GL_n(F) \to \GL_{n^2 - 1}(F)$ the representation obtained by letting $\GL_n(F)$ act on $\Sl_n(F)$ by conjugation.  The primary goal of this section is to prove that if $\rho_1, \rho_2 \colon G \to \GL_2(F)$ are semisimple representations such that $\ad^0 \rho_1 \cong \ad^0 \rho_2$, then $\rho_1 \cong \rho_2 \otimes \eta$ for some character $\eta\colon G\to F^\times$. 
This is done in \cref{isomorphic ad implies differ by twist}.  The easier case when the $\rho_i$ are not dihedral is treated first in \cref{the non-dihedral case}.  \cref{the dihedral case} is an analysis of dihedral representations that allows us to conclude \cref{isomorphic ad implies differ by twist} in full generality.  
The results of this section are probably well known to experts, but we give proofs for lack of a reference in the generality we need.  We were guided by Venkatarama's answer to \href{https://mathoverflow.net/questions/297746/}{MathOverflow question 297746}.  In the nondihedral case, this result can be found in \cite[Lemma 2.9]{MurtyPrasad2000}.  When the representations $\rho_1$ and $\rho_2$ arise from classical modular forms, the result can be found in \cite[Appendix]{Ramakrishnan}. 

\subsubsection{The nondihedral case}\label{the non-dihedral case}
Given a representation $\rho \colon G \to \GL_n(F)$, we write $\rho^*$ for its dual representation.  That is, if $V$ is the representation space of $\rho$, then $V^* \coloneqq \Hom(V, F)$ is the representation space of $\rho^*$ with $G$-action given by $(g\varphi)(v) \coloneqq g\varphi(g^{-1}v)$.  In terms of matrices, if we fix a basis for $V$ and take the dual basis for $V^*$, then $\rho^*(g)$ is the inverse transpose of $\rho(g)$.  

If $\rho$ is $2$-dimensional, then an explicit calculation shows that $\rho^* \cong \rho \otimes \Lambda^2 \rho^*$, where $\Lambda^2$ denotes the second exterior power of $\rho$.  (The conjugating matrix can be taken to be $\left(\begin{smallmatrix} 0 & -1\\
1 & 0
\end{smallmatrix}\right)$.)  We have that
\[
1 \oplus \ad^0 \rho \cong \rho \otimes \rho^*.
\]  
In particular, $\ad^0 \rho$ is self dual.  Furthermore, 
\[
1 \oplus \ad^0 \rho \cong \rho \otimes \rho^* \cong \rho \otimes \rho \otimes \Lambda^2 \rho^* \cong 1 \oplus (\Sym^2 \rho \otimes \Lambda^2 \rho^*),
\]
and so $\ad^0 \rho \cong \Sym^2 \rho \otimes \Lambda^2 \rho^* = \Sym^2 \rho \otimes \det \rho^{-1}$.

The following lemma is essentially a version of Schur's lemma that will be useful in what follows.

\begin{lemma}\label{schurs lemma}
If $\rho \colon G \to \GL_n(F)$ is a semisimple representation such that $\ad^0 \rho$ contains a copy of the trivial representation, then $\rho$ is reducible.
\end{lemma}

\begin{proof}
Let $V$ be the $F$-vector space on which $G$ acts via $\rho$.  Then $\End V$ is the representation space for $1 \oplus \ad^0 \rho$, where $1$ is the trivial representation, which corresponds to scalar endomorphisms of $V$.  If $\ad^0 \rho$ contains a copy of the trivial representation, then there is a nonscalar $\varphi \in \End V$ that commutes with the action of $G$.  By Schur's lemma, $\rho$ must be reducible. 
\end{proof}

\begin{lemma}\label{isomorphic ad implies differ by twist with hypothesis}
Let $\rho_1, \rho_2 \colon G \to \GL_2(F)$ be semisimple reducible representations.  
If $\ad^0 \rho_1 \cong \ad^0 \rho_2$, then there exists a character $\eta \colon G \to F^\times$ such that $\rho_1 \cong \eta \otimes \rho_2$.  
\end{lemma}

\begin{proof}
By assumption there exist for $i=1,2$ characters $\lambda_i, \mu_i \colon G \to F^\times$ such that $\rho_i \cong \lambda_i \oplus \mu_i$. 
It is straightforward to calculate 
\[
\lambda_1\mu_1^{-1} \oplus 1 \oplus \lambda_1^{-1}\mu_1 \cong \ad^0 \rho_1 \cong \ad^0 \rho_2 \cong \lambda_2\mu_2^{-1} \oplus 1 \oplus \lambda_2^{-1}\mu_2.
\]
Thus, up to switching $\lambda_2$ and $\mu_2$, we must have $\lambda_1\mu_1^{-1} = \lambda_2\mu_2^{-1}$. 
Let $\eta=\mu_1\mu_2^{-1}$. Then
\[
\rho_1 \cong \lambda_1 \oplus \mu_1 = \lambda_2\mu_1\mu_2^{-1} \oplus \mu_1 = (\mu_1\mu_2^{-1}) \otimes (\lambda_2 \oplus \mu_2) \cong \eta \otimes \rho_2.\qedhere
\]
\end{proof}

\begin{lemma}\label{isomorphic ad implies differ by twist with other hypothesis}
Let $\rho_1, \rho_2 \colon G \to \GL_2(F)$ be semisimple representations such that both $\ad^0 \rho_i$ are irreducible.  
If $\ad^0 \rho_1 \cong \ad^0 \rho_2$, then there exists a character $\eta \colon G \to F^\times$ such that $\rho_1 \cong \eta \otimes \rho_2$.  
\end{lemma}

\begin{proof}
We begin by showing that $\rho_1 \otimes \rho_2$ must be reducible (which does not make use of the assumption that $\ad^0 \rho_i$ is irreducible).  Indeed, by \cref{schurs lemma} if $\rho_1 \otimes \rho_2$ were irreducible then its endomorphism ring would contain a single copy of the trivial representation.  But 
\begin{align*}
\End(\rho_1 \otimes \rho_2) &= (\rho_1 \otimes \rho_2) \otimes (\rho_1 \otimes \rho_2)^* \cong (\rho_1 \otimes \rho_1^*) \otimes (\rho_2 \otimes \rho_2^*)\\
& \cong (1 \oplus \ad^0 \rho_1) \otimes (1 \oplus \ad^0 \rho_2)\\
& \cong 1 \oplus \ad^0 \rho_1 \oplus \ad^0 \rho_1 \oplus (\ad^0 \rho_1 \otimes \ad^0 \rho_1)\\
& \cong 1 \oplus \ad^0 \rho_1 \oplus \ad^0 \rho_1 \oplus (\ad^0 \rho_1 \otimes (\ad^0 \rho_1)^*),
\end{align*}
and $\ad^0 \rho_1 \otimes (\ad^0 \rho_1)^* \cong \End (\ad^0 \rho_1)$ contains a copy of the trivial representation, a contradiction.

Next we show that $\rho_1 \otimes \rho_2$ cannot be the sum of two $2$-dimensional representations.  Indeed, suppose that $\rho_1 \otimes \rho_2 \cong r_1 \oplus r_2$, where $r_1, r_2 \colon G \to \GL_2(F)$ are representations.  Take the second exterior product on both sides.  We have
\[
\Lambda^2 (\rho_1 \otimes \rho_2) \cong (\Lambda^2 \rho_1 \otimes \Sym^2 \rho_2) \oplus (\Sym^2 \rho_1 \otimes \Lambda^2 \rho_2)
\]
and
\[
\Lambda^2 (r_1 \oplus r_2) \cong \Lambda^2 r_1 \oplus \Lambda^2 r_2 \oplus (r_1 \otimes r_2).
\]

Since $\ad^0 \rho_i \cong \Sym^2 \rho_i \otimes \Lambda^2 \rho_i^*$, we have $\Sym^2 \rho_1 \otimes \Lambda^2 \rho_2 \cong \Sym^2 \rho_2 \otimes \Lambda^2 \rho_1$.  But if
\[
(\Lambda^2 \rho_1 \otimes \Sym^2 \rho_2)^{\oplus 2} \cong \Lambda^2 r_1 \oplus \Lambda^2 r_2 \oplus (r_1 \otimes r_2),
\]
then this contradicts irreducibility of $\ad^0 \rho_i$.  Thus $\rho_1 \otimes \rho_2$ must contain a $1$-dimensional representation; call it $\chi$.  Then we claim that $\rho_2 \cong \rho_1^* \otimes \chi \cong \rho_1 \otimes \det \rho_1^{-1} \otimes \chi$, and so $\rho_1$ and $\rho_2$ differ by a twist.  

To see that $\rho_2 \cong \rho_1^* \otimes \chi$, recall that $\rho_1 \otimes \rho_2 \cong \Hom(\rho_1^*, \rho_2)$.  Thus having a $1$-dimensional $G$-stable subspace corresponds to a nonzero linear map $\varphi \colon \rho_1^* \to \rho_2$ such that $g\varphi = \lambda(g)\varphi$ for some $\lambda(g) \in F^\times$ for all $g \in G$.  Define $f \colon \rho_1^* \to \rho_2 \otimes \chi^{-1}$ by $v \mapsto \varphi(v) \otimes e$, where $e$ is a basis for the $1$-dimensional vector space on which $G$ acts by $\chi$.  Note that $f \neq 0$ since $\varphi \neq 0$.  It is straightforward to check that $f(gv) = gf(v)$ for all $g \in G$.  Therefore $\Hom(\rho_1^*, \rho_2 \otimes \chi^{-1}) \neq 0$.  Since $\rho_1^*$ and $\rho_2 \otimes \chi^{-1}$ are irreducible, it follows that they must be isomorphic. 
\end{proof}

The following observation can be checked easily via a direct calculation on $2 \times 2$ matrices.

\begin{lemma}\label{trace of adjoint}
For any $g \in \GL_2(F)$ with (not necessarily distinct) eigenvalues $\lambda, \mu$, the eigenvalues of $\ad^0 g$ are $1, \lambda\mu^{-1}, \lambda^{-1}\mu$.  In particular, we have 
\[
\tr \ad^0 g = \frac{\tr(g)^2}{\det(g)} - 1.
\]
\end{lemma}

\subsubsection{The dihedral case}\label{the dihedral case}
In \cref{the dihedral case} we assume for simplicity that the characteristic of $F$ is not equal to 2.  The goal of \cref{the dihedral case} is to remove the assumption that both $\rho_i$ are reducible or both $\ad^0 \rho_i$ are irreducible from \cref{isomorphic ad implies differ by twist with hypothesis,isomorphic ad implies differ by twist with other hypothesis}.  We begin with a lemma that shows that, in light of \cref{isomorphic ad implies differ by twist with hypothesis,isomorphic ad implies differ by twist with other hypothesis}, we only need to consider the case when both $\rho_1$ and $\rho_2$ are dihedral representations.

\begin{lemma}\label{dihedral or ad irreducible}
If $\rho \colon G \to \GL_2(F)$ is irreducible but $\ad^0 \rho$ is reducible, then $\rho$ is dihedral.  
\end{lemma}

\begin{proof}
If $\ad^0 \rho$ is reducible, then so is $\Sym^2 \rho$ and $\Sym^2 \rho^*$ since $\ad^0 \rho \cong \Sym^2 \rho \otimes \det \rho^{-1}$.  But $\Sym^2 \rho^*$ can be identified with the action of $G$ on the $F$-vector space of quadratic forms on $F^2$.  Thus, there is a quadratic form $Q$ on which $G$ acts by a scalar.  Since $F$ is separably closed and $\mathrm{char} F\neq 2$, all quadratic forms are equivalent.  In particular, we may assume that $Q(x, y) = xy$.  But one checks immediately that the only matrices that preserve $Q$ up to scalars are diagonal and antidiagonal.  Thus $\rho$ must be dihedral. 
\end{proof}

The rest of this section is devoted to an analysis of dihedral representations.

\begin{lemma}\label{cm cst implies eta quadratic}
Assume that $\rho \colon G \to \GL_2(F)$ is a semisimple representation.  If $\rho \cong \eta \otimes \rho$ for some nontrivial character $\eta \colon G \to F^\times$, then the image of $\rho|_{\ker \eta}$ is abelian. 
\end{lemma}

\begin{proof}
This argument essentially comes from \cite[Proposition 4.4]{Ribet77}. 
Note that $\det \rho = \eta^2 \det \rho$ and so $\eta^2=1$.  Set $H \coloneqq \ker \eta$.  Thus $[G \colon H] = 2$ since $\eta$ is nontrivial.  By assumption, there is a matrix $M \in \GL_2(F)$ such that $M\rho(g)M^{-1} = \eta(g)\rho(g)$ for all $g \in G$.  In particular, $\rho(H)$ is contained in the commutant of $M$.  

We claim that $M$ is semisimple.  It suffices to show that $M$ has distinct eigenvalues. Up to 
a change of basis for $\rho$, we may assume that $M$ is upper triangular, say \mbox{$M = \bigl(\begin{smallmatrix}
a & b\\
0 & c
\end{smallmatrix}\bigr)$}.  The eigenvalues of $M$ acting on $M_2(F)$ by conjugation are $1, 1, ac^{-1}, a^{-1}c$ by \cref{trace of adjoint}.  Note that for any $g \in G \setminus H$, we have 
\[
M\rho(g)M^{-1} = -\rho(g).
\]
Thus $-1 = ac^{-1}$, which implies that $a \neq c$ and thus $M$ has distinct eigenvalues, as claimed.  Therefore $M$ is semisimple and so its commutant, and hence $\rho(H)$, is abelian.  
\end{proof}

If $H$ is a subgroup of $G$ of index 2, then we use $c$ to denote a fixed element in $G \setminus H$.  For a character $\chi \colon H \to F^\times$ and $g \in G$, we write $\chi^g \colon H \to F^\times$ for the character defined by $\chi^g(h) \coloneqq \chi(g^{-1}hg)$.  It is not difficult to check that $\chi^g$ depends only on the coset of $g$ in $G/H$.  Set $\chi^- \coloneqq \chi/\chi^c$.  We will write $\eta_H \colon G \to G/H \cong \{\pm1\}$ for the canonical projection map.  With this notation, we recall an explicit description of $\Ind_H^G \chi$.  Namely, $\Ind_H^G \chi$ is isomorphic to the representation
\begin{equation}\label{explicit induced rep}
g \mapsto \begin{cases} \begin{pmatrix} \chi(g) & 0 \\ 0 & \chi^c(g) \end{pmatrix} & \mbox{ if $g \in H$} \vspace{0.1cm} \\
\vspace{0.1cm} \begin{pmatrix} 0 & \chi(gc)\\  \chi^c(gc^{-1}) & 0 \end{pmatrix} & \mbox{ otherwise.} \end{cases}.
\end{equation}
Using Frobenius reciprocity it is easy to see that $\Ind_H^G \chi$ is irreducible if and only if $\chi \neq \chi^c$.

\begin{lemma}\label{resolve definitions of dihedral}
\begin{enumerate}
\item[]
\item If $\rho = \Ind_H^G \chi$ for a character $\chi \colon H \to F^\times$ and $[G \colon H] = 2$, then $\rho \cong \rho \otimes \eta_H$.
\item\label{second point} Conversely, if $\rho \colon G \to \GL_2(F)$ is a dihedral representation, then there is a subgroup $H$ of $G$ of index 2 and a character $\chi \colon H \to F^\times$ such that $\rho \cong \Ind_H^G \chi$ and $\chi \neq \chi^c$.  
\item Furthermore, $H$ as in \eqref{second point} is unique unless $\chi^2 = (\chi^c)^2$.
\item If $\chi^2 = (\chi^c)^2$ then there are exactly three index 2 subgroups $H_i$ of $G$ for $i = 1, 2, 3$ for which there exist characters $\chi_i \colon H_i \to F^\times$ such that $\rho \cong \Ind_{H_i}^G \chi_i$. 
\end{enumerate}
\end{lemma}

\begin{proof}
For the first point, note that $\chi$ is a constituent of $(\rho \otimes \eta_H)|_H = \rho|_H$.  By Frobenius reciprocity and dimension counting, it follows that $\Ind_H^G \chi \cong \rho \otimes \eta_H$.

If $\rho$ is dihedral, then there is a nontrivial character $\eta \colon G \to F^\times$ such that $\tr \rho = \eta \tr \rho$ and $\det \rho = \eta^2 \det \rho$.  In particular, $\eta^2 = 1$ and so $\eta$ is a quadratic character.  Let $H \coloneqq \ker \eta$.  Then $H$ is a subgroup of $G$ of index 2 and $\rho|_H$ is reducible by \cref{cm cst implies eta quadratic}.  Let $\chi \colon H \to F^\times$ be one of the constituents of $\rho|_H$.  By Frobenius reciprocity, $\Ind_H^G \chi$ is a constituent of $\rho$ and we deduce equality for dimension reasons.  Thus we have $\rho|_H = \chi \oplus \chi^c$.  Since $\rho$ is irreducible by the definition of being dihedral, it follows by Frobenius reciprocity that $\chi \neq \chi^c$.  This finishes the proof of the second point.

For the third point, suppose that $\rho = \Ind_{H'}^G \chi'$ for some character $\chi' \colon H' \to F^\times$ and $[G \colon H'] = 2$.  Let $c' \in G \setminus H'$.  Then by restricting to $H$ we have $\chi \oplus \chi^c = (\eta_{H'})|_H \cdot \chi \oplus (\eta_{H'})|_H \cdot \chi^c$.  Thus we either have $\chi = (\eta_{H'})|_H \cdot \chi$ or $\chi = (\eta_{H'})|_H\cdot \chi^c$.  In the first case, we see that $H = \ker \eta_{H'} = H'$.  In the second case we conclude that $\chi^2 = (\chi^c)^2$ since $\eta_{H'}$ is quadratic.  

Finally, suppose that $\chi^2 = (\chi^c)^2$.  Then $H_0 \coloneqq \ker(\chi/\chi^c)$ is a subgroup of index 2 in $H$.  We claim that $H_0$ is normal in $G$.  Recall that $\chi^c$ is independent of the choice of $c \in G \setminus H$.  If $h \in H_0$ and $g \in G \setminus H$ then 
\begin{align*}
\chi(g^{-1}hg)/\chi^c(g^{-1}hg) &= \chi(g^{-1}hg)/\chi^g(g^{-1}hg) = \chi^g(h)/\chi(h)\\
&= (\chi/\chi^g)(h)^{-1} = (\chi/\chi^c)(h)^{-1} = 1.
\end{align*}
Furthermore, the above calculation shows that the class of $c$ generates a subgroup of $G/H_0$ of order 2 distinct from $H$.  Thus $G/H_0$ is isomorphic to $(\Z/2\Z)^2$.  We claim that if $H'$ is any of the three subgroups of $G$ of index 2 containing $H_0$, then there is a character $\chi' \colon H' \to F^\times$ such that $\rho \cong \Ind_{H'}^G \chi'$.  By Frobenius reciprocity, it suffices to show that $\rho|_{H'}$ is reducible.  Since $\rho|_{H_0} = \chi|_{H_0} \oplus \chi^c|_{H_0}$, it follows from Frobenius reciprocity that $\rho|_{H'} = \Ind_{H_0}^{H'} \chi|_{H_0}$.  But $\chi|_{H_0} = \chi^c|_{H_0}$ and so it follows (again by Frobenius reciprocity) that $\rho|_{H'}$ is reducible. 
\end{proof}

Combining the following lemma with Frobenius reciprocity, we see that the irreducibility of $\Ind_H^G \chi$ is related to the question of whether the character $\chi \colon H \to F^\times$ extends to a character of $G$.

\begin{lemma}\label{extending characters}
Let $H$ be a subgroup of $G$ of index 2 and $\chi \colon H \to F^\times$ a character.  Then there exists an extension of $\chi$ to a character $G\to F^\times$ if and only if $\chi = \chi^c$.  If $\chi$ extends to a character of $G$, then there are exactly two different extensions, and they differ by $\eta_H$.
\end{lemma}

\begin{proof}
If such an $L$ and extension of $\chi$ exist, then certainly $\chi = \chi^c$. On the other hand, since $c^2 \in H$, we know that $\chi(c^2)$ is well defined.  Since $F$ is algebraically closed, we may choose a square root $r$ of $\chi(c^2)$ in $F$. Define a new character $\tilde{\chi} \colon G \to L^\times$ by 
\[
\tilde{\chi}(g) \coloneqq \begin{cases} 
\chi(g) & \text{if } g \in H\\
r \chi(c^{-1}g) & \text{if } g \not\in H.
\end{cases}
\] 
To see that $\tilde{\chi}$ is a character, it suffices to verify that it is multiplicative.  That is, one must check that $\tilde{\chi} (h) \tilde{\chi}(ch') = \tilde{\chi}(hch')$ and $\tilde{\chi}(ch) \tilde{\chi}(ch') = \tilde{\chi}(chch')$ for $h, h' \in H$. It is easy to see by direct computation that these are satisfied if $\chi = \chi^c$. 
\end{proof}

\begin{lemma}\label{adjoint of induced rep}
Let $\rho = \Ind_H^G \chi$ be a dihedral representation.  Then $\ad^0 \rho \cong \eta_H \oplus \Ind_H^G \chi^-$.  If $\ad^0 \rho$ is the sum of three characters, then $\chi^2 = (\chi^c)^2$ and $\ad^0 \rho \cong \eta_{H_1} \oplus \eta_{H_2} \oplus \eta_{H_3}$, where the $H_i$ are the index $2$ subgroups of $G$ given in \cref{resolve definitions of dihedral}.
\end{lemma}

\begin{proof}
The first claim is an explicit calculation. 
Let $e_1 \coloneqq \left(\begin{smallmatrix} 1 & 0\\
0 & -1
\end{smallmatrix}\right), e_2 \coloneqq \left(\begin{smallmatrix} 0 & 1\\
0 & 0
\end{smallmatrix}\right), e_3 \coloneqq \left(\begin{smallmatrix} 0 & 0\\
1 & 0
\end{smallmatrix}\right)$.  Assume that $\rho$ is given by \eqref{explicit induced rep}.  Then with respect to the basis $e_1, e_2, e_3$ we see that
\[
\ad^0(g) = \begin{cases} \left(\begin{smallmatrix} 1 & 0 & 0 \\ 0 & \chi^-(g) & 0 \\ 0 & 0 & \chi^-(g)^{-1} \end{smallmatrix}\right) & \mbox{ if $g \in H$} \vspace{0.1cm} \\
\vspace{0.1cm} \left(\begin{smallmatrix} -1 & 0 & 0 \\ 0 & 0 & \chi^-(gc) \\ 0 & \chi^-(gc^{-1})^{-1} & 0 \end{smallmatrix}\right) & \mbox{ otherwise.} \end{cases}
\]
We observe that $\eta_H$ appears in the upper left corner.  Furthermore, $(\chi^-)^c = (\chi^-)^{-1}$.  Therefore the lower right $2 \times 2$-matrix in $\ad^0 \rho$ is isomorphic to $\Ind_H^G \chi^-$ by \eqref{explicit induced rep}.  Thus $\ad^0 \rho \cong \eta_H \oplus \Ind_H^G \chi^-$.

If $\ad^0 \rho$ is the sum of three characters, then $\Ind_H^G \chi^-$ is reducible and thus $\chi^- = (\chi^-)^c$.  That is, $\chi^2 = (\chi^c)^2$.  By \cref{resolve definitions of dihedral}, it follows that there are exactly three subgroups $H_i$ of $G$ of index~2 for which $\rho \cong \Ind_{H_i}^G \chi_i$.  By the above calculation, each $\eta_{H_i}$ must be a constituent of $\ad^0 \rho$.  By counting dimensions, we find that $\ad^0 \rho \cong \eta_{H_1} \oplus \eta_{H_2} \oplus \eta_{H_3}$.  
\end{proof} 

\begin{theorem}\label{isomorphic ad implies differ by twist}
Let $F$ be a field whose characteristic is not 2. Let $\rho_1, \rho_2 \colon G \to \GL_2(F)$ be semisimple representations.  If $\ad^0 \rho_1 \cong \ad^0 \rho_2$ then there is a character $\eta \colon G \to L^\times$ such that $\rho_1 \cong \eta \otimes \rho_2$.
\end{theorem}

\begin{proof}
The case when either of $\rho_1$ or $\rho_2$ is not dihedral is settled by \cref{isomorphic ad implies differ by twist with hypothesis,isomorphic ad implies differ by twist with other hypothesis,dihedral or ad irreducible}. Therefore we may assume that both $\rho_1$ and $\rho_2$ are dihedral.  By \cref{adjoint of induced rep} there are index-2 subgroups $H_i$ of $G$ and characters $\chi_i \colon H_i \to F^\times$ such that $\rho_i \cong \Ind_{H_i}^G \chi_i$.  Note that the set of possible such $H_i$ can be read off from $\ad^0 \rho_i$ since $\eta_{H_i}$ is a constituent of $\ad^0 \rho_i$ by \cref{adjoint of induced rep} and $H_i = \ker \eta_{H_i}$.  In particular, since $\ad^0 \rho_1 \cong \ad^0 \rho_2$, we may assume that $H \coloneqq H_1 = H_2$.  By \cref{adjoint of induced rep} we have $\Ind_H^G \chi_1^- \cong \Ind_H^G \chi_2^-$. 
By restricting to $H$ it follows that $\chi_1^- \oplus (\chi_1^-)^c \cong \chi_2^- \oplus (\chi_2^-)^c$, and so up to replacing $\chi_2$ with $\chi_2^c$ (which is okay since $\Ind_H^G \chi_2 \cong \Ind_H^G \chi_2^c$), it follows that $\chi_1^- = \chi_2^-$.  That is, $\chi_1\chi_2^{-1} = (\chi_1\chi_2^{-1})^c$.  By \cref{extending characters} there is a character $\eta \colon G \to L^\times$ such that $\eta|_H = \chi_1\chi_2^{-1}$.  We claim that $\rho_1 \cong \eta \otimes \rho_2$.  Indeed, this is true upon restriction to $H$ since 
\[
\rho_1|_H = \chi_1 \oplus \chi_1^c = \eta|_H \otimes (\chi_2 \oplus \chi_2^c) = (\eta \otimes \rho_2)|_H.
\]
Therefore $\rho_1 \cong \eta \otimes \rho_2$ by Frobenius reciprocity since $\rho_1$ is irreducible and thus $\chi_1 \neq \chi_1^c$.
\end{proof}


\subsection{Trace field extensions} \label{tracefieldsec}
In this brief subsection, let $G$ be an abstract group, $F$ an abstract field, and $(t, d): G \to F$ a ({2-dimensional}) pseudorepresentation.

\begin{definition}  \label{tracefielddef} The \emph{trace field} of $(t, d)$ is the subfield of $F$ generated by the image $t(G)$ over the prime subfield of $F$. A pseudorepresentation $(t, d): G \to F$ is \emph{realizable} over an extension $L$ of $F$ if there {exists} a semisimple representation $\rho: G \to \GL_2(L)$ that \emph{carries} $(t, d)$ --- that is, with $t = \tr \rho$ and $d = \det \rho$. 
\end{definition}  

\begin{lemma}\label{atmostquad} 
Let $(t, d): G \to F$ is a pseudorepresentation with trace field $F$. If the characteristic of $F$ is not $2$, then $(t, d)$ realizable over an at-most quadratic extension of $F$. 
\end{lemma}

\begin{proof} To start with, $(t, d)$ is always realizable by a semisimple representation $V$ over $\ovl F$ \cite[Theorem~2.12]{ChenevierDet}. If we suppose that $V$ is irreducible, then the image of the associated $F$-algebra map \mbox{$\rho: F[G] \to \End_{\ovl F}(V)$} surjects onto the full matrix algebra $\End_{\ovl F}(V)$ after extending scalars to $\ovl F$, and is therefore a quaternion algebra $D$ over $F$. There are now two possibilities.
Either $D$ is split, in which case $D^\times \cong \GL_2(F)$ is a realization of $(t, d)$ over $F$. Or $D$ is an $F$-division algebra, in which case \mbox{$\left.\rho\right|_G: G \to D^\times$} carries $(t, d)$, in the sense that $\rho(g)$ has reduced trace $t(g)$ and reduced norm~$d(g)$, and any quadratic extension $L$ of $F$ that splits $D$ carries a realization of $(t, d)$ as an irreducible representation \mbox{$G \to \GL_2(L)$}. 

On the other hand, suppose $(t,d)$ splits into a sum of two characters \mbox{$\chi, \chi': G \to \ovl F^\times$}. The image of $\chi$ is a subgroup of $\ovl F^\times$ whose every element is contained in an at-most-quadratic extension of $F$. Suppose that $\alpha = \chi(a)$ and $\beta =\chi(b)$ for $a, b \in G$ generate \emph{different} quadratic extensions of $F$. Then on one hand, $\alpha \beta = \chi(ab)$ must generate the third quadratic subextension of $F(\alpha, \beta)$. But on the other hand, we claim that $\chi'(ab) = \chi'(a) \chi'(b)$ is equal to $\alpha \beta$: indeed,  let $c_\alpha$ be the generator of $\Gal\big(F(\alpha, \beta)/F(\beta)\big)$ viewed as an element of  $\Gal\big(F(\alpha, \beta)/F\big)$, so that $\chi'(a) = c_\alpha(\alpha)$; define $c_\beta$ similarly. Then $c_\alpha c_\beta$ generates $\Gal\big(F (\alpha, \beta)/F(\alpha \beta)\big)$ and hence fixes $\alpha \beta$. Therefore $\chi(ab) + \chi'(ab) = 2\alpha \beta$, which is not in $F$\footnote{The constraint on the characteristic is necessary: consider $G = \Z^2$ and $F = \F_2(x, y)$, with $(t, d)$ the pseudocharacter corresponding to the scalar $2$-dimensional representation $(1, 0)$ to $\sqrt{x}$ and $(0, 1)$ to $\sqrt{y}$.}  --- a contradiction.
\end{proof}

The following proposition is true in any dimension, but we state it here for dimension $2$. 
\begin{proposition}\label{tracefieldextension}
Let $(t,d): G \to F$ be a pseudorepresentation and $H \subseteq G$ a finite-index normal subgroup. Then the trace field of $(t, d)$ is a finite extension of the trace field of $(\left.t\right|_H, \left.d\right|_H)$.
\end{proposition}

\begin{proof}
Replace $F$ by the trace field of $(t, d)$, and let $E \subseteq F$ be the trace field of $(\left. t\right|_H, \left.d \right|_H)$.  
Let $\rho$ be a semisimple representation carrying $(t, d)$ over an extension of $F$. Note that $F$ is algebraic over $E$: every $g \in G$ satisfies $g^{[G:H]} \in H$, so that every eigenvalue of $\rho(g)$, and hence $t(g)$, is algebraic over the finite extension of $E$ containing the eigenvalues of $g^{[G:H]}$. Since $F$ is contained in the field generated by all these  eigenvalues, $F$ is algebraic over $E$, and in particular $\rho$ is realizable over an algebraic extension of $E$.

By Clifford's theorem \cite[Theorem 7.1.1]{Craven}, $\left.\rho\right|_H$ is still semisimple, so it is carried by a representation $V$ over a finite extension $E'$  of $E$. Since $F/E$ is algebraic, $V_{\ovl E} = V \otimes_{E'} \ovl E$ carries all of $\rho$. 
Let $g_1, \ldots, g_n$ be coset representatives for $H$ in $G$; write each $\rho(g_i)$ as a matrix in a fixed $E'$-basis of $V$ extended to $V_{\ovl E}$. Let $M$ be the subfield of $\ovl E$ generated over $E'$ by the matrix coefficients of all the $\rho(g_i)$. Then $M$ is a finite extension of $E'$, and hence of $E$, containing the values of $\tr \rho$.  Therefore $F/E$ is finite. 
\end{proof}

\subsection{Rings with involution}\label{commutative algebra lemmas}

Throughout \cref{commutative algebra lemmas}, let $A$ be a commutative noetherian ring equipped with an involution $*$.  Note that we will need to apply the results in this section to the universal constant-determinant pseudodeformation ring $\cA$, so we cannot assume that $A$ is a domain.  Let $A^\varepsilon = \{a \in A \colon a^* = \varepsilon a\}$ for $\varepsilon\in\{+,-\}$.  We will assume throughout that $*$ is not the identity on $A$ so that $A^- \neq 0$.  It is easy to see that $A^+$ is a subring of $A$ and $A^-$ is an $A^+$-module.  The following results have been adapted from \cite{Lanski1974} and \cite{ChuangLee1977}, where they are presented in the context when $A$ may be noncommutative.

\begin{definition}
We say that an $A$-ideal $\Aa$ is a \textit{$*$-ideal} if $\Aa^* = \Aa$.  We say that $A$ is \textit{$*$-prime} if whenever $\Aa$ and $\Bb$ are $*$-ideals such that $\Aa\Bb = 0$ then either $\Aa = 0$ or $\Bb = 0$.
\end{definition}

\begin{lemma}\label{* prime implies reduced}
If $A$ is $*$-prime then $A$ is reduced.
\end{lemma}

\begin{proof}
Let $0 \neq a\in A$ be nilpotent.  Then there is a smallest integer $n > 1$ such that $a^n = 0$.  Let $\Aa = aA$ and $\Bb = a^{n-1}A$.  Note that $\Aa \neq 0$ and $\Bb \neq 0$ by the minimality of $n$.  If $\Aa$ and $\Bb$ are $*$-ideals then we have reached a contradiction since $\Aa\Bb = a^nA = 0$.  In particular, if $a + a^* = 0$ then $a^* = -a \in aA$ and so $\Aa, \Bb$ are $*$-ideals.

If $a + a^* \neq 0$, then $a + a^*$ is still nilpotent since $A$ is commutative.  By replacing $a$ with $a + a^*$ in the above argument, we find that $\Aa$ and $\Bb$ are $*$-ideals and thus we reach a contradiction.  
\end{proof}

\begin{lemma}\label{R+ noetherian}
If $A$ is a noetherian commutative ring with $2 \in A^\times$, then $A^+$ is a noetherian ring.
\end{lemma}

\begin{proof}
The following argument comes from \cite[Lemma]{ChuangLee1977}.  Let $I_1 \subseteq I_2 \subseteq \cdots$ be an ascending chain of ideals in $A^+$.  Then $I_1A \subseteq I_2A \subseteq \cdots$ is an ascending chain of ideals in $A$.  Since $A$ is noetherian, there is some $n$ such that $I_nA = I_mA$ for all $m \geq n$.  

Fix $m \geq n$ and $a \in I_m \subseteq A^+$.  Since $a \in I_mA = I_nA$ we may write
\[
a = \sum_i b_ix_i
\]
with $b_i \in I_n$ and $x_i \in A$.  Applying the involution $*$ yields
\[
a = a^* = \sum_i b_ix_i^*.
\]
Thus
\[
2a = \sum_i b_i(x_i + x_i^*).
\]
Since $x_i + x_i^* \in A^+$ and $2 \in A^\times$ it follows that $a = \frac{1}{2}\sum_i b_i(x_i + x_i^*) \in I_n$.  In particular, $I_m = I_n$.
\end{proof}

We would like to show that $A$ is finitely generated as an $A^+$-module, which is equivalent to $A$ being a noetherian $A^+$-module since $A^+$ is a noetherian ring by \cref{R+ noetherian}.  The following lemma follows the proof of \cite[Lemma 6]{Lanski1974}.

\begin{lemma}\label{the domain case}
If there is an element $d \in A^-$ that is not a zero divisor in $A$, then $A$ is noetherian as an $A^+$-module.
\end{lemma}

\begin{proof}
Since $d$ is not a zero divisor, it follows that $A$ is isomorphic to $dA$ as an $A^+$-module.  On the other hand, for any $a \in A$ we can write
\[
da = \frac{1}{2}(d(a - a^*)) + \frac{1}{2}(d(a+a^*)) \in A^+ + dA^+.
\]
Thus $dA$ is a submodule of the finitely generated $A^+$-module $A^+ + dA^+$.  Since $A^+$ is noetherian by \cref{R+ noetherian}, it follows that $dA$, and hence $A$, is a finitely generated (and hence noetherian) $A^+$-module.
\end{proof}

\begin{proposition}\label{R fg as R+ module}
If $A$ is a commutative noetherian ring with $2 \in A^\times$, then $A$ is a noetherian $A^+$-module.
\end{proposition}

\begin{proof}
This proof combines elements of the proofs of \cite[Theorem]{ChuangLee1977} and \cite[Theorem 7]{Lanski1974}.

Suppose not.  Let $\Aa_0$ be the largest $*$-ideal of $A$ such that $A/\Aa_0$ is not a noetherian $A^+$-module, which exists since $A$ is a noetherian ring and is not noetherian as an $A^+$-module.  Thus, by replacing $A$ with $A/\Aa_0$, we may assume that $A/\Aa$ is a noetherian $A^+$-module for any $*$-ideal $\Aa \neq 0$.  

We claim that, under this assumption, $A$ is reduced.  It suffices to show that $A$ is $*$-prime by \cref{* prime implies reduced}.  Suppose that $\Aa$ and $\Bb$ are nonzero $*$-ideals of $A$ such that $\Aa\Bb = 0$.  Note that we can view $\Aa$ as an $A/\Bb$-module since $\Aa\Bb = 0$.  We know that $\Aa$ is noetherian as an $A/\Bb$-module since $\Aa$ is noetherian as an $A$-module.  Furthermore, $A/\Bb$ is a noetherian $A^+$-module since $\Bb \neq 0$.  Thus $\Aa$ is noetherian as an $A^+$-module.  We also know that $A/\Aa$ is a noetherian $A^+$-module since $\Aa \neq 0$.  Therefore $A$ is a noetherian $A^+$-module, a contradiction.  Thus $A$ is $*$-prime and hence reduced.

Since $A$ is a noetherian ring, it has only finitely many minimal prime ideals; call them $\p_1, \ldots, \p_n$.  Since $A$ is reduced, we have that
\[
\bigcap_{i = 1}^n \p_i = 0.
\]
Note that $n = 1$ corresponds to the case when $A$ is a domain, and in that case we have already seen that $A$ is a noetherian $A^+$-module by \cref{the domain case}.  Thus we assume henceforth that $n > 1$ and thus each $\p_i \neq 0$.  

If $\p_i^* \cap \p_i \neq 0$, then $\p_i \cap \p_i^*$ is a $*$-ideal and so $A/(\p_i \cap \p_i^*)$ is a noetherian $A^+$-module.  If every $\p_i$ satisfies $\p_i \cap \p_i^* \neq 0$ then we can view $A$ as a subring of
\[
\bigoplus_{i = 1}^n A/(\p_i \cap \p_i^*),
\]
which is noetherian as an $A^+$-module.  In particular, $A$ is a noetherian $A^+$-module, a contradiction, which proves the proposition.

Suppose there is some $k$ such that $\p_k \cap \p_k^* = 0$.  It is easy to check that $\p_k^*$ is another minimal prime ideal of $A$.  We claim that $n = 2$ in this case.  Indeed, if $\p$ is any minimal prime ideal of $A$, then we have $\p_k\p_k^* \subseteq \p_k \cap \p_k^* = 0$ and thus
$
\p_k\p_k^* = 0 \in \p.
$
Thus $\p = \p_k$ or $\p = \p_k^*$.  

Let us write $\p = \p_k$ henceforth.  We can embed $A$ into $A/\p \times A/\p^*$ by identifying $a \in A$ with $(a + \p, a + \p^*)$.  Note that $A^+ \cap \p = 0$ since if $a \in A^+ \cap \p$ then $a = a^* \in \p^* \cap \p = 0$.  Similarly, $A^+ \cap \p^* = 0$.  In particular, $A^+$ injects into $A/\p$ and is therefore a domain.  

Note that by \cref{the domain case}, we may assume that every element of $A^-$ is a zero divisor in $A$.  However, both $A/\p$ and $A/\p^*$ are domains, so the only zero divisors in $A/\p \times A/\p^*$ are elements of the form $(a + \p, \p^*)$ or $(\p, a + \p^*)$.  Recall that $(A^-)^2 \subseteq A^+$.  In particular, if $(a + \p, \p^*) \in A^-$, then $(a^2 + \p, \p^*) \in A^+$.  That is, there is some $a_+ \in A^+$ such that $a_+ - a^2 \in \p$ and $a_+ \in \p^*$.  But we have already seen that $A^+ \cap \p^* = 0$.  Similarly, any $(\p, a + \p^*) \in A^-$ must be trivial.  In other words, $A^- = 0$, a contradiction.  Therefore $A$ must be noetherian as an $A^+$-module.
\end{proof}

Given any ideal $\Aa$ of $A$, we define $\Aa^\varepsilon \coloneqq \Aa \cap A^\varepsilon$.  We call $\Aa$ a \textit{graded} ideal if $\Aa = \Aa^+ \oplus \Aa^-$.

\begin{proposition}\label{nongraded quotients}
Let $A$ be a commutative local noetherian ring such that $A$ and $A^+$ have the same residue field.  Assume that $2 \in A^\times$.  If $A'$ is the quotient of $A$ by a nongraded prime ideal, then $A'$ has the same field of fractions as the image of $A^+$ in $A'$.
\end{proposition}

\begin{proof}
Write $f \colon A \to A'$ for the quotient map.  It suffices to show that every element of $f(A^-)$ can be written as a quotient of elements in $f(A^+)$.  Since the prime ideal $\p = \ker f$ is assumed to be nongraded, it follows that there is some $a \in \p$ such that, if we decompose $a = a^+ + a^-$ with $a^+ \in A^+$ and $a^- \in A^-$, then neither $a^+$ nor $a^-$ is in $\p$.  It follows that $f(a^-) = -f(a^+)$, and so $f(a^-) \in f(A^+)$.  Note that $f(a^-) \neq 0$ since $a^- \not\in \p$.  For any $x \in A^-$ we have that $xa^- \in A^+$ since $(A^-)^2 \subseteq A^+$.  Thus $f(x) = f(xa^-)/f(a^-) \in Q(f(A^+))$, as desired.
\end{proof}

\subsection{Automorphisms and gradings}\label{autom grad}
We recall how ring automorphisms give rise to gradings. 

Let $A$ be a complete local ring and $X$ a finite abelian subgroup of the group of ring automorphisms of $A$. We write $\mu_n(A) \coloneqq \{a \in A^\times \colon a^n = 1\}$. 
Given a character $\varphi \colon X \to A^\times$, we define
\[
A^\varphi \coloneqq \{a \in A \colon \act{\sigma}{a} = \varphi(\sigma)a, \forall \sigma \in X\}.
\]

The following lemma is standard, so we leave the proof to the reader. 

\begin{lemma}\label{roots of unity in local rings}
Let $\F$ be a finite field of characteristic $p$ and $A$ a pro-$p$ local ring with residue field~$\F$.  If $p \nmid n$, then $\mu_n(A) = s(\mu_n(\F))$.
\end{lemma}

Assume the following:

\smallskip

$(\ast)$ for every positive integer $n$, if $X$ contains an element of order $n$, then $\#\mu_n(A) = n$. 

\smallskip

Then one has $\#X = \#\Hom(X, A^\times).$  (It is easily checked when $X$ is cyclic, and then for general~$X$ one applies the structure theorem of finite abelian groups.)

\begin{corollary}\label{using hypothesis ast}
Assume ($\ast)$.  If $p \nmid \#X$, then for any $1 \neq \sigma \in X$ we have
\[
\sum_{\varphi \in \Hom(X, A^\times)} \varphi(\sigma) = 0.
\]
\end{corollary}

\begin{proof}
First suppose that $X$ is cyclic of order $n$ and $\sigma$ is a generator for $X$.  Then $\Hom(X, A^\times)$ is cyclic, generated by any $\varphi_0$ such that $\varphi_0(\sigma)$ is a primitive $n^{\rm th}$ root of unity.  Let $H \coloneqq \langle \varphi_0^k \rangle$ be a nontrivial subgroup of $\Hom(X, A^\times)$.  Then by \cref{roots of unity in local rings} we have
\[
\sum_{\varphi \in H} \varphi(\sigma) = \sum_{i = 0}^{n/k} \varphi_0^{ki}(\sigma) = \sum_{\omega \in \mu_{n/k}(A)} \omega = \sum_{\omega \in \mu_{n/k}(\F)} s(\omega) = 0.
\]

Now we allow $X$ to be any finite abelian group such that $p \nmid \#X$ and $\sigma$ any nontrivial element of $X$.  Then we have an exact sequence
\[
0 \to \Hom(X/\langle \sigma \rangle, A^\times) \to \Hom(X, A^\times) \to \Hom(\langle \sigma \rangle, A^\times).
\]
Thus $\sum_{\varphi \in \Hom(X, A^\times)} \varphi(\sigma)$ is an integral multiple of 
$
\sum_{\varphi \in H} \varphi(\sigma),
$
where $H$ is the image of $\Hom(X, A^\times)$ in $\Hom(\langle \sigma \rangle, A^\times)$.  This sum is $0$ by the first paragraph. 
\end{proof}

\begin{lemma}\label{automorphisms to gradings}
Let $A$ and $X$ be as above.  Assume that $\#X \in A^\times$ and that condition $(\ast)$ holds.  Then $A$ admits a grading given by $A = \bigoplus_{\varphi \in \Hom(X, A^\times)} A^\varphi$.
Furthermore, for any $\Z[1/\#X][X]$-submodule $M \subseteq A$, letting $M^\varphi \coloneqq M \cap A^\varphi$, there is a decomposition
\[
M = \bigoplus_{\varphi \in \Hom(X, A^\times)} M^\varphi.
\]
\end{lemma}

\begin{proof}
For $\varphi \in \Hom(X, A^\times)$, define
$
e_\varphi \coloneqq \frac{1}{\#X}\sum_{\sigma \in X} \varphi(\sigma)\sigma^{-1} \in \Z[1/\#X][X].
$
A straightforward computation shows that $\{e_\varphi \colon \varphi \in \Hom(X, A^\times)\}$ is an orthogonal system of idempotents in $\Z[1/\#X][X]$.  (Note that \cref{using hypothesis ast} is needed to show that $\sum_{\varphi} e_\varphi = 1$.)  There is a natural ring homomorphism
$
\Z[1/\#X][X] \to \End A;
$
pushing forward the $e_\varphi$ to $\End A$ gives the result.
\end{proof}

\vspace{0.5cm}

\bibliography{CSTreferences}

\begin{thebibliography}{KMP00}

\bibitem[AB]{BoeckleArias}
Sara Aryas{-}de{-}Reina and Gebhard B{\"{o}}ckle.
\newblock Deformation rings and images of {G}alois representations.
\newblock ArXiv:
  {\href{https://arxiv.org/abs/2107.03114}{https://arxiv.org/abs/2107.03114}}.

\bibitem[AIP16]{AndreattaIovitaPilloniHilbert}
Fabrizio Andreatta, Adrian Iovita, and Vincent Pilloni.
\newblock On overconvergent {H}ilbert modular cusp forms.
\newblock {\em Ast\'{e}risque}, (382):163--193, 2016.

\bibitem[AM69]{AtiyahMacdonald}
M.~F. Atiyah and I.~G. Macdonald.
\newblock {\em Introduction to commutative algebra}.
\newblock Addison-Wesley Publishing Co., Reading, Mass.-London-Don Mills, Ont.,
  1969.

\bibitem[Amo21]{Amaros20}
Laia Amor\'{o}s.
\newblock Images of {G}alois representations in mod $p$ {H}ecke algebras.
\newblock {\em International Journal of Number Theory}, 17(5):1265--1285, 2021.

\bibitem[BC09]{Bellaiche-ChenevierBook}
Jo\"el Bella\"iche and Ga\"etan Chenevier.
\newblock Families of {G}alois representations and {S}elmer groups.
\newblock {\em Ast\'erisque}, 324:xii+314, 2009.

\bibitem[Bel19]{Bellaiche18}
Jo\"{e}l Bella\"{\i}che.
\newblock Image of pseudo-representations and coefficients of modular forms
  modulo {$p$}.
\newblock {\em Adv. Math.}, 353:647--721, 2019.

\bibitem[BH07]{BergerHarcos07}
Tobias Berger and Gergely Harcos.
\newblock {$l$}-{A}dic representations associated to modular forms over
  imaginary quadratic fields.
\newblock {\em Int. Math. Res. Not. IMRN}, (23):Art. ID rnm113, 16, 2007.

\bibitem[B{\"o}c13]{BoeckleDef}
Gebhard B{\"o}ckle.
\newblock Deforming {G}alois representations.
\newblock In {\em Elliptic curves, {H}ilbert modular forms and {G}alois
  deformations}, Adv. Courses Math. CRM Barcelona, pages 21--115.
  Birkh\"{a}user/Springer, Basel, 2013.

\bibitem[Buz07]{Buzzard07}
Kevin Buzzard.
\newblock Eigenvarieties.
\newblock In {\em {$L$}-functions and {G}alois representations}, volume 320 of
  {\em London Math. Soc. Lecture Note Ser.}, pages 59--120. Cambridge Univ.
  Press, Cambridge, 2007.

\bibitem[Car94]{Carayol94}
Henri Carayol.
\newblock Formes modulaires et repr{\'e}sentations galoisiennes {\`a} valeurs
  dans un anneau local complet.
\newblock {\em Contemp. Math.}, 165:213--237, 1994.

\bibitem[CEG]{CalegariEmertonGee}
Frank Calegari, Matthew Emerton, and Toby Gee.
\newblock Globally realizable components of local deformation rings.
\newblock To appear. ArXiv: \url{https://arxiv.org/abs/1807.03529}.

\bibitem[Che14]{ChenevierDet}
Ga{\"e}tan Chenevier.
\newblock The $p$-adic analytic space of pseudocharacters of a profinite group,
  and pseudorepresentations over arbitrary rings.
\newblock In {\em Automorphic forms and Galois representations. Proceedings of
  the LMS Durham Symposium 2011 (1)}, volume 414 of {\em London Math. Soc.
  Lecture Note Ser.}, pages 221--285. Cambridge Univ. Press, 2014.

\bibitem[CIT16]{CIT2016}
Andrea Conti, Adrian Iovita, and Jacques Tilouine.
\newblock Big image of {G}alois representations associated with finite slope
  {$p$}-adic families of modular forms.
\newblock In {\em Elliptic curves, modular forms and {I}wasawa theory}, volume
  188 of {\em Springer Proc. Math. Stat.}, pages 87--123. Springer, Cham, 2016.

\bibitem[CL77]{ChuangLee1977}
Cheng-Liang Chuang and Pjek-Hwee Lee.
\newblock Noetherian rings with involution.
\newblock {\em Chinese J. Math.}, 5(1):15--19, 1977.

\bibitem[CM98]{ColemanMazur98}
Robert Coleman and Barry Mazur.
\newblock The eigencurve.
\newblock In {\em Galois representations in arithmetic algebraic geometry
  ({D}urham, 1996)}, volume 254 of {\em London Math. Soc. Lecture Note Ser.},
  pages 1--113. Cambridge Univ. Press, Cambridge, 1998.

\bibitem[CM09]{CalegariMazur09}
Frank Calegari and Barry Mazur.
\newblock Nearly ordinary {G}alois deformations over arbitrary number fields.
\newblock {\em J. Inst. Math. Jussieu}, 8(1):99--177, 2009.

\bibitem[Cra19]{Craven}
David~A. Craven.
\newblock {\em Representation theory of finite groups: A guide}.
\newblock Springer, 2019.

\bibitem[Dic58]{Dickson58}
Leonard~Eugene Dickson.
\newblock {\em Linear groups: {W}ith an exposition of the {G}alois field
  theory}.
\newblock {W}ith an introduction by W. Magnus. Dover Publications, Inc., New
  York, 1958.

\bibitem[DK00]{Ramakrishnan}
William Duke and Emmanuel Kowalski.
\newblock A problem of {L}innik for elliptic curves and mean-value estimates
  for automorphic, with an appendix by {D}inakar {R}amakrishnan.
\newblock {\em Invent. Math.}, 139(1):1--39, 2000.

\bibitem[dSL97]{deSmitLenstra}
Bart de~Smit and Hendrik~W. Lenstra, Jr.
\newblock Explicit construction of universal deformation rings.
\newblock In {\em Modular forms and {F}ermat's last theorem ({B}oston, {MA},
  1995)}, pages 313--326. Springer, New York, 1997.

\bibitem[Eak86]{Eakin}
Paul~M. Eakin, Jr.
\newblock The converse to a well known theorem of {N}oetherian rings.
\newblock {\em Math. Ann.}, 177:278--282, 1986.

\bibitem[Eis95]{EisenbudCommAlg}
David Eisenbud.
\newblock {\em Commutative algebra}, volume 150 of {\em Graduate Texts in
  Mathematics}.
\newblock Springer-Verlag, New York, 1995.
\newblock With a view toward algebraic geometry.

\bibitem[Fis02]{Fischman02}
Ami Fischman.
\newblock On the image of {$\Lambda$}-adic {G}alois representations.
\newblock {\em Ann. Inst. Fourier (Grenoble)}, 52(2):351--378, 2002.

\bibitem[Hid12]{HidaGMF}
Haruzo Hida.
\newblock {\em Geometric modular forms and elliptic curves}.
\newblock World Scientific Publishing Co. Pte. Ltd., Hackensack, NJ, second
  edition, 2012.

\bibitem[Hid15]{Hida15}
Haruzo Hida.
\newblock Big {G}alois representations and {$p$}-adic {$L$}-functions.
\newblock {\em Compos. Math}, 151(4):603--664, 2015.

\bibitem[KMP00]{MurtyPrasad2000}
Vijaya Kumar~Murty and Dipendra Prasad.
\newblock Tate {C}ycles on a {P}roduct of {T}wo {H}ilbert {M}odular {S}urfaces.
\newblock {\em J. Number Theory}, 80(1):25--43, 2000.

\bibitem[Lan56]{Lang56}
Serge Lang.
\newblock Algebraic groups over finite fields.
\newblock {\em Amer. J. Math.}, 78:555--563, 1956.

\bibitem[Lan16]{Lang2016}
Jaclyn Lang.
\newblock On the image of the {G}alois representation associated to a non-{CM}
  {H}ida family.
\newblock {\em Algebra Number Theory}, 10(1):155--194, 2016.

\bibitem[Lan75]{Lanski1974}
Charles Lanski.
\newblock Chain conditions in rings with involution.
\newblock {\em J. London Math. Soc. (2)}, 9:93--102, 1974/75.

\bibitem[Man15]{Manohar15}
Jayanta Manoharmayum.
\newblock A structure theorem for subgroups of {$\GL_n$} over complete local
  noetherian rings with large residual image.
\newblock {\em Proceedings of the American Mathematical Society}, 7:2743--2758,
  2015.

\bibitem[Mat70]{matsumura}
Hideyuki Matsumura.
\newblock {\em Commutative algebra.}
\newblock Mathematics lecture note series. W. A. Benjamin, New York, 1970.

\bibitem[Maz89]{Mazur1989}
Barry Mazur.
\newblock Deforming {G}alois representations.
\newblock In {\em Galois groups over {${\bf Q}$} ({B}erkeley, {CA}, 1987)},
  volume~16 of {\em Math. Sci. Res. Inst. Publ.}, pages 385--437. Springer, New
  York, 1989.

\bibitem[Mok14]{Mok14}
Chung~Pang Mok.
\newblock Galois representations attached to automorphic forms on {${\rm
  GL}_2$} over {CM} fields.
\newblock {\em Compos. Math.}, 150(4):523--567, 2014.

\bibitem[Mom81]{Momose81}
Fumiyuki Momose.
\newblock On the {$l$}-adic representations attached to modular forms.
\newblock {\em J. Fac. Sci. Univ. Tokyo Sect. IA Math.}, 28(1):89--109, 1981.

\bibitem[MW86]{MazurWiles86}
Barry Mazur and Andrew Wiles.
\newblock On {$p$}-adic analytic families of {G}alois representations.
\newblock {\em Compositio Math.}, 59(2):231--264, 1986.

\bibitem[Nek12]{Nekovar2012}
Jan Nekov\'{a}\v{r}.
\newblock Level raising and anticyclotomic {S}elmer groups for {H}ilbert
  modular forms of weight two.
\newblock {\em Canad. J. Math.}, 64(3):588--668, 2012.

\bibitem[Nys96]{Nyssen96}
Louise Nyssen.
\newblock Pseudo-repr\'esentations.
\newblock {\em Math. Ann.}, 306(2):257--283, 1996.

\bibitem[Pin93]{Pink93}
Richard Pink.
\newblock Classification of pro-{$p$} subgroups of {$\text{SL}_2$} over a
  {$p$}-adic ring, where {$p$} is an odd prime.
\newblock {\em Compositio Math.}, 88(3):251--264, 1993.

\bibitem[Rib75]{RibetEllAdic}
Kenneth~A. Ribet.
\newblock On $\ell$-adic representations attached to modular forms.
\newblock {\em Invent. Math.}, 28:245--276, 1975.

\bibitem[Rib77]{Ribet77}
Kenneth~A. Ribet.
\newblock Galois representations attached to eigenforms with {N}ebentypus.
\newblock In {\em Modular functions of one variable, {V} ({P}roc. {S}econd
  {I}nternat. {C}onf., {U}niv. {B}onn, {B}onn, 1976)}, volume 601 of {\em
  Lecture Notes in Math.}, pages 17--51. Springer, Berlin, 1977.

\bibitem[Rib85]{Ribet85}
Kenneth~A. Ribet.
\newblock On {$l$}-adic representations attached to modular forms. {II}.
\newblock {\em Glasgow Math. J.}, 27:185--194, 1985.

\bibitem[Rou96]{Rouquier96}
Rapha{\"e}l Rouquier.
\newblock Caract\'erisation des caract\'eres et pseudo-caract{\`e}res.
\newblock {\em J. Algebra}, 180(2):571--586, 1996.

\bibitem[RV91]{FANT}
Dinakar Ramakrishnan and Robert~J. Valenza.
\newblock {\em Fourier Analysis on Number Fields}, volume 186 of {\em GTM}.
\newblock Springer, 1991.

\bibitem[Ser68]{Serre68}
Jean-Pierre Serre.
\newblock {\em Abelian {$l$}-adic representations and elliptic curves}.
\newblock McGill University lecture notes written with the collaboration of
  Willem Kuyk and John Labute. W. A. Benjamin, Inc., New York-Amsterdam, 1968.

\bibitem[Ser72]{SerrePropGalois}
Jean-Pierre Serre.
\newblock Propri\'et\'es galoisiennes des points d'ordre fini des courbes
  elliptiques.
\newblock {\em Invent. Math.}, 15:259--331, 1972.

\bibitem[Ser97]{SerreGalCoh}
Jean-Pierre Serre.
\newblock {\em Galois Cohomology}.
\newblock Springer, 1997.

\bibitem[Ser19]{Serban19}
Vlad Serban.
\newblock A finiteness result for {$p$}-adic families of {B}ianchi modular
  forms.
\newblock arXiv 1902.03217, 2019.

\bibitem[Tay91]{Taylor91}
Richard Taylor.
\newblock Galois representations associated to {S}iegel modular forms of low
  weight.
\newblock {\em Duke Math. J.}, 63:281--332, 1991.

\bibitem[Tay94]{Taylor94}
Richard Taylor.
\newblock {$l$}-adic representations associated to modular forms over imaginary
  quadratic fields. {II}.
\newblock {\em Invent. Math.}, 116(1-3):619--643, 1994.

\bibitem[Urb11]{Urban}
Eric Urban.
\newblock Eigenvarieties for reductive groups.
\newblock {\em Ann. of Math. (2)}, 174(3):1685--1784, 2011.

\bibitem[Wil90]{Wiles90}
Andrew Wiles.
\newblock The {I}wasawa conjecture for totally real fields.
\newblock {\em Ann. of Math. (2)}, 131(3):493--540, 1990.

\end{thebibliography}
\bibliographystyle{alpha}

\end{document}